\documentclass[11pt]{book}
\usepackage{array}
\usepackage{makeidx}
\usepackage{amsthm}
\usepackage{amsmath}
\usepackage{amsfonts}
\usepackage{eucal}
\usepackage{graphicx}

\def\cal{\CMcal}
\pdfoutput=1  

\usepackage[usenames,dvipsnames]{color}
\usepackage{color}
 \definecolor{programcode}{gray}{0.9}

\usepackage[colorlinks,%
linkcolor=BrickRed,%
filecolor =BrickRed,%
citecolor=RoyalPurple,%
]{hyperref}

\def\head#1{\medbreak\noindent\textbf{\textit{ #1}\quad }}

\def\Page#1{page~\pageref{#1}}

\def\FPE{\cite{frapowema91}}
\def\CHE{\cite{che84}}
\def\BRO{\cite{bro91}}
\def\ZD{\cite{zaddes69}}

\def\optchar{{\circ}}

\def\Span{\hbox{\rm Span}\,}
 
\def\ebox#1#2{%
\medskip
\begin{center} 
  \strut\epsfxsize=#1 \hsize\epsfbox{#2} 
\end{center}
\smallbreak}

\def\ebox#1#2{%
\smallbreak \centerline{
\includegraphics[width= #1\hsize]{#2}  
}
}

\def\eig{\hbox{eig}}

\def\defn#1{\textit{#1}\index{#1}}

\def\zero{\vartheta}

\def\inp#1,#2>{\langle #1 , #2 \rangle}

\def\underH{{ \cal H}}

\newcounter{l1}
\newcounter{l2}
\newcounter{l3}

\newcommand{\barablist}{\begin{list}{\arabic{l1}}{\usecounter{l1}}}
\newcommand{\balphlist}{\begin{list}{(\alph{l2})}{\usecounter{l2}} 
\setlength{\topsep}{0pt}
\setlength{\partopsep}{0pt}
\setlength{\itemsep}{0pt}
\setlength{\parsep}{5pt}
\setlength{\leftmargin}{-10pt}
\setlength{\rightmargin}{0pt}
\setlength{\itemindent}{-5pt}
}
\newcommand{\bromalist}{\begin{list}{(\roman{l3})}{\usecounter{l3}}}

\def\hwlabel#1{\label{hw-#1}}
\def\Exercise#1{Exercise~\ref{hw-#1}}

\newenvironment{matlab}{%
\subsection*{}
\begin{center}
\large
 \textit{\textbf{Matlab Commands}}
\end{center}
\begin{description}}{\end{description}}

\newenvironment{summary}{
\clearpage
\begin{center}
\large
{\textit{\textbf{Summary and References}}}
\end{center}
}{}

\newcounter{exercise}

\newenvironment{exercises}{%
\clearpage 
\section{Exercises} 
\setcounter{exercise}{0}%
\begin{list} 
{\textbf{\thesection.\textit{\arabic{exercise}}}}{\usecounter{exercise}}%
\setlength{\topsep}{0pt}
\setlength{\partopsep}{0pt}
\setlength{\itemsep}{0pt}
\setlength{\parsep}{5pt}
\setlength{\leftmargin}{-10pt}
\setlength{\rightmargin}{0pt}
\setlength{\itemindent}{-5pt} }
{\end{list}
\bigbreak}

\def\mindex#1{\index{#1}}

\def\notes#1{}

\def\clabel#1{\label{c:#1}}

\def\tlabel#1{\label{t:#1}}

\def\flabel#1{\label{f:#1}}

\def\elabel#1{\label{e:#1}}

\def\qed{\ifmmode\sq\else{\unskip\nobreak\hfil
\penalty50\hskip1em\null\nobreak\hfil\sq
\parfillskip=0pt\finalhyphendemerits=0\endgraf}\fi}

\long\def\defbox#1{\framebox[.9\hsize][c]{\parbox{.85\hsize}{%
\parindent=0pt
\baselineskip=12pt plus .1pt      
\parskip=6pt plus 1.5pt minus 1pt 
 #1}}}

\long\def\beginbox#1\endbox{\subsection*{}%
\hbox{\hspace{.05\hsize}\defbox{\medskip#1\bigskip}}%
\subsection*{}}

\def\endbox{}

\def\rank{{\rm rank\,}}

\def\sgn{{\rm sgn}}

\newsavebox{\junk}
\savebox{\junk}[1.6mm]{\hbox{$|\!|\!|$}}

\def\det{{\mathop{\rm det}}}

\def\argmin{\mathop{\rm arg\, min}}

\def\Co{\mathbb{C}}
\def\Re{\mathbb{R}}

\def\posRe{\mathbb{R}_+}

\newcommand{\hatV}{{\widehat{V}}}

\newcommand{\hatx}{{\widehat{x}}}
\newcommand{\haty}{{\widehat{y}}}

\def\tilx{{\widetilde{x}}}

\def\tilA{\widetilde A}
\def\tilB{\widetilde B}

\def\tilK{\widetilde K}

\def\til={{\widetilde =}}

\def\clA{{\cal A}}
\def\clB{{\cal B}}
\def\clC{{\cal C}}

\def\clF{{\cal F}}

\def\clL{{\cal L}}

\def\clN{{\cal N}}
\def\clO{{\cal O}}

\def\clR{{\cal R}}

\def\clU{{\cal U}}
\def\clV{{\cal V}}
\def\clW{{\cal W}}
\def\clX{{\cal X}}
\def\clY{{\cal Y}}

\def\proof{\head{Proof}}

\def\half{{\mathchoice{\textstyle \frac{1}{2}}%
{\textstyle{\frac{1}{2}}}%
{\hbox{\tiny $1\frac{1}{2}$}}%
{\hbox{\tiny $1\frac{1}{2}$}} }}

\def\eqdef{\mathbin{:=}}

 \def\eq#1/{(\ref{#1})}

\def\eye(#1){{\bf(#1)}\quad}
\def\epsy{\varepsilon}

\def\varble{\,\cdot\,}

\newtheorem{theorem}{Theorem}[chapter]
\newtheorem{lemma}[theorem]{Lemma}

\theoremstyle{definition}
\newtheorem{example}{Example}[section]

\newenvironment{ex}{\begin{example}}{\qed\end{example}}

\def\Lemma#1{Lemma~\ref{t:#1}}

\def\Theorem#1{Theorem~\ref{t:#1}}

\def\Section#1{Section~\ref{s:#1}}

\def\Figure#1{Figure~\ref{f:#1}}
\def\Chapter#1{Chapter~\ref{c:#1}}

\def\eq#1/{(\ref{e:#1})}

\def\beq{\begin{equation}}
\def\eeq{\end{equation}}

\def\beqa{\begin{eqnarray}}
\def\eeqa{\end{eqnarray}}

\def\bdes{\begin{description}}
\def\edes{\end{description}}

\newcommand{\dotP}{{\dot{P}}}

\newcommand{\dotX}{{\dot{X}}}
\newcommand{\dotY}{{\dot{Y}}}

\newcommand{\dote}{{\dot{e}}}

\newcommand{\dotp}{{\dot{p}}}

\newcommand{\dotx}{{\dot{x}}}
\newcommand{\doty}{{\dot{y}}}
\newcommand{\dotz}{{\dot{z}}}

\newcommand{\doteta}{{\dot{\eta}}}

\newcommand{\dotlambda}{{\dot{\lambda}}}

\newcommand{\barA}{{\bar{A}}}
\newcommand{\barB}{{\bar{B}}}
\newcommand{\barC}{{\bar{C}}}
\newcommand{\barD}{{\bar{D}}}

\newcommand{\barP}{{\bar{P}}}

\newcommand{\barW}{{\bar{W}}}

\newcommand{\barc}{{\bar{c}}}

\newcommand{\bare}{{\bar{e}}}

\newcommand{\baro}{{\bar{o}}}

\newcommand{\baru}{{\bar{u}}}

\newcommand{\barx}{{\bar{x}}}
\newcommand{\bary}{{\bar{y}}}

{\end{list}}

\def\ass(#1:#2){(#1\ref{#1:#2})}

\def\ritem#1{
\item[{\sf \ass(\current_model:#1)}]
}

\newenvironment{recall-ass}[1]{%
\begin{description}
\def\current_model{#1}}{
\end{description}
} 

\def\sq{\hbox{\rlap{$\sqcap$}$\sqcup$}}

\def\utw{\smash{\rlap{\lower5pt\hbox{$\sim$}}}}
\def\udtw{\smash{\rlap{\lower6pt\hbox{$\approx$}}}}

\def\bbbz{{\mathchoice {\hbox{$\sf\textstyle Z\kern-0.4em Z$}}
{\hbox{$\sf\textstyle Z\kern-0.4em Z$}}
{\hbox{$\sf\scriptstyle Z\kern-0.3em Z$}}
{\hbox{$\sf\scriptscriptstyle Z\kern-0.2em Z$}}}}
 
\def\diameter{{\ifmmode\oslash\else$\oslash$\fi}}

\def\FRAC#1#2#3{\genfrac{}{}{}{#1}{#2}{#3}}

\def\ddt{{\mathchoice{\FRAC{1}{d}{dt}}%
{\FRAC{1}{d}{dt}}%
{\FRAC{3}{d}{dt}}%
{\FRAC{3}{d}{dt}}}}

 \def\archive#1{}
\makeindex

\begin{document}

\pagenumbering{roman}
\thispagestyle{empty} 

\begin{center}
 {\bf\LARGE 
 ECE 515 \hfill Class Notes}
 \\[2cm] 
 \LARGE
Lecture Notes on 
\\[.5cm] 
\huge
  {\bf CONTROL SYSTEM THEORY}
\\[.2cm]
	{\bf AND DESIGN}  
\\[2cm]
{\Large Tamer Ba\c{s}ar, Sean P.\ Meyn, and William R.\ Perkins }
  \\[.2cm]
{\Large  University of Illinois at Urbana-Champaign}
  \\[5cm] 
  \normalsize
Originally prepared January 2007  (this version, July 2024)
 \end{center}

 \clearpage

\thispagestyle{empty} 

\chapter*{Preface}
\addcontentsline{toc}{chapter}{Preface}

This is a collection of the lecture notes of the three authors for
a first-year graduate course on control system theory and design
({\sl ECE~515}, formerly ECE~415) at the ECE Department of the
University of Illinois at Urbana-Champaign.  This is a fundamental
course on the modern theory of dynamical systems and their
control, and builds on a first-level course in control that
emphasizes frequency-domain methods (such as the course {\sl
ECE~486}, formerly ECE~386, at {\sl UIUC}).  The emphasis in this
graduate course is on state space techniques, and it encompasses
{\it modeling}, {\it analysis} (of structural properties of
systems, such as stability, controllability, and observability),
{\it synthesis} (of observers/compensators and controllers)
subject to design specifications, and {\it optimization}.
Accordingly, this set of lecture notes is organized in four parts,
with each part dealing with one of the issues identified above.
Concentration is on {\it linear systems}, with nonlinear systems
covered only in some specific contexts, such as stability and
dynamic optimization. Both continuous-time and discrete-time
systems are covered, with the former, however, in much greater
depth than the latter.

The notions of ``control'' and ``feedback'', in precisely the sense
they will be treated and interpreted in this text, pervade our
everyday operations, oftentimes without us being aware of it.  Leaving
aside the facts that the human body is a large (and very complex)
feedback mechanism, and an economy without its build-in (and
periodically fine-tuned) feedback loops would instantaneously turn to
chaos, the most common examples of control systems in the average
person's everyday life is the thermostat in one's living room, or the
cruise control in one's their automobile.  Control systems of a similar
nature can be found in almost any of today's industries.  The most
obvious examples are the aerospace industry, where control is required
in fly by wire positioning of ailerons, or in the optimal choice of
trajectories in space flight.  The chemical industry requires good
control designs to ensure safe and accurate production of specific
products, and the paper industry requires accurate control to produce
high quality paper.  Even in applications where control has not found
any use, this may be a result of inertia within industry rather than
lack of need!  For instance, in the manufacture of semiconductors
currently many of the fabrication steps are done without the use of
feedback, and only now are engineers seeing the difficulties that such
open-loop processing causes.

There are only a few fundamental ideas that are required to take this
course, other than a good background in linear algebra and
differential equations.  All of these ideas revolve around the concept
of {\it feedback}, which is simply the act of using measurements as
they become available to control the system of interest.  The main
objective of this course is to teach the student some fundamental
principles within a solid conceptual framework, that will enable
her/him to design feedback loops compatible with the information
available on the ``states'' of the system to be controlled, and by
taking into account considerations such as stability, performance,
energy conservation, and even robustness.  A second objective is to
familiarize her/him with the available modern computational,
simulation, and general software tools that facilitate the design of
effective feedback loops.\\ \\
\hfill TB, SPM, WRP \\
\hfill Urbana, January 2010

\vspace{2cm}

\noindent
\textit{Notes for the 2024 edition.}   The notes were revised to correct an  error in the definition of $y_t$ in \Figure{mod-ball}---many thanks to Professor Bruce Hajek at UIUC for bringing this to our attention.   

Some examples in the early chapters of these notes  appear in the recent monograph \textit{Control Systems and Reinforcement Learning} \cite{CSRL}.
  Part 1 aims to clarify the  relationship between concepts in these lecture notes 
  and reinforcement learning,  without relying on any theory from stochastic processes.

\tableofcontents
\newpage

\listoffigures

\pagenumbering{arabic}

\part{System Modeling and Analysis}

\pagenumbering{arabic}
\notes{BROGAN: Chap 3, 15.1}

\chapter{State Space Models}\clabel{chap1.1}

After a first course in control system design one learns that
intuition is a starting point in control design, but that intuition
may fail for complex systems such as those with multiple inputs and
outputs, or systems with nonlinearities.  In order to augment our
intuition to deal with problems as complex as high speed flight
control, or flow control for high speed communication networks, one
must start with a useful mathematical model of the system to be
controlled.  It is not necessary to ``take the system apart'' - to
model every screw, valve, and axle.  In fact, to use the methods to be
developed in this book it is frequently more useful to find a simple
model which gives a reasonably accurate description of system
behavior.  It may be difficult to verify a complex model, in which
case it will be difficult to trust that the model accurately describes
the system to be controlled.  More crucial in the context of this
course is that the control design methodology developed in this text
is model based.  Consequently, a complex model of the system will
result in a complex control solution, which is usually highly
undesirable in practice.

Although this text treats nonlinear models in some detail, the most
far reaching approximation made is \defn{linearity}.  It is likely
that no physical system is truly linear.  The reasons that control
theory is effective in practice even when this ubiquitous assumption
fails are that (i) physical systems can frequently be approximated by
linear models; and (ii) control systems are designed to be robust with
respect to inaccuracies in the system model.  It is not simply luck
that physical systems can be modeled with reasonable accuracy using
linear models.  Newton's laws, Kirchoff's voltage and current laws,
and many other laws of physics give rise to linear models.  Moreover,
in this chapter we show that a generalized Taylor series expansion
allows the approximation of even a grossly nonlinear model by a linear
one.

The linear models we consider are primarily of the state space form
\begin{equation}
\begin{aligned} 
\dot x &= A x + Bu 
\\ y &=Cx+ Du
\end{aligned} 
 \elabel{stateSpace}
\end{equation}
where $x$ is a vector signal in $\Re^n$, $y$ and $u$ are the output
and input of the system evolving in $\Re^p$ and $\Re^m$, respectively,
and $A, B,C,D$ are matrices of appropriate dimensions.  If these
matrices do not depend on the time variable, $t$, the corresponding
linear system will be referred to as \textit{ linear time invariant}
(LTI), whereas if any one of these matrices has time-varying entries,
then the underlying linear system will be called \textit{ linear time
varying} (LTV). If both the input and the output are scalar, then we
refer to the system as \textit{single input-single output} (SISO); if
either a control or output are of dimension higher than one, then the
system is
\textit{multi input-multi output} (MIMO).

Our first example now is a simple nonlinear system which we
approximate by a linear state space model of the form \eq stateSpace/
through linearization.

.

\begin{figure}[ht]
\ebox{.7}{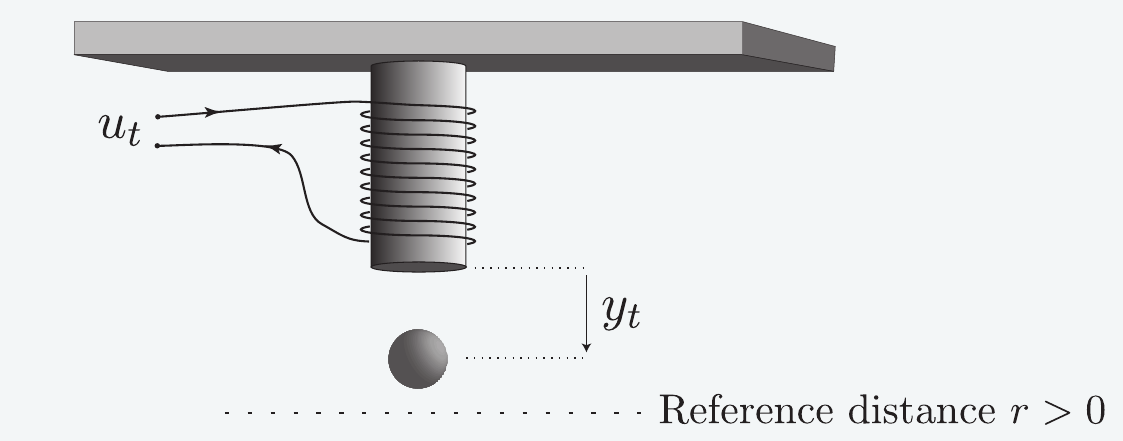}

\caption{Magnetically Suspended Ball}
\flabel{mod-ball}
\end{figure}

\section{An electromechanical system}
The magnetically suspended metallic ball illustrated in
\Figure{mod-ball} is a simple example which illustrates some of the
important modeling issues addressed in this chapter.
The input $u$ is
the current applied to the electro-magnet, and the output $y$ is the
distance between the center of the ball and some reference
height. Since positive and negative inputs are indistinguishable at
the output of this system, it follows that this cannot be a linear
system.  The upward force due to the current input is approximately
proportional to $u^2/y^2$, and hence from Newton's law for
translational motion we have
\begin{eqnarray*}
ma \ = \ m \ddot y\  = \ mg - c \frac{u^2}{y^2},
\end{eqnarray*}
where $g$ is the gravitational constant and $c$ is some constant
depending on the physical properties of the magnet and ball.  This
input-output model can be converted to (nonlinear) state space form
using $x_1 = y$ and $x_2 = \dot y$:
\begin{eqnarray*}
\dot x_1 = x_2, && \dot x_2 = g - \frac{c}{m }\frac{u^2}{x_1^2}
\end{eqnarray*}
where the latter equation follows from the formula $\dot x_2 =\ddot
y$. This pair of equations forms a two-dimensional state space model
\begin{eqnarray}
\dot x_1 & = & x_2  = f_1 ( x_1,  x_2, u) \elabel{eq1.1.1} \\
\dot x_2 & = & g - \frac{c  }{ m }\frac{u^2 }{ x^2_1} = f_2 (x_1, x_2,
u) \elabel{eq1.1.2}
\end{eqnarray}
It is nonlinear, since $f_2$ is a nonlinear function of
$\bigl(\begin{smallmatrix} x_1\\ x_2\end{smallmatrix}\bigr)$.  Letting
$x=\bigl(\begin{smallmatrix}x_1\\ x_2\end{smallmatrix}\bigr)$, and
$f=\bigl(\begin{smallmatrix}f_1\\ f_2\end{smallmatrix}\bigr)$, the
state equations may be written succinctly as
\begin{eqnarray*}
\dot x = f(x,u).
\end{eqnarray*}
The motion of a typical solution to a nonlinear state space model in
$\Re^2$ is illustrated in \Figure{mod-nonlinear}.
\begin{figure}[ht]
\ebox{.7}{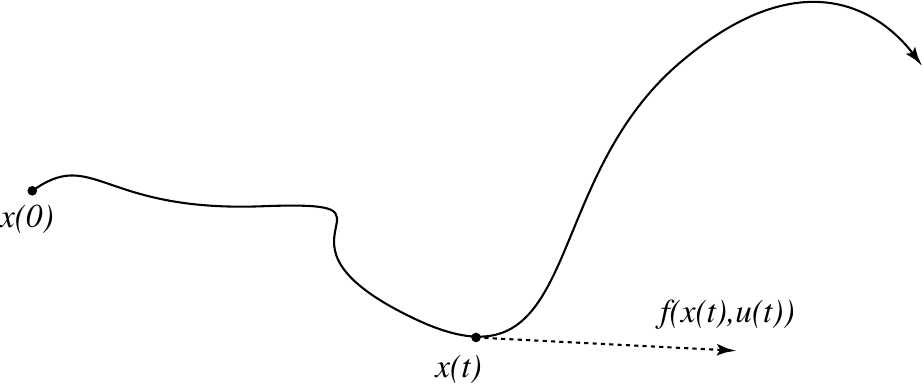}
\caption[Trajectory of a 2D nonlinear state space model]{Trajectory of a nonlinear state space model in two
dimensions: $\dot x = f(x,u)$} \flabel{mod-nonlinear}
\end{figure}

\section{Linearization about an equilibrium state}
\label{s:linearization}

Suppose that a fixed value $u_e$, say positive, is applied, and that
the state $x_e=\left(\begin{smallmatrix}x_{e1}\\
x_{e2}\end{smallmatrix}\right)$ has the property that
\begin{eqnarray*}
f(x_{e1},x_{e2},u_e) & = &
\left(\begin{smallmatrix}f_1(x_{e1},x_{e2},u_e) \\
f_2(x_{e1},x_{e2},u_e)\end{smallmatrix}\right) =\zero.
\end{eqnarray*}
From the definition of $f$ we must have $x_{e2}=0$, and
\begin{eqnarray*}
x_{e1} & = & \sqrt{\frac{c}{mg}} u_e
\end{eqnarray*}
which is unique when restricted to be positive.  The state $x_e$ is
called an \defn{equilibrium state} since the velocity vector $f(x,u)$
vanishes when $x=x_e$, and $u=u_e$.  If the signals $x_1(t)$, $x_2(t)$
and $u(t)$ remain close to the fixed point $(x_{e1},x_{e2},u_e)$, then
we may write
\begin{eqnarray*}
x_1 (t) & = & x_{e1} + \delta x_1 (t) \\ x_2 (t) & = & x_{e2} + \delta
x_2 (t) \\ u (t) & = & u_e + \delta u (t),
\end{eqnarray*}
where $\delta x_1(t)$, $\delta x_2(t)$, and $\delta u(t)$ are
small-amplitude signals.  From the state equations \eq eq1.1.1/ and
\eq eq1.1.2/ we then have
\begin{eqnarray*}
\delta \dot x_1 & = & x_{e2} + \delta x_2 (t) = \delta x_2 (t) \\
\delta \dot x_2 & = & f_2 (x_{e1} + \delta x_1, x_{e2} + \delta x_2,
u_e + \delta u)
\end{eqnarray*}
Applying a Taylor series expansion to the right hand side (RHS) of the
second equation above gives
\begin{eqnarray*}
\delta \dot x_2 & = & f_2 (x_{e1}, x_{e2}, u_e ) + \frac{\partial f_2
}{ \partial x_1} \Big|_{(x_{e1},x_{e2}, u_e )} \delta x_1 +
\frac{\partial f_2 }{ \partial x_2} \Big|_{(x_{e1},x_{e2},u_e )}
\delta x_2 \\ && + \frac{\partial f_2 }{ \partial u}
\Big|_{(x_{e1},x_{e2},u_e )} \delta u + \hbox{\it H.O.T.}
\end{eqnarray*}
After computing partial derivatives we obtain the formulae
\begin{eqnarray*}
\dot{\delta x_1 } & = & \delta x_2.  \\ \delta \dot x_2 & = & 2
\frac{c }{ m} \frac{u_e ^2 }{ x_{e1}^3} \delta x_1 - \frac{2c }{ m}
\frac{u_e }{ x_{e1}^2} \delta u + \hbox{\it H.O.T.}
\end{eqnarray*}
Letting $x$ denote the bivariate signal $x(t) =
\left(\begin{smallmatrix}x_1(t)\\ x_2(t)\end{smallmatrix}\right)$ we
may write the linearized system in state space form:
\begin{eqnarray*}
\delta \dot x & = & \left[ \begin{array}{cc} 0 & 1 \\ \alpha & 0
\end{array} \right] \delta x + \left[
\begin{array}{c} 0 \\ \beta \end{array}
\right] \delta u \\ \delta y & = & \delta x_1\,,
\end{eqnarray*}
where
\begin{eqnarray*}
\alpha = 2 \frac{c }{ m}\frac{u_e ^2 }{ x_{e1}^3}, && \beta = - 2
\frac{c }{ m} \frac{u_e }{ x_{e1}^2}.
\end{eqnarray*}
This linearized system is only an approximation, but one that is
reasonable and useful in control design as long as the state $\delta
x$ and the control $\delta u$ remain small. For example, we will see
in \Chapter{stability} that ``local'' stability of the nonlinear
system is guaranteed if the simpler linear system is stable.

\section{Linearization about a trajectory}

Consider now a general nonlinear model
\begin{eqnarray*}
\dot x\, (t) & = & f (x(t), u(t), t)
\end{eqnarray*}
where $f$ is a continuously differentiable ($C^1$) function of its
arguments.  Suppose that for an initial condition $x_0$ given at time
$t_0$ (i.e., $x(t_0) = x_0$), and a given input $u^n(t)$, the solution
of the nonlinear state equation above exists and is denoted by
$x^n(t)$, which we will refer to as the {\it nominal trajectory}
corresponding to the {\it nominal input} $u^n(t)$, for the given
initial state.  For example, in the control of a satellite orbiting
the earth, the nominal trajectory $x^n(t)$ might represent a desired
orbit (see \Exercise{satellite} below). We assume that the input and
state approximately follow these nominal trajectories, and we again
write
\begin{eqnarray*}
x (t) & = & x^n (t) + \delta x (t) \\ u (t) & = & u^n (t) + \delta u
(t)
\end{eqnarray*}
From a Taylor series expansion we then have
\begin{eqnarray*}
\dot x & = & \dot x^n + \delta \dot x = f (x^n, u^n, t) + \Bigl(
\underbrace{ \frac{\partial f }{ \partial x} \Big|_{(x^n,u^n,t)} }_{A
(t)} \Bigr) \delta x + \Bigl( \underbrace{ \frac{\partial f }{
\partial u} \Big|_{(x^n,u^n,t)}}_{B (t)} \Bigr) \delta u + H. O. T.
\end{eqnarray*}
Since we must have $\dot x^n = f(x^n, u^n,t)$, this gives the state
space description
\begin{eqnarray*}
\delta \dot x = A (t) \delta x + B (t) \delta u ,
\end{eqnarray*}
where the higher order terms have been ignored.

\section{A two link inverted pendulum}

Below is a photograph of the \textit{pendubot} found at the robotics
laboratory at the University of Illinois, and a sketch indicating its
component parts.  The pendubot consists of two rigid aluminum links:
link 1 is directly coupled to the shaft of a DC motor mounted to the
end of a table.  Link 1 also includes the bearing housing for the
second joint.  Two optical encoders provide position measurements: one
is attached at the elbow joint and the other is attached to the motor.
Note that no motor is directly connected to link 2 - this makes
vertical control of the system, as shown in the photograph, extremely
difficult!
\begin{figure}[ht]
\ebox{.95}{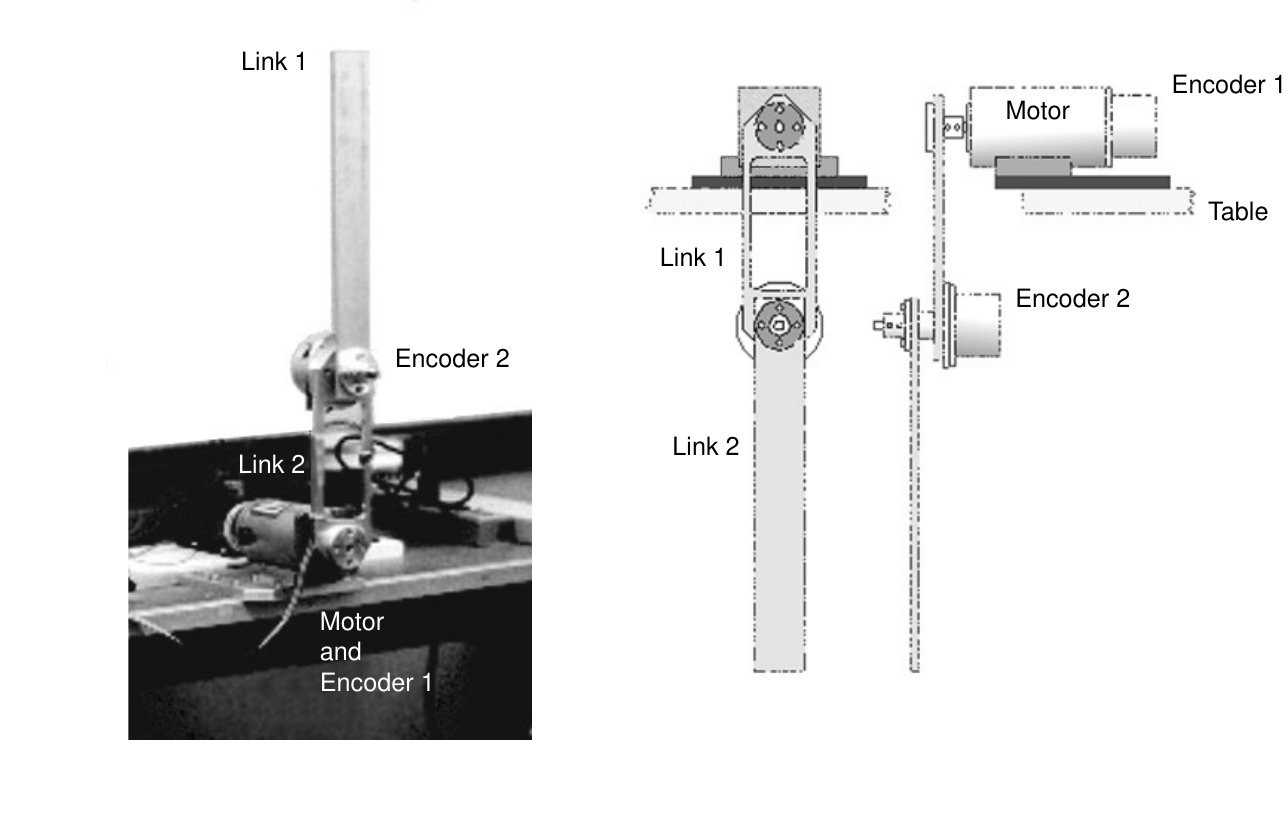}
\caption{The Pendubot}
\flabel{mod-pendubot}
\end{figure}

Since the pendubot is a two link robot with a single actuator, its
dynamic equations can be derived using the so-called Euler-Lagrange
equations found in numerous robotics textbooks \cite{spovid89}.
Referring to the figure, the equations of motion are
\begin{eqnarray}
d_{11} {\ddot q}_1 + d_{12} {\ddot q}_2 + h_1 + \phi_1 &=& \tau \elabel{tau}\\
d_{21} {\ddot q}_1 + d_{22} {\ddot q}_2 + h_2 + \phi_2 &=&  0 \elabel{pen}
\end{eqnarray}
where $q_1$, $q_2$ are the joint angles and
\begin{eqnarray*}
d_{11} &=&  m_1\ell_{c1}^2 + m_2 (\ell_1^2 + \ell_{c2}^2 + 2\ell_1\ell_{c2}
\cos (q_2)) + I_1 + I_2\\
d_{22} &=&  m_2\ell_{c2}^2 + I_2\\
d_{12} &=& d_{21} =  m_2 ( \ell_{c2}^2 + \ell_1\ell_{c2} \cos (q_2)) + I_2\\
h_1 &=&  -m_2\ell_1\ell_{c2} \sin (q_2){\dot q}_2^2 -2m_2\ell_1\ell_{c2}
\sin (q_2){\dot q}_2{\dot q}_1 \\
h_2 &=&  m_2\ell_1\ell_{c2} \sin (q_2){\dot q}_1^2\\
\phi_1 &=& (m_1\ell_{c1} + m_2\ell_1 )g\cos (q_1) + m_2\ell_{c2} g\cos (q_1+q_2)\\
\phi_2 &=& m_2\ell_{c2} g\cos (q_1+q_2)\\
\end{eqnarray*}
The definitions of the variables $q_i, \ell_1, \ell_{ci}$ can be deduced from
\Figure{mod-pendubot2}.

\begin{figure}[ht]
\ebox{.55}{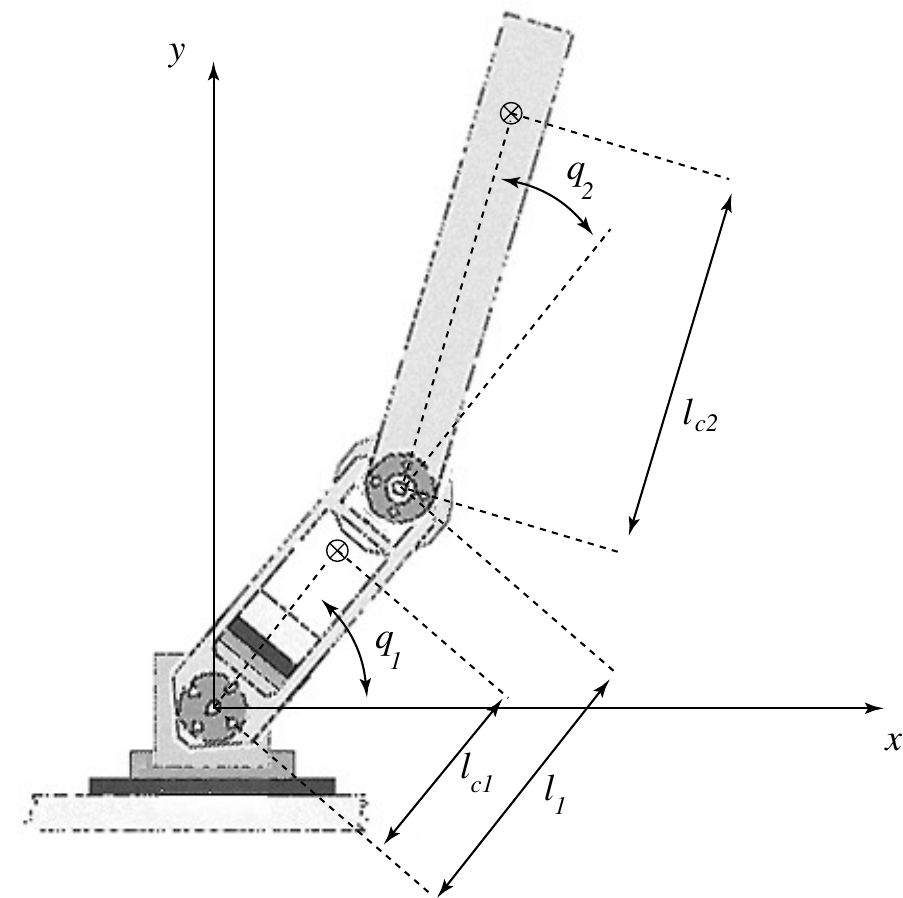}
\caption[Coordinate description of the pendubot]{Coordinate description of the pendubot:  $\ell_1$ is the length of the first
link, and $\ell_{c1},\ell_{c2}$ are the distances to the center of mass of the
respective links.  The variables $q_1, q_2$ are joint angles of the respective links.}
\flabel{mod-pendubot2}
\end{figure}

This model may be written in  state space form as a nonlinear vector differential
equation, $\dot x = f(x,u)$, where $x=(q_1,q_2,\dot q_1, \dot q_2)'$, and $f$ is
defined from the above equations. For a range of different constant torque inputs
$\tau$, this model admits various equilibria. For example, when $\tau =0$,   the
vertical downward position
$x_e=(-\pi/2,0,0,0)$  is an equilibrium, as illustrated in the right hand side of
\Figure{mod-pendubot}.  When
$\tau=0$ it follows from the equations of motion that the upright vertical position
$x_e=(+\pi/2,0,0,0)$  is  also an equilibrium.  It is clear from the photograph given in the
left hand side of \Figure{mod-pendubot}  that the upright equilibrium is strongly unstable in
the sense that with $\tau=0$, it is unlikely that the physical system will remain at rest.
Nevertheless, the velocity vector vanishes, $f(x_e,0) = 0$, so by definition the upright
position is an equilibrium when $\tau = 0$. Although complex, we may again linearize these
equations about the vertical equilibrium.  The control design techniques introduced
later in the book will provide the tools for stabilization of the pendubot in this
unstable vertical configuration via an appropriate controller.

\begin{figure}[ht]
\ebox{.35}{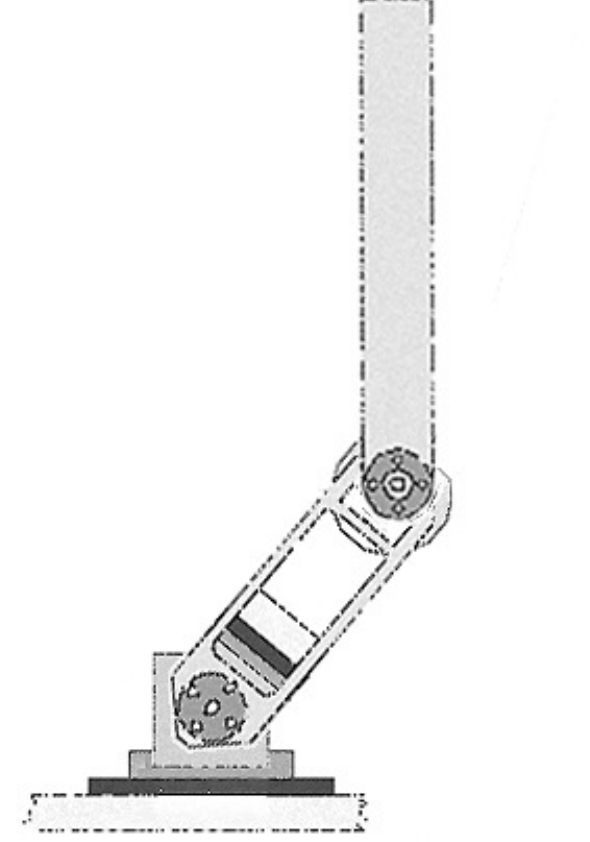}
\caption[Continuum of equilibrium positions for the Pendubot]{There is  a continuum of different equilibrium positions for the Pendubot
corresponding to different constant  torque inputs   $\tau$.}
\flabel{mod-pendubot3}
\end{figure}

\section{An electrical circuit}

A simple electrical circuit is illustrated in \Figure{mod-circuit1}.  Using Kirchoff's
voltage and current laws we may obtain a state space model in which the current through
the  inductor and the capacitor voltage become state variables.

\begin{figure}[ht]
\ebox{.6}{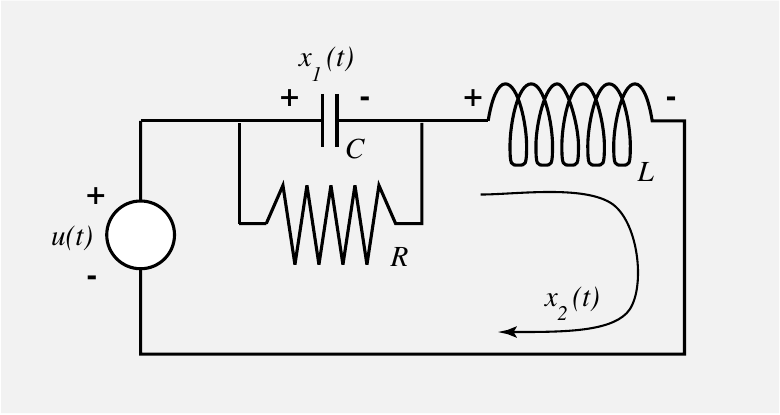}
\caption{A simple RLC circuit}
\flabel{mod-circuit1}

\ebox{.7}{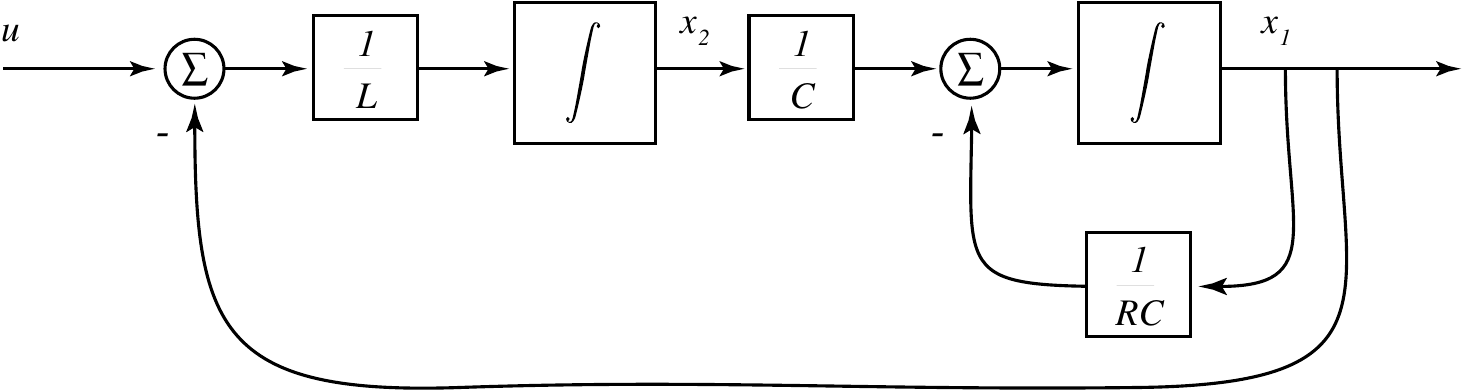}
\caption[Simulation diagram for the simple RLC circuit]{A simulation diagram for the corresponding state space model with $x_1=$
the voltage
across the capacitor, and $x_2=$ the current through the inductor.}
\flabel{mod-circuit2}
\end{figure}

Kirchoff's Current Law  (KCL) gives
\begin{eqnarray*}
x_2 = C \dot x_1 + \frac{ 1 }{ R } x_1,
\end{eqnarray*}
which may be written as
\begin{equation}
\dot x_1 = - \frac{1 }{ R  C} x_1 + \frac{1 }{ C} x_2
\elabel{eq1.3.1}
\end{equation}From Kirchoff's Voltage Law (KVL) we have
\begin{eqnarray*}
x_1+L \dot x_2=  + u,
\end{eqnarray*}
or,
\begin{equation}
  \dot x_2=  -\frac{1}{ L} x_1 + \frac{1}{ L}u.
\elabel{eq1.3.2}
\end{equation}Equations \eq eq1.3.2/ and \eq eq1.3.1/ thus give the system
of state space equations
\begin{eqnarray*}
\dot x = \left[
\begin{matrix} -\frac{1}{ RC} & \frac{1}{ C} \\
 -\frac{1}{ L} &   0 \end{matrix}
\right]x
+
\left[
\begin{matrix} 0 \\
  \frac{1}{ L}   \end{matrix}
\right] u
\end{eqnarray*}

\section{Transfer Functions \&\ State Space Models}
 \clabel{TransferFns}

In the previous examples we began with a physical description of the
system to be controlled, and through physical laws obtained a state
space model of the system. In the first and second examples, a linear
model could be obtained through linearization, while in the last
example, a linear model was directly obtained through the KVL and KCL
circuit laws. For a given LTI system, a state space model which is
constructed from a given transfer function $G(s)$ is called a
\defn{realization} of $G$.  Realization of a given $G$ is not unique,
and in fact there are infinitely many such realizations, all of which
are equivalent however from an input-output point of view.  In this
section we show how to obtain some of the most common realizations
starting from a transfer function model of a system.

Consider the third-order model
\begin{eqnarray}
\dddot y + a_2 \ddot y  + a_1 \dot y + a_0 y =  b_2 \ddot u  + b_1 \dot u +b_0 u
  \elabel{eq1.3.3}
\end{eqnarray}
where the coefficients $\{a_i, b_i\}$ are arbitrary real numbers.  The
corresponding transfer function description is
\begin{eqnarray*}
G (s) = \frac{Y (s) }{ U (s)} = \frac{B(s)}{A(s)}=\frac{b_2 s^2 + b_1 s +  b_0 }{ s^3 +
a_2 s^2 + a_1 s + a_0}
\end{eqnarray*}
where $Y$ and $U$ are the Laplace transforms of the signals $y$ and $u$, respectively.  The
numerator term $B(s) = b_2 s^2 + b_1 s +  b_0$ can create   complexity in determining a state
space model.  So, as a starting point, we consider the ``zero-free system'' where
$B(s)\equiv 1$, which results in the model
\begin{eqnarray*}
\dddot{ w} + a_2 \ddot w + a_1 \dot w + a_0 w = u
\end{eqnarray*}
An equivalent simulation diagram description for this system is:
            \ebox{.7}{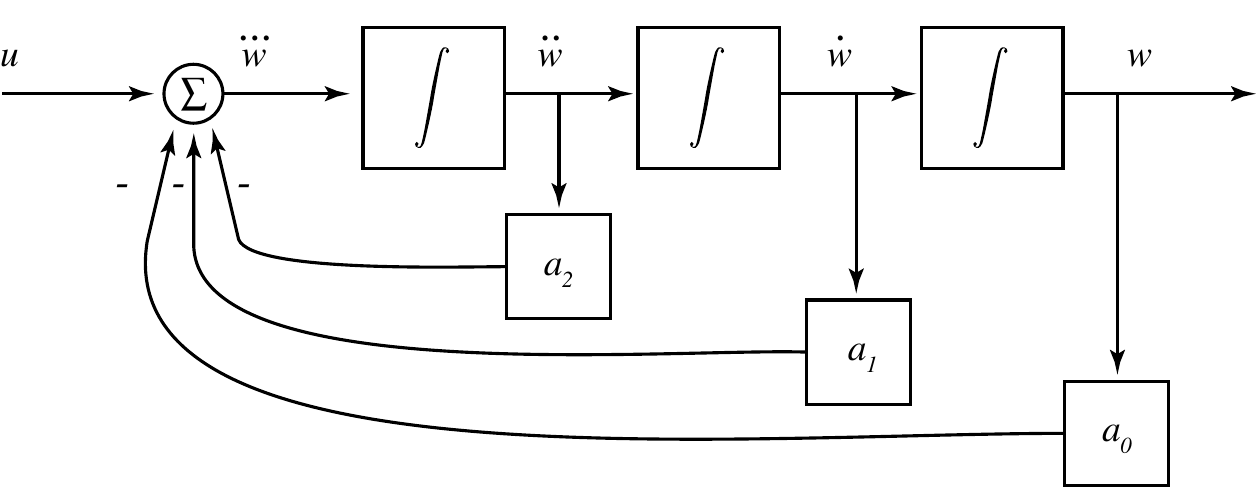}
With zero initial conditions one has the relation $Y(s) = B(s) W(s)$,
so that the signal $y$ can be obtained by adding several interconnections to the
above simulation diagram, yielding the simulation diagram depicted in \Figure{ccf}.
Letting the outputs of the integrators   form states for
the system we obtain the state space model
\begin{eqnarray*}
x_1 = w && \dot x_1 = x_2 \\
x_2 = \dot w && \dot x_2 = x_3 \\
x_3 = \ddot w && \dot x_3 = -a_2 x_3 - a_1 x_2 - a_0 x_1 + u
\end{eqnarray*}
and
\begin{eqnarray*}
y = b_0 x_1 + b_1 x_2 + b_2 x_3
\end{eqnarray*}
This may be written in matrix form as
\begin{eqnarray*}
\dot x & = & \left[ \begin{array}{ccc} 0 & 1 & 0 \\ 0 & 0 & 1 \\ -a_0 &
-a_1 & -a_2 \end{array} \right] x + \left[ \begin{array}{c} 0 \\ 0 \\
1 \end{array} \right]  u \\
y & = & [b_0\ \ b_1 \ \ b_2 \ ] x + [0] u.
\end{eqnarray*}
This final state space model is called the  \defn{controllable canonical form}
(CCF) - one of the most important system descriptions from the point of view of
analysis.

\begin{figure}

\ebox{.9}{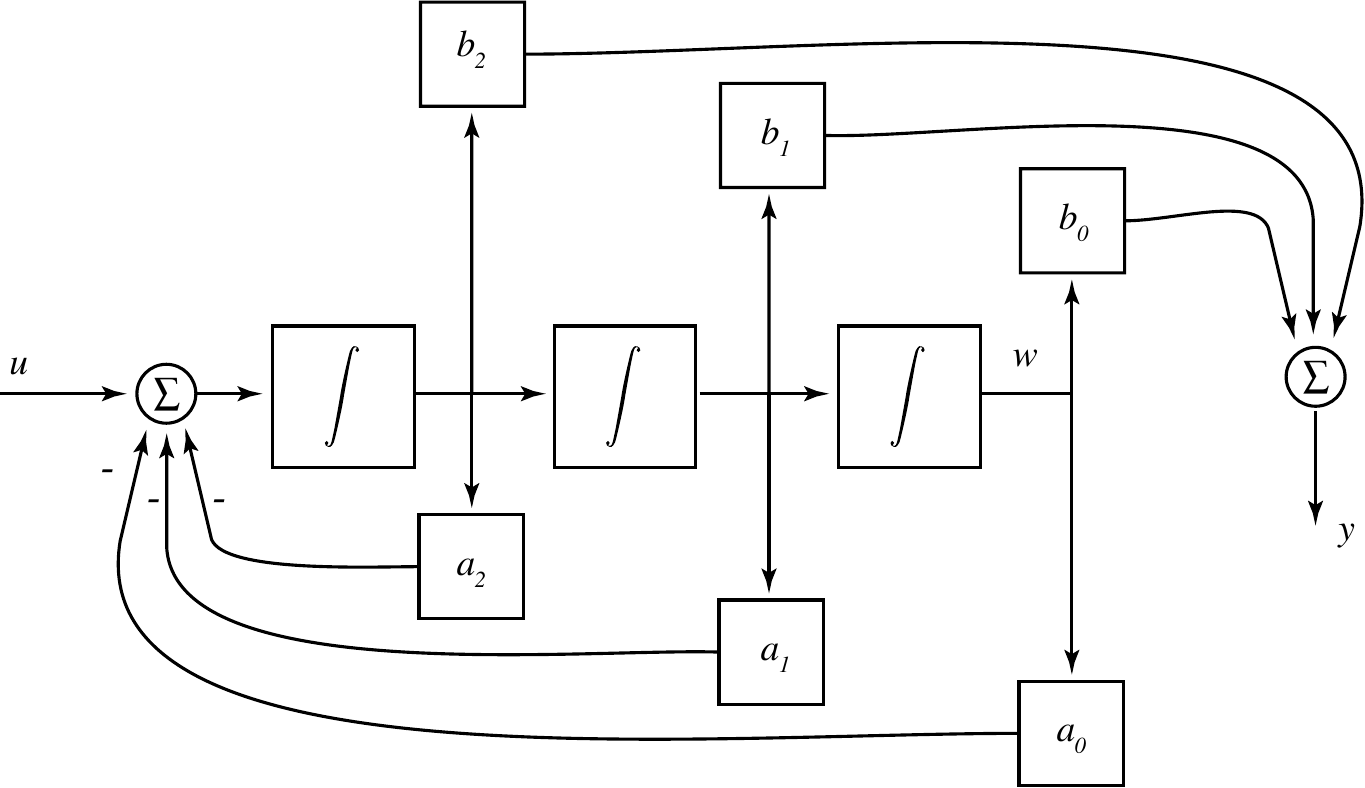}
\caption{Controllable Canonical Form}
\flabel{ccf}
\ebox{.9}{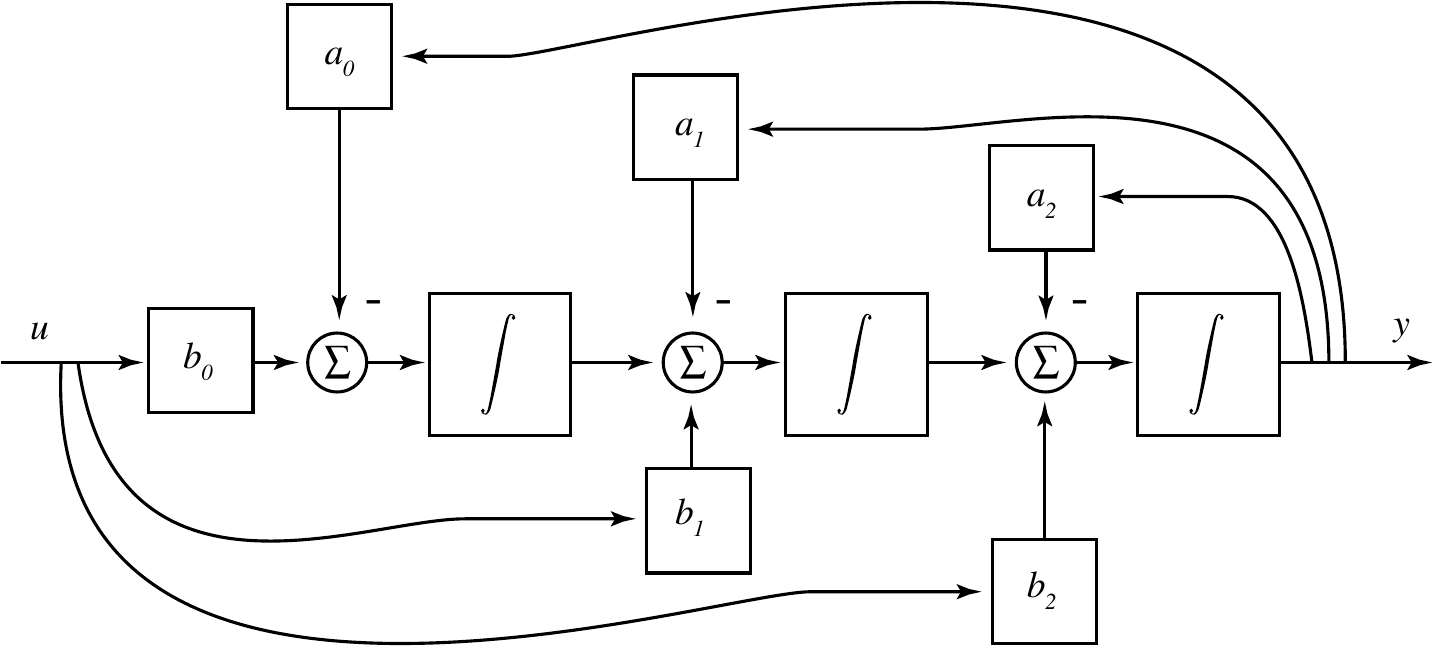}
\caption{Observable Canonical Form}
\flabel{ocf}
\end{figure}

Several alternative descriptions can be obtained by manipulating the
transfer function $G(s)$ before converting to the time domain, or by
defining states in different ways.  One possibility is to take the
description
\begin{eqnarray}
(s^3 + a_2 s^2 + a_1 s + a_0) Y (s) = (b_0 + b_1 s + b_2 s^2) U (s)
\elabel{eq1.4.4}
\end{eqnarray}
and divide throughout by $s^3$ to obtain
\begin{eqnarray*}
\left( 1 + \frac{a_2 }{ s} + \frac{a_1 }{ s^2} + \frac{a_0 }{ s^3} \right) Y
(s) = \left( \frac{b_0 }{ s^3} + \frac{b_1 }{ s^2} + \frac{b_2 }{ s} \right)
U (s)
\end{eqnarray*}
Rearranging terms then gives
\begin{eqnarray*}
Y & = & \frac{1 }{ s^3} (b_0 U - a_0 Y) + \frac{1 }{ s^2} (b_1 U - a_1 Y)
\\
&& + \frac{1 }{ s} (b_2 U - a_2 Y)
\end{eqnarray*}
We may again describe this equation using a simulation diagram, as given in \Figure{ocf}.
As before, by letting $x_1, x_2, x_3$ denote  the outputs of the integrators we
obtain a state space model which now takes the form
\begin{eqnarray*}
\dot x_1 & = & x_2 - a_2 x_1 + b_2 u \\
\dot x_2 & = & x_3 - a_1 x_1 + b_1 u \\
\dot x_3 & = & - a_0 x_1 + b_0 u\\
y & = & x_1
\end{eqnarray*}
or in matrix form
\begin{eqnarray*}
\dot x & = & \left[ \begin{array}{ccc} -a_2 & 1 & 0 \\ -a_1 & 0 & 1 \\
-a_0 & 0 & 0 \end{array} \right] x + \left[ \begin{array}{c} b_2 \\
b_1 \\ b_0 \end{array} \right] u \\ y & = & [1\ \ 0\ \ 0] x + [\ 0\ ]
u.
\end{eqnarray*}
This final form is called the \defn{observable canonical
form} (OCF).

In the example above, the degree $n_0$ of the denominator of $G(s)$ is $3$, and the
degree
$m_0$ of the numerator is $2$,  so that $n_0>m_0$. In this case the model is called
\defn{strictly proper}.  In the case where
$n_0 = m_0$, the ``$D$'' matrix in \eq stateSpace/ will be non-zero in
any state space realization.  To see this, try adding the term $b_3
s^3$ to the right hand side of \eq eq1.4.4/, or solve
\Exercise{not-strict} of this chapter.

Both controllable and observable canonical forms admit natural generalizations to the
$n$-dimensional case. For a SISO LTI system, let the input-output transfer function be
given by
\begin{eqnarray*}
G (s) = \frac{Y (s) }{ U (s)} = \frac{B(s)}{A(s)}=\frac{b_{n-1} s^{n-1} + \dots + b_1 s
+  b_0 }{ s^n + a_{n-1}s^{n-1} + \dots + a_2 s^2 + a_1 s + a_0}
\end{eqnarray*}
Then, the $n$-dimensional state space realization in \textit{controllable canonical
form} is identified by the following $A$, $B$, and $C$ matrices:
\begin{eqnarray*}
A =
\left[
\begin{array}{cccccc}
0 & 1 & 0 &  0 &\cdots & 0 \\
0 & 0 & 1 & 0 & \cdots & 0 \\
\vdots &\vdots  &\vdots && \ddots &  \\
-a_0 & -a_1 & -a_2 & a_3 &\cdots & a_{n-1} \end{array} \right]\, ;\quad
B= \left[ \begin{array}{c}
0 \\ \vdots\\ 0 \\ 1 \end{array} \right]\, ;\quad
C= \left[ \begin{array}{c}
b_0 \\ b_1 \\ \vdots\\ b_{n-1}  \end{array} \right]^T\, .
\end{eqnarray*}
The \textit{observable canonical form}, on the other hand, will have the
following $A$, $B$, and $C$ matrices:
\begin{eqnarray*}
A =
\left[
\begin{array}{cccccc}
-a_{n-1} & 1 & 0 &  0 &\cdots & 0 \\
-a_{n-2} & 0 & 1 & 0 & \cdots & 0 \\
\vdots &\vdots  &\vdots && \ddots &  \\
-a_0 & 0 & 0 & 0 &\cdots & 0 \end{array} \right]\, ;\quad
B= \left[ \begin{array}{c}
b_{n-1} \\ \vdots\\ b_1 \\ b_0 \end{array} \right]\, ;\quad
C= \left[ \begin{array}{c}
1 \\ 0 \\ \vdots\\ 0  \end{array} \right]^T\, .
\end{eqnarray*}

Another alternative is obtained by applying a partial fraction
expansion\index{Partial fraction expansion} to the transfer function
$G$:
\begin{eqnarray*}
G (s) & = & \frac {b(s)}{ a(s)} = d + \sum_{i=1}^n \frac{k_i }{
s-p_i},
\end{eqnarray*}
where $\{p_i : 1\le i\le n\}$ are the poles of $G$, which are simply
the roots of $a$.  A partial expansion of this form is always possible
if all of the poles are \textit{distinct}.  In general, a more complex
partial fraction expansion must be employed.  When a simple partial
fraction expansion is possible, as above, the system may be viewed as
a parallel network of simple first order simulation diagrams; see
\Figure{mod-parallel}.

The significance of this form is that it yields a strikingly simple
system description:
\begin{eqnarray*}
\left. \begin{array}{c} \dot x_1 = p_1 x_1 + k_1 u \\ \vdots \\ \dot
x_n = p_n x_n + k_n u
\end{array} \right\} \mbox{decoupled dynamic equations.}
\end{eqnarray*}
This gives the state space model
\begin{eqnarray*}
\dot x & = & \left[
\begin{array}{cccc}
p_1 & 0 & \cdots & 0 \\ 0 & p_2 & \cdots & 0 \\ \vdots & & \ddots & \\
0 & 0 & \cdots & p_n \end{array} \right] x + \left[ \begin{array}{c}
k_1 \\ k_2\\ \vdots\\ k_n \end{array} \right] u \\ y &=& [1,\dots,1] x
+ [d] u.
\end{eqnarray*}
This is often called the \defn{modal form}, and the states $x_i(t)$
are then called \defn{modes}.  It is important to note that this form
is not always possible if the roots of the denominator polynomial are
not distinct.  Exercises~\ref{hw-modal-repeat} and
\ref{hw-modal-complex} below address some generalizations of the modal
form.
\begin{figure}[bt]
\ebox{.9}{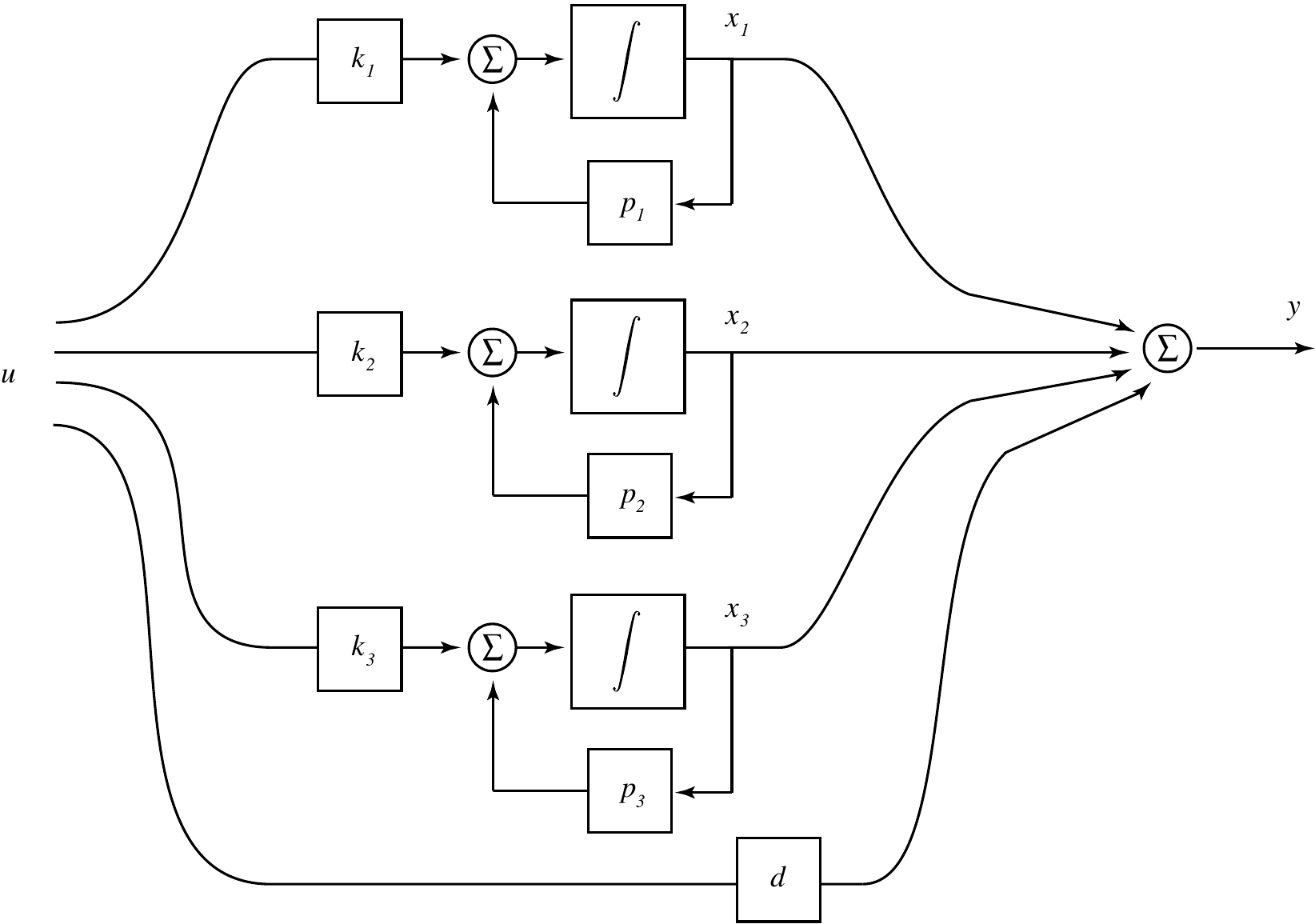}
\caption{Partial fraction expansion of a transfer function}
\flabel{mod-parallel}
\end{figure}

\begin{matlab}
\item
\textit{Matlab} is not well suited to nonlinear problems.  However,
the Matlab program \textit{Simulink} can be used for simulation of
both linear and nonlinear models.  Some useful Matlab commands for
system analysis of linear systems are

\item[RLOCUS] calculates the root locus. e.g.\ {\tt rlocus(num,den)},
or {\tt rlocus(A,B,C,D)}.

\item[STEP] computes the step response of a linear system.

\item[BODE] computes the Bode plot.

\item[NYQUIST] produces the Nyquist plot.

\item[TF2SS] gives the \textit{CCF} state space representation of a
transfer function model, but in a different form than given here.

\item[SS2TF] computes the transfer function of a model, given any
state space representation.

\item[RESIDUE] may be used to obtain the partial fraction expansion of
a transfer function.
\end{matlab}

\begin{summary}
State space models of the form
\begin{eqnarray*}
\dot x &=& f(x,u) \\ y &=& g(x,u)
\end{eqnarray*}
occur naturally in the mathematical description of many physical
systems, where $u$ is the input to the system, $y$ is the output, and
$x$ is the state variable.  In such a state space model, each of these
signals evolves in Euclidean space. The functions $f$ and $g$ are
often linear, so that the model takes the form
\begin{eqnarray*}
\dot x &=& A x + B u \\ y &=& C x + Du.
\end{eqnarray*}
If not, the functions $f,g$ may be approximated using a Taylor series
expansion which again yields a linear state space model.  The reader
is referred to \FPE\ for more details on modeling from a control
perspective.

A given linear model will have many different state space
representations.  Three methods for constructing a state space model
from its transfer function have been illustrated in this chapter:
\balphlist
\item
Construct a simulation diagram description, and define the outputs of
the integrators as states.  This approach was used to obtain the
controllable canonical form, but it is also more generally applicable.

\item
Manipulate the transfer function to obtain a description which is more
easily described in state space form.  For instance, a simple division
by $s^n$ led to the observable canonical form.

\item
Express the transfer function as a sum of simpler transfer functions
-- This approach was used to obtain a modal canonical form.
\end{list}
These three system descriptions, the modal, controllable, and
observable canonical forms, will be applied in control analysis in
later chapters.

Much more detail on the synthesis of state space models for linear
systems may be found in Chapter~6 of \CHE, and Chapter~3 of \BRO.
\end{summary}

\begin{exercises}
\item \hwlabel{simple-linear} You are given a nonlinear input-output
system which satisfies the nonlinear differential equation:
\begin{eqnarray*}
\ddot y(t) & = & 2y -(y^2+1)(\dot y +1) + u.
\end{eqnarray*}
\balphlist
\item Obtain a nonlinear state-space representation.

\item Linearize this system of equations around its equilibrium output
trajectory when $u(\cdot)\equiv 0$, and write it in state space form.
\end{list}

\item
Repeat \Exercise{simple-linear} with the new system
\begin{eqnarray*}
\ddot y(t)= 2y -(y^2+1)(\dot y +1) + u +2\dot u .
\end{eqnarray*}

\item Obtain state equations for the following circuit.  For the
states, use the voltage across the capacitor, and the current through
the inductor.  \ebox{.65}{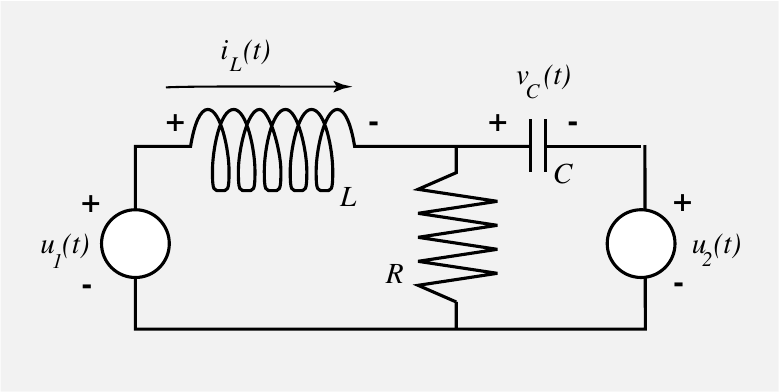}

\item In each circuit below, \balphlist
\item
Obtain a transfer function and a state space realization.

\item
Sketch a frequency response.

\item
Use the {\tt step} command in \textit{Matlab} to obtain a step
response.
\end{list}
\ebox{.65}{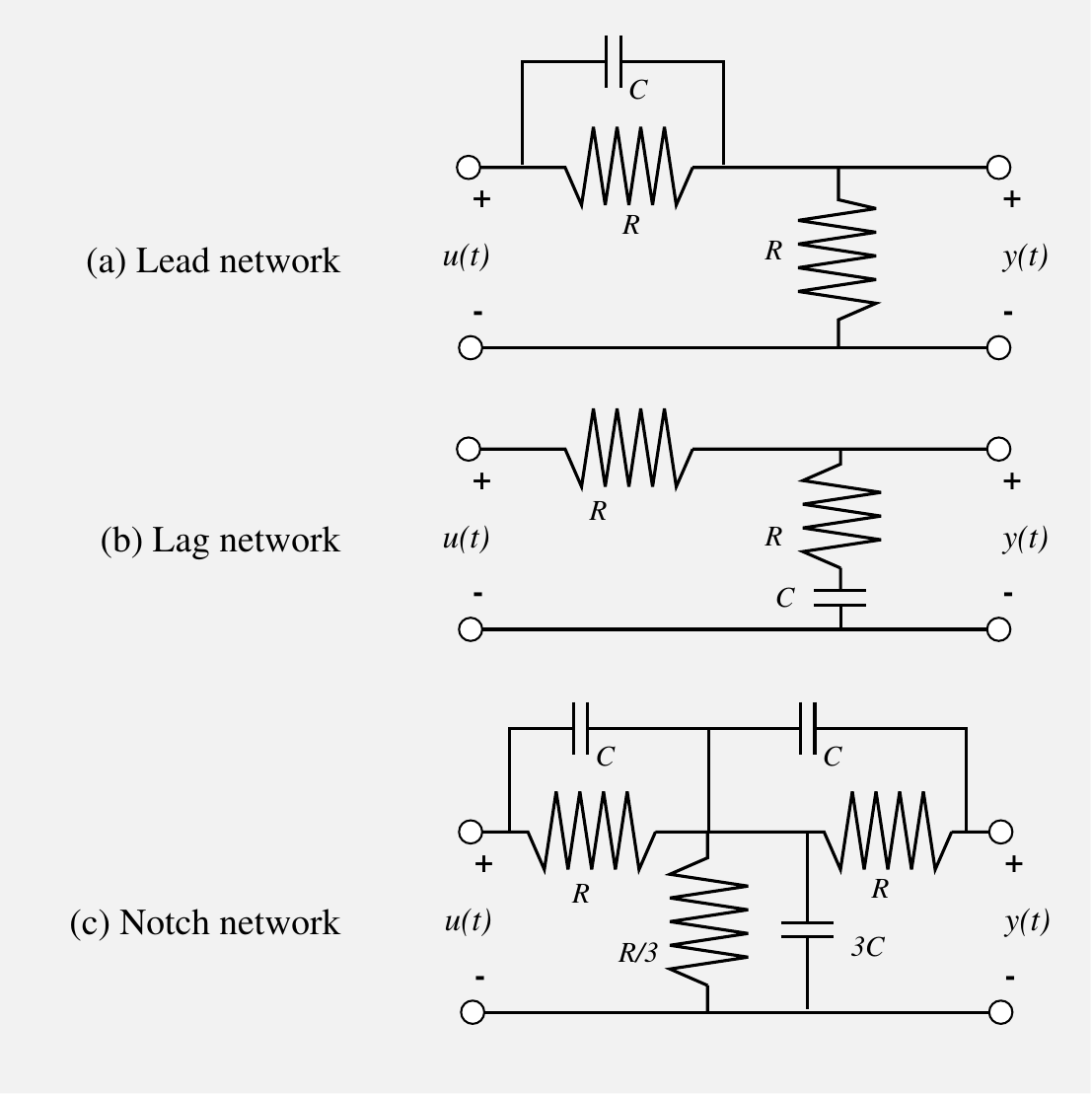}

\item
Consider the mass-spring system shown below \ebox{.65}{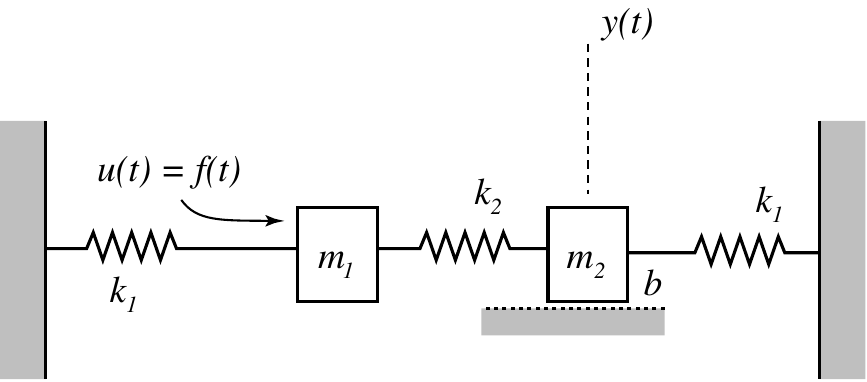}
Assume that a force is acting on $m_1$, and let the horizontal
position of $m_2$ represent the output of this system.

\balphlist
\item Derive a set of differential equations which describes this
input-output system.  To solve this problem you will require Newton's
law of translational motion, and the following facts: (i) The force
exerted by a spring is proportional to its displacement, and (ii) the
force exerted by a frictional source is proportional to the relative
speed of the source and mass.

\item Find the transfer function for the system.

\item Obtain a state space description of the system.
\end{list}

\item
The $n$-dimensional nonlinear vector differential equation $\dot
x=f(x)$ has a unique solution from any $x\in\Re^n$ if the function $f$
has continuous partial derivatives.  To see that just continuity of
$f$ is not sufficient for uniqueness, and that some additional
conditions are needed, consider the scalar differential equation
\begin{eqnarray*}
\dot x & = & \sqrt{1-x^2},\qquad x(0) = 1.
\end{eqnarray*}
Show that this differential equation with the given initial condition
has at least two solutions: One is $x(t)\equiv 1$, and another one is
$x(t) = \cos(t)$.

\item \hwlabel{satellite} Consider a satellite in {\it planar} orbit
about the earth.  The situation is modeled as a point mass $m$ in an
inverse square law force field, as sketched below.  The satellite is
capable of thrusting (using gas jets, for example) with a radial
thrust $u_1$ and a tangential ($\theta$ direction) thrust $u_2$.
Recalling that acceleration in polar coordinates has a radial
component $(\ddot r - r \dot\theta^2)$, and a tangential component
$(r \ddot\theta + 2\dot r \dot\theta)$, Newton's Law gives
\begin{eqnarray*}
m(\ddot r - r \dot \theta^2) = -\frac{k}{ r^2} + u_1\\ m(r\ddot\theta
+ 2\dot r \dot\theta) &=& u_2,
\end{eqnarray*}
where $k$ is a gravitational constant.

\balphlist \item Convert these equations to (nonlinear) state space
form using $x_1 = r$, $x_2 = \dot r$, $x_3 = \theta$, $x_4 = \dot
\theta$.

\item
Consider a nominal {\it circular} trajectory $r(t) = r_0$; $\theta(t)
= \omega_0 t$, where $r_0$ and $\omega_0$ are constants.  Using
$u_1(t) = u_2(t) = 0$, obtain expressions for the nominal state
variables corresponding to the circular trajectory.  How are $k$,
$r_0$, $m$, and $\omega_0$ related?

\item
Linearize the state equation in (a) about the state trajectory in (b).
Express the equations in matrix form in terms of $r_0$, $\omega_0$ and
$m$ (eliminate $k$).
\end{list}
\ebox{.4}{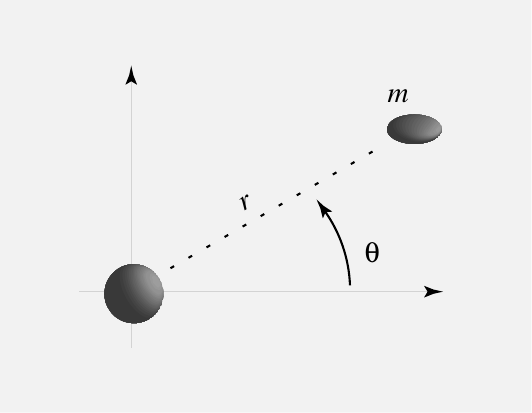}

\item
Using \textit{Matlab} or \textit{Simulink}, simulate the nonlinear
model for the magnetically suspended ball.  \balphlist
\item
Using proportional feedback $u = -k_1 y + k_2 r $, can you stabilize
the ball to a given reference height $r$?  Interpret your results by
examining a root locus diagram for the linearized system.

\item
Can you think of a better control law?  Experiment with other designs.
\end{list}

\vspace{.25in}\centerline{\textbf{\textit{Transfer Functions \&\ State
Space Models}}}

\item
A SISO LTI system is described by the transfer function
\begin{eqnarray*}
G(s) & = & \frac{s+4}{(s+1)(s+2)(s+3)}
\end{eqnarray*}

\balphlist
\item Obtain a state space representation in the controllable
canonical form;

\item Now obtain one in the observable canonical form;

\item Use partial fractions to obtain a representation with a diagonal
state matrix $A$ (modal form).
\end{list}
In each of (a)--(c) draw an all-integrator simulation diagram.

\item \hwlabel{not-strict} A SISO LTI system is described by the
transfer function
\begin{eqnarray*}
G(s) & = & \frac{s^3+2}{(s+1)(s+3)(s+4)}
\end{eqnarray*}
\balphlist

\item Obtain a state space representation in the controllable
canonical form;

\item Now obtain one in the observable canonical form;

\item Use partial fractions to obtain a representation of this model
with a diagonal state matrix $A$.
\end{list}
In each of (a)--(c) draw an all-integrator simulation diagram.

\item
For the multiple input-multiple output (MIMO) system described by the
pair of differential equations
\begin{eqnarray*}
\dddot y_1 +2\dot y_1 +3y_2&=&u_1+ \ddot u_1 + \dot u_2 \\ \ddot y_2
-3 \dot y_2 + \dot y_1 +y_2+y_1 &=&u_2+ \dot u_3 -u_3
\end{eqnarray*}
obtain a state space realization by choosing $y_1$ and $y_2$ as state
variables.  Draw the corresponding simulation diagram.

\item \hwlabel{modal-repeat} This exercise generalizes modal form to
the case where some eigenvalues are repeated.  For each of the
following transfer functions, obtain a state-space realization for the
corresponding LTI system by breaking $H$ into {\it simple} additive
terms, and drawing the corresponding simulation diagrams for the sum.
Choose the outputs of the integrators as state variables.  \balphlist
\item
$\displaystyle H(s)=\frac{2s^2}{s^3-s^2+s-1} $

\item
$\displaystyle H(s)=\frac{s^2+s+1} { s^3+4s^2+5s+2} $
\end{list}

\item \hwlabel{modal-complex} This exercise indicates that a useful
generalization of the modal form may be constructed when some
eigenvalues are complex.  \balphlist
\item
Use partial fractions to obtain a diagonal state space representation
for a SISO LTI system with transfer function
\begin{eqnarray*}
G(s) & = & \frac{s+6}{s^2+2s+2}.
\end{eqnarray*}
Note that complex gains appear in the corresponding all-integrator
diagram.

\item
Given a transfer function in the form
\begin{eqnarray*}
G(s) & = & \frac{s + \beta}{(s-\lambda_1)(s-\lambda_2)}
\end{eqnarray*}
and a corresponding state space realization with A matrix Compute the
eigenvalues $\lambda_1,\lambda_2$ of the matrix
\begin{eqnarray*}
A & = & \left(\begin{matrix}\sigma \quad \omega\\ -\omega\quad
\sigma\end{matrix}\right)
\end{eqnarray*}
where $\sigma\ge 0$, $\omega>0$, find the relationships between
$\lambda_1,\lambda_2, \beta$ and $\sigma , \omega$.  In view of this,
complete the state space realization by obtaining the $B$ and $C$
matrices, and draw the corresponding simulation diagram.

\item
For the transfer function $\displaystyle H(s)=\frac{s^2+s+1}
{s^3+4s^2+5s+2}$, obtain a state-space realization for the
corresponding LTI system by breaking $H$ into simple additive terms,
and drawing the corresponding simulation diagrams for the sum.

\item
Apply your answer in (b) to obtain a state space realization of $G(s)$
in (a) with only \textit{real} coefficients.
\end{list}
\end{exercises}

\chapter{Vector Spaces}
\clabel{chap1.5}

Vectors and matrices, and the spaces where they belong, are
fundamental to the analysis and synthesis of multivariable control
systems. The importance of the theory of vector spaces in fact goes
well beyond the subject of vectors in finite-dimensional spaces, such
as $\Re^n$.  Input and output signals may be viewed as vectors lying
in infinite dimensional function spaces.  A system is then a mapping
from one such vector to another, much like a matrix maps one vector in
a Euclidean space to another.  Although abstract, this point of view
greatly simplifies the analysis of state space models and the
synthesis of control laws, and is the basis of much of current optimal
control theory.  In this chapter we review the theory of vector spaces
and matrices, and extend this theory to the infinite-dimensional
setting.

\section{Fields}
\label{s:chap1.5.1}

A \defn{field} is any set of elements for which the operations of
addition, subtraction, multiplication, and division are defined.  It
is also assumed that the following axioms hold for any
$\alpha,\beta,\gamma\in\clF$
\balphlist
\item    $\alpha+\beta \in\clF$ and $\alpha\cdot\beta \in\clF$.

\item Addition and multiplication are commutative:
\begin{eqnarray*}
\alpha+\beta=\beta+\alpha,\qquad\alpha \cdot \beta=\beta\cdot\alpha.
\end{eqnarray*}

\item Addition and multiplication are associative:
\begin{eqnarray*}
(\alpha+\beta) + \gamma =\alpha+(\beta+\gamma), \qquad (\alpha \cdot
\beta)\cdot \gamma= \alpha \cdot (\beta \cdot \gamma).
\end{eqnarray*}

\item
Multiplication is distributive with respect to addition:
\begin{eqnarray*}
(\alpha+\beta) \cdot \gamma =\alpha\cdot\gamma + \beta\cdot\gamma
\end{eqnarray*}

\item There exists a unique null element $0$ such that $0 \cdot \alpha
= 0$ and $0+\alpha = \alpha$.

\item There exists a unique identity element $1$ such that $1\cdot
\alpha = \alpha$.

\item For every $\alpha \in\clF$ there exists a unique element $\beta
\in\clF$ such that $\alpha + \beta =0$; this unique element is
sometimes referred to as the \textit{additive inverse} or
\textit{negative} of $\alpha$, and is denoted as $-\alpha$.
\item
To every $\alpha \in\clF$ which is not the element $0$ (i.e., $\alpha
\not= 0$), there corresponds an element $\gamma$ such that
$\alpha\cdot \gamma = 1$; this element is referred to as the
\textit{multiplicative inverse} of $\alpha$, and is sometimes written
as $\alpha^{-1}$.
\end{list}

Fields are a generalization of $\Re$, the set of all real numbers.
The next  example is the set of all complex numbers,
denoted $\Co$.  These are the only examples of fields that will be
used in the text, although we will identify others in the exercises at
the end of the chapter.

\section{Vector Space}
\label{s:chap1.5.2}

A set of vectors $\clX$ is a set on which vector addition, and scalar
multiplication are defined, generalizing the concept of the vector
space $\Re^n$. In this abstract setting a \defn{vector space} over a
field $\clF$, denoted $(\clX,\clF)$, is defined as follows: \balphlist
\item
For every $x_1,x_2\in\clX$, the vector sum $x_1+x_2\in\clX$.

\item
Addition is commutative: For every $x_1,x_2\in\clX$, the sum
$x_1+x_2=x_2+x_1$.

\item
Addition is associative: For every $x_1,x_2,x_3\in\clX$,
\begin{eqnarray*}
(x_1+x_2) +x_3 & = & x_1+ (x_2 +x_3)
\end{eqnarray*}

\item
The set $\clX$ contains a vector $\zero$ such that $\zero + x=x$ for
all $x\in\clX$.

\item
For every $x\in\clX$, there is a vector $y\in\clX$ such that
$x+y=\zero$.

\item
For every $x\in\clX$, $\alpha\in\clF$, the scalar product $\alpha
\cdot x \in \clX$.

\item
Scalar multiplication is associative: for every $\alpha,\beta\in\clF$,
and $x\in\clX$,
\begin{eqnarray*}
\alpha(\beta x) & = & (\alpha \beta) x.
\end{eqnarray*}
\end{list}

Below is a list of some of the vector spaces which are most important
in applications
\begin{description}
\item
$(\Re^n,\Re)$ -- the real vector space of $n$ dimensional real-valued
vectors.

\item
$(\Co^n,\Co)$ -- the complex vector space of $n$ dimensional
complex-valued vectors.

\item
$(C^n[a,b],\Re)$ -- the vector space of real-valued continuous
functions on the interval $[a,b]$, taking values in $\Re^n$.

\item
$(D^n[a,b],\Re)$ -- the vector space of real-valued
piecewise-continuous functions on the interval $[a,b]$, taking values
in $\Re^n$.

\item
$(L_p^n [a,b],\Re)$ -- the vector space of functions on the interval
$[a,b]$, taking values in $\Re^n$, which satisfy the bound
\begin{eqnarray*}
\int_a^b | f(t) |^p\, dt < \infty, && f \in L_p[a,b].
\end{eqnarray*}

\item
$(\clR(\Co),\Re)$ -- the vector space of rational functions
$\frac{b(s)}{ a(s)}$ of a complex variable $s$, with real
coefficients.
\end{description}

A \defn{subspace} $\clY$ of a vector space $\clX$ is a subset of
$\clX$ which is itself a vector space with respect to the operations
of vector addition and scalar multiplication.  For example, the set of
complex $n$-dimensional vectors whose first component is zero is a
subspace of $(\Co^n,\Co)$, but $\Re^n$ is not a subspace of
$(\Co^n,\Co)$.

\section{Bases}
\label{s:chap1.5.4}

A set of vectors $S=(x^1, \ldots, x^n)$ in $(\clX, \clF)$ is said to
be \defn{linearly independent} if the following equality
\begin{eqnarray*}
\alpha_1 x^1 + \alpha_2 x^2 + \cdots + \alpha_n x^n & = & 0
\end{eqnarray*}
holds for a set of $n$ elements $\{ \alpha_i:1\le i \le
n\}\subset\clF$, then $\alpha_1 = \alpha_2 = \cdots \alpha_n = 0$. If
the set $S$ contains an infinite number of vectors, then we call $S$
linearly independent if every finite subset of $S$ is linearly
independent, as defined above.

In the case where $(\clX, \clF) = (\Re^n,\Re)$, we will have linear
independence of $\{x^1,\dots, x^n\} $ if and only if
\begin{eqnarray*}
\det [x^1 x^2 \cdots x^n] & = & \det \left[
\begin{array}{ccc} x_1^1 & \quad & x_1^n \\ \vdots & \cdots &  \vdots\\ x_n^1
&& x_n^n \end{array} \right] \neq 0.
\end{eqnarray*}

The maximum number of linearly independent vectors in $(\clX, \clF)$
is called the \defn{dimension} of $(\clX, \clF)$. For example,
$(\Re^n, \Re)$ and $(\Co^n, \Co)$ both have dimension $n$.
\textit{What is the dimension of $(\Co^n, \Re)$?} (see
\Exercise{dimCn}).

More interesting examples can be found in function spaces.  If for
example $(\clX, \clF)=(C[0,1],\Re)$, where $C[0,1]$ is the set of
real-valued continuous functions on $[0,1]$, then we can easily find a
set $S$ of infinite size which is linearly independent.  One such set
is the collection of simple polynomials $S= \{t, t^2, t^3, \ldots \}$.
To see that $S$ is linearly independent, note that for any $n$,
\begin{eqnarray*}
\sum_{i=1}^n \alpha_i t^i & = & \zero \qquad \mbox{ only if } \qquad
\alpha_i = 0, \ 1\le i\le n,
\end{eqnarray*}
where $\zero\in C[0,1]$ is the function which is identically zero on
$[0,1]$.  We have thus shown that the dimension of $(C[0,1],\Re)$ is
infinite.\medskip

A set of linearly independent vectors $S=\{e^1,\dots,e^n\}$ in $(\clX,
\clF)$ is said to be a \defn{basis} of $\clX$ if every vector in
$\clX$ can be expressed as a unique linear combination of these
vectors.  That is, for any $x \in \clX$, one can find $\{\beta_i, 1
\leq i \leq n\}$ such that
\begin{eqnarray*}
x & = & \beta_1 e^1 + \beta_2 e^2 + \cdots+ \beta_n e^n .
\end{eqnarray*}
Because the set $S$ is linearly independent, one can show that for any
vector $x$, the scalars $\{\beta_i\}$ are uniquely specified in
$\clF$.  The $n$-tuple $\{\beta_1, \ldots, \beta_n \}$ is often called
the \textit{representation} of $x$ with respect to the basis $\{e^1,
\dots ,e^n\}$.

We typically denote a vector $x \in \Re^n$ by
\begin{eqnarray*}
x & = & \left[ \begin{array}{c} x_1 \\ \vdots \\ x_n \end{array}
\right].
\end{eqnarray*}
There are two interpretations of this equation: \balphlist
\item
$x$ is a vector (in $\Re^n$), independent of basis.
\item
$x$ is a representation of a vector with respect to the \defn{natural
basis}:
\begin{eqnarray*}
x &=& x_1 \left[ \begin{array}{c} 1 \\ 0 \\ \vdots \\ 0 \\ 0
\end{array} \right] + x_2 \left[ \begin{array}{c} 0 \\ 1 \\ 0 \\
\vdots \\ 0
\end{array} \right] + \cdots + x_n \left[ \begin{array}{c} 0 \\ 0 \\ \vdots
\\ 0 \\ 1 \end{array} \right]  .
\end{eqnarray*}
\end{list}

The following theorem is easily proven in $\Re^n$ using matrix
manipulations, and the general proof is similar:
\begin{theorem}
In any $n$-dimensional vector space, {\it any} set of $n$ linearly
independent vectors qualifies as a basis.  \qed
\end{theorem}

In the case of Euclidean space $(\clX ,\clF) = (\Re^n, \Re)$, with
$\{e^1,\dots,e^n\}$ a given basis, any vector $x\in\Re^n$ may be
expressed as
\begin{eqnarray*}
x & = & \beta_1 e^1 + \cdots + \beta_n e^n,
\end{eqnarray*}
where $\{ \beta_i \}$ are all real scalars.  This expression may be
equivalently written as $ x = E \beta$, where
\begin{eqnarray*}
x & = & \left[ \begin{matrix}x_1 \\ \vdots \\ x_n \end{matrix}
\right], \quad E=\left[ \begin{matrix}e_{11} & \cdots & e_{1n}\\
\vdots & & \vdots\\ e_{n1} & \cdots & e_{nn}\end{matrix}\right], \quad
\beta = \left[ \begin{matrix}\beta_1 \\ \vdots \\ \beta_n \end{matrix}
\right],
\end{eqnarray*}
which through inversion shows that $\beta$ is uniquely given by
\begin{equation}
\beta = E^{-1} x \elabel{repFormula}
\end{equation}
Here $E^{-1}$ stands for the matrix inverse of $E$ (i.e., $E^{-1}E = E
E^{-1} = I$, where $I$ is the identity matrix), which exists since
$e^i$'s are linearly independent.

Consider the numerical example with $e^1=\left(\begin{smallmatrix}1 \\
1\end{smallmatrix}\right)$, $e^2 = \left(\begin{smallmatrix}0\\
1\end{smallmatrix}\right)$, and $x=\left(\begin{smallmatrix}2\\
5\end{smallmatrix}\right)$.  Then we have $x = \beta_1 e^1 + \beta_2
e^2$, with
\begin{eqnarray*}
\left(\begin{matrix}\beta_1\\\beta_2\end{matrix}\right) =
\left(\begin{matrix} 1 & 0\\ 1 & 1\end{matrix}\right)^{-1}
\left(\begin{matrix}2\\ 5\end{matrix}\right) =
\left(\begin{matrix} 1 & 0\\ -1 & 1\end{matrix}\right)
\left(\begin{matrix}2\\ 5\end{matrix}\right) =
\left(\begin{matrix}2\\ 3\end{matrix}\right)
\end{eqnarray*}

\section{Change of basis}
Suppose now we are given two sets of basis vectors:
\begin{eqnarray*} \{e^1 , \cdots, e^n \} \qquad
\{\bar{e}^1, \cdots , \bar e^n \}
\end{eqnarray*}
A vector $x$ in $\clX$ can be represented in two possible ways,
depending on which basis is chosen:
\begin{equation}
x= \beta_1 e^1 + \cdots \beta_n e^n = \sum_{k=1}^n \beta_k e^k
\elabel{eq1.5.1}
\end{equation}
or,
\begin{equation}
x = \bar\beta_1 \bar e^1 + \cdots \bar\beta_n \bar e^n = \sum_{k=1}^n
\bar\beta_k \bar e^k. \elabel{eq1.5.2}
\end{equation}
Since $\{e^i\} \subset \clX$, there exist scalars $\{ p_{ki} : 1\le
k,i \le n\}$ such that for any $i$,
\begin{eqnarray*}
e^i = p_{1i} \bar e^1 + \cdots + p_{ni} \bar e^n = \sum_{k=1}^n p_{ki}
\bar e^k.
\end{eqnarray*}
From \eq eq1.5.1/ it then follows that a vector $x$ may be represented
as
\begin{equation}
x = \sum_{i=1}^n \beta_i \left( \sum_{k=1}^n p_{ki} \bar e^k \right) =
\sum_{k=1}^n \left( \sum_{i=1}^n p_{ki} \beta_i \right)\bar e^k.
\elabel{eq1.5.1b}
\end{equation}
In view of \eq eq1.5.2/ and \eq eq1.5.1b/ we have by subtraction
\begin{eqnarray*}
\sum_{k=1}^n \left[ \sum_{i=1}^n p_{ki} \beta_i - \bar\beta_k \right]
\bar e^k = \zero.
\end{eqnarray*}
By linear independence of $\{ \bare^k\}$, this implies that each
coefficient in brackets is $0$. This gives a matrix relation between
the coefficients $\{\beta_i\}$ and $\{\bar\beta_i\}$:
\begin{eqnarray*}
\bar\beta_k & = & \sum_{i=1}^n p_{ki} \beta_i, \quad k = 1, \ldots, n
\end{eqnarray*}
or using compact matrix notation,
\begin{eqnarray*}
\bar\beta & = & \left[ \begin{array}{c} \bar\beta_1 \\ \vdots \\
\bar\beta_n \end{array} \right] = \left[ \begin{array}{ccc} p_{11} &
\ldots & p_{1n} \\ \vdots && \\ p_{n1} & \ldots & p_{nn} \end{array}
\right] \left[ \begin{array}{c} \beta_1 \\ \vdots \\ \beta_n
\end{array} \right] = P \beta.
\end{eqnarray*}
The transformation $P$ maps $\clF^n\to\clF^n$, and is one to one, and
onto.  It therefore has an inverse $P^{-1}$, so that $\beta$ can also
be computed through $\bar \beta$:
\begin{eqnarray*}
\beta & = & P^{-1} \bar\beta.
\end{eqnarray*}

For the special case where $(\clX,\clF) = (\Re^n,\Re)$, the vectors
$\{e_i\}$ can be stacked to form a matrix to obtain as in \eq
repFormula/,
\begin{eqnarray*}
x& = & E \beta = \bar E \bar\beta.
\end{eqnarray*}
Hence the transformation $P$ can be computed explicitly: $\bar\beta =
\bar E^{-1} E \beta$, so that $ P=\bar E^{-1} E$. The inverse $\bar
E^{-1}$ again exists by linear independence.

\section{Linear Operators}
\label{s:chap1.5.7}

A \defn{linear operator} $\clA$ is simply a function from one vector
space $(\clX,\clF)$ to another $(\clY,\clF)$, which is \defn{linear}.
This means that for any scalars $\alpha_1,\alpha_2$, and any vectors
$x_1, x_2$,
\begin{eqnarray*}
\clA(\alpha_1 x_1 + \alpha_2 x_2) & = & \alpha_1\clA(x_1) + \alpha_2
\clA(x_2).
\end{eqnarray*}
For a linear operator $\clA$ or a general function from a set $\clX$
into a set $\clY$ we adopt the terminology
\begin{eqnarray*}
\clX & : & \mbox{ Domain of the mapping $\clA$} \\ \clY & : & \mbox{
Co-Domain of the mapping $\clA$}
\end{eqnarray*}

When $\clA$ is applied to every $x \in \clX$, the resulting set of
vectors in $\clY$ is called the \defn{range} (or \defn{image}) of
$\clA$, and is denoted by $\clR (\clA)$:
\begin{eqnarray*}
\clR (\clA) & \eqdef & \bigcup_{x\in\clX} \clA(x)
\end{eqnarray*}
Pictorially, these notions are illustrated as follows:
\ebox{.65}{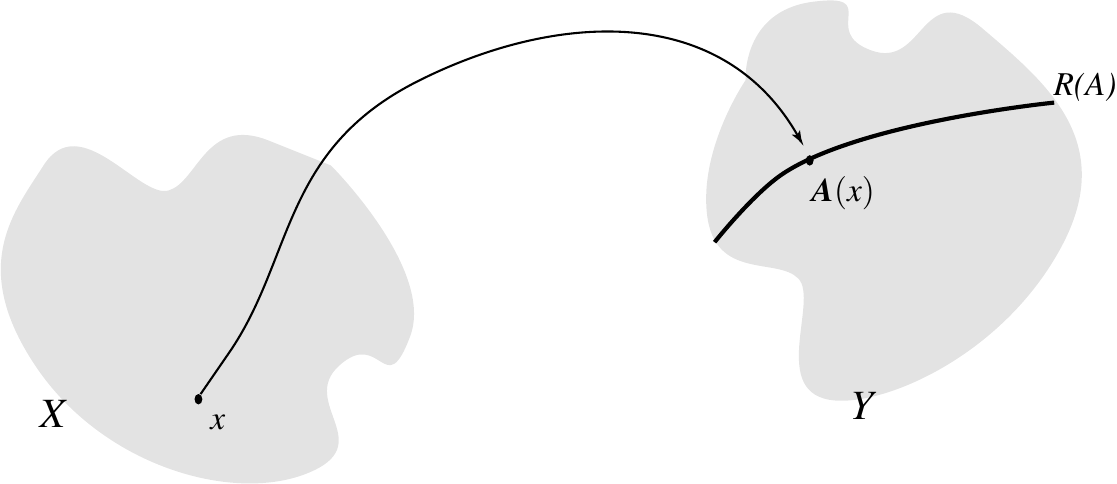}

The \defn{rank} of $\clA$ is defined to be the dimension of
$\clR(\clA)$.

In the special case of a linear operator $\clA : \Re^n \to \Re^m$
defined as $\clA(x) = A x$ for an $m\times n$ matrix $A$,
\begin{eqnarray*}
y & = & A x = [a^1 \ldots a^n] \left[ \begin{array}{c} x_1 \\ \vdots
\\ x_n \end{array} \right] \\ & = & x_1 a^1 + \cdots + x_n a^n,
\end{eqnarray*}
the range of $A$ is the set of all possible linear combinations of the
columns $\{a_i\}$ of $A$.  That is, the space spanned by the columns
of $A$.  The dimension of $\clR (A)$ is then the maximum number of
linearly independent columns of $A$.

\begin{theorem}
For a linear operator $\clA$, the set $ \clR (\clA)$ is a subspace of
$\clY$.
\end{theorem}

\proof To prove the theorem it is enough to check closure under
addition and scalar multiplication.  Suppose that $y^1,
y^2\in\clR(\clA)$, and that $\alpha_1,\alpha_2\in\clF$.  Then by
definition of $\clR(\clA)$ there are vectors $x_1,x_2$ such that
\begin{eqnarray*}
\clA (x^1) = y^1, && \clA (x^2) = y^2,
\end{eqnarray*}
and then by linearity of the mapping $\clA$,
\begin{eqnarray*}
\clA (\alpha_1 x^1 + \alpha_2x^2) = \alpha_1\clA (x^1) + \alpha_2\clA
(x^2) =\alpha_1 y^1 + \alpha_2 y^2.
\end{eqnarray*}
Hence, $\alpha_1 y^1 + \alpha_2 y^2\in\clR(\clA)$, which establishes
the desired closure property.  \qed

By assumption, the set $\clY$ contains a zero vector, which could be
mapped from numerous $x \in \clX$.  The set of all such $x \in \clX$
is called the \defn{nullspace} of $\clA$, denoted $\clN(\clA)$:
\begin{eqnarray*}
\clN (\clA) & \eqdef & \left\{ x \in \clX \mbox{ such that } \clA (x)
= 0 \right\}.
\end{eqnarray*}

\begin{theorem}
For any linear operator $\clA \colon \clX \to \clY$, the nullspace
$\clN (\clA)$ is a subspace of $\clX$.  \end{theorem}

\proof Again we check that $\clN(\clA)$ is closed under addition and
scalar multiplication.  Suppose that $x^1, x^2\in\clN(\clA)$, and that
$\alpha_1,\alpha_2\in\clF$.  Then by definition,
\begin{eqnarray*}
\clA(x^1) & = & \clA(x^2) =\zero.
\end{eqnarray*}
Again by linearity of $\clA$ it is clear that $\clA(\alpha_1 x^1 +
\alpha_2 x^2) = \zero$, which proves the theorem.  \qed

\section{Linear operators and matrices}

Suppose that $\clV$ and $\clW$ are two finite-dimensional vector
spaces over the field $\clF$, with bases $\{ v^1, \ldots, v^n \}$ and
$\{w^1, \ldots, w^m\}$ respectively.  If $ \clA : \clV \to \clW $,
then $\clA$ may be represented by a matrix.  To see this, take any
$v\in\clV$, and write
\begin{eqnarray*}
v & = & \sum_{i=1}^n \alpha_i v^i
\end{eqnarray*}
for scalars $\{\alpha_i\}\subset\clF$.  By linearity we then have,
\begin{eqnarray*}
\clA (v) & = & \clA \Bigl( \sum_{i=1}^n \alpha_i v^i \Bigr) =
\sum_{i=1}^n \alpha_i \clA (v_i)
\end{eqnarray*}
But for any $i$ we have that $\clA (v^i) \in \clW$, which implies that
for some scalars $\{a_{ji}\}$,
\begin{eqnarray*}
\clA (v^i) & = & \sum_{j=1}^m a_{ji} w^j,\qquad 1\le i\le n.
\end{eqnarray*}
From the form of $v$ we must therefore have
\begin{eqnarray}
\clA (v) & = & \sum_{i=1}^n \alpha_i \sum_{j=1}^m a_{ji} w^j \nonumber
\\ & = & \sum_{j=1}^n \Bigl( \sum_{i=1}^m a_{ji} \alpha_i \Bigr) w^j
\elabel{eqref}
\end{eqnarray}
Recall that the vector $w=\clA(v)$ in $\clW$ has a \textit{unique}
representation
\begin{eqnarray*}
\clA (v) & = & \sum_{j=1}^m \beta_j w^j
\end{eqnarray*}
Consequently, the terms in parentheses in \eq eqref/ are identical to
the $\{\beta_j\}$:
\begin{eqnarray*}
\beta_j & = & \sum_{i=1}^n a_{ji} \alpha_i \quad j = 1, \ldots, m.
\end{eqnarray*}
From this we see how the representations of $v $ and $ w$ are
transformed through the linear operator $\clA$:
\begin{eqnarray*}
\beta & = & \left[ \begin{array}{c} \beta_1 \\ \vdots \\ \beta_m
\end{array}\right] = \left[\begin{array}{ccc} a_{11} & \ldots & a_{in}
\\ \vdots && \vdots \\ a_{mi} & \ldots & a_{mn}
\end{array} \right]  \left[\begin{array}{c} \alpha_1 \\ \vdots \\
\alpha_n \end{array} \right] = A \alpha
\end{eqnarray*}
so that $A$ is the representation of $\clA$ with respect to $\{v^i\}$
and $\{w^i\}$.

The special case where $\clA$ is a mapping of $(\clX,\clF)$ into
itself is of particular interest.  A question that frequently arises
is ``if the linear operator is represented by a matrix $A$ of elements
in $\clF$, how is $A$ affected by a change of basis?''. Let $x^b= \clA
(x^a)$, and write
\begin{eqnarray*}
x^a &=& \sum_{i=1}^n \alpha_i e^i = \sum_{i=1}^n \bar \alpha_i \bar
e^i \\ x^b&=& \sum_{i=1}^m \beta_i e^i = \sum_{i=1}^m \bar\beta_i \bar
e^i
\end{eqnarray*}
where the $\alpha$ and $\beta$ are related by
\begin{eqnarray*}
\beta & = & A \alpha \qquad \bar\beta =\bar A \bar\alpha.
\end{eqnarray*}
To see how $A$ and $\bar A$ are related, recall that there is a matrix
$P$ such that
\begin{eqnarray*}
\bar \alpha = P \alpha \qquad \bar\beta= P \beta
\end{eqnarray*}
Combining these four equations gives
\begin{eqnarray*}
P A \alpha = P\beta =\bar \beta = \bar A\bar\alpha = \bar A P \alpha.
\end{eqnarray*}
Since $\alpha\in\clX$ is arbitrary, we conclude that $PA=\bar A P$,
and hence we can also conclude that
\begin{eqnarray*}
\bar A & = & PAP^{-1} \\ A & = & P^{-1} \bar AP
\end{eqnarray*}
When these relationships hold, we say that the matrices $A$ and $\bar
A$ are \defn{similar}.

\section{Eigenvalues and Eigenvectors}
\label{s:chap1.5.9}

For a linear operator $\clA\colon\clX\to\clX$ on an arbitrary vector
space $(\clX,\clF)$, a scalar $\lambda$ is called an \defn{eigenvalue}
of $\clA$ if there exists a non-zero vector $x$ for which
\begin{equation}
\clA(x) = \lambda x. \elabel{e-value}
\end{equation}
The vector $x$ in \eq e-value/ is then called an \defn{eigenvector}.

Let $\clX=\Co^n$, $\clF=\Co$, and $A$ be a matrix representation for
$\clA$.
If an eigenvalue $\lambda$ of $\clA$ does
exist, then one may infer from the equation
\begin{eqnarray*}
(A - \lambda I) x & = & 0,
\end{eqnarray*}
that the matrix $ A - \lambda I$ is singular.  For nontrivial
solutions, we must then have
\begin{equation}
\Delta(\lambda) \eqdef \det (\lambda I -A)= 0. \elabel{charPoly}
\end{equation}
The function $\Delta(\varble)$ is called the \defn{characteristic
polynomial} of the matrix $A$, and \eq charPoly/ is known as the
\defn{characteristic equation}.  The characteristic polynomial is a
polynomial of degree $n$, which must therefore have $n$ roots. Any
root of $\Delta$ is an eigenvalue, so at least in the case of
operators on $(\Co^n,\Co)$, eigenvalues always exist.
Note that if $\bar{A}$ is some other matrix representation for $\clA$,
since $A$ and $\bar{A}$ are necessarily {\it similar}, $\bar{A}$ has
the same characteristic polynomial as $A$. Hence, the eigenvalues do not
depend on the specific representation picked.

If the roots of the characteristic polynomial are distinct, then the
vector space $\Co^n$ admits a basis consisting entirely of
eigenvectors:
\begin{theorem}
\tlabel{distinct-evalues} Suppose that $\lambda_1, \ldots, \lambda_n$
are the \textit{distinct} eigenvalues of the $n \times n$ matrix $A$,
and let $v^1, \ldots ,v^n$ be the associated eigenvectors.  Then the
set $\{v^i, i = 1 \ldots n\}$ is linearly independent over $\Co$.
\end{theorem}

When the eigenvalues of $A$ are distinct, the modal matrix defined as
\begin{eqnarray*}
M\eqdef [v^1 \ldots v^n]
\end{eqnarray*}
is nonsingular.  It satisfies the equation $AM = M \Lambda$, where
\begin{eqnarray*}
\Lambda = \left[\begin{matrix} \lambda_1 && 0 \\ & \ddots & \\ 0 &&
\lambda_n
\end{matrix}\right]
\end{eqnarray*}
and therefore $A$ is similar to the diagonal matrix $\Lambda$:
\begin{equation}
\Lambda = M^{-1} AM. \elabel{modal-prop}
\end{equation}
Unfortunately, not every matrix can be diagonalized, and hence the
eigenvectors do not span $\Co^n$ in general. Consider the matrix
\begin{eqnarray*}
A & = & \left[ \begin{array}{rrr} 1 & 1 & 2 \\ 0 & 1 & 3 \\ 0 & 0 & 2
\end{array} \right]
\end{eqnarray*}
The characteristic polynomial for $A$ is
\begin{eqnarray*}
\Delta(\lambda) & = & \det (\lambda I -A) = (\lambda - 1) (\lambda - 1)
(\lambda - 2)
\end{eqnarray*}
So, eigenvalues are $\lambda_1 = 1$, $\lambda_2 = 1$, $\lambda_3 = 2$.

Considering the eigenvector equation
\begin{eqnarray*}
(A - \lambda_1 I) x & = & \left[ \begin{array}{rrr} 0 & 1 & 2 \\ 0 & 0
& 3 \\ 0 & 0 & 1
\end{array} \right] x = \zero,
\end{eqnarray*}
we see that $x = \left( \begin{array}{r} 1 \\ 0 \\ 0 \end{array}
\right)$ and its constant multiples are the only eigenvectors
associated with $\lambda_1$. It then follows that there does not exist
a state transformation $P$ for which
\begin{eqnarray*}
\Lambda & = & \left( \begin{array}{rrr} 1 & 0 & 0 \\ 0 & 1 & 0 \\ 0 &
0 & 2 \end{array} \right) = PAP^{-1}.
\end{eqnarray*}

To obtain a nearly diagonal representation of $A$ we use
\defn{generalized eigenvectors}, defined as follows.  Search for a
solution to
\begin{eqnarray*}
(A - \lambda_1 I)^k x & = & 0
\end{eqnarray*}
such that $(A - \lambda_1 I)^{k-1} x \neq 0$.  Then $x$ is called a
generalized eigenvector of grade $k$.

Note that it then follows that the vector $(A - \lambda_1 I) x$ is a
generalized eigenvector of grade $k-1$.  In the example,
\begin{eqnarray*}
(A - \lambda_1 I)^2 x & = & \left[ \begin{array}{rrr} 0 & 1 & 2 \\ 0 &
0 & 3 \\ 0 & 0 & 1 \end{array} \right]
\left[ \begin{array}{rrr} 0 & 1 & 2 \\ 0 &
0 & 3 \\ 0 & 0 & 1 \end{array} \right] x \\
& = & \left[ \begin{array}{rrr} 0 & 0 & 5 \\ 0 &
0 & 3 \\ 0 & 0 & 1 \end{array} \right] x.
\end{eqnarray*}
So $x = \left(  \begin{array}{c} 0 \\ 1 \\ 0 \end{array} \right)$ is a
generalized eigenvector of grade $k$.

Letting $y = (A - \lambda_1 I) x$, we have
\begin{eqnarray*}
(A - \lambda_1 I) y & = & (A - \lambda_1 I)^2 x = 0.
\end{eqnarray*}
So $y$ is an eigenvector to $A$, and $x,y$ are linearly independent.  In fact, $y = (1,0,0)'$ is
the eigenvector computed earlier.  To obtain an approximately diagonal form, let $x_1 = y$, $x_2
=x$,  $x_3$ any eigenvector corresponding to
$\lambda_3 = 2$.  We then have
\begin{eqnarray*}
A x_1 &=& \lambda_1 x_1 \\
A x_2 & = & Ax  = y + \lambda_1 x = x_1 + \lambda_1 x_2  \\
A x_3 & = & \lambda_3 x_3
\end{eqnarray*}
Letting $M=[x_1|x_2|x_3]$ it follows that
\begin{eqnarray*}
A M = M \left[ \begin{array}{ccc} \lambda_1 & 1 & 0 \\ 0 &
\lambda_1 & 0 \\ 0 & 0 & \lambda_3 \end{array} \right] & = &
\left[
M \left(
\begin{array}{c} \lambda_1 \\ 0 \\ 0 \end{array}
\right) \left|
M \left(
\begin{array}{c} 1 \\ \lambda_1 \\ 0 \end{array}
\right) \right|
 M \left(
\begin{array}{c} 0 \\ 0 \\ \lambda_3 \end{array} \right)   \right] \\
& = & M J
\end{eqnarray*}
where
\begin{eqnarray*}
J = \left[ \begin{array}{ccc} \lambda_1 & 1 & 0 \\ 0 & \lambda_1 & 0
\\ 0 & 0 & \lambda_3 \end{array} \right] & = & M^{-1} A M
\end{eqnarray*}
This representation of $A$ with respect to a basis of generalized
eigenvectors is known as the \defn{Jordan form}.

\section{Inner Products}
\label{s:chap1.5.8}

Inner products are frequently applied in the solution of optimization problems
because they  give a natural measure of distance between vectors. This abstract notion
of distance can often be interpreted as a \textit{cost} in applications to finance, or
as \textit{energy} in mechanical systems. Inner products also provide a way of defining
angles between vectors, and thereby introduce geometry to  even infinite dimensional
models where any geometric structure is far from obvious at first glance.

To define an inner product we restrict our attention to a vector space $\clX$ over the
complex field $\Co$.
An inner product is then a complex-valued function of two vectors, denoted
$\inp \varble,\varble>$, such that the following three properties hold:
\balphlist
\item
$\inp x,y> = \overline{\inp y,x>}$ (complex conjugate).
\item
$\inp x, \alpha_1 y^1 + \alpha_2 y^2> = \alpha_1 \inp  x, y^1> + \alpha_2 \inp  x,
y^2>$,\qquad   for all $x,y^1,y^2\in\clX$,  $\alpha_1,\alpha_2\in\Co$.
\item
$\inp x, x> \geq 0$ for all $x\in\clX$,
 and $\inp x,x> = 0$ if and only if $x = 0$.
\end{list}

In the special case where $\clX=\Co^n$ we typically define
\begin{eqnarray*}
\inp  x,y> = x^* y
\end{eqnarray*}
where $x^*$ denotes the complex conjugate transpose of the vector $x$.  Another important
example is the function space  $L_p[a,b]$ with $p=2$. It can be shown that the formula
\begin{eqnarray*}
\langle f, g\rangle \eqdef \int_a^b f(t) g(t)\, dt,\qquad f,g\in L_2[a,b],
\end{eqnarray*}
defines on inner product on $L_2[a,b]$.

The most obvious application of the inner product is in the formulation of a
\defn{norm}. In general, the norm of a vector $x$ in a vector space $(\clX,\Co)$, denoted
$\| x\|$,  is a real-valued function of a vector $x$ such that
\begin{description}
\item[1.] $\| x \| \geq 0$, and $\| x \| = 0$ if and only if $x = 0$.

\item[2.] $\| \alpha x \| = | \alpha | \ \| x \|$, for any $\alpha\in\Co$.

\item[3.] $\| x + y \| \leq \| x \| + \| y \|$.
\end{description}
The third defining property is known as the \textit{triangle inequality}.

In  $\Re^n$ we usually define the norm as the usual Euclidean norm,
\begin{eqnarray*}
\| x \| = \sqrt{x^T x} = \sqrt{\sum_{i=1}^n x_i^2} \,,
\end{eqnarray*}
which we will henceforth write as $|x|$, and reserve the notation $\|\cdot \|$ for
norm of an infinite-dimensional vector.
This Euclidean norm can also be defined using the inner product:
\begin{equation}
| x | = \sqrt{\inp x,x>}
\elabel{inner-norm}
\end{equation}
In fact, one can show that  the expression \eq inner-norm/ defines a norm
in an arbitrary (finite- or infinite-dimensional) inner product space.

We define the norm of a vector $f\in L_p[a,b]$ as
\begin{eqnarray*}
\| f\|_{L_p} \eqdef \Bigl( \int_a^b |f(t)|^p \, dt\Bigr)^{1/p}.
\end{eqnarray*}
In the case $p=2$, this norm is derived from the inner product on $L_2[a,b]$, but for
general $p$ this norm is not consistent with \textit{any} inner product.

\section{Orthogonal vectors and reciprocal  basis vectors}
 \label{s:chap1.5.8.1}

Two vectors $x,y$ in an inner product space $(\clX,\Co)$ are said to be
\defn{orthogonal} if $\inp  x,y> = 0$.   This concept has many applications in optimization,
and orthogonality is also valuable in computing representations of vectors.   To see the
latter point, write
\begin{eqnarray*}
x = \sum_{i=1}^n \alpha_i v^i ,
\end{eqnarray*}
where $\{v^i, i = 1 \ldots n\}$ is a  basis for $(\clX,\Co)$.  By
orthogonality we then have
\begin{eqnarray*}
\inp  v^j, x> = \inp  v^j, \sum_{i=1}^n \alpha_i v^i> = \sum_{i=1}^n
\alpha_i \inp  v^j, v^i> \quad j = 1 \ldots n
\end{eqnarray*}
This may be written explicitly as
\begin{eqnarray*}
\inp  v^1, x> &=& \alpha_1 \inp   v^1, v^1> + \alpha_2 \inp   v^1 , v^2> +
\cdots  + \alpha_n \inp  v^1, v^n>
\\
&\vdots&
\\
\inp  v^n, x> &=&  \alpha_1 \inp  v^n, v^1> + \alpha_2 \inp v^1, v^2> +
\cdots + \alpha_n \inp v^n, v^n>
\end{eqnarray*}
or in matrix form
\begin{eqnarray*}
\left[ \begin{array}{c} \inp v^1, x> \\ \vdots \\ \inp v^n, x>
\end{array} \right] & = & \underbrace{  \left[ \begin{array}{ccc} \inp v^1,v^1>
& \cdots & \inp v^1,v^n> \\ \vdots && \\
\inp v^n,v^1> &\cdots & \inp v^n, v^n> \end{array} \right] }_{G}
\left[
\begin{array}{c} \alpha_1 \\ \vdots \\ \alpha_n \end{array} \right]
\end{eqnarray*}
The $n\times n$ matrix $G$ is called the \defn{Grammian}.  Its inverse gives a formula
for the representation $\alpha$:
\begin{eqnarray*}
\alpha = G^{-1} \left[ \begin{array}{c} \inp v^1, x> \\ \vdots \\ \inp v^n, x>
\end{array} \right]
\end{eqnarray*}
If the basis is orthogonal, then $G$ is diagonal, in which case the computation of the
inverse $G^{-1}$ is straightforward.

A basis is said to be \defn{orthonormal} if
\begin{eqnarray*}
\inp  v^j , v^i> = \delta_{ij} = \left\{ \begin{array}{ll} 1, & i = j \\
0, & i \neq j \end{array} \right.
\end{eqnarray*}
In this case $G=I$ (the identity matrix), so that $ G^{-1} = G$.

A basis $\{r^i\}$ is said to be \defn{reciprocal} to (or \defn{dual} to) the basis
$\{v^i\}$ if
\begin{equation}
\inp r^i, v^j> = \delta_{ij}, \qquad i = 1 ,\ldots, n, \quad j = 1, \ldots, n.
\elabel{mod-dual}
\end{equation}If a dual basis is available, then again the representation of $x$ with respect to the
basis $\{v^i\}$ is easily computed.  For suppose that $x$ is represented by $\alpha$:
\begin{eqnarray*}
x = \sum_{i=1}^n \alpha_i v^i
\end{eqnarray*}
Then, by the dual property and linearity we have
\begin{eqnarray*}
\inp  r^j, x> = \inp r^j, \sum_{i=1}^n \alpha_i v^i> = \sum_{i=1}^n \alpha_i
 \inp r^j, v^i> .
\end{eqnarray*}
Since $ \inp r^j, v^i> = \delta_{ij}$, this shows that $
x = \sum_{i=1}^n \inp r^i, x> v^i$.
Of course, to to apply this formula we must have the reciprocal basis  $\{ r^i\}$,
which may be as difficult to find as the inverse Grammian.

In the vector space $(\Co ^n,\Co)$ define the matrices
\begin{eqnarray*}
M = [v^1\cdots v^n] \qquad
R = \left[ \begin{array}{c} r^{1*} \\ \vdots \\ r^{n*} \end{array}
\right]
\end{eqnarray*}
From the defining property \eq mod-dual/ of the dual basis, we must have $RM=I$, so that
$R = M^{-1}$.

\section{Adjoint transformations}
 \label{s:chap1.5.8.4}

Suppose that $\clA$ is a linear operator from the vector space $(\clX,\Co)$ to another
vector space $(\clY,\Co)$, that is, $\clA : \clX \to \clY$.  Then the \defn{adjoint} is
a linear operator working in the reverse direction:
\begin{eqnarray*}
\clA^* : \clY \to \clX.
\end{eqnarray*}
Its definition is subtle since it is not directly defined through $\clA$.  We say that
$\clA^*$ is the adjoint of $\clA$, if for any $x\in\clX$ and any $y\in \clY$,
\begin{eqnarray*}
\inp  \clA (x), y> = \inp x, \clA^* (y)>.
\end{eqnarray*}

To illustrate this concept, let us begin with the finite dimensional case
\begin{eqnarray*}
\clX = \Co^n,\qquad \clY = \Co^m
\end{eqnarray*}
and suppose that $\clA$ is defined through an $m\times n$ matrix $A$, so that $\clA (x)
=  A x$,  $ x\in\clX$.  We may then compute the adjoint using the definition of the
inner products on $\clX$ and $\clY$ as follows
\begin{eqnarray*}
\inp \clA(x), y> = (Ax)^* y = \bar  x^T \bar A^T y  = \inp x, \bar A^T y>
\end{eqnarray*}
Thus,  the adjoint of $\clA$ is defined through the complex conjugate transpose of $A$:
\begin{eqnarray*}
\clA^* (y)=\bar A^T y = A^* y
\end{eqnarray*}

\begin{matlab}

\item
Matlab is virtually designed to deal with the numerical aspects of the vector space concepts
described in this chapter.  Some useful commands are

\item[INV] to compute the inverse of a matrix.
\item[DET] to compute the determinant.
\item[EIG] finds eigenvalues and eigenvectors.
\item[RANK] computes the rank of a matrix.

\end{matlab}

\begin{summary}

In this chapter we have provided a brief background on several different topics in
linear algebra, including
\balphlist
\item
Fields and vector spaces.

\item
Linear independence and bases.

\item
Representations of vectors, and how these representations change under a change of
basis.

\item
Linear operators and matrices.

\item
Inner products and norms.

\item
Adjoint operators.

\item
Eigenvectors.
\end{list}
Good surveys on linear algebra and matrix theory may be found in Chapter~2 of \CHE, or
Chapters~4 and 5 of \BRO.

\end{summary}

\begin{exercises}

\item
Determine conclusively which of the following are fields:
\balphlist
\item The set of integers.
\item The set of rational numbers.
\item The set of polynomials of degree less than 3 with real coefficients.
\item The set of all $n \times n$ nonsingular matrices.
\item The set $\{ 0 , 1\}$ with addition being
                binary ``exclusive-or'' and multiplication being binary ``and''.

\end{list}

\item
Define rules of addition and multiplication such that the   set consisting of three
elements
$\{ a,b,c\}$ forms a field.  Be sure to define the zero and unit elements.

\item
Let $(X,\clF)$ be a vector space, and $Y\subset X$ a subset of $X$.
If $Y$ satisfies the closure property $y_1,y_2\in Y,\ \alpha_1,\alpha_2\in \clF
\Longrightarrow \alpha_1 y_1+\alpha_2y_2\in Y$,
show carefully using the definitions that $(Y,\clF)$ is a subspace of $(X,\clF)$.

\vspace{.25in}\centerline{\textbf{\textit{Linear independence and bases}}}

\item
Which of the following sets are linearly independent?
\balphlist
\item
$\displaystyle
\left[
\begin{matrix}1 \\ 4 \\ 2\end{matrix}
\right],
\left[
\begin{matrix} 0 \\ 2 \\ 3\end{matrix}
\right],\left[
\begin{matrix} 0 \\ 0 \\ 1\end{matrix}
\right]
$
in $(\Re^3,\Re)$.

\item$\displaystyle
\left[\begin{matrix}1 \\ 2 \\ 3\end{matrix}\right] , \left[\begin{matrix}4 \\ 5 \\ 6\end{matrix}\right] , \left[\begin{matrix}7 \\ 8 \\ 9\end{matrix}\right] ,
        \left[\begin{matrix}7912 \\ -314 \\ 0.098\end{matrix}\right]$
in $(\Re^3,\Re)$.

\item
$\displaystyle
\left[
\begin{matrix} 1 \\ 4j \\ 2\end{matrix}
\right],
\left[
\begin{matrix} 0 \\ -2 \\ j\end{matrix}
\right],\left[
\begin{matrix} j \\ 0 \\ 0\end{matrix}
\right]
$
in $(\Co^3,\Co)$.

\item
$\displaystyle
 \sin(t), \cos(t), t$ in $(C[0,1],\Re)$ - the set of real-valued
continuous functions on $[0,1]$ over the real field.

\item
$\displaystyle
e^{jt}, \sin(t), \cos(t), t$ in $(C[0,1],\Co)$ - the set of complex-valued continuous
functions on $[0,1]$ over the complex field.

\end{list}

\item
Determine which of the following sets of vectors are linearly independent in $\Re^3$
by computing the  determinant of an appropriate matrix.
\balphlist
\item
$\displaystyle
\left[\begin{matrix}1 \\ 0 \\ 2\end{matrix}\right] , \left[\begin{matrix}2 \\ 0 \\ 1\end{matrix}\right] , \left[\begin{matrix}0 \\ 5 \\ 1\end{matrix}\right]
$.

\item
$\displaystyle
\left[\begin{matrix}4 \\ 5 \\ 1\end{matrix}\right] , \left[\begin{matrix} 1 \\ 2 \\ -1\end{matrix}\right] , \left[\begin{matrix}2 \\ 1 \\ 3\end{matrix}\right]
$.

\end{list}

\item
\hwlabel{dimCn}
For the vector space $(\Co^n, \Re)$,
\balphlist
\item Verify that this is indeed a vector space.
\item What is its dimension?
 \item Find a basis for $(\Co^n, \Re)$.
\end{list}

\item
 Let $\Re^{2 \times 2}$ be the set of all $2 \times 2$ real matrices.
\balphlist
\item Briefly verify that $\Re^{2 \times 2}$ is a vector space under the usual
        matrix addition and scalar multiplication.
\item What is the dimension of $\Re^{2 \times 2}$ ?
\item Find a basis for $\Re^{2 \times 2}$.
\end{list}

\item\hwlabel{pre-CH}
Is the set $\left\{ I, A ,A^2 \right\}$ linearly dependent or
independent in $(\Re^{2 \times 2}, \Re)$, with $
A = \left[\begin{matrix}1& 1 \\ 0 & 2\end{matrix}\right]$?

 \vspace{.25in}\centerline{\textbf{\textit{Representations of vectors}}}

\item
Given the basis
$\left\{ \displaystyle
v^1=\left[\begin{matrix}1\\ 1\\ 1\end{matrix}\right],
\
v^2=\left[\begin{matrix}2\\ 0\\ 0\end{matrix}\right],
\
v^3=\left[\begin{matrix}1\\ 0\\ 1\end{matrix}\right] \right\}
$
and the vector
\begin{eqnarray*}
x = \left[\begin{matrix}3\\ 2\\ 1\end{matrix}\right]  = \alpha_1 v^1 +\alpha_2 v^2 +\alpha_3 v^3
\end{eqnarray*}
\balphlist
\item
Compute the reciprocal basis.

\item
Compute the Grammian.

\item
Compute the representation of $x$ with respect to $\{v^i\}$ using your answer to (a).

\item
Compute the representation of $x$ with respect to $\{v^i\}$ using your answer to (b).

\end{list}

\vspace{.25in}\centerline{\textbf{\textit{Linear operators and matrices}}}

\item Compute the null space,   range space, and rank of the following matrices.
\balphlist
\item
 $ \displaystyle
A = \left[\begin{matrix}1&1&2\\
0&2&2\\0&3&3\\\end{matrix}\right].$

\item
$ \displaystyle
A = \left[\begin{matrix}1 & 3 & 2 & 1 \\ 2 & 0 & 1 & -1 \\ -1 & 1 & 0 & 1\\\end{matrix}\right].
$

\end{list}

\item Let $b \in \Re^n$ and $A \in \Re^{n \times m}$.  Give necessary and
sufficient conditions on $A$ and $b$ in order that the linear system of
equations $Ax=b$ has a solution $x \in \Re^m$.

\item
Let ${\clA} : \Re^3 \rightarrow \Re^3$ be a linear operator. Consider the
two sets $B=\{b_1,b_2,b_3\}$ and $C=\{c_1,c_2,c_3\}$ below
\begin{eqnarray*}
B=\left\{ \left[\begin{matrix}1\\0\\0\end{matrix}\right], \left[\begin{matrix}0\\1\\0\end{matrix}\right] , \left[\begin{matrix}0\\0\\1\end{matrix}\right] \right\} , \
C=\left\{ \left[\begin{matrix}1\\1\\0\end{matrix}\right], \left[\begin{matrix}0\\1\\1\end{matrix}\right] , \left[\begin{matrix}1\\0\\1\end{matrix}\right] \right\}
\end{eqnarray*}
It should be clear  that these are bases for $\Re^3$.

\balphlist

\item
Find the transformation $P$ relating the two bases.

\item
Suppose the linear operator $\clA$ maps
\begin{eqnarray*}
{\clA}(b_1) = \left[\begin{matrix}2\\-1\\0\end{matrix}\right] , \
{\clA}(b_1) = \left[\begin{matrix}0\\0\\0\end{matrix}\right] , \
{\clA}(b_1) = \left[\begin{matrix}0\\4\\2\end{matrix}\right]
\end{eqnarray*}
Write down the matrix representation of $\clA$ with
respect to the basis $B$ and also with respect to the basis $C$.
\end{list}

\item Find the inverse of the matrix $A$, where $B$ is a matrix.
\begin{eqnarray*}
A = \left[\begin{matrix}I&B\\0&I\end{matrix}\right]
\end{eqnarray*}

\item Consider the set $P_n$ of all polynomials of degree {\it strictly
less than} $n$, with real coefficients, where $x\in P_n$ may be written
$x = a_0 + a_1t + \cdots + a_{n-1} t^{n-1}$.

\balphlist
\item
Verify that $P_n$ is a vector space, with the usual definitions of polynomial
addition, and scalar multiplication.

\item
Explain why $\{ 1,t,\dots,t^{n-1}\}$ is a basis, and thus why $P_n$ is $n$-dimensional.

\item
Suppose $x = 10 - 2t + 2t^2-3t^3$.  Find the representation $\alpha$ of $x$ with respect to the basis in (b) for $n=4$.

\item
Consider {\it differentiation}, $\frac{d}{dt}$, as an operator $\clA\colon P_n\to P_{n-1}$.  That is, $\clA(x) = \frac{d}{dt}x$.  Show that $\clA$ is a linear operator, and compute its null space and range space.

\item
Find $A$, the matrix representation of $\clA$ for $n=4$, using the basis in $(b)$.
Use your $A$ and $\alpha$ to compute the derivative of $x$ in (c).

\end{list}

\vspace{.25in}\centerline{\textbf{\textit{Inner products and norms}}}

\item
Let $(\clV,\Co)$ be an inner product space.
\balphlist
\item
Let $x,y \in \clV$ with $x$
orthogonal to $y$. Prove the Pythagorean theorem:
\begin{eqnarray*}
\| x+ y\|^2 = \|x\|^2 + \|y\|^2
\end{eqnarray*}

\item
 Prove that in an inner product space,
\begin{eqnarray*}
\|x+y\|^2  + \|x-y\|^2 = 2 \|x\|^2 +2\|y\|^2
\end{eqnarray*}
where $\| \cdot \|$ is the norm induced by the inner product. This
is called the {\em Parallelogram law}. Can you give a geometric
interpretation of this law in $\Re^2$?

\end{list}

\vspace{.25in}\centerline{\textbf{\textit{Adjoint operators}}}

\item\hwlabel{grammian}
If $\clA\colon \clX\to\clY$ where $\clX$ and $\clY$ are inner product spaces,
the adjoint $\clA^*$ is a mapping $\clA^* \colon\clY\to\clX$.  Hence, the composition $ Z =
\clA^* \circ \clA$ is a mapping from $\clX$ to itself.  Prove that $\clN(Z)
=\clN(\clA)$.

\item\hwlabel{pre-adjointSystem}
For $1\le p < \infty$,          let $L_p$ denote functions
$f\colon(-\infty,\infty)\to \Co$ such that $\int_{-\infty}^\infty |f(s)|^p\,
ds<\infty$. For $p=\infty$, $L_p$ denotes bounded functions
$f\colon(-\infty,\infty)\to \Co$. The set $L_p$ is a vector space over the
complex field.

Define the function $\clA \colon L_p \to L_p$ as $\clA(f) = a*f$, where  ``$*$''
denotes convolution.  We assume that  for some constants
$C<\infty$, $c>0$,  we have the bound  $|a(t)| \le Ce^{-c|t|}$ for all
$t\in\Re$.  This is sufficient to ensure that $\clA\colon L_p\to L_p$ for any $p$.
\balphlist

\item
First consider the case where $p=\infty$, and let $f_\omega(t) = e^{j\omega t}$,
where $\omega\in\Re$.  Verify that  $f_\omega\in L_\infty$ is an eigenvector of
$\clA$. What is the corresponding eigenvalue?

\item
In the special case $p=2$, $L_p$ is an inner product space with
\begin{eqnarray*}
<f,g>_{L_2} = \int f^*(s) g(s) \, ds,\qquad f,g\in L_2.
\end{eqnarray*}
Compute the adjoint $\clA^*\colon L_2\to  L_2$, and find conditions on $a$
under which $\clA$ is self adjoint.

\end{list}

\item
Let $X=\Re^n$ with the usual inner product.  Let $Y=L_2^n[0,\infty)$,
the set of functions $f\colon[0,\infty)\to\Re^n$ with $\int_0^\infty
|f(s)|^2\,ds <\infty$. We define the inner product as before:
\begin{eqnarray*}
<f,g>_Y =\int f^\top(s) g(s) \, ds\qquad f,g\in Y
\end{eqnarray*}

For an $n\times n$ Hurwitz matrix $A$, consider the differential equation $\dot x = Ax$.
By stability, for each initial condition $x_0\in X$, there exists a unique
solution $x\in Y$.  Define $\clA\colon X\to Y$ to be the map which takes
the initial condition $x_0\in X$ to the solution $x\in Y$.
\balphlist

\item
Explain why $\clA$ is a linear operator.

\item
What is the  null space $N(\clA)$? What is the rank of $\clA$?

\item
Compute the adjoint $\clA^*$.
\end{list}

 \vspace{.25in}\centerline{\textbf{\textit{Eigenvectors}}}

\item Find the eigenvalues of the matrix
\begin{eqnarray*}
A = \left[\begin{matrix}1&1&0&0\\0&1&0&0\\4&5&1&5\\1&2&0&1\end{matrix}\right]
\end{eqnarray*}
\item
An $n\times n$ matrix $A$ is called \defn{positive definite} if it is symmetric,
\begin{eqnarray*}
A = A^* = \bar A^T,
\end{eqnarray*}
and if for any $x\neq \zero $, $x\in\Re^n$,
\begin{eqnarray*}
x^* A x > 0.
\end{eqnarray*}
The matrix $A$ is \defn{positive semi-definite} if the strict inequality in the above
equation is replaced by ``$\ge$''.   Show that for a positive definite matrix,
\balphlist
\item
Every eigenvalue is real and strictly positive.

\item
If $v^1$ and $v^2$ are eigenvectors corresponding to different eigenvalues $\lambda_1$
and $\lambda_2$, then  $v^1$ and $v^2$ are orthogonal.

\end{list}

\item
For a square matrix $X$ suppose that (i) all of the eigenvalues of $X$ are strictly
positive, and (ii)  the domain of $X$ possesses an orthogonal basis consisting
entirely of eigenvectors of $X$. Show that $X$ is a positive definite matrix (and
hence that these two properties completely characterize positive definite matrices).

\textit{Hint}: Make a change of basis using the modal matrix $M=[v^1\cdots v^n]$, where
$\{v_i\}$ is an orthonormal basis of eigenvectors.

\item\hwlabel{left-eig}
Left eigenvectors $\{\omega^i\}$ of an $n\times n$ matrix $A$ are defined by
$\omega^i A = \lambda_i \omega^i$, where  $\{\omega^i\}$ are row vectors
($1\times n$ matrices) and the $\{\lambda_i\}$ are the eigenvalues of $A$.
Assume that the eigenvalues are distinct.
\balphlist
\item
How are the $\{\omega^i\}$ related to the ordinary (right) eigenvectors of $A^*$?
\item
How are the $\{\omega^i\}$ related to the reciprocal (dual) eigenvectors of $A$?

\end{list}
\end{exercises}

\chapter{Solutions of State Equations}
\clabel{chap1.6}

In this chapter we investigate solutions to the linear state space model \eq stateSpace/ in the
general setting with $ x\in\Re^n$, $y\in \Re^p$ and $u\in\Re^m$.  We begin with the linear,
time-invariant (LTI) case where the matrices $(A,B,C,D)$ do not depend upon time.

\section{LTI state space models}

Assume that an initial condition $x_0$ is given at time $t=0$.  The LTI case is then
particularly simple to treat since we can take
Laplace transforms in \eq stateSpace/ to obtain the transformed equations
\begin{eqnarray*}
sX(s) - x_0 &=& A X (s) + BU (s)
\\
Y (s)  &=& C X (s) + D U (s)
\end{eqnarray*}
Solving for $X(s)$ gives $[sI -A] X (s) = x_0 + B U (s)$, and therefore
\begin{eqnarray}
X(s) = [ sI - A]^{-1} x_0 + [sI - A]^{-1} B U (s)
\elabel{eq1.6.1}
\end{eqnarray}
Letting
\begin{eqnarray*}
\Phi (s) \eqdef [sI - A]^{-1}
\end{eqnarray*}
we finally obtain
\begin{equation}
X (s) = \Phi (s) x_0 + \Phi (s) B U (s).
\elabel{soln-laplace}
\end{equation}To convert this into the time domain, we will take inverse Laplace transforms.  First
define
\begin{equation}
\phi(t)=\clL^{-1} (\Phi(s))=\clL^{-1} ([sI - A]^{-1})
\elabel{stm-def}
\end{equation}where the inverse Laplace transform of the matrix-valued function $\Phi(s)$ is taken
term by term.  Since products in the transform domain are mapped to convolutions in
the time domain, the formula \eq soln-laplace/ may be inverted to give
\begin{eqnarray*}
x (t) = \phi (t) x_0 + \int_0^t \phi (t- \tau) Bu (\tau)\, d \tau.
\end{eqnarray*}
The function $\phi$ defined in \eq stm-def/ is evidently important in
understanding the solutions of \eq stateSpace/.  Since in the absence of inputs the
state evolves according to the equation
\begin{equation}
x (t) = \phi (t) x_0,
\elabel{phi}
\end{equation}the function $\phi$ is known as the   \defn{state transition matrix}.

As for the output, in the transform domain we have
\begin{eqnarray*}
Y (s) = CX (s) + D U (s)
\end{eqnarray*}
Assuming for simplicity that $x (0) = 0$, this gives
\begin{eqnarray*}
Y (s) = CX (s) + D U (s) = C \Phi (s) B U (s) + D U (s)
\end{eqnarray*}
Therefore the transfer function is a $p\times m$ matrix-valued function of  $s$ which
takes the form
\begin{eqnarray*}
G(s) = C \Phi (s) B + D ,
\end{eqnarray*}
and the impulse response is the inverse Laplace transform of the transfer function:
\begin{eqnarray*}
g(t) \eqdef \clL^{-1} \{G (s)\}  =  C \phi (t) B   + D\delta (t).
\end{eqnarray*}

In the time domain, we conclude that for zero initial conditions,
\begin{eqnarray*}
y (t) = g* u\, (t) \eqdef  \int_0^t g (t - \tau) u (\tau) d \tau
     = \int_0^t C \phi (t - \tau) Bu (\tau) d \tau + Du (t)
\end{eqnarray*}

\section{Other descriptions of the state transition matrix}
 \label{s:chap1.6.2.1}

We have seen that the solution of the differential equation
\begin{equation}
\dot x = Ax, \qquad x (0) = x_0
\elabel{control-free}
\end{equation}is given by the formula \eq phi/.
We also know that in the scalar case where $n=1$,  the solution of \eq control-free/ is given
by  the exponential $x(t) = e^{At}x_0$, so that $\phi(t) = e^{At}$ in this special case.  In the
multivariable case, define the matrix exponential through the power series
\begin{eqnarray*}
e^{At} \eqdef \sum_{k=0}^\infty \frac{t^k }{ k!} A^k = I + tA + \frac{t^2}{ 2} A^2 \cdots,
\end{eqnarray*}
so that $e^{At}$ is  a \textit{matrix-valued function} of the matrix $A$.  We now show that
we do indeed have $\phi(t) = e^{At}$
with this definition of the matrix exponential.  To see that $x(t) = e^{At} x_0$ is a
solution to \eq control-free/, first note that the specified initial conditions are
satisfied:
\begin{eqnarray*}
x(0) =  e^{A\, 0} x_0  = I x_0 = x_0.
\end{eqnarray*}
To see that the differential equation is satisfied, take derivatives term by term
to obtain the formula
\begin{eqnarray*}
\frac{d}{dt} e^{At} = A + tA^2 + \cdots = Ae^{At}
\end{eqnarray*}
This implies that $ \frac{d}{dt}( e^{At}x_0) = A( e^{At}x_0)$, so that the equation \eq
control-free/ is satisfied,  and hence it is correct to assert that $\phi (t) =
e^{At}$.

The matrix exponential satisfies the following three properties, for any $t, s$,
\begin{enumerate}
\item
$e^ {A 0} = I$.

\item
$\left( e^{At} \right)^{-1} =  e^{-At}  $.

\item
$e^{At}  e^{As} = e^{A(t+s)}  $.
\end{enumerate}
The third property, which is known as the \defn{semigroup property}, is analogous to the
equation $e^A e^B = e^{A + B}$ when $A$ and $B$ are scalars.  However, this equation \textit{does not}
hold in the matrix case unless the matrices $A$ and $B$ commute.

To summarize, we now have three descriptions of the state transition matrix $\phi (t)$:
\balphlist
\item
$  \phi (t) = \clL^{-1} \{ (sI - A)^{-1} \} $
\item
$\displaystyle
	\phi (t) = e^{At} = \sum_{k=0}^\infty  \frac{t^k }{  k!} A^k$
\item
$\dot\phi (t) = A \phi (t)$, $\phi_0 = I$
\end{list}

\section{Change of state variables}

We now consider the effect of a change of basis on the state transition matrix $\phi$.
Suppose that the $n\times n$ matrix $P$ is used to make a change of basis, and define
the   variable $\bar x$ as
\begin{eqnarray*}
\bar  x = Px \qquad x = P^{-1} \bar  x
\end{eqnarray*}
After pre-multiplying by $P$, the state equation
\begin{eqnarray*}
\dot x = Ax + Bu
\end{eqnarray*}
 becomes a state equation for $\bar x$:
\begin{eqnarray*}
\dot {\bar  x} = \underbrace{PAP^{-1}}_{\bar   A} \bar  x + \underbrace{PB}_{\bar  B} u
\end{eqnarray*}
The output equation $ y = Cx + Du$ can be similarly   transformed:
\begin{eqnarray*}
  y = \underbrace{CP^{-1}}_{\bar  C} \bar  x + \underbrace{D}_{\bar  D}u
\end{eqnarray*}
Hence for any   invertible matrix
$P$, we can form a new state space model
\begin{eqnarray}
\dot {\bar x} &=& \bar A \bar x + \bar Bu
\nonumber\\
y &=& \bar C\bar x+ \bar Du,
\elabel{bar-stateSpace}
\end{eqnarray}
with $\bar A = PAP^{-1}$, $\bar B = PB$,
$\bar C = C P^{-1}$, and $\bar D = D$.  Moreover,  if $\bar x (0) = P x(0)$ and identical
inputs are applied to both systems, then the outputs of \eq stateSpace/ and \eq
bar-stateSpace/ will be identical.

The state transition matrix of the system \eq bar-stateSpace/ can be expressed in terms
of the state transition matrix of \eq stateSpace/ as follows.  Starting with the
defining formula $\bar  \Phi (s) = (sI - \bar A)^{-1} = (sI - PAP^{-1})^{-1}$ and
noting that $I = P P^{-1}$, we have
\begin{eqnarray*}
\bar  \Phi (s) = (P [sI - A] P^{-1})^{-1} = P [sI - A]^{-1} P^{-1}.
\end{eqnarray*}
It follows that $\bar  \Phi (s) = P \Phi (s) P^{-1}$, and hence that
\begin{equation}
\bar  \phi (t) = P \phi(t) P^{-1}
\elabel{phi-cob}
\end{equation}To show how \eq phi-cob/ may be applied, suppose that the eigenvalues of $A$ are
distinct, so that the modal matrix $M= [v^1 \cdots v^n]$ exists and is non-singular.
Defining $P= M^{-1}$, we have  seen in equation \eq modal-prop/  that
\begin{eqnarray*}
\bar A = PAP^{-1} = \Lambda
\end{eqnarray*}
The point is, when $\bar A$ is diagonal, the state transition matrix is given by the
direct formula
\begin{eqnarray*}
\bar  \phi (t) = e^{\Lambda t} = \left[ \begin{array}{lll} e^{\lambda_1
t} && 0 \\ & \ddots & \\ 0 && e^{\lambda_n t} \end{array} \right]
\end{eqnarray*}
and from \eq phi-cob/ we   have
\begin{eqnarray*}
\phi (t) = M \bar  \phi (t) M^{-1},
\end{eqnarray*}
so that $\phi(t)$ is quickly computed via the modal matrix.

\section{Numerical examples}

Take the two-dimensional case, with $A$ given by
\begin{eqnarray*}
A = \left[ \begin{array}{rr} 0 & 1 \\ 8 & - 2 \end{array} \right]
\end{eqnarray*}
We will compute $\phi(t)$ for the system \eq control-free/ using the various methods
described in the previous sections.  The simple-minded approach is to write down the
formula
\begin{eqnarray*}
\phi (t) = e^{At} \eqdef \left[ \begin{array}{rr} 1 & 0 \\ 0 & 1 \end{array} \right]
+ t \left[ \begin{array}{rr} 0 & 1 \\ 8 & -2 \end{array} \right] + \frac{t^2
}{ 2} \left[ \begin{array}{rr} 8 & -2 \\ -16 & 12 \end{array}
\right] + \cdots ,
\end{eqnarray*}
but this does not give a closed form expression unless the infinite series can be
computed.

In the transform domain, we have
\begin{eqnarray*}
\Phi(s) \eqdef
(sI - A)^{-1} = \left[ \begin{array}{rr} s & -1 \\ -8 & s + 2
\end{array} \right]^{-1} =   \frac{\left[ \begin{array}{cc} s+2 & 1 \\ 8 & s \end{array}
\right] }{ s^2 + 2s - 8}
\end{eqnarray*}
Taking inverse Laplace transforms then gives
\begin{equation}
\phi (t) = \left[
 \begin{array}{ll}
\displaystyle
\frac{1 }{ 3} e^{-4t} + \frac{2 }{ 3}e^{2t},
 &
\displaystyle - \frac{1 }{ 6} e^{-4t} + \frac{1 }{ 6} e^{2t}
 \\
\displaystyle
- \frac{4 }{ 3} e^{-4t} + \frac{4 }{ 3} e^{2t} ,
&
\displaystyle\frac {2 }{ 3} e^{-4t} + \frac{1}{ 3} e^{2t}
\end{array}
\right]
\elabel{example-phi}
\end{equation}Observe that $\phi (0) = I$, as it should be.

A third approach is to convert the state space model to modal form.
To begin, we must compute the eigenvalues of the matrix $A$.  The characteristic
polynomial is given by
\begin{eqnarray*}
\Delta(s) \eqdef \det \left[ \begin{array}{rr} s & -1 \\ -8 & s + 2 \end{array} \right]
= (s+4)(s-2)
\end{eqnarray*}
The roots of $\Delta$, which are also the eigenvalues of $A$, are $\lambda = -4, 2$.
Hence the eigenvalues are distinct, and a modal form does exist.  Solving the equation
$Ax =
\lambda x$
for the eigenvector $x$ in this example gives for $\lambda = -4$,
\begin{eqnarray*}
\left[ \begin{array}{rr} 0 & 1 \\ 8 & -2 \end{array} \right] \left[
\begin{array}{c} x_1 \\ x_2 \end{array} \right] = -4 \left[
\begin{array}{c} x_1 \\ x_2 \end{array} \right]
\end{eqnarray*}
which is solved by $x_2 = -4 x_1$.  It follows that the vector
\begin{eqnarray*}
  v^1 =  \left(\begin{matrix} 1 \\ -4\end{matrix}\right)
\end{eqnarray*}
is an eigenvector for the eigenvalue $\lambda_1 = -4$.  Similarly, the vector $v^2
=\left(\begin{smallmatrix} 1 \\ 2\end{smallmatrix}\right)$ can be shown to be an
eigenvector for the eigenvalue $\lambda_2 =2$.

The modal matrix becomes $M=[v^1\mid v^2] = \left[\begin{matrix}1 & 1 \\ -4 & 2\end{matrix}\right]$, and the state transition
matrix may be expressed as
\begin{eqnarray*}
\phi (t) = M \bar  \phi (t) M^{-1} = \left[\begin{matrix}1 & 1 \\ -4 &
2\end{matrix}\right]\left(\begin{matrix}e^{-4t} \ 0\\0\ e^{2t}
\end{matrix}\right)
\left[\begin{matrix}1 & 1
\\ -4 & 2\end{matrix}\right]^{-1}.
\end{eqnarray*}
After computing the matrix inverse and products, one again obtains the
expression \eq example-phi/.

\section{Cayley-Hamilton Theorem}

Another useful tool in systems theory which also gives a method for computing $\phi$ is
the following theorem attributed to Cayley and Hamilton.   For an $n\times n$ matrix
$A$  with real coefficients, we recall that the characteristic polynomial   defined
as
$\Delta(x) = \det (x I - A)
$ is a polynomial of degree $n$, which therefore may also be expressed as
\begin{eqnarray*}
\Delta(x) =
x^n +
\alpha_1 x^{n-1} + \cdots + \alpha_n
\end{eqnarray*}
for real coefficients $\{\alpha_i\}$.

 \begin{theorem}[Cayley-Hamilton Theorem]
\tlabel{CH}
 The matrix $A$ satisfies the matrix characteristic equation
\begin{equation}
\Delta (A) = A^n + \alpha_1 A^{n-1} + \cdots + \alpha_n I = 0_{n\times n}
\elabel{CH}
\end{equation}\end{theorem}

\proof  We provide a proof for the case of distinct
eigenvalues only - for the general case, see \BRO, pages  286, 294-5.

When the eigenvalues are distinct, we can work with the transformed matrix
\begin{eqnarray*}
\Lambda
=\left[ \begin{array}{rrr} \lambda_1 && 0 \\ & \ddots & \\ 0 &&
\lambda_m \end{array} \right]
=
 M^{-1} AM
\end{eqnarray*}
where $M$ is the modal matrix, $ M = [v^1 \cdots v^n] $.
It follows that
\begin{eqnarray*}
A =  M \Lambda M^{-1},
\end{eqnarray*}
and that for any $k$,
\begin{eqnarray*}
A^k = M \Lambda^k M^{-1}
\end{eqnarray*}
Applying this to \eq CH/ gives
\begin{eqnarray*}
\Delta (A) & = &
M (\Lambda^n + \alpha_1 \Lambda^{n-1} \cdots + \alpha_nI) M^{-1}
\\
& = &
M \left( \left[ \begin{array}{ccc} \lambda_1^n & & 0 \\ & \ddots
& \\ 0 && \lambda_n^n \end{array} \right]
+
\alpha_1 \left[
\begin{array}{ccc} \lambda_1^{n-1} && 0 \\ & \ddots & \\ 0 &&
\lambda_n^{n-1} \end{array} \right]
+ \cdots  +
\alpha_n  \left[
\begin{array}{ccc} 1 && 0 \\ & \ddots & \\ 0 & & 1\end{array} \right]
\right) M^{-1}
\\
& = &
 M \left[ \begin{array}{cccc}
\lambda_1^n + \alpha_1 \lambda_1^{n-1} \ldots + \alpha_n &&& 0 \\
& \lambda_2^n + \alpha_1 \lambda_2^{n-1} + \cdots + \alpha_n \\
&& \ddots & \\
0 &&& \lambda_n^n + \alpha_1 \lambda_n^{n-1} + \cdots + \alpha_n
\end{array} \right] M^{-1}
\end{eqnarray*}
For each $i$, the eigenvalue $\lambda_i$ satisfies the   characteristic equation
$\Delta(\lambda_i)=0$. So,   $
\Delta (A) = 0$, as claimed.
\qed

\begin{ex}
Consider the $2\times 2$ matrix
\begin{eqnarray*}
A = \left[ \begin{array}{rr} 0 & 1 \\ -2 & -3 \end{array} \right]
\end{eqnarray*}
The characteristic polynomial in this example is
\begin{eqnarray*}
\Delta(\lambda)
=
\det (\lambda I - A) =
\det \left[ \begin{array}{rr} \lambda & -1 \\
2 & \lambda + 3 \end{array} \right] = \lambda^2 + 3 \lambda + 2
\end{eqnarray*}
We therefore have
\begin{eqnarray*}
\Delta(A) = A^2 + 3 A + 2 I
=
\left[ \begin{array}{rr} -2 & -3 \\ 6 & 7 \end{array} \right]
+
3\left[ \begin{array}{rr} 0 & 1 \\ -2 & -3 \end{array} \right]
+
2
\left[ \begin{array}{rr} 1 & 0 \\ 0 & 1 \end{array} \right]
\end{eqnarray*}
which is evidently equal to  $0$ (the $2\times2$ null matrix).
\end{ex}

In view of \Theorem{CH}, the  infinite series used to define the transition matrix
may be written as a finite sum
\begin{eqnarray*}
\phi (t) = e^{At} = \sum_{k=0}^\infty \frac{t^k }{ k!} A^k =\sum_{k=0}^{n-1} \beta_k (t)
A^k
\end{eqnarray*}
To see this, apply  the Cayley-Hamilton theorem to replace any $A^k$ for $k \geq
n$ by a linear combination of $I, A, \ldots, A^{n-1}$.

To compute $\phi(t)$, we must then find $\beta_k (t)$.  We illustrate this with a
numerical example, using the matrix $A$ defined above.

\begin{ex}
Our goal is to compute
\begin{equation}
\phi(t)
=
e^{At} = \beta_0 (t) I + \beta_1 (t) A
\elabel{phi-ex}
\end{equation}by computing the scalar functions $\beta_0, \beta_1$.
Solving the characteristic equation $\lambda^2 + 3 \lambda + 2 = 0$ then gives
the distinct eigenvalues
\begin{eqnarray*}
\lambda_{1} = -1,\quad
\lambda_2 = -2
\end{eqnarray*}

From the definition of the matrix exponential we have for $i=1,2$,
\begin{eqnarray*}
e^{At} v^i
=   I v^i   + tA v^i + \frac{t^2}{ 2} A^2v^i   \cdots
=    v^i   + t\lambda_i v^i + \frac{t^2}{ 2}\lambda_i^2   v^i \cdots
=    e^{\lambda_i t} v^i .
\end{eqnarray*}
That is, the eigenvectors $\{v^i\}$ for the matrix $A$ are also eigenvectors for
the matrix $e^{At}$. We also have
\begin{eqnarray*}
[\beta_0(t) I + \beta_1 (t) A ] v^i = [\beta_0(t)   + \beta_1 (t) \lambda_i ]v^i
\end{eqnarray*}
and since $e^{At} = [\beta_0(t) I + \beta_1 (t) A ]$, this gives the system of equations
\begin{eqnarray*}
 e^{\lambda_1 t} = e^{-t}  & = & \beta_0(t) + \beta_1 (t)(-1)
\\
 e^{\lambda_2 t} = e^{-2t} & = & \beta_0(t) + \beta_1(t) (-2)
\end{eqnarray*}
Rewriting this as
\begin{eqnarray*}
e^{-t} & = & \beta_0(t) - \beta_1(t) \\
e^{-2t} & = & \beta_0(t) - 2 \beta_1(t)
\end{eqnarray*}
we may solve to obtain
\begin{eqnarray*}
\beta_0 (t)& = & 2e^{-t} -e^{-2t} \\
\beta_1(t) & = & e^{-t} - e^{-2t}
\end{eqnarray*}
Applying the formula \eq phi-ex/ then gives
\begin{eqnarray*}
\phi (t) & = & e^{At} = \left[ \begin{array}{cc} \beta_0 & 0 \\ 0 &
\beta_0 \end{array} \right] + \left[ \begin{array}{cc} 0 & \beta_1 \\
-2 \beta_1 & - 3 \beta_1 \end{array} \right] \\
& = & \left[ \begin{array}{cc} \beta_0 & \beta_1 \\ -2 \beta_1 &
\beta_0 - 3 \beta_1 \end{array} \right] \\
 & = & \left[
\begin{array}{cc} 2e^{-t} -e^{-2t} & e^{-t} -e^{-2t} \\ -2 e^{-t} +
2e^{-2t} & - e^{-t} + 2e^{-2t} \end{array} \right]
\end{eqnarray*}
\end{ex}

\section{Linear Time-Varying Systems}
\label{s:chap1.6.3}

We conclude this chapter with a development of the linear time-varying (LTV) case
\begin{equation}
\dot x=A (t) x, \qquad x (t_0) = x_0
\elabel{LTVcontrol-free}
\end{equation}Given our experience with the LTI case, one might presume that the response to \eq
LTVcontrol-free/ may be expressed as $x(t) = \exp( \int_{t_0}^t A(s) \, ds ) x_0$.
While this expression is valid in the scalar case, it is unfortunately \textit{not}
correct in general.  This is a more complicated situation since we can no longer take
Laplace transforms, and in general we do not even know if \eq LTVcontrol-free/ has a
solution.  However, we do have the following result.

\begin{theorem}
\tlabel{exist}
If $A (t)$ is a piecewise-continuous function of $t$, then:
\balphlist
\item
For any  $x_0 \in \Re^n$, the
differential equation \eq LTVcontrol-free/ possesses a unique solution $x(t)$ defined on
$(-\infty,\infty)$;

\item
The set of all possible solutions to $\dot x = A (t) x$ forms an $n$-dimensional vector
space.
\end{list}
\end{theorem}

\proof
The proof of (a) is outside of the scope of this book.  To see (b), let $\{x_0^i \}$
be a basis for $\Re^n$, and define for each $1\le i\le n$ the function $\psi^i$ to be the
 solution to \eq LTVcontrol-free/ with $x_0=x_0^i$.    For any $x_0=\sum \alpha_i x_0^i
\in\Re^n$ one may easily verify that $x(t) = \sum \alpha_i \psi^i(t)$ is the unique
solution to \eq LTVcontrol-free/.  Hence the span of $\{\psi^i\}$ contains all solutions
to \eq LTVcontrol-free/.  We now show that the $\{\psi^i\}$ are linearly independent,
and hence provide a basis for the solution set.

To show that the $\{\psi^i (t) \}$ are linearly independent time
functions, assume   that   there exist $\alpha_1, \ldots,
\alpha_n$   such that
\begin{eqnarray*}
\alpha_1 \psi^1 (t) + \cdots + \alpha_n \psi^n (t) = \zero, \qquad t\in\Re.
\end{eqnarray*}
Setting $t=t_0$ gives $\sum \alpha_i x_0^i = \zero$, and since $\{x_0^i\}$ was assumed
to be a basis, we must have $\alpha_i=0$ for all $i$.  By definition, the set
$\{\psi^i (t) , 1\le i\le n\}$ is linearly independent,  and hence forms a basis for
the solution set of \eq LTVcontrol-free/.
\qed

\section{Fundamental matrices}
\label{s:fundamental}

Using \Theorem{exist} we may define an analog of the state transition matrix $\phi(t) =
e^{At}$ constructed earlier in the LTI case.  Let $\{\psi^i (t) , 1\le i\le n\}$ be any
set of solutions to \eq LTVcontrol-free/ which are linearly independent as vectors in
$C^n(-\infty,\infty)$, and define the $n\times n$  \defn{fundamental matrix} $U$ as
\begin{equation}
U (t) = [ \psi^1 (t) \cdots \psi^n (t)]
\elabel{fund-mat}
\end{equation}Since each $\psi_i$ is a solution to \eq LTVcontrol-free/, the fundamental matrix
satisfies the matrix differential equation
\begin{equation}
\dot U(t) = A(t) U (t),
\elabel{fund-eqn}
\end{equation}just as the state transition matrix does in the LTI case.  To proceed, we need to
invert the fundamental matrix.
\begin{lemma}
  \tlabel{fund-invert}
The inverse $U^{-1} (t)$ exists for all $ - \infty < t < \infty$.
\end{lemma}

\proof
Fix $t_0$, and suppose that for  some vector
$x\in\Re^n$,
\begin{eqnarray*}
U(t_0) x = \zero\in\Re^n .
\end{eqnarray*}
We will show that $x$ is necessarily equal to $\zero$, and from this it will follow
that $U(t_0)$ is invertible as claimed.

The function $x(t) = U(t) x = \sum x_i \psi^i(t)$ is evidently a solution to \eq
LTVcontrol-free/ with $x(t_0) = \zero$.  By uniqueness of solutions, we must then have
$x(t)=\zero\in\Re^n$ for all $t$, and hence  $\sum x_i \psi^i=\zero$ as a vector in
$C^n(-\infty,\infty)$.  Since the $\{\psi^i (t) , 1\le i\le
n\}$ are linearly independent, it follows that $x=\zero$.
\qed

Given \Lemma{fund-invert}, we may define
\begin{eqnarray*}
\phi (t, \tau) \eqdef U (t) U^{-1} (\tau) ,
\qquad
 t, \tau \in (- \infty, \infty)
\end{eqnarray*}
The matrix valued function $\phi$ is again called   the \defn{state transition matrix}
of the differential equation $\dot x = A (t) x$.  The state transition matrix shares
many of the properties of the matrix $e^{At}$ defined earlier, although we do
\textit{not have} $\phi(t) = e^{A(t) t}$ in the LTV case!  We do have, however,
the following properties, which hold even in the LTV
case:
\label{properties-phi}  For any $t_0, t_1, t_2$,
\begin{enumerate}
\item $\phi (t_0,t_0) = I$.
\item $\phi^{-1} (t_0,t_1) = \phi (t_1,t_0)$.
\item$\phi (t_2, t_0) = \phi (t_2, t_1) \phi (t_1, t_0) $.

\end{enumerate}
The third property is again known as the \defn{semigroup property}.

The solution to \eq LTVcontrol-free/ can now be conveniently expressed in terms of the
state transition matrix
\begin{eqnarray*}
x (t) = \phi (t , t_0) x (t_0).
\end{eqnarray*}
To see this, use property (1) above, and the matrix differential equation \eq fund-eqn/.
The state transition matrix also gives the solution to the controlled system \eq stateSpace/
in the time varying case.  We have that for any control $u$, initial time $t_0$, and any
initial condition $x(t_0)=x_0$,
\begin{eqnarray*}
x (t) = \phi (t, t_0) x_0 + \int_{t_0}^{t} \phi (t, \tau) B (\tau)u
(\tau) d \tau
\end{eqnarray*}
To see this, check the equation $\dot x = Ax + Bu$ by  differentiating both sides of  the
equation above.

\section{Peano-Baker Series}
 \label{s:chap1.6.5}

There are few tools available for computing the state transition matrix in the LTV case.
One approach is by iteration, which is related to the Taylor series expansion used in the
LTI case.  From the differential equation
\begin{eqnarray*}
\frac{d}{dt} \phi (t, t_0) = A (t) \phi (t, t_0), \qquad \phi (t_0,t_0) = I.
\end{eqnarray*}
and the fundamental theorem of calculus, we have the integral equation
\begin{eqnarray*}
\phi (t, t_0) = I + \int_{t_0}^{t} A (\sigma) \phi  (\sigma, t_0)\, d
\sigma
\end{eqnarray*}
This suggests the following iterative  scheme for finding $\phi$:  Set $\phi_0(t,t_0)\equiv I$,
and
\begin{eqnarray*}
\phi_{k+1} (t, t_0) = I + \int_{t_0}^{t} A (\sigma) \phi_k (\sigma,t_0)\, d \sigma,
\qquad k\ge 0.
\end{eqnarray*}
It can be shown that if $A$ is piecewise continuous, then this \defn{Peano-Baker
series} does converge to the state transition matrix.

We can show this explicitly in the LTI case where $A(t) \equiv A$. The Peano-Baker
series then becomes
\begin{eqnarray*}
\phi_0 (t, t_0) & = & I\\
\phi_1 (t, t_0) & = & I + A \int_{t_0}^{t} I  \, d \sigma   = I + A (t - t_0) \\
\phi_2 (t, t_0) & = & I + A \int_{t_0}^{t} (I + A (\sigma - t_0))\, d
\sigma \\
& = & I + (t - t_0) A + \frac{(t - t_0)^2 }{ 2} A^2
\end{eqnarray*}
Thus, we see again that in the LTI case,
\begin{eqnarray*}
\phi (t , t_0) = \sum_{i=0}^{k} \frac{(t - t_0)^k  }{ k!} A^k.
\end{eqnarray*}

\begin{matlab}

\item[EXPM] computes the matrix exponential.

\end{matlab}
\begin{summary}

We have provided in this chapter a detailed description of solutions to \eq
stateSpace/  through the construction of the state transition matrix.
Chapter~9 of \BRO\ contains a
similar exposition on the solution of state equations; for a  more
detailed exposition, the reader is referred to \ZD.
\end{summary}

\begin{exercises}

\item
Use the (matrix) exponential series to evaluate $e^{At}$ for
\begin{eqnarray*}
\hbox{(a)}\quad
 A=\left[\begin{matrix}0&1\\ 0&0\\\end{matrix}\right]
\qquad
\hbox{(b)}\quad
 A=\left[\begin{matrix}1&0\\ 1&1\\\end{matrix}\right]
\qquad
\hbox{(c)}\quad
  A=\left[\begin{matrix}0&\pi\\  -  \pi &0\\\end{matrix}\right]\;.
\end{eqnarray*}

\item
If $A$ is an $n\times n$ matrix of full rank, show that
\begin{eqnarray*}
\int_0^t e^{A\sigma}d\sigma = A^{-1}[e^{At}-I] \;.
\end{eqnarray*}
 \textit{Hint}: This is easy once you realize that both functions of $t$ solve
the same differential equation.

Using this result, obtain the solution to the linear time-invariant equation
\begin{eqnarray*}
{\dot x} = Ax+B\bar u \;,\quad x(0)=x_0\;,
\end{eqnarray*}
where $\bar u$ is a constant $m$-dimensional vector and $B$ is an
$(n\times m)$-dimensional matrix.

\item
  Given the single-input single-output system
\begin{eqnarray*}
\dot x =\left[\begin{matrix}0&1\\ 1&0\\\end{matrix}\right]x +\left[\begin{matrix}0\\ 1\\\end{matrix}\right]\;;\quad
x(0)=\left[\begin{matrix}1\\ 0\\\end{matrix}\right]\;,\quad y=\left[\begin{matrix}1&0\\\end{matrix}\right]x\;,
\end{eqnarray*}
let $\phi(t,s)$ denote the corresponding state transition matrix.
\balphlist
\item
 Evaluate $\phi(t,s)$ using the transfer function approach
({\it i.e.} $(sI-A)^{-1}$).
\item
 Evaluate $\phi(t,s)$ by transforming the system into
diagonal form ({\it i.e.} using modal decomposition).

\item
With $u(t)\equiv 1$, compute $y(1)$.
\end{list}

\item
Given the second-order system
\begin{eqnarray*}
\dot x =\left[\begin{matrix}0&1\\  -  4& -  5\\\end{matrix}\right]x\qquad
;\quad x(0)=x_0\;,
\end{eqnarray*}
find the conditions on the initial state vector $x_0$ such that only the
mode corresponding to the smaller (in absolute value) eigenvalue is
excited.

\item  Given the matrix
\begin{eqnarray*}A=\left[\begin{matrix} -  1&2&0\\ 1&1&0\\ 2& - 1&2\\\end{matrix}\right]\;,
\end{eqnarray*}
\balphlist
\item
 Compute its inverse
using the Cayley-Hamilton Theorem.

\item
Compute $A^6$.
\end{list}

\item
 Revisit \Exercise{pre-CH} in \Chapter{chap1.5},
and using the Cayley-Hamilton Theorem, explain why these three matrices must be
linearly dependent.

\item
Suppose that $A$ and $B$ are constant square matrices.
Show that the
state transition matrix for the time-varying system described by
\begin{eqnarray*}
{\dot x}(t) = e^{-At}Be^{At}x(t)
\end{eqnarray*}
is
\begin{eqnarray*}
\Phi(t,s)=e^{-At}e^{(A+B)(t-s)}e^{As}\;.
\end{eqnarray*}

\item\hwlabel{pre-Lyapunov}
Suppose that $A$ and $Q$ are $n\times n$ matrices, and consider the differential equation
\begin{equation}
\dot Z = A Z + Z A^*,\qquad Z(0) = Q.
\elabel{L-ode}
\end{equation}\balphlist
\item
Show using the product rule that the unique solution to \eq L-ode/ is
\begin{eqnarray*}
Z(t) =  e^{At} Q e^{A^* t}
\end{eqnarray*}
\item
Integrating both sides of \eq L-ode/ from $0$ to $t_f$, and applying the fundamental
theorem of calculus, conclude   that if $ e^{At}\to 0$, then
\begin{eqnarray*}
P = \lim_{t_f\to\infty} \int_0^{t_f} Z(t) \, dt
\end{eqnarray*}
is a solution to the
\defn{Lyapunov equation}
\begin{eqnarray*}
 A P +P A^* + Q  = 0.
\end{eqnarray*}
This fact will be very useful when we study stability, and it will come up again in our study
of optimization.
\end{list}

\item  If $\phi$ denotes the state transition matrix for the system $\dot x\, (t) = A(t) x(t)$,
show that   $V(t) =\phi(t,t_0) Q \phi^*(t,t_0)$ is a solution to the equation
\begin{eqnarray*}
\frac{d}{dt}  V(t) = A(t) V(t) + V(t) A^*(t),\qquad  V(t_0) = Q.
\end{eqnarray*}

\item \hwlabel{adjoint-stm}
Associated with the control free differential equation \eq LTVcontrol-free/ is the
so-called \defn{adjoint equation}
\begin{eqnarray*}
\dot z = -A^* z,
\end{eqnarray*}
where $A^* = \bar A^T$.  Show that the state transition matrix $\phi_a$ for the adjoint
equation satisfies
\begin{eqnarray*}
\phi_a(t_1,t_0) = \phi^*(t_0,t_1) \,,
\end{eqnarray*} where $\phi (t , t_0)$ is the state transition matrix for \eq LTVcontrol-free/.

\item
The LTV differential equation \eq stateSpace/, with $D=0$, has the \textit{adjoint
equations}
\begin{eqnarray*}
\dot z &=& -A^*z  + C^* v\\
w &=& B^* z\,.
\end{eqnarray*}
\begin{description}
\item[(i)]
Show that the impulse response matrix $g_a$ for the adjoint model satisfies
\begin{eqnarray*}
g_a(t_1,t_0) = g^*(t_0,t_1) .
\end{eqnarray*}
\item[(ii)]
In the linear case, this relation becomes
\begin{eqnarray*}
g_a(s) = g^*(-s) .
\end{eqnarray*}
Comment on this result, in view of \Exercise{pre-adjointSystem} of the previous chapter.
\end{description}

\notes{What is the adjoint system for CCF?}

\item
Derive the three properties of the state transition matrix given on
page~\pageref{properties-phi},
\balphlist
\item
In the LTI case, using properties of the matrix exponential.

\item
In the LTV case, using the definition of $\phi$ in terms of the fundamental matrix.

\end{list}
\end{exercises}


\part{System Structural Properties}

\chapter{Stability}
\clabel{stability} 
 
Depicted in \Figure{trajectories} are results from several
simulations conducted using \textit{Matlab}.  Three different linear
state space models of the form \eq stateSpace/ have been simulated
using three different initial conditions.
\begin{figure}[ht]
\ebox{.95}{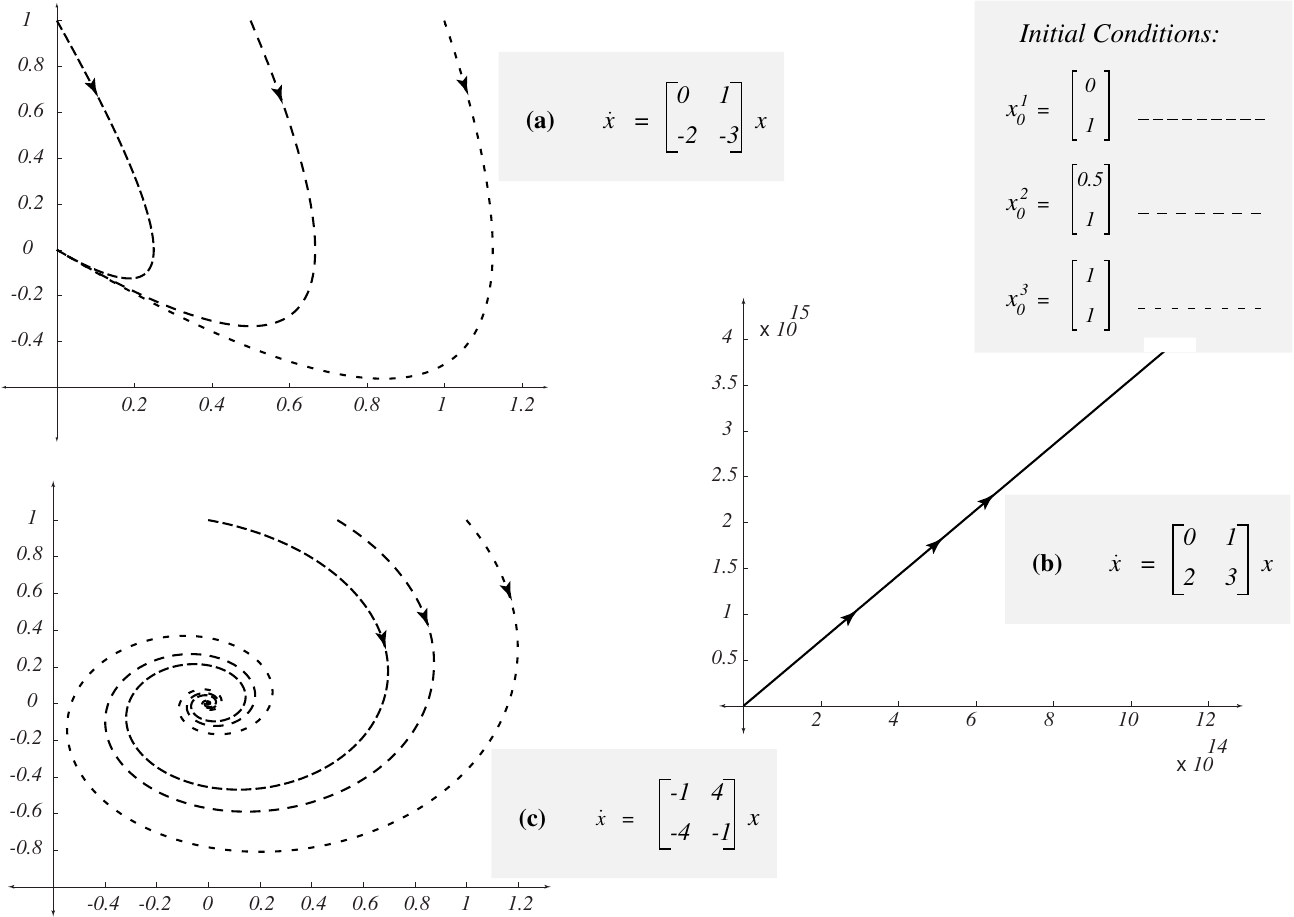}
\caption{Typical solutions of the linear state space model}
\flabel{trajectories}
\end{figure}

One can see in these plots three substantially different modes of
behavior.  In (a) the solutions are attracted to the origin along the
line $x_2 = - x_1$, while in (c) solutions are also attracted to the
origin, but they move there in a cyclical fashion.  In (b) solutions
explode to infinity approximately along a straight line.

These simulations suggest that some solidarity among different states
exists, but this is not true in general.  Consider for example the
interconnection shown in
\Figure{structure-unstableSim}.  From the state equations
\begin{eqnarray*}
\dot x_1 & = & x_2 \\
\dot x_2 & = & 2x_1 + x_2 + u
\end{eqnarray*}
we have, in the Laplace transform domain,
\begin{eqnarray*}
X_1(s) = \frac{1}{ s^2 - s - 2} U(s) \qquad X_2(s) =  \frac{s}{ s^2 - s - 2} U(s) 
\end{eqnarray*}
Thus the output may be expressed in terms of the input as
\begin{eqnarray*}
Y(s) = X_2(s) - 2 X_1(s) =   \frac{s-2 }{ s^2 - s - 2} U(s) 
            =   \frac{1 }{ s + 1} U(s) 
\end{eqnarray*}
Since the single pole $p=-1$ is in the left hand plane, this is a BIBO stable model, as
defined in any undergraduate signals and systems course.  However, the state equations also
give 
\begin{eqnarray*}
A = \left[ \begin{array}{rr} 0 & 1 \\ 2 & 1 \end{array} \right]
\end{eqnarray*}
Solving for the eigenvalues of $A$ yields
\begin{eqnarray*}
\det (sI - A) = 0 \Longrightarrow \lambda_1 = -1, \lambda_2 = 2
\end{eqnarray*}
Thus, we have an unstable root $\lambda_2 = 2$, but it is cancelled by
the numerator of the transfer function due to the feed forward of
$x_2$.  For certain initial conditions one would expect the state
to diverge to infinity, even if the input and output remain bounded.  

\begin{figure}[ht]
\ebox{.5}{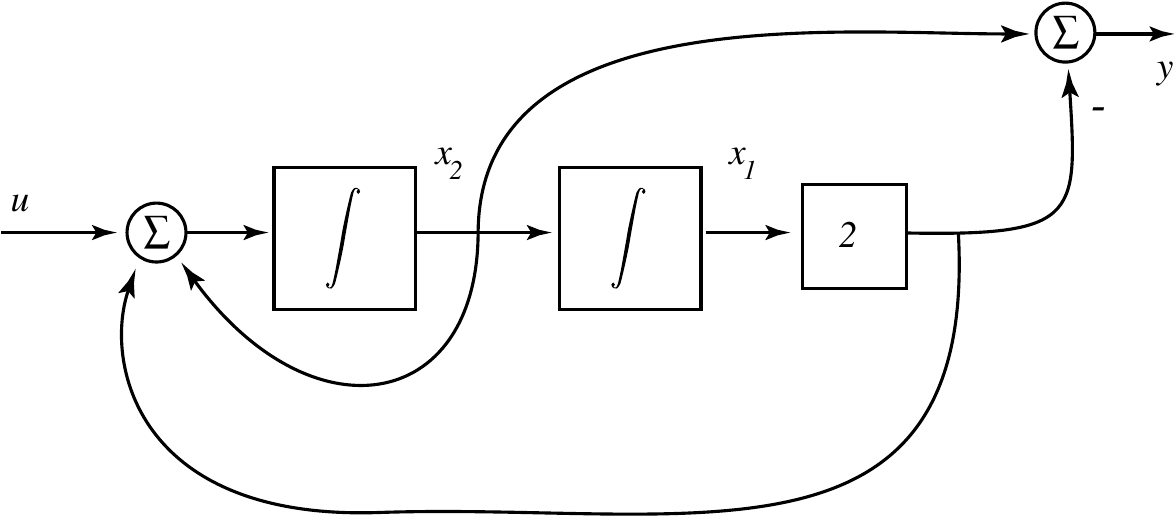}
\caption{Is this model stable, or unstable?}
\flabel{structure-unstableSim}
\end{figure}

In this chapter we  present an introduction to stability theory for nonlinear models, 
and in the process we  give a complete characterization of the types of behavior
possible for a linear state space model.   
 
\begin{figure}[ht] 
\ebox{.5}{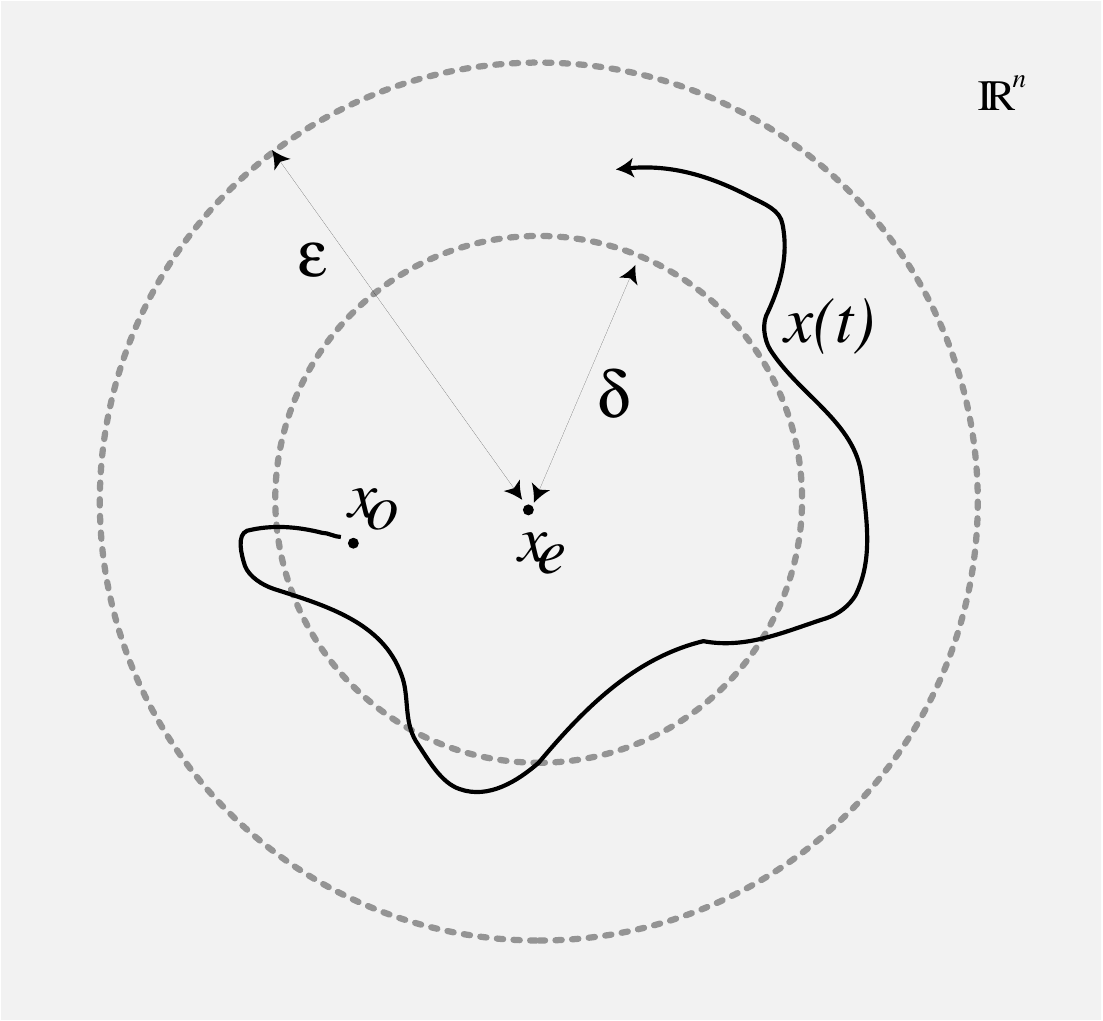}
\caption{ If $x_0 \in B_\delta (x_e)$,  
then $ x(t)=\varphi (t; x_0) \in B_\epsilon (x_0)$ for
all $t\ge 0$.}  
\flabel{Trajectory}
\end{figure}

\section{Stability in the sense of Lyapunov}

Once a feedback control law has been specified, it is frequently possible to express the
closed-loop system in the form of a nonlinear state space model
\begin{equation}
\dot x = f (x); \qquad x (0) = x_0,
\elabel{nonlinear-ss}
\end{equation}where $x(t)$ evolves in $\Re^n$.  We assume throughout this chapter that  the function
$f\colon\Re^n\to\Re^n$ has continuous first and second partial derivatives  (denoted $f\in
C^2(\Re^n)$).   The solution is denoted by
$x(t)$, or
$\varphi (t; x_0)$ when we wish to emphasize the dependence of the state on its initial
condition.    Recall that
$x_e$ is an
\defn{equilibrium} of \eq nonlinear-ss/ if $f (x_e) =\zero $, in which case $\varphi (t; x_e) =
x_e$ for all $t$.  Our most basic, and also the weakest form of stability considered in
this book is the following:
The equilibrium $x_e$ is \defn{stable in the sense of Lyapunov} if for all
$\epsilon > 0$, there exists $\delta > 0$ such that  if  $\| x_0 - x_e \| < \delta$,
then
\begin{eqnarray*}
\| \varphi (t; x_0) - \varphi (t; x_e) \| < \epsilon
\qquad \hbox{for all $t \geq 0$.}
\end{eqnarray*} 
In words, if an initial condition is close
to the equilibrium, then it will stay close forever.  An illustration is
provided in \Figure{Trajectory}.

\begin{figure}[ht] 
\ebox{.7}{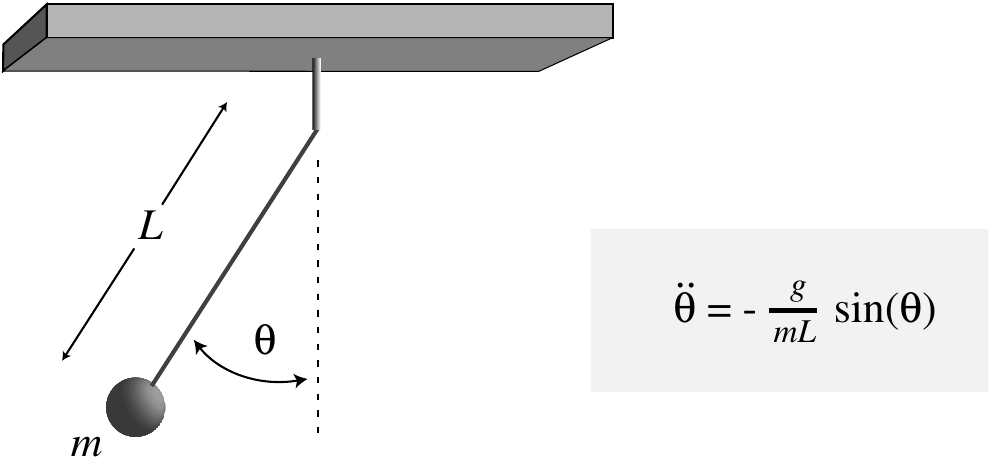}
\caption{Frictionless pendulum}
\flabel{Pendulum}
\end{figure} 

\begin{ex}

Consider the frictionless pendulum illustrated in \Figure{Pendulum}.
Letting $x_1 = \theta$, $x_2 = \dot\theta $ we obtain the state space model
\begin{eqnarray*}
\left(\begin{matrix}\dot x_1 \\ \dot x_2\end{matrix}\right) = 
\left(\begin{matrix}x_2\\ - \tfrac{g}{  m L}\sin (x_1)   \end{matrix}\right).
\end{eqnarray*}
In \Figure{PendFlow} some sample trajectories of $x(t)$ are illustrated.   As
\Figure{PendFlow} shows, the equilibrium $x_e =\left(\begin{smallmatrix}\pi \\
0\end{smallmatrix}\right)$ is not stable in any sense, which agrees with physical
intuition.  Trajectories which begin near the equilibrium
$x_e=\zero$ will remain near this equilibrium thereafter. Consequently, the origin is stable
in the sense of Lyapunov.  \end{ex}

Although the equilibrium $x_e=\zero$ is stable, there is no reason to
expect that   trajectories $x(t)$ will converge to $x_e$.  Since it is frequently a goal in
control design to make the model fall into a  specific regime starting from any initial
condition, we often require the following stronger forms of stability.
\begin{description}
\item[(a)]
An equilibrium $x_e$ is said to be \defn{asymptotically stable} if (i)
$x_e$ is stable in the sense of Lyapunov; and (ii) for  some $\delta_0 > 0$,
whenever $\| x_0 - x_e \| < \delta_0 $,
\begin{eqnarray*}
\lim_{t \to \infty} \varphi (t; x_0) = x_e.
\end{eqnarray*}

\item[(b)]
An equilibrium $x_e$ is \defn{globally asymptotically stable} if  it is
asymptotically stable with $\delta_0=\infty$.  Hence $x(t)\to x_e$ from any
initial condition.
\end{description}

\begin{figure}[ht]
\ebox{.8}{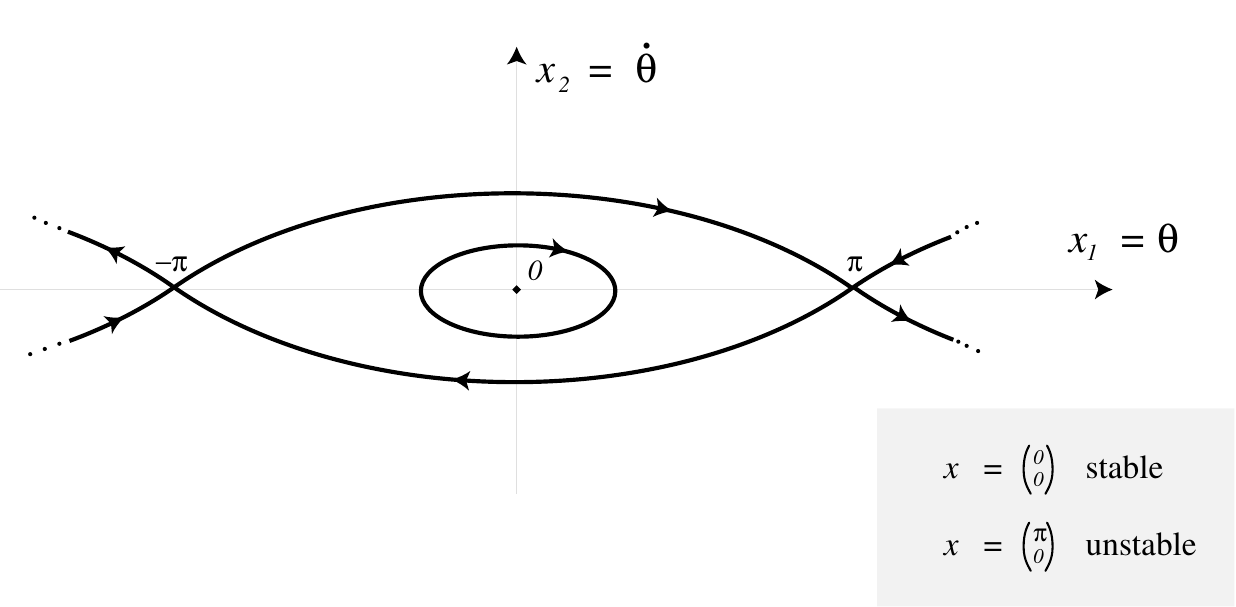}
\caption{Stable and unstable equilibria.}
\flabel{PendFlow}
\end{figure}

Consider now the case where  $f$ is linear, so that $\dot x = Ax$. When is this LTI
model asymptotically stable? First note that by the
\textit{homogeneity} property of linear systems, asymptotic   stability and global
asymptotic  stability are equivalent concepts.  In the linear case,  $x_e$ is an
equilibrium if and only if $ x_e\in {\cal N}(A)$.  Hence the origin is always an
equilibrium.  For asymptotic stability we must in fact have $\clN (A) = \{0\}$ since, by
global asymptotic stability, we cannot have more than one equilibrium.  Unstable modes are
also ruled out. The following result can be demonstrated using these facts:

\begin{theorem}
\tlabel{a-stable}
An LTI model is asymptotically stable if and only if ${\rm Re\,}
(\lambda) < 0$ for every eigenvalue $\lambda$ of $A$. 
\qed
\end{theorem}

A matrix $A$ satisfying the conditions of \Theorem{a-stable} is said to be \defn{Hurwitz}.
 
Stability in the sense of Lyapunov is more subtle. For example,
consider the two matrices
\begin{eqnarray*}
A_1= \left(\begin{matrix} 0 &  0\\ 0 &  0\end{matrix}\right)
\qquad
A_2= \left(\begin{matrix} 0 &  1\\ 0 & 0\end{matrix}\right) 
\end{eqnarray*}
The eigenvalues of both matrices are the same: $\lambda_1 = \lambda_2
= 0$. In the first case,
\begin{eqnarray*}
(sI-A_1)^{-1} = \left[ \begin{array}{rr} 1/s & 0 \\ 0 & 1/s
\end{array},
 \right] 
\end{eqnarray*}
so that the state transition matrix becomes $ \phi_1 (t) = I$.  Since state trajectories are
constant, it follows that every state is an equilibrium, and that every equilibrium is stable
in the sense of Lyapunov.

However for $A_2$,
\begin{eqnarray*}
(sI - A_2)^{-1} = \left[ \begin{array}{rr} 1/s & 1/s^2 \\ 0 & 1/s
\end{array} \right] 
\end{eqnarray*}
Therefore,
\begin{eqnarray*}
\phi_2 (t) = \left[ \begin{array}{rr} 1 & t \\ 0 & 1 \end{array} \right]
\end{eqnarray*}
which is growing with time.  This implies that $x(t)\to\infty$ for certain initial
conditions, and so this system is not stable  in any sense.

Thus, the eigenvalues alone do not foretell stability in the sense of
Lyapunov.  Looking at the eigenvectors, we see in the first instance from
the equation
\begin{eqnarray*}
\left[ \begin{array}{rr} 0 & 0 \\ 0 & 0 \end{array} \right] \left[
\begin{array}{c} x_1 \\ x_2 \end{array} \right] = 0 \left[
\begin{array}{c} x_1 \\ x_2 \end{array} \right]
\end{eqnarray*}
that we may take as eigenvectors
\begin{eqnarray*}
x^1= \left[ \begin{array}{c} 1 \\ 0 \end{array} \right] \qquad
x^2 = \left[ \begin{array}{c} 0 \\ 1 \end{array} \right] 
\end{eqnarray*}
However, in the second case, the eigenvector equation becomes
\begin{eqnarray*}
\left[ \begin{array}{rr} 0 & 1 \\ 0 & 0 \end{array} \right] \left[
\begin{array}{c} x_1 \\ x_2 \end{array} \right] = 0 \left[
\begin{array}{c} x_1 \\ x_2 \end{array} \right]
\end{eqnarray*}
This has a solution
\begin{eqnarray*}
x^1 = \left[ \begin{array}{c} 1 \\ 0 \end{array} \right],
\end{eqnarray*}
but in this example it is not possible to  choose another linearly independent
eigenvector.  This is a special case of the more general result:
 
 \begin{theorem}
The LTI model is stable in the sense of Lyapunov if and only if
\balphlist
\item
${\rm Re\,} (\lambda) \leq 0$ for every eigenvalue $\lambda$.
\item 
If ${\rm Re\,} (\lambda) = 0$, where $\lambda$ is an eigenvalue of multiplicity
$m$, then there exist $m$ linearly independent corresponding eigenvectors.
\end{list}
\end{theorem}

\proof
This follows from Theorem~8-14 of \CHE.
\qed

\section{Lyapunov's direct method}
\label{s:LyapunovDirect}

Lyapunov's direct method is a general approach to verifying stability for linear
or nonlinear models. This method allows us to determine stability without
solving an ordinary differential equation.  The idea is to search for a
``bowl-shaped'' function $V$ on the state space
$\Re^n$. If $V(x(t))$ is decreasing in time, then  we can qualitatively describe
the behavior of the trajectory $x(t)$, and   determine if the model is stable
in the sense of Lyapunov.
 
\begin{figure}[ht]
  \ebox{.8}{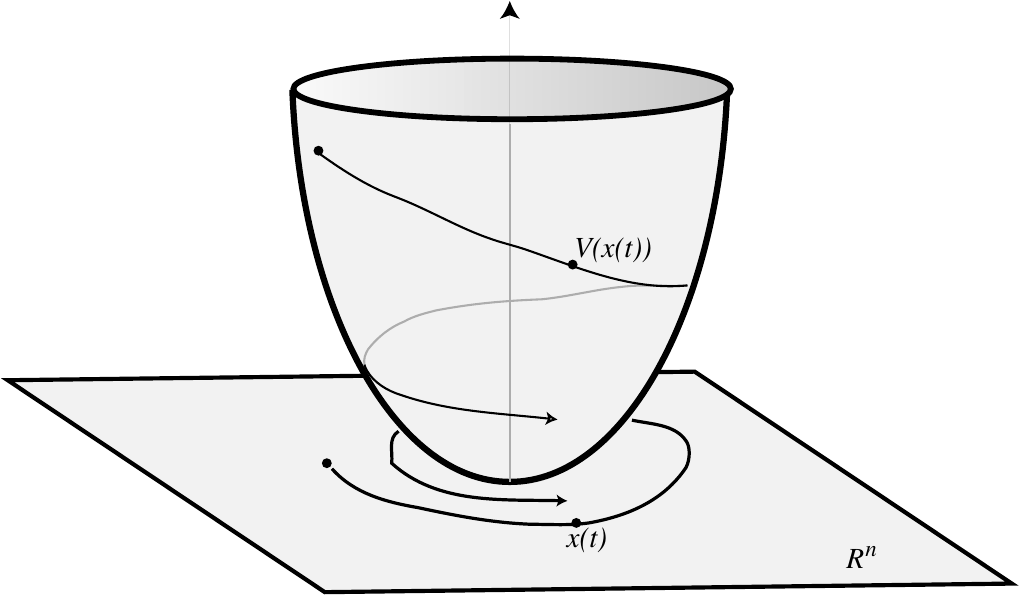}
\caption[$V (x(t))$ represents the height of the `lifted' trajectory]{$V (x(t))$ represents the height of the `lifted' trajectory from the state space
$\Re^n$}
 \flabel{Bowl} 
\ebox{.8}{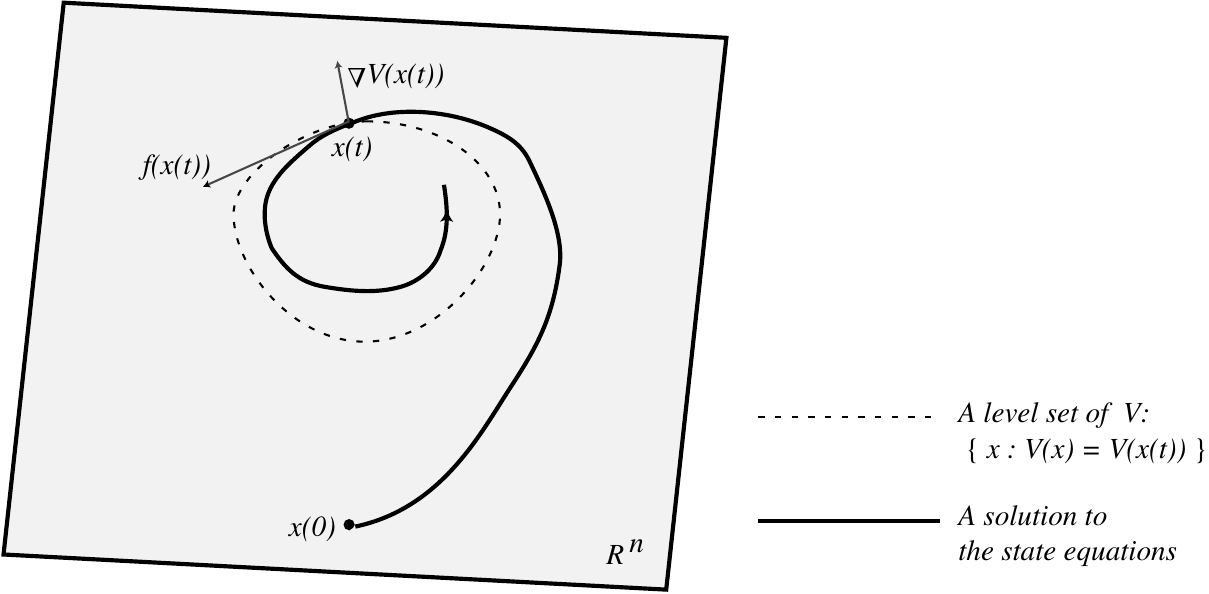}   
\caption{$V(x(t))$ is decreasing with time if $ \nabla V(x)\, f(x) \le 0$}
\flabel{Level} 
\end{figure}

To capture what we mean by ``bowl-shaped'', we give the following definition.
A scalar-valued function $V : \Re^n \to \Re$ with continuous partial derivatives
is called \defn{positive definite} on an (open) region $\Omega$ about the
equilibrium $x_e$ if  $V (x_e) = 0$, and $V (x) > 0$ for all $x \in \Omega$, 
$x\neq x_e$. 
 
\begin{theorem}
\tlabel{LyapunovStabilityCriteria}
Suppose $x_e\in\Re^n$ is an equilibrium, and that $\Omega\subset\Re^n$ is an open set
containing $x_e$.
\balphlist
\item 
$x_e$ is stable in the sense of Lyapunov if there exists a positive
definite function $V$ such that 
\begin{eqnarray*}
\frac{d}{dt}  V (x(t)) \leq 0 \qquad \mbox{ whenever } x (t) \in \Omega
\end{eqnarray*}

\item 
$x_e$ is asymptotically stable if there exists a positive definite function $V$
such that
\begin{eqnarray*}
\frac{d}{dt}  V (x(t)) < 0 \qquad \mbox{ whenever } x (t) \in \Omega,\ x (t)
\neq x_e
\end{eqnarray*} 
\item 
$x_e$ is globally asymptotically stable if there exists a positive definite function
$V$ on $\Re^n$ such that  $V(x) \to \infty$ as $x \to \infty$, and 
\begin{eqnarray*}
\frac{d}{dt}  V (x(t)) < 0\qquad \mbox{ for  all } x (t), \; x (t) \neq  x_e. 
\end{eqnarray*}

\end{list} 
\end{theorem}

\proof
See for example \cite{laslef61}.
\qed

\clearpage

Using the formula 
\begin{eqnarray*}
\frac{d}{dt}  V(x(t)) = \nabla V (x (t))\, f (x(t)),
\end{eqnarray*}
the condition for asymptotic stability becomes 
\begin{eqnarray*}
\nabla V (x) \, f (x) < 0\quad  \mbox{ for } x \in \Omega, x \neq x_e.
\end{eqnarray*}
This is illustrated geometrically in \Figure{Level}.

As an application of \Theorem{LyapunovStabilityCriteria} consider the  Van der Pol
oscillator described by the pair of equations
\begin{equation}
\begin{array}{rcl}
\dot x_1 & = & x_2 \\
\dot x_2 & = & - (1 - x_1^2) x_2 - x_1.
\end{array}
\elabel{VanDerPol}
\end{equation}
We first search for equilibria:  If $f(x) = \zero$, it follows that    $x_2 = 0 $ and  $-(1 -
x_1^2) x_2 - x_1 = 0$.  The unique equilibrium is therefore $x_e = \zero$.

To see if $x_e$ is stable in any sense, try the Lyapunov function 
$V (x) = \half | x |^2$.  We then have
\begin{eqnarray*}
 \nabla V (x) = x^T , \qquad f(x) = \left(\begin{matrix} x_2 \\ - (1 - x_1^2) x_2
- x_1\end{matrix}\right)
\end{eqnarray*}
and so
\begin{eqnarray*}
\nabla V (x)\, f (x) = x_1 x_2 - (1 - x_1^2) x_2^2 - x_1 x_2
= - (1 - x_1^2) x_2^2 \le  0  \mbox{\quad If    $ | x_1 | < 1$.}
\end{eqnarray*}
The set $\Omega = \{ x \in \Re^2 : | x_1 | < 1 \}$ contains $x_e$ and is open.
Hence, $ x_e$ is stable in the sense of Lyapunov.

\begin{figure}[ht]
\ebox{.9}{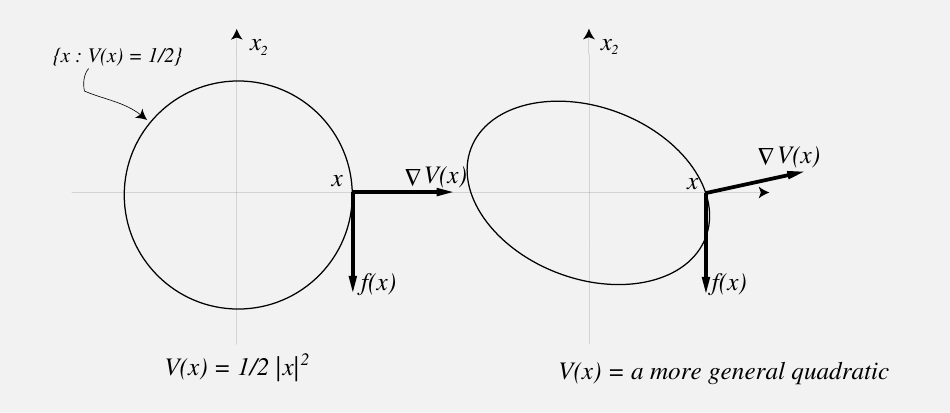}
\caption[The vectors $\nabla V(x) $ and $f(x)$ meet at ninety degrees]{The vectors $\nabla V(x) $ and $f(x)$ meet at ninety degrees
when $x=(\gamma, \ 0 )^T$, and $V(x) =|x|^2$.  By distorting $V$, and
hence also its level set $\{x : V(x) = 1/2\}$, one can increase this
angle beyond ninety degrees.}  \flabel{CircLevel}
\end{figure}

Is the origin asymptotically stable?  Considering \Figure{CircLevel}, it is
apparent that we need to be more careful in choosing
$V (x)$. To obtain elliptical level sets of the form shown in the second
illustration of \Figure{CircLevel}, set
\begin{eqnarray*}
V (x) = \half  (x_1^2 + x_2^2) + \epsilon x_1 x_2
\end{eqnarray*}
where $\epsilon$ is small.  If $| \epsilon | < 1$,   the function $V$ is then
positive definite. The derivative is now
\begin{eqnarray*}
\nabla V (x) = \left(\begin{matrix} x_1 + \epsilon x_2 
\\ x_2+ \epsilon x_1\end{matrix}\right)^T
\end{eqnarray*}
and the equation to be analyzed is
\begin{eqnarray*}
\nabla V (x) \, f (x) & = & x_1 x_2 + \epsilon x_2^2  - \left[
(1 - x_1^2) x_2^2 + x_1 x_2 + \epsilon (1 - x_1^2) x_1 x_2 + \epsilon
x_1^2 \right]\\ 
& = & - \left[ (1 - x_1^2) - \epsilon \right] x_2^2  -
\epsilon (1 - x_1^2) x_1 x_2 -  \epsilon x_1^2.
\end{eqnarray*}
If we take $x_1^2<\half$ this then   gives  
\begin{eqnarray*}
\nabla V (x) \, f (x) \leq - \left( \half  - \epsilon \right)
x_2^2 + \epsilon | x_1 x_2 | - \epsilon x_1^2
\end{eqnarray*}
Now use the bound  $|x_1 x_2| \leq \half  (x_1^2 + x_2^2)$:
\begin{eqnarray*}
\nabla V (x) \, f (x) & \leq & - \left( \half  - \epsilon
\right) x_2^2 + \half  \epsilon (x_1^2 + x_2^2) - \epsilon x_1^2
\\ & = & - \left( \half  - \epsilon - \half  \epsilon
\right) x_2^2 -  \frac{\epsilon }{ 2} x_1^2.
\end{eqnarray*}
Finally, take $\epsilon =  \frac{1 }{ 6}$.  Then for any $x\neq \zero$,
\begin{eqnarray*}
\nabla V (x) \, f (x) & \leq & - \left( \half  -  \frac{1 }{ 6} -
 \frac{1 }{ 12} \right) x_2^2 - \frac {1 }{ 12} x_1^2 \\
& = & -  \frac{1 }{ 4} x_2^2 -  \frac{1 }{ 12} x_1^2 < 0.
\end{eqnarray*}
With $\Omega = \{ x : x_1^2 < \half  \}$, we see that
\begin{eqnarray*}
\nabla V (x) \, f (x) < 0 \mbox{ for } x \in \Omega, x \neq x_e = 0.
\end{eqnarray*}
It follows that the equilibrium
$ x_e$ is asymptotically stable.

\section{Region of asymptotic stability}
 \label{s:chap2.2.6}

Let $V (x)$ be a positive-valued function on $\Re^n$ with continuous first partial
derivatives, and let $\Omega_\ell$ designate the  open region within
which $V (x) $ is less than $ \ell$:
\begin{eqnarray*}
\Omega_\ell \eqdef \{ x : V (x) <  \ell\}.
\end{eqnarray*}
Suppose that, within $\Omega_\ell$, $V(x(t))$ is always decreasing.  This will be the case  
if and only if
\begin{equation}
\nabla V (x)\, f(x) \le 0 \qquad \hbox{for all } x \in \Omega_\ell. 
\elabel{omega-invariant}
\end{equation}Since under \eq omega-invariant/ we have $V(x(t))< \ell$ for all $t$ if  $V(x(0))< \ell$,
it follows that the
set $\Omega_\ell$ is \defn{invariant}:  if
$x_0\in\Omega_\ell$, then
$x(t)\in\Omega_\ell$ for all $t>0$.  Under a stronger set of conditions, the set
$\Omega_\ell$ becomes a region of asymptotic stability.

\begin{theorem}
\tlabel{region-stability}
Suppose that $V\colon\Re^n\to \Re$ is $C^1$, and that the nonlinear model \eq nonlinear-ss/
satisfies for some $\ell>0$, 
\begin{eqnarray}
 &&\hbox{The set $\Omega_\ell$ is bounded.}
\elabel{omega-invariant0}\\
V (x) > 0     &&\qquad \hbox{for all } x \neq \zero, x \in \Omega_\ell.
\elabel{omega-invariant2}\\
\nabla V (x)\, f(x) < 0 &&\qquad \hbox{for all } x \neq \zero, x \in \Omega_\ell.
\elabel{omega-decreasing}
\end{eqnarray} 
Then $x_e =  \zero $ is an asymptotically stable equilibrium, and every solution beginning in
$\Omega_\ell$ tends to $\zero$ as $t \to \infty$.  That is, $\Omega_\ell$ is a region
of asymptotic stability. 
\end{theorem}
 
\proof Since $V(x(t))$ is decreasing in $t$ and is bounded from below,
it must have a limit $v_0 \geq 0$.  If $v_0=0$, then by \eq
omega-invariant2/ it follows that $x(t)\to \zero$, as claimed.  To
prove the theorem we will show that it is impossible to have $v_0>0$.

We proceed by contradiction:  If $v_0>0$ then for all $t$, $x(t)$ lies in the
closed and bounded set
\begin{eqnarray*}
S=\{ x : v_0\le V(x)\le \ell  \}.
\end{eqnarray*}
Since $\nabla V(x)\, f(x)$ is strictly negative on $S$, by continuity it follows that
for some
$\epsy>0$,
\begin{eqnarray*}
\nabla V(x)\, f(x) \le -\epsy ,\qquad x\in S.
\end{eqnarray*}
Hence, from this bound, the bound $V(x(0))\le \ell$,  and the fundamental theorem of
calculus, for any
$T>0$,
\begin{eqnarray*}
-\ell\le V(x(t_f))-V(x(0)) = \int_0^{t_f} \nabla V(x(t))\, f(x(t))\, dt \le -\epsy t_f.
\end{eqnarray*}
Letting $t_f \to\infty$, we obtain a contradiction, and this proves the theorem.
\qed
   
\begin{ex}

Consider the nonlinear state space model
\begin{eqnarray*}
\dot x_1 & = & - x_1- x_2 \\
\dot x_2 & = & x_1 - x_2 +   x_2^3
\end{eqnarray*} 
With the Lyapunov function $V (x) =  x_1^2 + x_2^2$,
we have  for all $x \neq \zero$
\begin{eqnarray*}
\nabla V(x)\, f(x) &=& 2 (x_1,x_2)\cdot(- x_1- x_2 , x_1 - x_2 +   x_2^3)
\\
&=&  
-2[x_1^2 + x_2^2(1-x_2^2) ]
\end{eqnarray*}
which will be negative  for $x \neq \zero$ so long as $x_2^2 < 1$.  Considering
\Figure{region-stability}, it follows from \Theorem{region-stability} that $\Omega_1=\{ x :
x_1^2 + x_2^2 <1\}$ is a region of asymptotic stability.
\end{ex}

\begin{figure}[ht]
\ebox{.95}{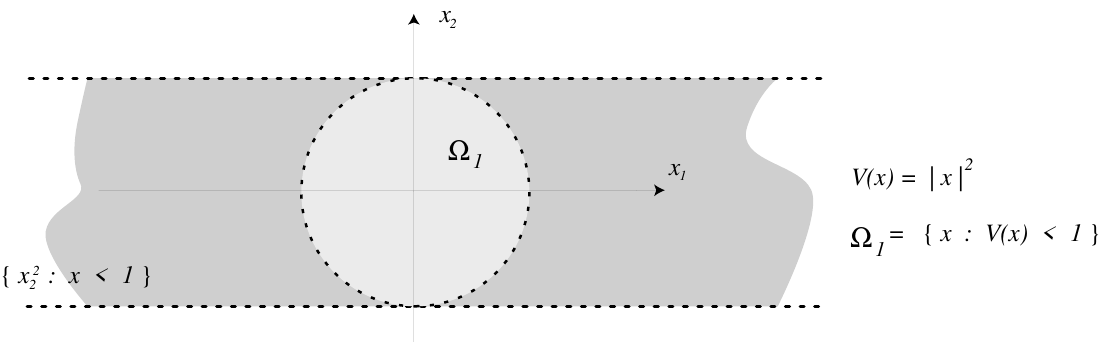}
\caption[The sub-level set of $V$]{The sub-level set of $V$ shown as a grey disk in the figure is a region of asymptotic
stability}
\flabel{region-stability} 
\end{figure}

\section{Stability of linear state space models}

We have already seen in \Theorem{a-stable} that to determine asymptotic stability of a linear
state space model, it is sufficient to compute the eigenvalues of an $n\times n$ matrix. 
By specializing \Theorem{LyapunovStabilityCriteria} to the LTI case, we derive here an
algebraic characterization of stability which is simpler to apply than the
eigenvalue test. Before we proceed,  we  need a better understanding of positive definite
\textit{quadratic} functions.  An  $n\times n$  matrix $M$ is called
\defn{Hermitian} if $M = M^*$ (complex-conjugate transpose).  A matrix $M$ is
said to be \defn{positive definite} if $M$ is Hermitian,  and
\begin{eqnarray*}
x^* M x > 0 \qquad \hbox{for all } x\neq \zero, x\in\Co.
\end{eqnarray*}
If the strict inequality above is replaced by $ x^* M x \ge 0$, then the matrix $M$ is
called \defn{positive semi-definite}.  We adopt the notation $M>0$,  $M\ge 0$,
respectively.  If $V(x) = x^T M x$ for a real matrix $M>0$, then the function $V$ is positive
definite on $\Re^n$.

There are several tests for positive definiteness which avoid computation of  $ x^* M x$ for
every $x\in\Re^n$.  Suppose for instance that $M$ is Hermitian, and is partitioned 
as shown below:
\ebox{.95}{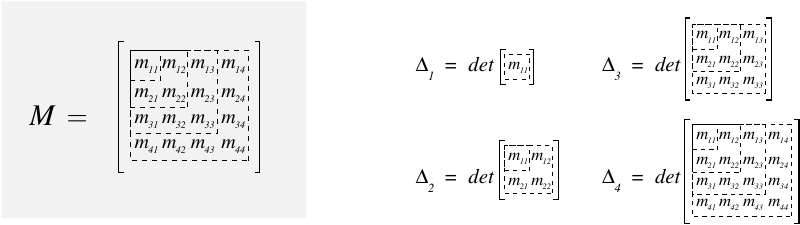}
In this $4\times 4$ example, we have taken four submatrices found on the diagonal, and
indicated  their respective determinants, $\Delta_1,\Delta_2,\Delta_3,\Delta_4 $.  These are
called the (leading) \defn{principle minors} of $M$.  One of the most readily verified tests
for positivity is   given in the following theorem.

 \begin{theorem}[Sylvester-Principal Minors Test]
\tlabel{sylvester}
A Hermitian matrix $M$ is positive definite if and only if each of its leading principle
minors is strictly positive.
\end{theorem}
\proof See \cite{bel70}.  \qed

Consider now the LTI model
\begin{eqnarray*}
\dot x = Ax
\end{eqnarray*}
Choose $V (x) = x^T P x$, where we can assume without loss of generality that
$P$ is real and   Hermitian  $(P = P^T=P^*)$.  If  $P $ is a positive definite matrix, then
$V$ is also positive definite.

The gradient $\nabla V (x)$, which is a row vector, may be computed using the product
rule,
\begin{eqnarray*}
\nabla V (x)^T = Px + P^T x = 2 P x.
\end{eqnarray*}

Also, $f (x) = Ax$, so the condition for asymptotic stability is
\begin{eqnarray*}
\nabla V (x) \, f (x) = 2 x^T P^T A x < 0, \quad x \neq 0.
\end{eqnarray*}
Since $x^T P^T Ax = x^T A^T Px$, we have that
$2 x^T P^T Ax =   x^T [PA + A^T P]x$.  
Thus, we need to find $P > 0$ such that the \defn{Lyapunov equation} 
 \begin{equation}
A^T P + PA = - Q 
\elabel{LyapunovEqn}
\end{equation}
is satisfied for some positive definite matrix $Q$. This approach
characterizes stability for LTI models:

\begin{theorem}
\tlabel{LyapunovEquation}
The LTI model
$\dot x = Ax$ is asymptotically stable if and only if for some $Q > 0$, there exists
a solution $P > 0$ to the Lyapunov equation. 
\end{theorem}

\proof
If for \textit{one} $Q > 0$ we can find a $P > 0$ which solves the
Lyapunov equation \eq LyapunovEqn/, then we have seen that that $V (x) = x^T Px$ is a
valid Lyapunov function.  Hence the model is globally asymptotically stable by
\Theorem{LyapunovStabilityCriteria}.

Conversely, if the model is asymptotically stable, then ${\rm Re\,} (\lambda) < 0$ for
any $\lambda$, so $e^{At} \to 0$ as $t\to\infty$ at an exponential rate.  Let
\begin{eqnarray*}
P =
\int_0^\infty e^{A^Tt} Q e^{At}\, dt
\end{eqnarray*}
where $Q > 0$ is given.  Given the form of the state transition matrix $\phi$, it follows
that the matrix $P$ has the following interpretation:
\begin{eqnarray*} x^T P x = \int_0^\infty x (t)^T Q x (t)\, dt
\qquad \hbox{when $x (0) = x$.}
\end{eqnarray*}

To see that $P$ solves the Lyapunov equation, let $F (t) =e^{A^T t} Q e^{At}$, and apply the
product rule for differentiation:
\begin{eqnarray*}
\frac{d}{dt}  F (t)= A^T e^{A^T t} Qe^{At} + e^{A^T t} Qe^{At} A
=   A^T F (t) + F (t) A.
\end{eqnarray*}
So, by the Fundamental Theorem of Calculus,
\begin{eqnarray*}
F (s) - F (0)= \int_0^s \frac{d}{dt}  F (t)\,  d t =
A^T \Bigl(  \int_0^s F (t) dt\Bigr) + \Bigl(\int_0^s F (t)\, dt\Bigr) A.
\end{eqnarray*}
We have $F(0) = Q$, and letting $s \to \infty$ we have $F (s) \to 0$, and $\int_0^s F (t) dt
\to P$. In the limit we thus see that the identity $-Q=  A^T P + PA$ holds, so the
Lyapunov equation is satisfied.
\qed

\begin{ex}
Consider the two-dimensional model
\begin{eqnarray*}
\dot x = A x =  \left[ \begin{array}{cc} 0 & 1 \\ -2 & -2 \end{array} \right] x.
\end{eqnarray*}
To see if this model is stable we will attempt to solve the Lyapunov
equation with $Q =\left(\begin{smallmatrix}1 & 0 \\ 0 & 1 \end{smallmatrix}\right)$:
\begin{eqnarray*}
\left[ \begin{array}{cc} 0 & -2 \\ 1 & -2 \end{array} \right]
\left[ \begin{array}{cc} p_{11} & p_{12} \\ p_{12} & p_{22}
\end{array} \right] + 
\left[ \begin{array}{cc} p_{11} & p_{12} \\ p_{12} & p_{22}
\end{array} \right]
\left[ \begin{array}{cc} 0 & 1 \\ -2 & -2 \end{array} \right] = \left[
\begin{array}{cc} -1 & 0 \\ 0 & -1 \end{array} \right]  
\end{eqnarray*}
Solving the resulting three equations gives:
\begin{eqnarray*}
-2 p_{11} - 2 p_{12} & = & - 1 \\
- 2 p_{22} + p_{11} - 2 p_{12} & = & 0 \\
2 p_{12} - 4 p_{22} & = & - 1,
\end{eqnarray*}
which has the solution  $P = \left[ \begin{array}{rr} 5/4 & 1/4 \\ 1/4 & 3/8 \end{array}
\right]$.

To see if the system is stable, we must now check to see if the matrix $P$ is positive
definite.  The principle minors can be computed as
\begin{eqnarray*}
 5/4 > 0,   \quad (5/4)(3/8) - (1/4)(1/4) > 0
\end{eqnarray*}
Hence the Sylvester-Principal Minors  test is positive.  From this we can conclude that
$P > 0$,  and that the model is globally asymptotically stable. 
\end{ex}

\section{Stability subspaces}
\label{s:chap2.3.1} 

The \defn{stable subspace} $\Sigma_s$ of $\dot x = Ax$ is the set of all $x
\in \Re^n$ such that
\begin{eqnarray*}
e^{At} x \to 0 \mbox{ as } t \to \infty.
\end{eqnarray*} 
In the special case where the matrix $A$ has distinct eigenvalues, the stable subspace 
is easily characterized using the modal form.  Letting $\{\lambda_i\}$ denote the distinct
eigenvalues, and $\{ v^i\}$ the corresponding eigenvectors, we may write any solution
$x$ as
\begin{eqnarray*}
x (t) = \sum_{i=1}^n \alpha_i e^{\lambda_i t} v^i  
\end{eqnarray*}
where $\{\alpha _i\}$ is the representation of $x(0)$  with respect to the basis $\{ v^i\}$.
From this equation it is evident that
\begin{eqnarray*}
\Sigma_s = \hbox{Span}\, \{ v^i : {\rm Re}(\lambda_i) < 0 \}.
\end{eqnarray*}
In this case we may also define an unstable subspace as
\begin{eqnarray*}
\Sigma_u = \hbox{Span}\, \{ v^i : {\rm Re}(\lambda_i) \ge 0 \}.
\end{eqnarray*}

Strictly speaking, if the stable subspace is to be real, then these formulae will not
be valid if there are complex eigenvalues.  Consider for example the three-dimensional
model with
\begin{eqnarray*}
A = \left[ \begin{array}{rrr} -2 & 1 & -1 \\ -2 & -5 & 6 \\ -1 & -3 & 4 \end{array} \right]
\end{eqnarray*}
The eigenvalues of $A$ are distinct, $\lambda = 1, -2 \pm j $, and hence a modal form of the
state equations is possible with a complex matrix $\Lambda$.

The unstable subspace is the span of the eigenvector $v^1$ corresponding to
$\lambda_1=1$:
\begin{eqnarray*}
\Sigma_u : \{x: (\lambda_1 I - A) x = 0\},
\end{eqnarray*}
from which we compute $x_1 = 0$, $x_2 = x_3$.  Therefore,
\begin{eqnarray*}
\Sigma_u = \hbox{Span}\, \left\{ \left[ \begin{array}{r} 0 \\ 1 \\ 1 \end{array}
\right]\right\}
\end{eqnarray*}

The following vectors are eigenvectors corresponding to $\lambda_2$, $\lambda_3$:
\begin{eqnarray*}
v^2 = \left[ \begin{array}{r} 1 \\ 2j  \\ j  \end{array} \right],
\qquad 
v^3 = \left[ \begin{array}{r} 1 \\ -2j  \\ -j  \end{array} \right].
\end{eqnarray*}
Since the stable subspace is real, it follows that
\begin{eqnarray*}
\Sigma_s : x = k  \left[ \begin{array}{r} 1 \\ j2 \\ j1 \end{array}
\right] + \bar k \left[ \begin{array}{r} 1 \\ -j2 \\ -j1 \end{array}
\right], \qquad k \in \Co.
\end{eqnarray*}
where $\bar k$ denotes the complex conjugate of $k$.
The formula for the stable subspace may be equivalently expressed as 
\begin{eqnarray*}
\Sigma_s = \hbox{Span}\, \left\{ \left[ \begin{array}{r} 1 \\ 0 \\ 0 \end{array}
\right], \left[ \begin{array}{r} 0 \\ 2 \\ 1 \end{array} \right] \right\}.
\end{eqnarray*}

\section{Linearization of nonlinear models and stability}

We now return to the nonlinear case, and show that the  simply stated characterization of
stability given in \Theorem{LyapunovEquation} for linear models has remarkably strong
implications to general nonlinear state equations.   Consider again the nonlinear model \eq
nonlinear-ss/ with $x_e = \zero$.    With $A = \frac{d}{dx} f \, (\zero)$,   the linearized
model may be expressed
 \begin{equation}
\elabel{linearizedStability}
\dot{\delta  x} = A\, \delta x.
\end{equation}
\begin{theorem}
\tlabel{linearizedStability} 
For the nonlinear state space model with  $x_e = \zero$,
\balphlist
\item 
If $x_e$ is asymptotically stable for \eq linearizedStability/, then $x_e$ is also
asymptotically stable for  \eq nonlinear-ss/.
\item 
If $A$ has an eigenvalue with ${\rm Re\,}( \lambda) > 0$, then $x_e$ is {\it not}
  stable in the sense of Lyapunov for  \eq nonlinear-ss/.
\item 
If ${\rm Re\,} (\lambda) \leq 0$ for all $\lambda$, but ${\rm Re\,} (\lambda) = 0$
for some $\lambda$, nothing can be said about the stability of the nonlinear model.
\end{list}  
\end{theorem}

\proof (sketch)
For (a), the proof proceeds by constructing a Lyapunov function for the linear model, and
then applying this to the nonlinear model.  The result then follows from
\Theorem{LyapunovStabilityCriteria}. 

 The proof of (b) is similar:  If the linear model is
strictly unstable in this sense, then there exists an initial condition $x_0$ for which the
corresponding solution to the linearization tends to infinity.  One may then construct  
$\epsy>0,\delta>0$ such that the solution to the nonlinear equations with initial condition
$x(0) = c x_0$ will eventually leave the ball  $B_\epsy(\zero)$ centered at the origin with
radius $\epsy$, whenever $c < \delta$. This follows from the Taylor series expansion. 
Hence, the origin cannot be stable in the sense of Lyapunov.

For (c), there is nothing to prove, but we may give an example.  Consider the two state
space models on $\Re^1$:
\begin{eqnarray*}
\hbox{(a)}\quad
\dot x = -x^3
\qquad
\hbox{(b)}\quad
\dot x = +x^3
\end{eqnarray*}
In the first instance, using the Lyapunov function $V(x) = x^2$ we see that 
$\nabla V(x)\, f_1(x) = -2 x^4$, and we conclude that the model is globally
asymptotically stable.  In the second case, we have $\nabla V(x)\, f_2(x) = +2 x^4$,
which shows that trajectories explode to infinity from any initial condition. 
Although the behavior of these two models is very different, both have the same
linearization with $A= 0$.
\qed

We can now revisit the Van der Pol equations defined in \eq VanDerPol/.  The
corresponding matrix $A$ can be computed as
\begin{eqnarray*}
A =\left.  \frac{d}{dx} f (x) \right|_{x=\zero} = \left[
\begin{array}{cc} 0 & 1 \\ -1 & -1 \end{array} \right]
\end{eqnarray*}
Taking determinants,
$
|Is - A| = \det \left(\begin{smallmatrix} s &  -1 \\ 1&  s+1 \end{smallmatrix}\right) =
s^2 + s + 1$, which gives $
\lambda_i =  \frac{-1 \pm \sqrt{1-4} }{ 2} = -\half  \pm  \frac{\sqrt{3 }}{ 2} j$.
So, the origin is asymptotically stable. Previously we found that the Lyapunov function
\begin{eqnarray*}
V (x) = \half  x_1^2 + \half  x_2^2 + \epsilon x_1 x_2
\end{eqnarray*}
is   effective   for the nonlinear model.  The function $V$   may
be written
\begin{eqnarray*}
V (x) = x^T Px \mbox{ where } P = \left[
\begin{array}{cc}
\half  & \half  \epsilon \\
\half  \epsilon & \half 
\end{array}
\right]
\end{eqnarray*}
This also works for the linearized model, as can be seen from the following Lyapunov
equation:
\begin{eqnarray*}
A^T P + PA & = & \half  \left( \begin{array}{cc} 0 & -1 \\ 1 & -1
\end{array} \right) \left( \begin{array}{cc} 1 & \epsilon \\ \epsilon
& 1 \end{array} \right) + \half  \left( \begin{array}{cc} 1 &
\epsilon \\ \epsilon & 1 \end{array} \right) \left( \begin{array}{cc}
0 & 1 \\ -1 & -1 \end{array} \right) \\
& = & \half  \left( \begin{array}{cc} - \epsilon & -1 \\ 1 -
\epsilon & \epsilon -1 \end{array} \right) + \half  \left(
\begin{array}{cc} - \epsilon & 1 - \epsilon \\ -1 & \epsilon - 1
\end{array} \right) \\
& = & -  \left( \begin{array}{cc}  \epsilon &   \half  \epsilon \\
  \half  \epsilon &  1 -\epsilon \end{array} \right) = - Q.
\end{eqnarray*}
Since the matrix $Q$ is positive definite for $\epsilon < 4/5$, we again see that the
origin is an asymptotically stable for  the linearized model. 

Although this is obviously a powerful method, it does not in general tell the complete
story about the stability of a nonlinear system.  In some applications  we are in the
situation of part (c) where the linearized model is not asymptotically stable, and in
this case the nonlinear equations must be tackled directly.  Another shortcoming of
this method is that when using a linearized model it is not possible to estimate a
region of asymptotic stability for the nonlinear model.

\begin{figure}[ht] 
\ebox{.7}{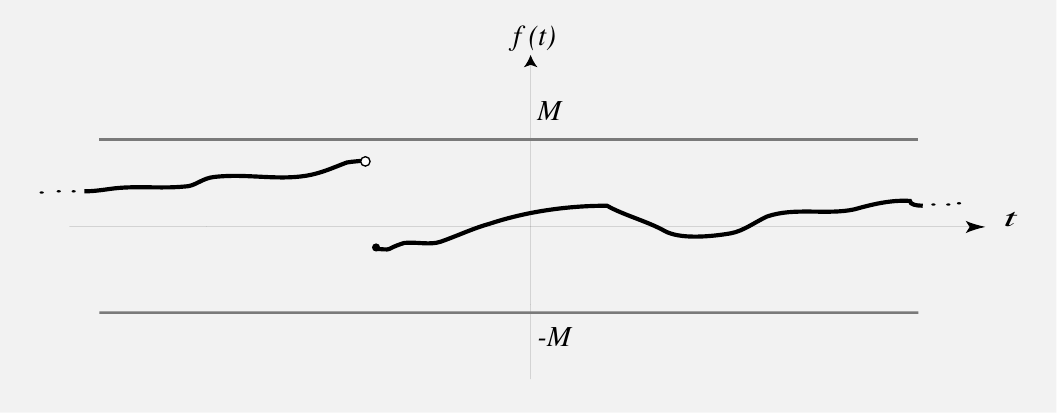}

\caption[The graph of $f$ remains in a tube of radius $M$]{A typical bounded function:  The graph of $f$ remains in a tube of radius $M$.}
\flabel{Bounded}

\end{figure}  
\section{Input-output stability}

A function $f : [0, \infty) \to \Re^n$ is called \defn{bounded}
if there exists a constant $M$ such that  $| f (t) | \leq M$ for all $t$.
A typical bounded function is illustrated in \Figure{Bounded}.  
A state space model with input
$u$ and output $y$ is  
\defn{bounded input/bounded output stable}  (BIBO stable) if for any $t_0$,
and any bounded input $u (t)$, $t\geq t_0$, the output is bounded on
$t_0 \leq t < \infty$ \textit{if } $x(t_0) = \zero$.

Consider for example the LTI model  whose state space and  transfer function descriptions
are expressed as
\begin{eqnarray*} 
\begin{array}{rcl}
\dot x & = & Ax + Bu\\
y & = & Cx + Du ,
\end{array}
\qquad
G (s) = C [ Is - A]^{-1} B + D.
\end{eqnarray*}
This model is BIBO stable if and only if each pole of every entry of
$G (s)$ has a \textit{strictly} negative real part.

For an LTI model, we see that asymptotic stability of $\dot x = A x$
implies BIBO stability. The converse is \textit{false} since pole-zero
cancellations may occur.  The connections between BIBO stability and
asymptotic stability may be strengthened, but first we must develop
controllability and observability for LTI models.

\begin{matlab} 
\item[EIG] finds eigenvalues and eigenvectors.
\item[ROOTS] finds the roots of a polynomial. 
\item[LYAP]   solves the   Lyapunov matrix equation \eq LyapunovEqn/. 
\end{matlab}
 
\begin{summary}
In this chapter we have introduced a powerful approach to
understanding the qualitative behavior of dynamical systems.  The main
tool is the \defn{Lyapunov function}, which is interpreted
geometrically as ``lifting'' the state trajectory onto the surface
defined by the Lyapunov function, and checking to see if the
trajectory travels ``down hill''.  Starting with these ideas, we then
developed \balphlist
\item
Criteria for stability and instability of nonlinear state space
models.
\item
Strong connections between stability of a nonlinear state space model,
and stability of its linearization.
\item
An algebraic characterization of stability for linear systems, based
on the matrix \defn{Lyapunov equation}.
\end{list}
There are many excellent texts on stability and the Lyapunov function
approach.  See in particular \cite{laslef61}. A historical development
of stability and the Lyapunov equation is contained in
\cite{kalber60}.
\end{summary}

\begin{exercises}
\item
Explain the behavior exhibited in \Figure{trajectories}.

\balphlist
\item
In this first simulation, the trajectories appear to converge to a
line within the two dimensional state space.  Show that this is the
case for the given initial conditions.  Will this occur for
\textit{every} nonzero initial condition?  Is the model stable in any
sense?

\item
Explain the behavior in simulation (b).  Will the trajectories
approximately follow a line from each initial condition?

\item
Explain the behavior in the last simulation.  Will the trajectories
cycle to zero from each nonzero initial condition?
\end{list}
 
\item
This exercise shows how to use state space methods to solve numerical
equations and optimization problems.  \balphlist

\item
Suppose that one wishes to compute numerically a solution $x^\star$ to
$f(x) = g(x)$, where $f,g\colon\Re^n\to\Re^n$ are $C^1$ functions.
Find conditions on $f$ and $g$ which ensure that solutions to the
differential equation
\begin{eqnarray*}
\dot x = A(f(x) - g(x))
\end{eqnarray*}
will converge to $x^\star$ for initial conditions near $x^\star$.  In
this equation $A$ is an $n\times n$ matrix which can be chosen by the
user. Experiment with the algorithm using \textit{Simulink} with
$f(x)=e^x$, and $g(x) = \sin(x)$, for a suitable scalar $A$.

\item
Suppose that one wants to minimize a $C^2$ function $V\colon\Re^n\to
\posRe$.  A necessary condition for a point $x^\star\in\Re^n$ to be a
minimum is that it be a \defn{stationary point}; that is, $\nabla V
(x^\star) =\zero$.
 
Consider the differential equation
\begin{eqnarray*}
\dot x = - P \nabla V (x)^T,
\end{eqnarray*}
where $P$ is a positive definite matrix.  Find conditions on the
function $V$ to ensure that a given stationary point $x^\star$ will be
asymptotically stable for this equation. \textit{Hint}: Find
conditions under which the function $V$ is a Lyapunov function.  Try
this algorithm out using \textit{Simulink} with a function of your
choice with $n=2$.  For example, try $V(x) = x_1^6 -2 x_1^2 +x_1 x_2 +
6 + |x_2|^3$.
\end{list}

\item
Consider the pendulum described in \Figure{Pendulum}, with the
addition of an applied torque $u$.  This is a nonlinear system which
is described by the differential equation $\ddot\theta = - \sin
(\theta) +u$.
 
Find a PD controller (state feedback) which stabilizes the system
about $\theta = \pi$. That is, for controlled state space model with
$x=\left(\begin{smallmatrix}\theta\\
\dot\theta\end{smallmatrix}\right) $, the equilibrium
$\left(\begin{smallmatrix}\pi\\ 0\end{smallmatrix}\right)$ is
asymptotically stable.

Is $\left(\begin{smallmatrix}\pi\\ 0\end{smallmatrix}\right)$ globally
asymptotically stable with your control law?

\item
Suppose that $V(x) = x^T A x$ for a real matrix $A$, and suppose that
the function $V$ is positive definite.  Prove that the matrix $M=
\half (A + A^T)$ is positive definite.  Can you find an example in
which the matrix $A$ is \textit{not} symmetric?

\item
Study the stability of the equilibrium state $x^e=0$
($x\eqdef(y\quad\dot y)^T$) for the nonlinear system
\begin{eqnarray*}
\ddot y+\dot y+y^3=0
\end{eqnarray*}
using the function
\begin{eqnarray*}
V(y,\dot y)=y^4+y^2+2y\dot y+2{\dot y}^2
\end{eqnarray*}
as a candidate Lyapunov function.
 
\item For the LTI system described by
\begin{eqnarray*}
\dot x=\left[\begin{matrix}-1&1\\ -2&3\\\end{matrix}\right]x\;,
\end{eqnarray*}
investigate asymptotic stability using the Lyapunov equation
\begin{eqnarray*}
A^TP+PA=-Q\qquad {\rm with\;}\;\; Q=I.
\end{eqnarray*}
Can you arrive at any definitive conclusion?

\item
Consider the linear model
\begin{eqnarray*}
\dot x = \left(\begin{matrix} 1 &1 \\ -1 - 2 \alpha + \alpha^2 &
-1 - \alpha + \alpha^2 \end{matrix}\right) x,
\end{eqnarray*}
where $\alpha$ is a scalar parameter.  For what values of $\alpha$ is
the origin (i) asymptotically stable?  (ii) stable in the sense of
Lyapunov?

\item
Consider the simple scalar model
\begin{eqnarray*}
\dot y = \alpha_0 u.
\end{eqnarray*}
We would like to use feedback of the form $u=-Ky + N r$ so that the
closed-loop pole is at $-1$, where $K$ and $N$ are scalar constants.
However, the parameter $\alpha_0 $ is not known!

To solve this problem, we estimate $\alpha_0$ using the ``observer'':
\begin{eqnarray*}
\frac{d}{dt}{\hat\alpha} = \ell (\dot y-\hat\alpha u)u.
\end{eqnarray*} 
Here $\ell$ is a positive constant.

\balphlist
\item
Verify that $\frac{d}{dt} \, (\hat\alpha-\alpha_0)^2 \le 0$ for all
time. Hence the absolute identification error $ |\hat\alpha-\alpha_0|
$ can never increase, regardless of the control. Find a control signal
$u(t)$ which forces the error to go to zero.
\end{list}

Assuming that $\hat\alpha(t)$ does converge to $\alpha_0$, a natural
control law is
\begin{eqnarray*}
u(t)= \frac{1}{\hat\alpha(t)}(r(t)-y(t))
\end{eqnarray*}
We now consider the {\it regulation problem}, where $r(t)=r_0\in\Re$
for all $t$.

\balphlist \setcounter{l2}{1}
\item
Obtain a nonlinear state space model for the {\it closed-loop }
system, with states $y$ and $\hat\alpha$.

\item
Show that $y=r_0$; $\hat\alpha = \beta$ is an equilibrium for {\it
any} $\beta\neq 0$. Deduce that the equilibrium
$\left(\begin{smallmatrix}r_0\\ \alpha_0\end{smallmatrix}\right)$ {\it
cannot} be asymptotically stable.

\item
Linearize the nonlinear model about an equilibrium.  What can you
conclude about the stability of an equilibrium
$\left(\begin{smallmatrix}r_0\\ \beta\end{smallmatrix}\right)$?

\item
Show that $V(y,\hat\alpha) = \half (y-r_0)^2 + \half
(\hat\alpha-\alpha_0)^2$ is a Lyapunov function for this model.  Can
you conclude that the equilibrium $(r_0,\alpha_0)$ is stable in the
sense of Lyapunov?
\end{list} 
These ideas are the starting point of {\it adaptive control}; an
approach based on simultaneous system identification and control.

\item
We have seen that stability of an LTI model may be characterized in
terms of the existence of a Lyapunov function, and we found that the
Lyapunov function in the stable case can be taken as
\begin{eqnarray*}
V(x) = x^T\int_0^\infty \phi(\tau)^T Q \phi(\tau)\, d\tau \, x = x^T P
x \qquad\hbox{with $Q>0$.}
\end{eqnarray*}
A similar result is possible in the nonlinear case under suitable
stability conditions, although the function $V$ will not necessarily
be a quadratic if the model is not linear.

Suppose that for the nonlinear model $\dot x = f(x)$, the origin
$x_e=0$ is globally asymptotically stable, and suppose moreover that
the stability is such that the solutions $\phi(t;x)$ are square
integrable over time.  Hence for any $Q>0$ the function $V(x) =
\int_0^\infty \phi(\tau;x)^T Q \phi(\tau;x)\, d\tau $ is finite
valued.

\begin{list}{(\alph{l2})}{\usecounter{l2}}
\item
Compute $V(x(t))$ for $t>0$, and simplify using properties of the
state transition function, e.g.,
\begin{eqnarray*}
\phi(0;x)=x,\qquad \phi(t;\phi(s;x)) = \phi(t+s;x)
\end{eqnarray*}

\item
Using (a), show that $\frac{d}{dt} V(x(t))|_{t=0} = -x(0)^T Q x(0)$.
Conclude that $V$ is an effective Lyapunov function for this model.
\end{list} 
\end{exercises}

\chapter{Controllability}
\clabel{Controllability}
 
In this chapter we consider exclusively the LTV state space model
\begin{eqnarray}
\dot x = A (t) x + B (t) u, \qquad x (t_0) = x_0 \,,
\elabel{LTV-controlled}
\end{eqnarray}
and as a special case its LTI version.  Our interest lies in the
effectiveness of the input $u$ to control the state $x$ to arbitrary
points in $\Re^n$.  Consider for example the two-dimensional LTI model
illustrated in \Figure{controlExample}.

\begin{figure}[ht]
\ebox{.9}{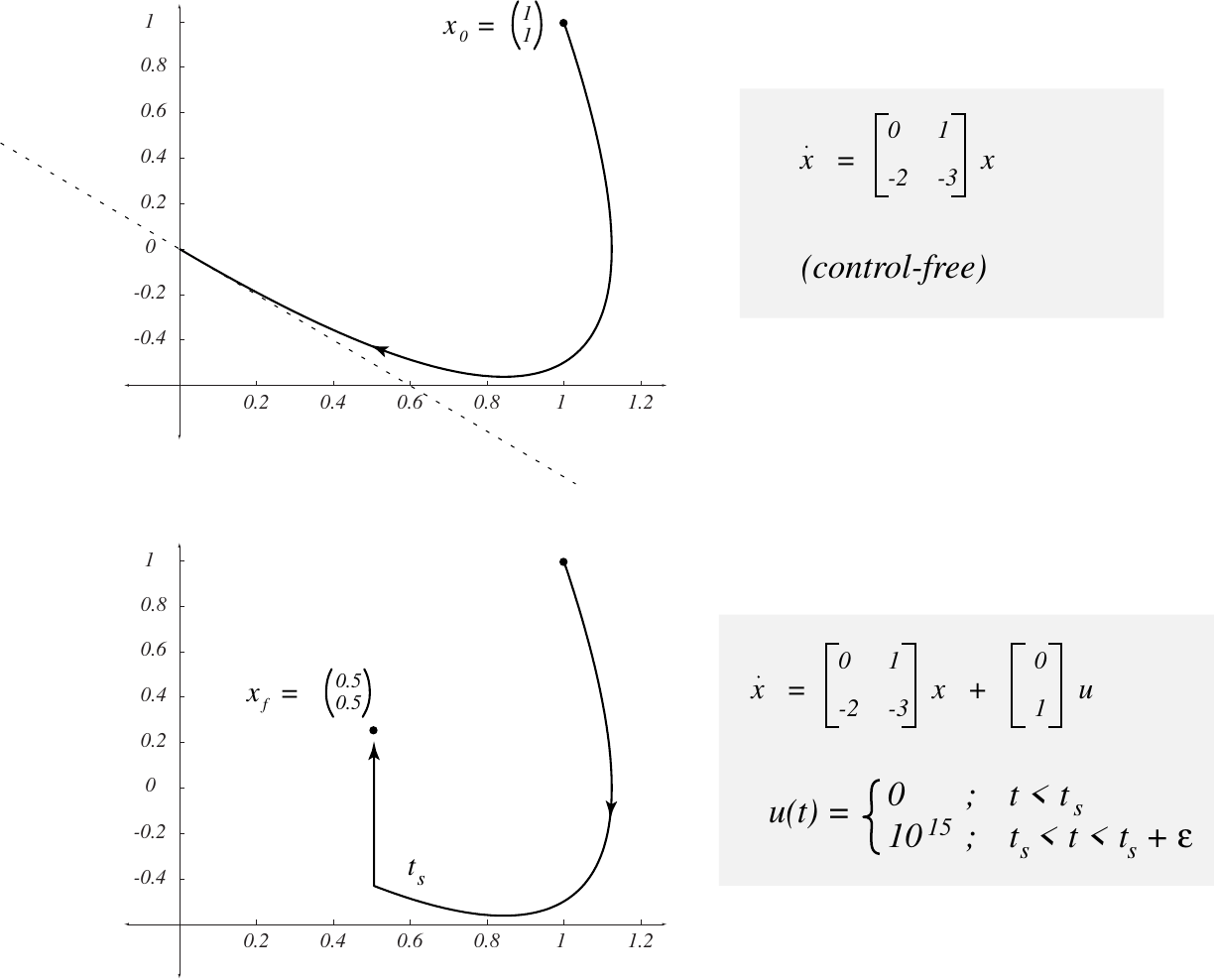}
\caption[The controlled model can be steered to any state]{The controlled model can be steered to any state in $\Re^n$
using impulsive control.}  \flabel{controlExample}
\end{figure}

From the state equation
\begin{eqnarray*}
\dot x = \left[\begin{matrix}0&1\\ -2 & -3\end{matrix}\right] x +
\left[\begin{matrix}b_1\\ b_2\end{matrix}\right] u.
\end{eqnarray*}
we see that when $b_1 = b_2= 0$, this is the example illustrated in
\Figure{trajectories}~(a) at the beginning of the previous chapter.
Obviously, in this case the control cannot affect the state, so when
$b=\zero$ the model is \textit{uncontrollable} in the strongest sense.
Consider now $b=\left(\begin{smallmatrix}0\\
1\end{smallmatrix}\right)$, and suppose that it is desired to choose a
control which takes the initial state
$x_0=\left(\begin{smallmatrix}1\\ 1\end{smallmatrix}\right)$ to a
final state $x_f=\left(\begin{smallmatrix}0.5\\
0.5\end{smallmatrix}\right)$ at some time $t_f$.  One approach is
illustrated in \Figure{controlExample}: Simply follow the control-free
trajectory until a time $t_s$ at which the state is nearly ``below''
the desired final state $x_f$.  Then, apply a large positive value so
that the vector $b u$ dominates the state derivative $Ax + b u$. This
procedure will place $x(t_f)$ at $x_f$ for \textit{some} $t_f$, and a
similar strategy can be used to drive any state to any other state.

A substantially different situation occurs with $b =
\left(\begin{smallmatrix}-1\\ 1\end{smallmatrix}\right)$.  In this
case the vector $b$ is parallel to the subspace indicated by the
dotted line $\hbox{Span}\, \{v^1 \} = \hbox{Span}\, \{b \}$, where
$v^1$ is an eigenvector for $A$.  It can be shown in this case that no
control will take a state from one side of this line to the other, and
that if the initial state $x_0$ lies on this line, then it will stay
there forever, regardless of the control.

The general LTV model \eq LTV-controlled/ is called
\defn{controllable} at time $t_0$ if there exists a finite time,
$t_f$, such that for \textit{any} initial state $x_0$, and
\textit{any} final state $x_f$, there is a control $u$ defined on
$[t_0, t_f]$ such that $x (t_f) = x_f$.  This concept is not obviously
directly useful in control design, except possibly in a ballistics
problem where $x_f$ denotes a target!  However, we will see that this
is a central idea in the theory of state space models, and that
controllability will be an underlying assumption in many of the
control design methods to be developed.  The purpose of this chapter
is to characterize controllability to better understand examples such
as these simple models.  We will see in \Theorem{Hautus} that the
pathology illustrated in \Figure{uncontrolExample} is the only way
that a lack of controllability can occur for an LTI model.

\begin{figure}[ht]
\ebox{.8}{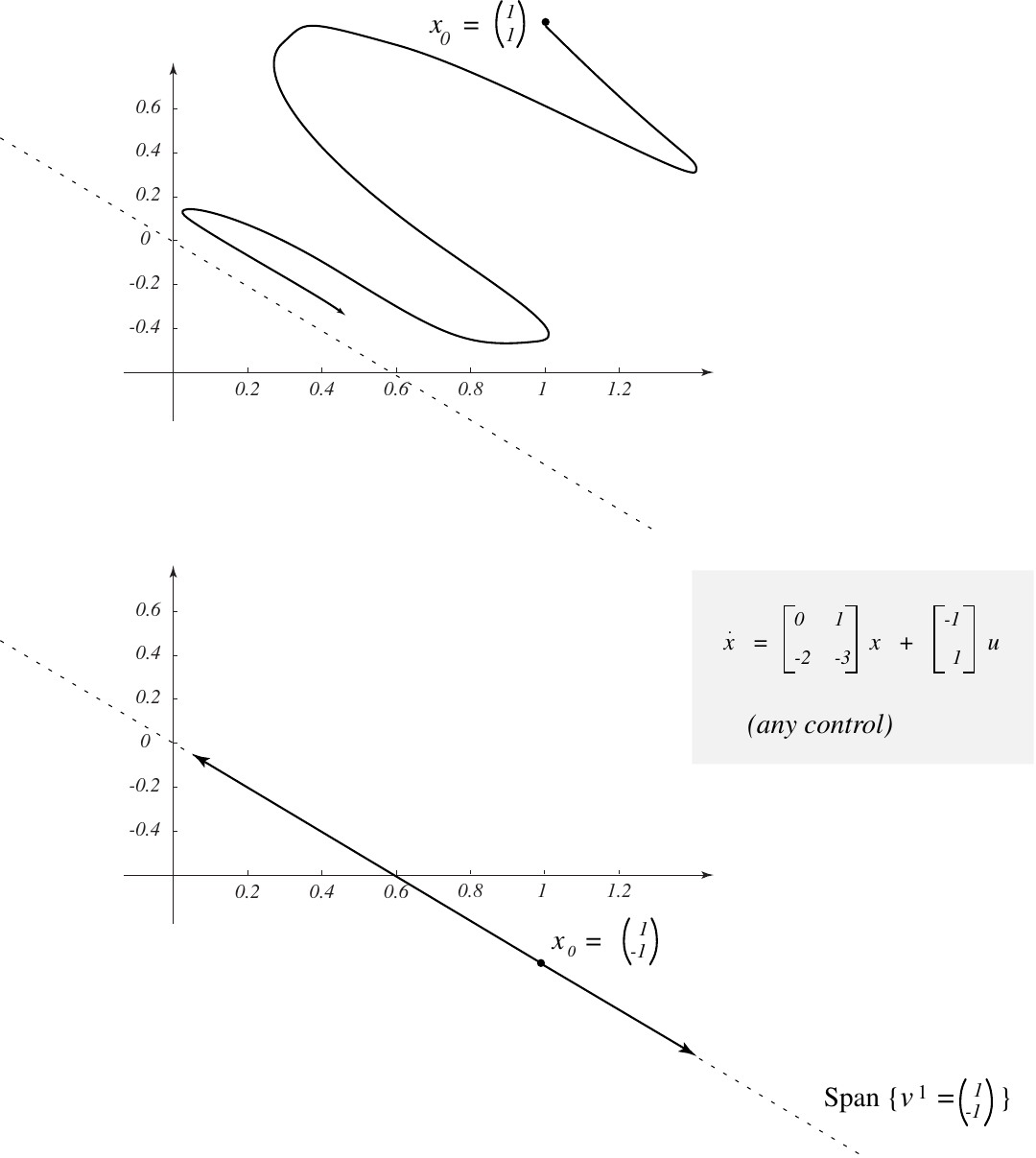}
\caption[The state can never
leave this line, regardless of the control]{If the initial condition satisfies $x_1 = -x_2$, so that the
state lies on the indicated diagonal line, then the state can never
leave this line, regardless of the control.  If the initial condition
does not start on the diagonal, then it may never reach it, though it
will come arbitrarily close to it for large $t$.}
\flabel{uncontrolExample}
\end{figure}
\clearpage

\section{A preview: The LTI discrete-time case}
\label{s:chap2.4.1}

Controllability is most easily understood in the LTI discrete-time
case
\begin{eqnarray*}
x (k+1) = A x (k) + B u (k), \qquad k= 0, 1, 2, 3, \ldots,
\end{eqnarray*}
where again $x$ is of dimension $n$ and $u$ is of dimension $m$.
Assume that we start at $x (0) =\zero $, and that it is desired to
move to another state $x_f$ in finite time.  We have $x (1) = Bu (0)$,
and in general
\begin{eqnarray*}
k = 2\qquad x (2) & = & A x (1) + B u (1) = A Bu (0) + Bu (1) \\ k = 3
\qquad x (3) & = & A x (2) + B u (2) \\ & = & A^2 Bu (0) + A B u (1) +
Bu (2) \\ && \qquad \qquad \vdots \\ k = r \qquad x (r) & = & A^{r-1}
Bu (0) + A^{r-2} Bu (1) + \ldots Bu (r-1).
\end{eqnarray*}
This can be written compactly as
\begin{eqnarray*}
x (r) = [B\ AB \cdots A^{r-1} B] \left[
\begin{array}{c}  u (r-1) \\ u (r-2) \\ \vdots \\ u (0) \end{array}
\right]
\end{eqnarray*}
which shows that we can reach arbitrary $x_f\in \Re^n$ at time $t_f =
r$ if and only if $\clR([B\ AB \cdots A^{r-1} B]) =\Re^n $, which
holds if and only if the rank of $[B\ AB \cdots A^{r-1} B] $ is equal
to $n$.

How do we know when to stop iterating?  The Cayley-Hamilton
\Theorem{CH} tells us that $A^k$ can be written as a linear
combination of $A^i$, $i<n$.  Hence the rank of the matrix $[B \mid AB
\mid \cdots \mid A^{r-1} B]$ cannot increase for $r\ge n$.  The
\defn{controllability matrix} is the $n\times (nm) $ matrix
\begin{eqnarray*} 
\clC \eqdef [B \mid AB \mid \cdots \mid A^{n-1} B].
\end{eqnarray*}
From this discussion it follows that in the discrete time model it is
possible to steer from the origin to any other state in $\Re^n$ if and
only if $\rank (\clC) = n$.

\section{The general LTV continuous-time case}
We saw in \Section{fundamental} that the solution to the state space
model \eq LTV-controlled/ may be expressed
\begin{eqnarray*}
x (t) = \phi (t, t_0) x_0 + \int_{t_0}^t \phi (t, \tau) B (\tau)
u(\tau) \, d\tau,
\end{eqnarray*} 
Controllability is thus determined by the range of the linear operator
$\clA\colon D[t_0,t_f]\to\Re^n$ defined as
\begin{equation}
\clA(u) = \int_{t_0}^{t_f} \phi (t_f, \tau) B (\tau) u (\tau) \, d
\tau. \elabel{controlOp}
\end{equation}
The input space $D[t_0,t_f]$ is the set of all piecewise continuous
time functions on $[t_0, t_f]$.  We leave the proof of the following
as an exercise.

\begin{theorem}\tlabel{rangeControllable}
The following are equivalent for the LTV model \eq LTV-controlled/
\balphlist
\item The model is controllable at time $t_0$;
\item $\clR(\clA)=\Re^n$ for some $t_f>t_0$.
\item There exists a time $t_f>t_0$ such that for any $x_0\in\Re^n$,
there is a control $u$ defined on $[t_0,t_f]$ such that $x(t_f) =
\zero$ when $x(t_0) = x_0$.
\end{list}
\end{theorem}
 
The \defn{controllability grammian} $W (t_0, t_f)$ for the LTV model
\eq LTV-controlled/ is the positive definite matrix
\begin{equation}
W (t_0, t_f) \eqdef \int_{t_0}^{t_f} \phi (t_0, \tau) B (\tau) B^*
(\tau) \phi^* (t_0, \tau) \, d\tau.  \elabel{grammian}
\end{equation}The following deeper characterization of controllability gives a finite dimensional test 
in terms of the $n\times n$ matrix $W$. The essence of the proof is
this: the grammian is full rank if and only if the $n$ rows of the
function $ \phi (t_0, \varble) B (\varble) $ are linearly independent
as functions in $D[t_0, t_f]$. This then is equivalent to
controllability of the model.
 
\begin{theorem}\tlabel{grammian}
The LTV model \eq LTV-controlled/ is controllable at time $t_0$ if and
only if there exists a finite time $t_f>t_0$ such that $W (t_0, t_f)$
is nonsingular.
\end{theorem}

\proof \head{(i)} Assume $W (t_0, t_f)$ is nonsingular, so that $W
(t_0, t_f)^{-1}$ exists.  We will show that the model is controllable
by showing that it is possible to control any $x_0$ to $\zero$ at time
$t_f$.  Consider the control
\begin{eqnarray*}
u (t) = - B^* (t) \phi^* (t_0, t) W^{-1} (t_0, t_f) x_0.
\end{eqnarray*}
Given the formula
\begin{eqnarray*}
x (t_f) = \phi (t_f, t_0) x_0 + \int_{t_0}^{t_f} \phi (t_f, \tau) B
(\tau) u (\tau) \, d\tau
\end{eqnarray*}
for the final value of the state, we may substitute the formula for
the control to give
\begin{eqnarray*}
x (t_f) & = & \phi (t_f, t_0) x_0 + \int_{t_0}^{t_f} \phi (t_f, \tau)
B (\tau) \left[ - B^* (\tau) \phi^* (t_0, \tau) W^{-1} (t_0, t_f) x_0
\right] \, d\tau \\ & = & \phi (t_f, t_0) x_0 - \Bigl(
\int_{t_0}^{t_f} \phi (t_f, \tau) B (\tau)B^* (\tau) \phi^* (t_0,
\tau) \, d\tau \Bigr) W^{-1} (t_0, t_f) x_0.
\end{eqnarray*}
The term $W^{-1} (t_0, t_f) x_0$ is taken outside of the integral
since it does not depend on $\tau$.  By the semigroup property
\begin{eqnarray*}
\phi (t_f, t_0) \phi (t_0, \tau) = \phi (t_f, \tau),
\end{eqnarray*}
the integral becomes $\phi (t_f, t_0) W (t_0, t_f)$.  Therefore
\begin{eqnarray*}
x (t_f) = \phi (t_f, t_0) x_0 - \Bigl(\phi (t_f, t_0) W
(t_0,t_f)\Bigr) W^{-1} (t_0, t_f) x_0
\end{eqnarray*}
which gives $x (t_f) = \zero$. Hence by \Theorem{rangeControllable}
the model is controllable at time $t_0$.

\head{(ii)} We now show that controllability of the LTV model implies
that $W(t_0, t_f)^{-1}$ exists.  We prove the contrapositive: Assuming
$W (t_0, t_f)$ is singular, we show that the model is not
controllable.

If $W$ is singular, there exists an $\alpha \neq \zero$ such that
\begin{eqnarray}
W (t_0, t_f) \alpha = \zero \qquad \hbox{and hence also} \qquad
\alpha^* W (t_0, t_f) \alpha = 0.  \elabel{eq2a}
\end{eqnarray}
From the definition of $W$ it then follows that
\begin{eqnarray*}
\int_{t_0}^{t_f} \alpha^* \phi (t_0, \tau) B (\tau) B^* (\tau) \phi^*
(t_0, \tau) \alpha \, d\tau = 0
\end{eqnarray*}
This can be equivalently expressed as
\begin{eqnarray*}
\int_{t_0}^{t_f} \left\langle B^* (\tau) \phi^* (t_0, \tau) \alpha,
B^* (\tau) \phi^* (t_0, \tau) \alpha \right\rangle \, d\tau = 0,
\end{eqnarray*}
or
\begin{eqnarray*}
\int_{t_0}^{t_f} | B^* (\tau) \phi^* (t_0, \tau) \alpha |^2 \, d\tau =
0.
\end{eqnarray*}
Since the integral of the square norm is zero, and all of these
functions are assumed to be piecewise continuous, we must have $B^*
(\tau) \phi^* (t_0, \tau) \alpha=\zero$ and hence
\begin{equation}
\elabel{alphaNull}
\alpha^* \phi  (t_0, \tau) B (\tau)  =\zero^*
\end{equation}
for all $\tau\in [t_0,t_f]$.

If the model is controllable at time $t_0$, then starting at the
particular state $x_0 = \alpha$, there is a control $u$ which makes
$x(t_f)=\zero$:
\begin{eqnarray*}
x (t_f) = \zero = \phi (t_f, t_0) \alpha + \int_{t_0}^{t_f} \phi (t_f,
\tau) B (\tau) u (\tau) \, d\tau.
\end{eqnarray*}
From the property $\phi^{-1} (t_f, t_0) = \phi (t_0, t_f)$, we can
multiply this equation through by by $\phi^{-1}$ to obtain
\begin{eqnarray*}
\zero = \alpha + \int_{t_0}^{t_f} \phi (t_0, t_f) \phi (t_f, \tau) B
(\tau) u (\tau) \, d\tau.
\end{eqnarray*}
From the semigroup property we must therefore have
\begin{eqnarray*}
\alpha = - \int_{t_0}^{t_f} \phi (t_0, \tau) B (\tau) u (\tau) \,
d\tau.
\end{eqnarray*}
Multiplying both sides on the left by $\alpha^*$ gives
\begin{eqnarray*}
| \alpha |^2 = - \int_{t_0}^{t_f} \alpha^* \phi (t_0, \tau) B(\tau) u
(\tau)\, d \tau
\end{eqnarray*}
which by \eq alphaNull/ must be zero.  However, $| \alpha |^2 \neq 0$
by assumption.  Hence the model cannot be controllable if $W$ is
singular, and the proof is complete.  \qed

\section{Controllability using the controllability matrix}

In the LTI discrete-time case we have shown that the controllability
matrix is key to understanding controllability.  To show that this is
still true for the continuous-time LTI model
\begin{equation}
\begin{array}{rcl}
\dot x & = & Ax + Bu\\ y & = & Cx + Du,
\end{array}
\elabel{structure-LTI}
\end{equation}
we consider how the controllability grammian $W$ is related to the
controllability matrix
\begin{eqnarray*}
\clC = [ B \mid AB \mid A^2 B \mid \cdots \mid A^{n-1} B].
\end{eqnarray*}  
Using this approach we will prove the following theorem.
  
\begin{theorem}\tlabel{control-matrix}
For the LTI model \eq structure-LTI/, the model is controllable if and
only if the $n \times (nm)$ matrix $\clC$ has rank $n$.
\end{theorem}

\proof To begin, note that in the LTI case we have for $t_0=0$,
\begin{eqnarray}
W (0, t_f) & = & \int_{0}^{t_f} \phi (- \tau) B B^* \phi^* (- \tau) \,
d\tau \nonumber \\ & = & \int_0^{t_f} e^{-A \tau} BB^* e^{- A^*\tau}
\, d\tau. \elabel{Wlti}
\end{eqnarray}

To proceed, we show that
\begin{eqnarray*}
\rank (\clC) < n \Longleftrightarrow \hbox{$W$ is singular}.
\end{eqnarray*}

\head{(i)} We first show that assuming $\rank (\clC) < n$, we must
also have that $W$ is singular.

If $\rank (\clC) < n$, there exists $\alpha \neq \zero$ such that
\begin{eqnarray*}
\alpha^* [B\ AB \cdots A^{n-1}B] = \zero^*
\end{eqnarray*}
which can be equivalently expressed as
\begin{eqnarray*}
\alpha^* A^k B = \zero^*, \quad k = 0, \ldots, n-1.
\end{eqnarray*}
From the Cayley-Hamilton \Theorem{CH}, the state transition matrix can
be expressed as
\begin{eqnarray*}
\phi (- \tau) = e^{-A\tau}= \sum_{k=0}^{n-1} \beta_k (- \tau) A^k
\end{eqnarray*}
from which it follows that $\alpha^* \phi (- \tau) B =\zero^*$ for all
$\tau$.  From \eq Wlti/ we see that $\alpha^* W (0, t_f) =\zero^*$,
which shows that $W$ is singular.

\head{(ii)} Assume now that $W$ is singular.  We will show that $\rank
(\clC) < n$.
  
We have seen in equation \eq alphaNull/ that if $W(0,t_f)$ is
singular, then there exists an $\alpha\in\Re^n$ such that
\begin{eqnarray*}
\alpha^* \phi (0, \tau) B (\tau) = \zero^*,\qquad 0\le\tau\le t_f.
\end{eqnarray*}
In the LTI case, it follows that $\alpha^* e^{-A\tau} B=\zero^*$ for
all $\tau$.  By setting $t=0$ this shows that $\alpha^* B=\zero^*$.
Moreover, if a function is zero over an interval, so is its
derivative:
\begin{eqnarray*}
\left. \frac{d}{dt} (\alpha^* e^{-At} B) = \alpha^T (-Ae^{-At} B),
\right|_{t=0} = 0
\end{eqnarray*}
which shows that also $- \alpha^* AB = \zero^*$.  Continuing this
procedure, we see that
\begin{eqnarray*}
\left.  \frac{d^2}{dt^2} (\alpha^* e^{-At} B) \right|_{t=0} = \alpha^*
A^2 B =\zero ^*,
\end{eqnarray*}
and for any $k$
\begin{eqnarray*}
(-1)^{k-1} \alpha^* A^{k-1} B =\zero^*.
\end{eqnarray*}
It then follows that
\begin{eqnarray*}
\alpha^* [B|AB|A^2B \mid \cdots \mid A^{n-1} B] =\zero ^*,
\end{eqnarray*}
which establishes the desired result that $\rank (\clC) < n$.  \qed

A careful examination of the proof shows that we in fact have
established the stronger result
\begin{eqnarray*}
\clR(\clC) =\clR(\clA).
\end{eqnarray*}
Hence, the subspace of $\Re^n$
\begin{eqnarray*}
\Sigma_c \eqdef \clR(\clC) =\{\clC z : z\in\Re^{mn} \}
\end{eqnarray*}
is equal to the set of all states which can be reached from the
origin.  The set $\Sigma_c$ is called the \defn{controllable
subspace}.

\section{Other tests for controllability} 
The controllability of an LTI state space model is much more apparent
when the model is in modal form.  To see this, we must construct the
modal form using a similarity transformation, so we begin with the
following:

\begin{theorem}\tlabel{control-invariance}
Controllability of the LTI model
\begin{eqnarray*}
\dot x = Ax + Bu
\end{eqnarray*}
is invariant under any equivalence transformation $\barx = P x$
\end{theorem}

\proof For the transformed model $\dot{\bar x} = \bar A\bar x + \bar
Bu$ we have
\begin{eqnarray*}
\bar\clC & = & [\barB \mid \barA \barB \mid \cdots \mid \barA^{n-1}
\barB] \\ & = & [PB \mid (PAP^{-1}) PB \mid (PAP^{-1})(PAP^{-1}) PB
\mid \cdots \mid (PAP^{-1})\cdots(PAP^{-1}) PB].
\end{eqnarray*}
Therefore, the controllability matrix for the transformed model is
\begin{eqnarray*}
\bar{\clC} = [PB \mid PAB \mid PA^2 B \mid \cdots \mid PA^{n-1} B],
\end{eqnarray*}
or $\bar{\clC} = P\clC $.  Since the matrix $P$ is nonsingular, it
follows that $\rank \bar\clC= \rank \clC$.  \qed

\begin{figure}
\ebox{0.4}{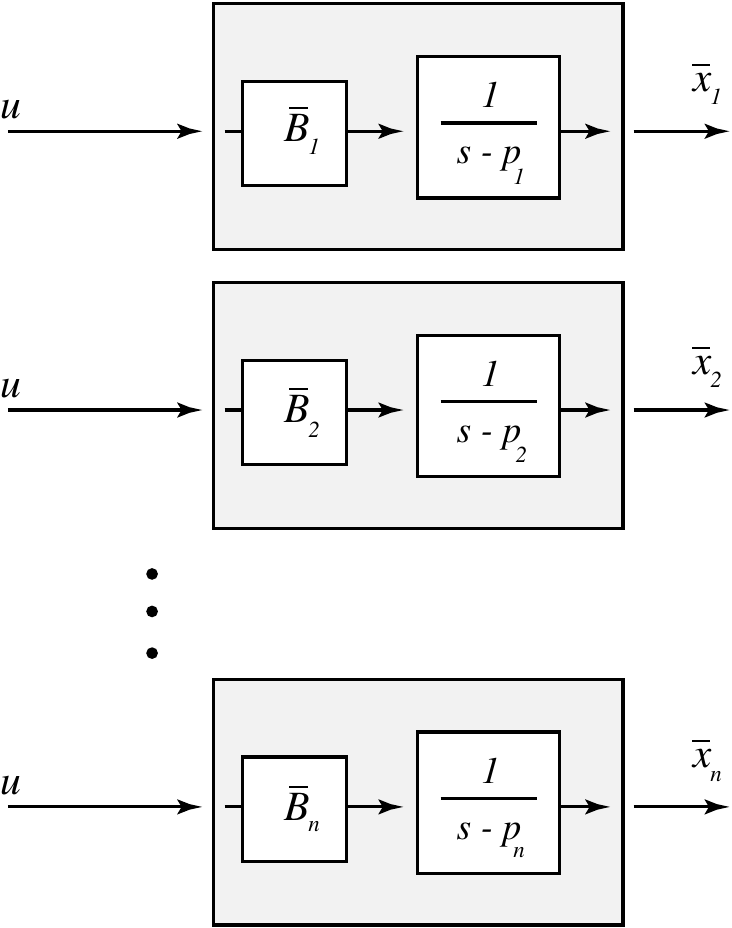}
\caption[The modes of the model are decoupled]{The modes of the model are decoupled, and hence the model
cannot be controllable if some $\bar B_i$ is zero.}
\flabel{structure-modal}
\end{figure}

Suppose that the eigenvalues of $A$ are distinct, so that $A$ is
similar to a diagonal matrix $\Lambda$:
\begin{eqnarray*}
\barA = \Lambda, \qquad P = M^{-1},
\end{eqnarray*}
where $M = [v^1 \cdots v^n]$.  This gives the modal form
\begin{eqnarray*}
\dot{\barx} = \Lambda \barx + \barB u,
\end{eqnarray*}
where the modes of the model are decoupled.  \Figure{structure-modal}
illustrates the structure of the transformed model. If some $\barB_i$
is zero, which means that a row of the matrix $\barB$ is zero, it is
apparent that the model cannot be controllable.  It turns out that the
converse is also true:

\begin{theorem}\tlabel{modal-control} 
Suppose that the matrix $A$ has distinct eigenvalues.  Then, $(A,B)$
is controllable if and only if $\barB = M^{-1} B$ has no zero rows.
\end{theorem}

\proof We prove the theorem in two steps.

\head{(i)} If row $i$ of of $\barB = M^{-1} B$ is zero, it follows
that if the initial state satisfies $x_i(0)=0$, then $x_i(t)=0$ for
all $t$.  Obviously then the model is not controllable.

\head{(ii)} Conversely, if the model is not controllable, there exists
an $\alpha\in\Re^n$, $\alpha\neq 0$, such that
\begin{eqnarray*}
0= \alpha^* \barW \alpha = \int_0^{t_f} \alpha^*e^{-\Lambda \tau}
\barB\, \barB^* e^{- \Lambda^*\tau} \alpha\, d\tau.
\end{eqnarray*}
Equivalently, we must have
\begin{eqnarray*}
\sum_{i=1}^n \alpha_i \barB_i^* e^{-\lambda_i \tau} \equiv \zero^*.
\end{eqnarray*}
If none of the $\{\barB_i \}$ are zero, choose $k$ such that
$\alpha_k\neq 0$.  We then have
\begin{eqnarray*}
e^{-\lambda_k\tau} = \frac{1}{\alpha_k \barB_k\barB_k^*} \sum_{        
\substack{
 i=1
\\
 i \not = k}
}^n \alpha_i \barB_k\barB_i^* e^{-\lambda_i \tau}.
\end{eqnarray*}
This is impossible, since the functions $\{e^{\lambda_i t}\}$ are
linearly independent in $C[t_0,t_f]$.  Hence, some row $B_k$ must be
zero, which thus proves the theorem.  \qed

\Theorem{modal-control} may also be understood in terms of reciprocal
eigenvectors.  Recall that the matrix $\barB$ may be expressed as
\begin{eqnarray*}
\barB = M^{-1} B = \left[ \begin{array}{c} r^{1*} \\ \vdots \\ r^{n*}
\end{array} \right] B
\end{eqnarray*}
The $j$th row of $\barB$ is $r^{j*} B$.  From \Theorem{modal-control}
we see that $(A,B)$ is controllable if and only if there is no $r^j$
within the null space of $B^*$.

Since $r^{j*}$ is also a left eigenvector (see \Exercise{left-eig} of
\Chapter{chap1.5}), we always have
\begin{eqnarray*}
r^{j*} (\lambda_j I - A) = \zero^*.
\end{eqnarray*}
Then, $(A,B)$ is uncontrollable if and only if there is a $j$ such
that
\begin{eqnarray*}
r^{j*} (\lambda_j I - A) = \zero^* \qquad \hbox{and}\qquad r^{j*} B
=\zero^*,
\end{eqnarray*}
which can be written compactly as
\begin{eqnarray*}
r^{j*} \left[ \left. \lambda_j I - A \right| B \right] =\zero^*.
\end{eqnarray*}
To summarize, we see that the pair $(A,B)$ in uncontrollable if and
only if for some $j$,
\begin{eqnarray*}
\rank [\lambda_j I - A \mid B ] < n.
\end{eqnarray*} 
This final test was derived using a modal form, but it turns out to be
true in general:

\begin{theorem}[Hautus-Rosenbrock test]
\tlabel{Hautus} The pair $(A,B)$ is controllable if and only if for
every $s\in\Co$,
\begin{equation}
\elabel{Hautus} \rank [ s I - A \mid B] = n.
\end{equation} 
\end{theorem} 

\proof One direction is easy: If $r^* [ \lambda I - A \mid B ] =\zero
^*$, then we must have for any input $u$,
\begin{eqnarray*}
\frac{d}{dt} (r^* x\, (t)) = r^* A x(t) + r^* B u(t) = \lambda r^*
x(t),
\end{eqnarray*}
so that $ r^* x (t) = e^{\lambda t} r^* x (0)$, regardless of the
control!  Hence, if \eq Hautus/ is violated for $s=\lambda$, then the
model is not controllable.

We have seen that the converse is a direct consequence of
\Theorem{modal-control} in the case of distinct eigenvalues, but the
proof is more subtle in full generality. The reader is referred to
page~184 of \CHE\ for details.  \qed

Note that when $s$ is not an eigenvalue of $A$ the matrix $sI-A$
already has rank $n$, so that we also have $\rank ([s I - A \mid B]) =
n$.  Hence the rank of the matrix \eq Hautus/ need only be checked for
$s$ equal to the eigenvalues of $A$.

\begin{ex}
We illustrate some of the controllability tests developed so far using
the simple model
\begin{eqnarray*}
\dot x = \left[ \begin{array}{rr} -2 & 0 \\ -1 & -1 \end{array}
\right] x + \left[ \begin{array}{r} 1 \\ 1 \end{array} \right] u.
\end{eqnarray*}

\head{(i)} The controllability matrix for this model is
\begin{eqnarray*}
\clC = [ B \mid AB] = \left[ \begin{array}{rrr} 1 & -2 \\ 1 &
-2\end{array} \right].
\end{eqnarray*}
Since the rank of $\clC$ is 1, the model is not controllable.  The
controllable subspace $\Sigma_c$ is equal to the range space of
$\clC$, which in this case is the diagonal in $\Re^2$:
\begin{eqnarray*}
\Sigma_c = \{\alpha \left(\begin{smallmatrix}1\\
1\end{smallmatrix}\right) : \alpha \in\Re\}.
\end{eqnarray*}

\head{(ii)} The eigenvalues of $A$ are $\lambda_{1} = -1$,
$\lambda_{2} -2$, and the modal matrix can be taken as $M = \left[
\begin{array}{rr} 0 & 1 \\ 1 & 1 \end{array} \right]$, $M^{-1} = \left[ \begin{array}{rr} -1 &
1 \\ 1 & 0 \end{array} \right]$.  We can then compute
\begin{eqnarray*}
\Lambda = \left[ \begin{array}{rr} -1 & 0 \\ 0 & -2 \end{array}
\right] \qquad \barB = M^{-1} B= \left[
\begin{array}{r} 0 \\ 1 \end{array} \right], 
\end{eqnarray*}
and again we conclude that the model is not controllable.  Since the
first row is zero, it follows that the first mode $(\lambda_1 = - 1)$
is not connected to the input in modal form.
 
\head{(iii)} We have
\begin{eqnarray*}
[sI - A \mid B] = \left[ \begin{array}{rrr} s+2 & 0 & 1 \\ 1 & s+1 & 1
\end{array} \right].
\end{eqnarray*}
Evaluating this matrix with $s=-1$ gives
\begin{eqnarray*}
\left[ \begin{array}{rrr} 1 & 0 & -2 \\ 1 & 0 & -2
\end{array} \right].
\end{eqnarray*}
Since the rank of this matrix is one, we again see from the
Hautus-Rosenbrock test that the model is not controllable.\end{ex}

\begin{exercises} 
\item
Investigate the controllability properties of the LTI model $\dot x =
Ax+Bu$, for the three pairs of $(A,B)$ matrices given below.
\balphlist
\item
$A = \left[\begin{matrix}-5&1\\ 0&4\\\end{matrix}\right]\;,\quad B =
\left[\begin{matrix}1\\ 1\\\end{matrix}\right]$.

\item
$A = \left[\begin{matrix}3&3&6\\ 1&1&2\\
2&2&4\\\end{matrix}\right]\;,\quad B = \left[\begin{matrix}0\\ 0\\
1\\\end{matrix}\right]$

\item
$A = \left[\begin{matrix} 0&1&0\\ 0&0&1\\
0&0&0\\\end{matrix}\right]\;, \quad B = \left[\begin{matrix}0\\ 0\\
1\\\end{matrix}\right]$
\end{list}

\item
For the system model corresponding to $A=\left[\begin{matrix}-5&1\\
0&4\\\end{matrix}\right]$, $B=\left[\begin{matrix}1\\
1\\\end{matrix}\right]$, obtain a control that drives the system from
$x=(1\quad 0)^T$ at $t=0$, to $x=(0\quad 1)^T$ at $t=1$.  Simulate
your control using \textit{Simulink}.

\item
Given the linear time-invariant model
\begin{eqnarray*}
\dot x &=& \left[\begin{matrix}-7&-2&6\\ 2&-3&-2\\
-2&-2&1\\\end{matrix}\right]x+ \left[\begin{matrix}1&1\\ 1&-1\\
1&0\\\end{matrix}\right] u\quad = Ax+Bu \\
y&=&\left[\begin{matrix}-1&-1&2\\ 1&1&-1\\\end{matrix}\right]x = Cx\,,
\end{eqnarray*}
check controllability using \balphlist \item the controllability
matrix

\item
the rows of $\barB =M^{-1}B$, where $M$ is chosen such that $M^{-1}AM$
is diagonal

\item The   Hautus-Rosenbrock condition.
\end{list}

\item
Transform the state space model below into CCF, and from the resulting
equations compute its transfer function.
\begin{eqnarray*}
\dot x = \left(\begin{matrix}0&0&-2\\ 1&0&1\\
0&1&2\\\end{matrix}\right) x + \left(\begin{matrix}1\\1\\
1\\\end{matrix}\right) u \qquad y = (0\ 0\ 1) x
\end{eqnarray*}

\item \hwlabel{CtbleMode} Consider the LTI model
\begin{eqnarray*}
\dot x &=& \left[\begin{matrix}-3&1&0\\ 0&-3&0\\
0&0&4\\\end{matrix}\right] x + \left[\begin{matrix}0\\ 1 \\ 0
\\\end{matrix}\right] u \\ y &=& \left[1 \ 0 \ 1\right] x.
\end{eqnarray*}
For what initial states $x(0)=x_0$ is it possible to choose a control
$u$ so that $y(t) = t e^{-3t}$ for $t>1$?

\item Let $A$ be an $n\times n$ matrix and $B$ be an $n\times r$
matrix, both with real entries. Assume that the pair $(A,B)$ is
controllable.  Prove or disprove the following statements. (If the
statement is false, then producing a counterexample will suffice.)
\balphlist
\item
The pair $(A^2,B)$ is controllable.
\item
Let $k(\cdot)$ be an $r$-dimensional (known) function, piecewise
continuous in $t$. Then, the model described by
\begin{eqnarray*}
\dot x =Ax+Bu+k(t)
\end{eqnarray*}
is completely controllable, in the sense that any state can be reached
from any other state using some control.

\item Given that the model $\dot x =Ax+Bu$ has the initial condition
$x(0)=x_0\not =\zero$, it is possible to find a piecewise continuous
control, defined on the interval $[0,\infty)$, so that the model is
brought to {\it rest} at $t=1$ (that is, $x(t)=\zero$ for all $t\geq
1$).

\item
Assume now that the model above is initially (at $t=0$) at rest, and
we wish to find a piecewise continuous control which will bring the
state to ${\bar x}\in \Re^n$ by time $t=1$ and maintain that value for
all $t\geq 1$. Such a control can always be found.
\end{list}
 
\item
Consider the LTI model $\dot x = A x + B u$.  Suppose that $A$ has $n$
distinct eigenvalues, and corresponding eigenvectors $\{
v^1,\dots,v^n\}$.

Show that if $\clR(B) \subset \Span
(v^1,\dots,v^{j-1},v^{j+1},\dots,v^n)$ for some $j$ then the model is
not controllable.

{\it Hint: Look at the reciprocal basis vectors: How are $r^j$ and
$\clR(B)$ related}?

\item  
Compute the adjoint $\clA^*$ of the operator $\clA$ defined in \eq
controlOp/, and then compute the composition
\begin{eqnarray*}
V=\clA \circ \clA^* \colon \Re^n\to\Re^n
\end{eqnarray*}
as an $n\times n$ matrix.  How are the controllability grammian $W$
and the operator $V$ related?  Comment on this result in view of
\Theorem{grammian} and \Exercise{grammian} of
\Chapter{chap1.5}.

\item
Consider the $m$-input, $p$-output LTI model described by the
convolution
\begin{eqnarray*}
y(t)=\int_{-\infty}^\infty g(t-\tau) u(\tau)\, d\tau,
\end{eqnarray*}
where $g$ is a $p\times m$ matrix-valued function.

Call the model {\it output controllable} if for any $y_0\in \Re^p$,
$t_0\in\Re$, there exists a control $u_0$ such that $y(t_0)=y_0$.
\balphlist
\item
Formulate output controllability as a range space dimension problem
for a specific linear operator.

\item
Compute the adjoint of the linear operator determined in (a).

\item
Derive a finite dimensional test for output controllability, based on
the rank of some $p\times p$ matrix (see the previous problem).
\end{list}

\item
Prove \Theorem{rangeControllable}.
\end{exercises}
 
\chapter{Observability, Duality and Minimality}\clabel{Observability}

In many control problems we obtain measurements $y(t)$, but we may not
be able to measure directly all of the system states.  For a complex
system with a large number of states, it may be expensive to position
a sensor to measure every state, and it may even be impossible to do
so.  In the light of this, the question we address in this chapter is
the following: \textit{Is it possible to determine the state using
only the input-output information $\{u(t), y(t) \}$}?  Informally, we
say that some states are \textit{observable} if one can determine them
unambiguously based on these input-output measurements.  This is
obviously a very important property for a physical system.

\begin{figure}[ht]
\ebox{0.6}{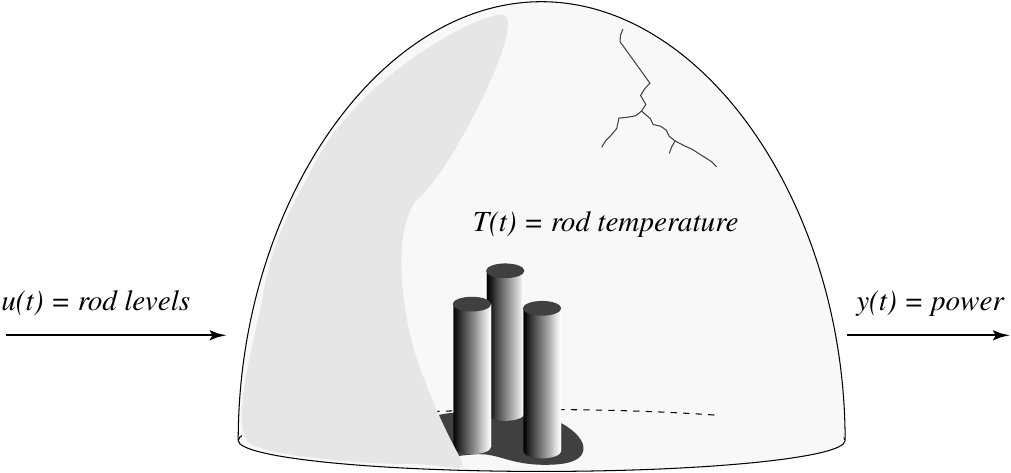}
\caption[If the input does not affect the internal temperature...]{If the input to this plant does not affect the internal
temperature, another actuator is needed.  If the temperature cannot be
estimated based on input-output measurements, another sensor is needed
to directly measure the temperature. } \flabel{structure-reactor}
\end{figure}

To give a simple, if artificial example, consider the academic model
of a nuclear reactor given in \Figure{structure-reactor}.  Suppose
that there is a state space model for this plant with the input equal
to the rod levels, the output equal to the power generated, and one of
the states equal to the rod temperature.  Suppose first that using
this control it is not possible to influence the temperature $T$.
Then the system is not controllable, and for safe operation it will be
necessary to add another \textit{actuator} to raise or lower the
internal temperature.  This situation seems highly unlikely for a real
reactor.  Suppose now that given the input-output measurements, it is
not possible to estimate the internal temperature $T$. We then say
that the system is not observable.  This again is a highly undesirable
situation, though it can be resolved by adding another measurement - a
\textit{sensor} to directly measure $T$.  Hence, an uncontrollable or
unobservable plant can be ``repaired'' by expanding the number of
inputs or outputs for the plant.

To take a more realistic example, consider the circuit given in
\Figure{structure-circuit}.  It is obvious here that if the input
voltage is zero, then the output voltage is zero, regardless of the
voltage across the capacitor.  Hence, the voltage $x_c$ is not
observable.

\begin{figure}
\ebox{0.55}{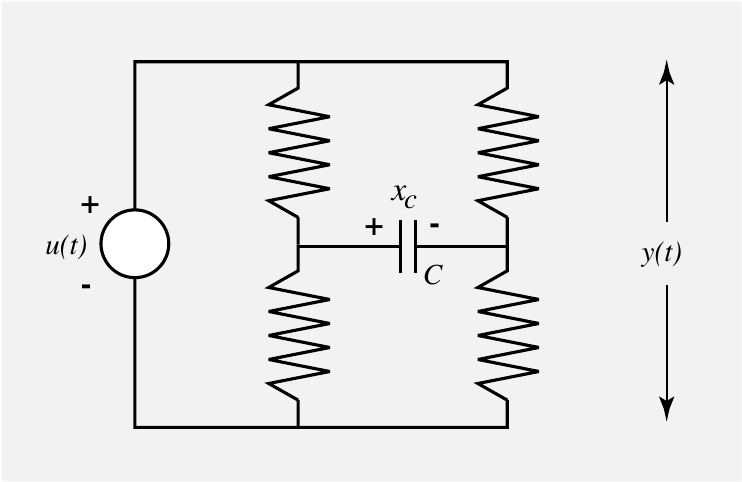}
\caption[If the input to this circuit is zero...]{If the input to this circuit is zero, then the output is
zero.  Hence the charge across the capacitor is \textit{not}
observable. } \flabel{structure-circuit}
\end{figure}

To give a formal definition of observability, consider the LTV model
with an output equation $y$ of the form
\begin{equation}
\begin{array}{rcl}
\dot x (t) & = & A (t) x + B (t) u; \\ y (t) & = & C (t) x + D (t) u.
\end{array} \qquad x (t_0) = x_0.  \elabel{structure-LTV}
\end{equation}
Since we are assuming that the system is defined by a linear state
space model, to determine $x(t)$ for all $t$ it is sufficient to
compute the initial state $x_0=x(t_0)$.

The model \eq structure-LTV/ is said to be \defn{observable} at time
$t_0$ if there exists a finite time $t_1 > t_0$ such that for any
initial $x_0$ at $t_0$, knowledge of the input $u(t)$ and the output
$y(t)$ for $t_0 \le t \le t_1$ suffices to determine $x_0$. Given the
expression for the state
\begin{eqnarray*}
x (t) = \phi (t, t_0) x_0 + \int_{t_0}^t \phi (t, \tau) B (\tau) u
(\tau) \, d\tau
\end{eqnarray*}
the output may be expressed as
\begin{eqnarray*}
y (t) = \underbrace{\vphantom{\int_{t_0}^t} C(t) \phi (t, t_0)
x_0}_{\hbox{\tiny unknown}} + \underbrace{C(t) \int_{t_0}^t \phi (t,
\tau) B (\tau) u (\tau) \, d\tau + D (t) u (t)}_{\baru (t)\
\hbox{\tiny known}}.
\end{eqnarray*}
Since $\baru$ is known, we can simplify matters by considering the
known quantity
\begin{equation}
\bary (t) = y (t) - \baru (t) = C(t) \phi (t, t_0) x_0.
\elabel{bary-def}
\end{equation}
Equivalently, when analyzing observability we may assume without loss
of generality that the input is zero.  Defining the linear operator
$\clB\colon \Re^n\to D[t_0,t_1]$ as
\begin{eqnarray*}
\clB(x) = f, \qquad \hbox{where}\ f(t) = C(t) \phi (t, t_0) x,
\end{eqnarray*}
we see that the model is observable if, for some $t_1$, the operator
$\clB$ is one to one.

\section{The observability matrix}

In the LTI case we have
\begin{eqnarray*}
\bary (t) = C e^{At} x_0 \qquad t \in [0, t_1],
\end{eqnarray*}
where in this case we have taken $t_0 = 0$.  We have some information
about $x_0$ through the equation $y (0) = Cx_0$. However, typically
$C$ is not square, so inversion cannot be used directly to solve this
equation for $x_0$.  Taking derivatives yields
\begin{eqnarray*}
\dot{\bary} (t) & = & CAe^{At} x_0 \\ \dot{\bary} (0) & = & CAx_0,
\end{eqnarray*}
and this process can be continued until sufficient information has
been generated to find $x_0$:
\begin{eqnarray*}
\left[ \begin{array}{c} \bary (0) \\ \dot{\bary} (0) \\ \ddot{\bary}
(0) \\ \vdots \\ \bary^{(n-1)}(0)
\end{array} \right] = \left[ \begin{array}{c} C \\ CA \\ CA^2 \\
\vdots \\ CA^{n-1} \end{array} \right] x_0.
\end{eqnarray*}

Defining the \defn{observability matrix} as the $(np)\times n$-matrix
\begin{eqnarray*}
\clO \eqdef \left[ \begin{array}{c} C \\ CA \\ CA^2 \\\cdots \\
CA^{n-1}
\end{array} \right],
\end{eqnarray*}
where $p$ is the dimension of $y$,
we can now prove the following theorem.

\begin{theorem}
\tlabel{obsMatrixTest} The LTI model \eq structure-LTI/ is observable
if and only if
\begin{eqnarray*}
\rank \clO = n.
\end{eqnarray*}
\end{theorem}

\proof If $\clO$ is full rank, then we have seen that $x_0$ can be
computed using $y$ and its derivatives.

Conversely, suppose that there is an $x_0\in \Re^n$, $x\neq \zero $,
such that $\clO x = \zero$. We then have $C A^k x_0 = \zero$ for any
$k$, and hence also $C e^{At} x_0 =\zero $ for all $t$. Since $x_0$ is
indistinguishable from the state $\zero$, it follows that the model is
not observable.  \qed

The {\it unobservable subspace} $\Sigma_{\bar o}$ is defined to be the null
space of the observability matrix.  From the proof above, we see that
for an initial condition $x(0)\in \Sigma_{\bar o}$, the output $y$ is
identically zero if the input is equal to zero.

\section{LTV models and the observability grammian}

In the LTV case we define the \defn{observability grammian} as
\begin{equation}
H (t_1, t_0) \eqdef \int_{t_0}^{t_1} \phi^* (\tau, t_0) C^* (\tau)
C(\tau) \phi (\tau, t_0) \, d \tau. \elabel{obsGram}
\end{equation}

\begin{theorem}\tlabel{obs-grammian}
The model
\begin{eqnarray*}
\dot x = A (t) x, \quad y = C(t) x
\end{eqnarray*}
is observable at time $t_0$ if and only if there exists a finite time
$t_1$ such that the observability grammian $H(t_1,t_0)$ has rank $n$.
\end{theorem}

\proof We first prove that non-singularity of the observability
grammian is \textit{sufficient} for observability:

\head{(i)} From \eq bary-def/ we have
\begin{eqnarray*}
\int_{t_0}^t \phi^* (\tau, t_0)C^* (\tau) \bary (\tau) \, d\tau =
\int_{t_0}^{t} \phi^* (\tau, t_0)C^* (\tau) C (\tau) \phi (\tau, t_0)
\, d\tau x_0.
\end{eqnarray*}
Thus
\begin{equation}
\underbrace{\int_{t_0}^{t_1} \phi^* (\tau, t_0)C^* (\tau) \bary (\tau)
\, d\tau}_{\hbox{\tiny known $n \times 1$ vector}} \ = \
\underbrace{\vphantom{\int_{t_0}^t} H (t_1,
t_0)}_{\parbox{1.4cm}{\tiny known $n \times n$ \newline matrix}} x_0.
\elabel{eq2.3.2}
\end{equation}
If $H$ is nonsingular, $x_0$ is uniquely determined by \eq eq2.3.2/,
and we conclude that the model is observable.  We now prove the
converse:

\head{(ii)} Assuming the model is observable, we prove that $H$ is
nonsingular by contradiction. If $H$ is singular, then there exists
$\alpha \neq 0$ such that
\begin{eqnarray*}
H (t_1, t_0) \alpha= \zero \qquad \hbox{and hence }\quad \alpha^* H
(t_1, t_0) \alpha = \zero.
\end{eqnarray*}
This may also be written as
\begin{eqnarray*}
\alpha^*\left(\int_{t_0}^{t_1} \phi^*(\tau)C^*(\tau) C(\tau)
\phi(\tau) \, d\tau \right) \alpha = 0
\end{eqnarray*}
or equivalently
\begin{eqnarray*}
\langle C \phi \alpha, C \phi \alpha\rangle_{L_2} = \| C \phi \alpha
\| ^2_{L_2} = 0.
\end{eqnarray*}
This implies that $C (t) \phi (t) \alpha = 0$ for all $t_0\le t\le
t_1$.  Since the two initial conditions $x (t_0) = x_0 + \alpha$ and
$x(t_0) = x_0$ yield the same output $\bary (t)$, we conclude that the
model is not observable.  \qed

\section{Duality}
By exploiting the similarities between controllability and
observability, we can generate further tests for observability.  This
approach is known as duality: observability and controllability are
dual concepts.  This can be made precise by defining a model which is
``dual'' to the model under study.

Given two models
\begin{quote}
\begin{description}
\item[(I)] $\displaystyle\begin{array}{rcl} \dot x &=& A (t) x + B (t)
u \\ y &=& C (t) x + D (t) u, \qquad x(t_0)\in \Re^n. \end{array}$
\medskip

\item[(II)] $\displaystyle\begin{array}{rcl} \dotz &=& \hbox to 0pt
{\hss $-\,$}A^* (t) z +C^* (t) v \\ \gamma &=& B^* (t) z + D^* (t) v
,\qquad z(t_0)\in \Re^n.\end{array}$
\end{description}
\end{quote} 
The second model is called the \defn{dual} of the first.

In the LTI case we know that $\rank (\clC_{II}) = n$ if and only if
(II) is controllable, where
\begin{eqnarray*}
\clC_{II} = [C^* \mid A^* C^* \mid A^{*2} C^* \mid \cdots \mid
A^{*n-1} C^*].
\end{eqnarray*}
Since we have
\begin{eqnarray*}
\clC_{II}^* = \left[ \begin{array}{c} C\\ CA \\ CA^2 \\ \vdots
\end{array} \right]  =  \mbox{observability matrix for model I},
\end{eqnarray*}
it is obvious that the controllability test for model II is equivalent
to the observability test for I in the LTI case.  This is a special
case of the following result.

\begin{theorem}
\tlabel{duality} For the general LTV model, \balphlist
\item
The LTV model (I) is controllable at $t_0$ if and only if (II) is
observable at $t_0$.
\item
The model (I) is observable at $t_0$ if and only if (II) is
controllable at $t_0$.
\end{list}
\end{theorem}

\proof We just prove (a) since the proof of (b) is identical.

\Exercise{adjoint-stm} of \Chapter{chap1.6} gives the relation
\begin{eqnarray*}
\phi_{II}(t_1,t_0) = \phi_I(t_0,t_1)^*.
\end{eqnarray*}
Hence the controllability grammian for I and the observability
grammian for II are identical:
\begin{eqnarray*}
W_I &\eqdef & \int_{t_0}^{t_f} \phi_I (t_0, \tau) B_I (\tau) B_I^*
(\tau) \phi_I^* (t_0, \tau)\, d\tau \\ &=& \int_{t_0}^{t_f}
\phi_{II}^* (\tau, t_0) C^*_{II} (\tau) C_{II} (\tau) \phi_{II}
(\tau,t_0)\, d\tau \\ &=& H_{II}.
\end{eqnarray*}
In view of \Theorem{grammian} and \Theorem{obs-grammian} the proof of
(a) is complete.  \qed
 
As a direct application of \Theorem{duality} and \Theorem{Hautus} we
have the following dual of the Hautus-Rosenbrock test.

\begin{theorem}
\tlabel{dualHautus} The LTI model is observable if and only if
\begin{eqnarray*}
\rank [sI - A^* \mid C^*] = n \mbox{ for any } s\in\Co.
\end{eqnarray*}
\qed
\end{theorem}

We then call a mode $\lambda$ unobservable if $\rank [\lambda I - A^*
\mid C^*]< n$.

\begin{figure}
\ebox{0.8}{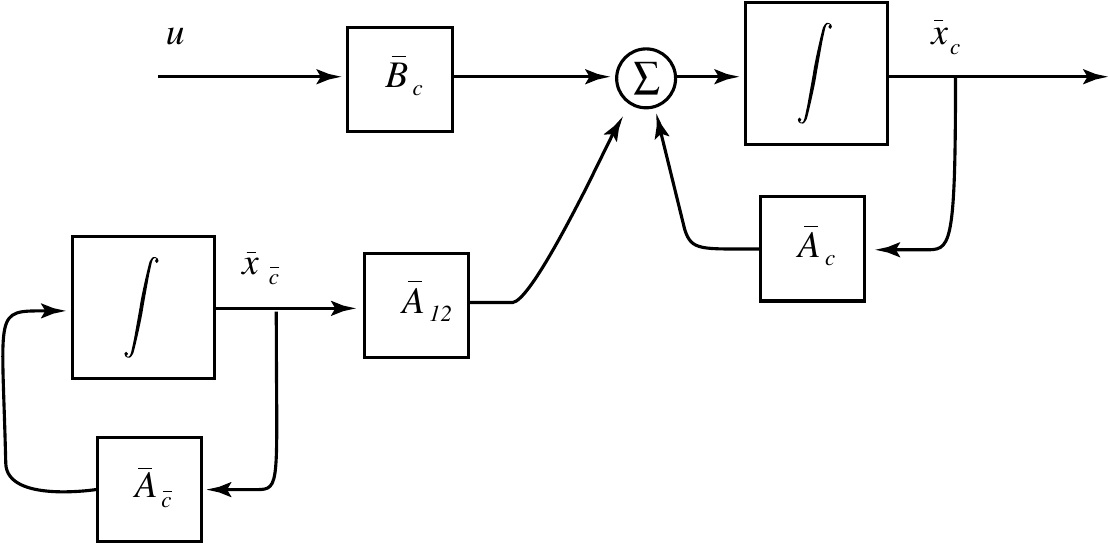}
\caption{The Kalman Controllability Canonical Form}
\flabel{structure-KalmanCCF}
\end{figure}

\section{Kalman canonical forms}
 If the model is not controllable, then we can construct a state
transformation $P$ in which the state $\barx = Px$ has the form $\barx
= \left(\begin{smallmatrix}\barx_c\\ \barx_{\bar
c}\end{smallmatrix}\right)$, as shown in \Figure{structure-KalmanCCF}.
After this transformation, the state space model takes the form
\begin{eqnarray*}
\dot{\barx} & = & \left[ \begin{array}{cc} A_c & A_{12} \\ 0 &
A_{\barc} \end{array} \right] \left[ \begin{array}{c} \barx^c \\
\barx^{\barc} \end{array} \right] + \left[ \begin{array}{c} B_c \\ 0
\end{array} \right] u
\\ y& = & [C_c \ \ C_{\barc} ] \left[ \begin{array}{c} \barx^c \\
\barx^{\barc} \end{array} \right].
\end{eqnarray*}
This is called the \defn{Kalman Controllability Canonical Form}
(KCCF).  In the case of distinct eigenvalues, the eigenvalues of
$(A_c)$ are precisely those for which $\rank [sI - A\mid B] = n$.  The
matrix $P$ can be defined as
\begin{eqnarray*}
P^{-1} = \left[ \left. \begin{array}{c} \leftarrow n_1 \rightarrow \\
\mbox{linear independent} \\ \mbox{columns} \\ \mbox{of}\ \clC
\end{array} \right| \begin{array}{c} \leftarrow (n - n_1) \rightarrow
\\ \mbox{set of vectors} \\ \mbox{linearly independent of } \\
\mbox{columns of} \ \clC
\end{array} \right] 
\end{eqnarray*}
where $\rank (\clC) = n_1 < n$ (see Brogan, \S 11.7).

\begin{ex}
Consider the LTI model with
\begin{eqnarray*}
A = \left[ \begin{array}{rr} -2 & 0 \\ -1 & -1 \end{array} \right],
\quad B = \left[ \begin{array}{r} 1 \\ 1 \end{array} \right]
\end{eqnarray*}
The controllability matrix is $\clC = \left[\begin{matrix}1 & -2\\ 1&
-2\end{matrix}\right]$, which has rank one.  Hence the model is not
controllable.  Define the matrix $P$ by
\begin{eqnarray*}
P^{-1} = \left[ \left. \begin{array}{c} 1 \\ 1 \end{array} \right|
\begin{array}{c} 1 \\ 0 \end{array} \right]
\end{eqnarray*}
where the column $\left(\begin{smallmatrix}1\\
1\end{smallmatrix}\right)$ was taken as the first column of the
controllability matrix, and the column $\left(\begin{smallmatrix}1\\
0\end{smallmatrix}\right)$ was chosen to make $P^{-1}$ full
rank. Thus,
\begin{eqnarray*}
P = \left[ \begin{array}{rr} 0 & 1 \\ 1 & - 1 \end{array} \right]
\end{eqnarray*}
The Kalman controllability canonical form can be written as
\begin{eqnarray*}
\barA = PAP^{-1} = \left[ \begin{array}{rr} -2 & -1 \\ 0 & -1
\end{array} \right] 
\qquad \barB = PB = \left[ \begin{array}{c} 1 \\ 0 \end{array} \right]
\end{eqnarray*}  
\end{ex}

For general LTI models, the observability canonical form can be
obtained using duality: Using a state transformation, the transformed
state may be written $\barx=\left(\begin{smallmatrix} {\barx}_o \\
\barx_\baro  \end{smallmatrix}\right)$ with
\begin{eqnarray*}
\dot x & = & \left[ \begin{array}{ll} A_O & 0 \\ A_{21} & A_{\baro}
\end{array} \right] \left[ \begin{array}{c} \barx_o \\ \barx^{\baro}
\end{array} \right] + \left[
\begin{array}{c} B_o \\ B_{\baro} \end{array} \right] u \\
y & = & [C_o\ \ 0] \left[ \begin{array}{c} \barx^o \\ \barx^{\baro}
\end{array} \right]
\end{eqnarray*}
This \defn{Kalman Observability Canonical Form} (KOCF) is illustrated
pictorially in \Figure{structure-KalmanOCF}.

\begin{figure}
\ebox{0.8}{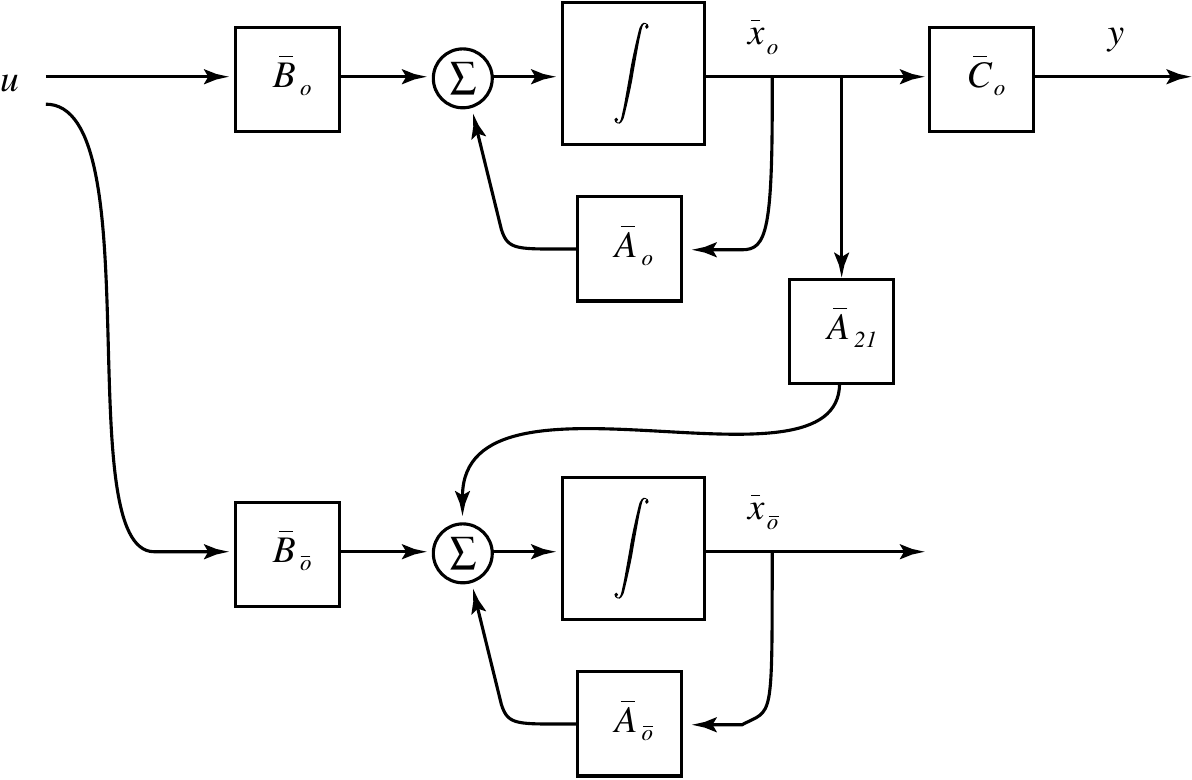}
\caption{The Kalman Observability Canonical Form}
\flabel{structure-KalmanOCF}
\end{figure}
 
\section{State space models and their transfer functions}
It should now be clear that given a transfer function, there are many
possible state space descriptions.  Moreover, the dimension of the
realization is not unique since we can always append uncontrollable or
unobservable dynamics.  A state space realization of a transfer
function $P(s)$ is called \defn{minimal} if it has minimum dimension
over all state space models for which
\begin{eqnarray*}
C (sI - A)^{-1} B + D = P (s).
\end{eqnarray*}

\begin{theorem}
\tlabel{min-controllable-observable} A realization $(A, B, C, D)$ of a
transfer function $P(s)$ is minimal if and only if it is both
controllable and observable.
\end{theorem} 

\proof For a general proof, see \BRO, \#12.17, page. 439. We prove the
theorem here in the special case where the eigenvalues of $A$ are
distinct.  When the eigenvalues are distinct, we have the partial
fraction expansion
\begin{equation}
P (s) = D + \sum_{i=1}^u \frac{k_i}{s -\lambda_i} \elabel{PFE0}
\end{equation}
where in the MIMO case, the $\{k_i\}$ are $p\times m$ matrices.  By
considering the modal form, the transfer function can also be written
as
\begin{equation}
P (s) = \barD + \barC (sI - \barA)^{-1} \barB \elabel{PFE}
\end{equation}where
\begin{eqnarray*}
\begin{array}{rclrcl}
\barA & = & \Lambda \qquad & \barB & = & M^{-1} B \\ \barC & = & CM
\qquad &\barD & = & D
\end{array}
\end{eqnarray*}
This latter form will allow us to compute the gains $\{k_i\}$ in \eq
PFE0/.  First break up the matrices $\barC$, $\barB$ into their
component columns and rows, respectively:
\begin{eqnarray*}
\barC & = & [ \gamma^1 \cdots \gamma^n] \\ \barB & = & [ \beta^1
\cdots \beta^n]^T = \left[ \begin{array}{c} \beta^{1T} \\ \vdots
\\\beta^{nT}
\end{array} \right] 
\end{eqnarray*}
Then the expansion \eq PFE/ can be written explicitly as
\begin{eqnarray*}
P (s) = \barD + \sum_{i=1}^n \frac {\gamma^i \beta^{iT}}{s-\lambda_i}.
\end{eqnarray*}
The model is controllable if and only if each $\beta_i$ is non-zero,
and it is observable if and only if each $\gamma_i$ is non-zero.  From
the equation above, it follows that the modal form is controllable and
observable if and only if none of the poles $\lambda_i$ are cancelled
in this sum.  That is, there are no pole/zero cancellations in the
expression $C (sI -A)^{-1} B + D$.  Thus, the transfer function has
$n$ poles, and hence any realization must have dimension no less than
$n$.  \qed

Consider for example
\begin{eqnarray*}
P (s) = \frac{s+1 }{s^2 + 6s + 5} = \frac{s+1 }{(s+5)(s+1)}.
\end{eqnarray*}
The second order modal realization is not minimal, but the first order
realization is.
 
\section{Realization of MIMO transfer functions}

For a SISO model with transfer function $P$ it is straightforward to
obtain a minimal realization. First, cancel any common pole-zero
pairs, count the number of poles remaining, say $n$, and construct any
$n$-dimensional state space model. One can for example use the CCF.
For a MIMO transfer function the situation is more complex.  We refer
the reader to \CHE\ for details - here we provide an example to
illustrate some of the issues involved.

Consider the transfer function $P$ which describes a $2$-input,
$2$-output model ($m=p=2$):
\begin{eqnarray*}
P(s)=\left[\begin{array}{cc} \frac{1}{s+1} & \frac{2}{s+1} \\
\frac{-1}{(s+1)(s+2)} & \frac{1}{s+2}
\end{array}\right].
\end{eqnarray*} 
By examining each entry of the matrix $P$ we see that there are no
direct pole-zero cancellations.  Since there are two poles, one might
expect that a realization can be constructed of degree two.  On the
other hand, since $P$ consists of three first-order, and one
second-order transfer functions, one might suspect that a degree four
state space model will be required.  It turns out that for this model
both predictions are incorrect.  The degree of any minimal realization
is \textit{three}, as we now show.

\begin{figure}[ht]
\ebox{.85}{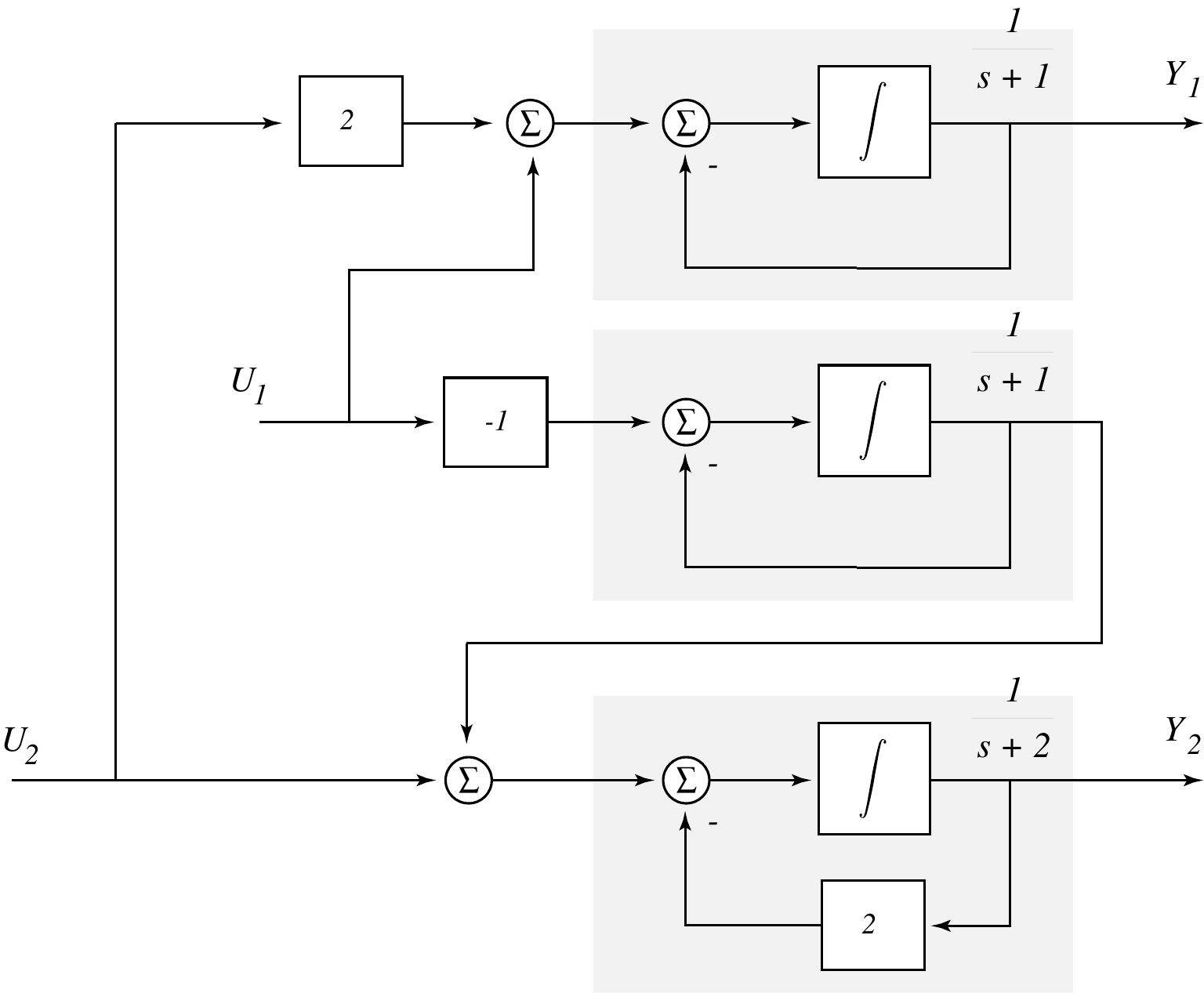}
\caption{A block diagram realization of $P$.}
\flabel{structure-realization}
\end{figure}

The input-output description $Y(s) = P(s) U(s)$ may be written as
\begin{eqnarray*}
Y_1(s) &=& \frac {1}{s+1} [U_1(s) + 2 U_2(s)] \\ Y_2(s) &=& \frac{1}{
s+2} [- \frac{1}{s+1}U_1(s) + U_2(s)]
\end{eqnarray*}
Hence, using just three first order filters we may simulate this model
using the block diagram description shown in
\Figure{structure-realization}.  Taking the outputs of the integrators
as states, we obtain the following state space realization:
\begin{eqnarray*}
\dot x &=& \left[\begin{array}{rrr} -1 & 0 & 0 \\ 0 & -1 & 0 \\ 0 & 1
& -2 \end{array}\right] x + \left[\begin{array}{rr} 1 & 2 \\ -1 & 0 \\ 0
& 1 \end{array}\right] u \\ y &=& \left[\begin{array}{rrr} 1&0&0 \\
0&0&1 \end{array}\right] x.
\end{eqnarray*}
By computing the controllability and observability matrices one can
show that this model is both controllable and observable.  Hence, this
is a minimal realization.

A general approach to obtaining a minimal realization for a MIMO
transfer function is to mimic the partial fractions expansion used in
the SISO case to form the modal form.  This will again yield a
diagonal, minimal realization.  For the transfer function $P$ given
above we may write
\begin{eqnarray*}
P(s) = \frac{1}{(s+1)(s+2)}\left[\begin{array}{rr} s+2 & 2(s+2) \\ -1
& s+1
\end{array}\right].
\end{eqnarray*}
Performing a partial fractions expansion on each term we obtain
\begin{eqnarray*}
P(s) = \frac{1}{(s+1)} \left[\begin{array}{rr} 1 & 2 \\ -1 & 0
\end{array}\right]
+ \frac{1}{(s+2)} \left[\begin{array}{rr} 0 & 0 \\ 1 &1
\end{array}\right] 
\eqdef \frac{1}{(s+1)} R_1 + \frac{1}{(s+2)} R_2.
\end{eqnarray*}
The rank of $R_1$ is two, and the rank of $R_2$ is one.  The sum of
these ranks gives the degree $3$ of the minimal realization.

To compute the diagonal realization let
\begin{eqnarray*}
R_1 = \underbrace{C_1}_{2\times 2}\underbrace{B_1}_{2\times 2} \qquad
R_2 = \underbrace{C_2}_{2\times 1}\underbrace{B_2}_{1\times 2}
\end{eqnarray*}
For instance, we can take
\begin{eqnarray*}
R_1 &=& \left[\begin{array}{rr} 1 & 2 \\ -1 & 0
\end{array}\right] = \left[\begin{array}{rr} 1 & 2 \\ -1 & 0
\end{array}\right] \left[\begin{array}{rr} 1 & 0 \\ 0 & 1
\end{array}\right] \\ R_2  &=& \left[\begin{array}{rr}
0 & 0 \\ 1 &1 \end{array}\right] = \left[\begin{array}{r} 0 \\ 1
\end{array}\right] \left[\begin{array}{rr} 1 & 1 \end{array}\right] 
\end{eqnarray*}
With this factorization, the transfer function $P$ may be expressed
\begin{eqnarray*}
P(s) &=& \frac{1}{(s+1)} C_1B_1 + \frac{1}{(s+2)} C_2 B_2 \\ &=& C_1
\left[\begin{array}{cc} \frac{1}{(s+1)} & 0 \\ 0 &\frac{1}{(s+1)}
\end{array}\right]  B_1 + C_2\left [ \frac{1}{(s+2)}\right]B_2
\\ \\ &=& \left[ C_1 \mid C_2\right] \left[
\begin{array}{ccc} \frac{1}{(s+1)} & 0 & 0 \\ 0 &\frac{1}{(s+1)} & 0
\\ 0 & 0 & \frac{1}{(s+2)} \end{array} \right] \left[\begin{array}{c}
B_1 \\ B_2 \end{array}\right] \\ \\ &=& \left[\begin{array}{ccc} 1 & 2
& 0 \\ -1 & 0 & 1 \end{array}\right] \left[\begin{array}{ccc}
\frac{1}{(s+1)} & 0 & 0 \\ 0 &\frac{1}{(s+1)} & 0 \\ 0 & 0 &
\frac{1}{(s+2)} \end{array}\right] \left[\begin{array}{cc} 1 & 0 \\ 0
& 1 \\ 1 & 1 \end{array}\right] \\ &=& C(Is - A)^{-1} B.
\end{eqnarray*}
These values of $A,B$ and $C$ define a minimal realization of the
transfer function.

In summary, when a simple partial fractions expansion is possible,
then a minimal realization is obtained by a factorization of the
matrix gains found in this expansion.  In general, this approach
might fail if one of the denominator polynomials in $P$ possesses a
repeated root. However, a minimal realization may still be found by a
more general factorization.  Write the transfer function $P$ as
\begin{eqnarray*}
P(s) = N_R(s) D_R(s)^{-1}
\end{eqnarray*}
where $N_R,D_R$ are matrices whose entries are polynomial functions of
$s$.  This is known as a right factorization.  A left factorization
can be defined analogously: $ P(s) = D_L(s)^{-1} N_L(s)$.  Using the
right factorization one may construct a generalization of the
controllable canonical form, and using the left factorization an
observable canonical form may be constructed.  To ensure minimality,
one must impose a co-primeness condition on the matrices $N$ and $D$ -
for details, the reader is referred to \CHE.

\begin{exercises} 
\item
Consider the transfer function
\begin{eqnarray*}
P(s)= \frac{s+1}{s^2+3s +2}
\end{eqnarray*}
\balphlist
\item Obtain a second-order state space realization in controllable
canonical form.  Is your realization controllable?  Observable?

\item Obtain a second-order state space realization in observable
canonical form.  Is your realization controllable?  Observable?

\end{list}
\item
Given the linear time-invariant model
\begin{eqnarray*}
\dot x &=& \left[\begin{matrix}-7&-2&6\\ 2&-3&-2\\
-2&-2&1\\\end{matrix}\right]x+ \left[\begin{matrix}1&1\\ 1&-1\\
1&0\\\end{matrix}\right] u\quad = Ax+Bu \\
y&=&\left[\begin{matrix}-1&-1&2\\ 1&1&-1\\\end{matrix}\right]x = Cx
\end{eqnarray*}
check observability using \balphlist
\item The observability matrix
\item 
columns of $\barC=CM$, where $M$ is chosen such that $M^{-1}AM$ is
diagonal

\item
the Hautus-Rosenbrock test.
\end{list}

\item 
Prove the following assertion: If the state space model
\begin{eqnarray*}
\dot x = A x + Bu,\qquad y = Cx + Du
\end{eqnarray*}
is a minimal realization of the transfer function $G(s)$, and if the
model is BIBO stable, then the state space model $\dot x = A x$ is
asymptotically stable.

\textit{For example}: Take the BIBO stable model $Y(s) =
(s-1)/[(s-1)(s+1)] U(s)$.  Any second order realization will have
characteristic polynomial $s^2 -1$, and hence the second order model
$\dot x = A x$ cannot be asymptotically stable. A \textit{minimal}
realization will be first order with $A=-1$, so that $\dot x = A x$ is
asymptotically stable.

\item
Show that the observability grammian given by \eq obsGram/ is equal to
$\clB^*\circ\clB$, where $ \clB(x) = C(t) \phi (t, t_0) x$.
\end{exercises}

\part{Feedback}

\chapter{Pole Placement}

In this chapter, we discuss a number of feedback control issues for the  LTI state
space model
\begin{equation}
\elabel{feedback-ss}
\begin{array}{rcl}
\dot x & = & Ax + Bu \\ y & = & Cx + Du.
\end{array}
\end{equation}
We consider initially the regulation problem: By translating the state
space, we may assume that we are regulating to the value $x=\zero$.
To perform this translation we must assume that the DC gain of the
plant is known exactly, and that there are no disturbances.  We will
develop more ``robust'' designs which work in spite of disturbances later in
\Section{IMP}.

The following example serves to illustrate most of the concepts
introduced in this chapter.  Consider the special case of \eq
feedback-ss/,
\begin{equation}
\begin{array}{rcl}
\dotx&=&Ax+Bu=\left[\begin{matrix}1&0&0\\ 0&2&0\\
0&0&-1\\\end{matrix}\right]x +\left[\begin{matrix}1\\ 1\\
0\\\end{matrix}\right] u \\
y&=&Cx=\left[\begin{matrix}1&0&1\end{matrix}\right]x
\end{array}
\elabel{star1}
\end{equation}

The system has three modes:
\begin{center}
\begin{tabular}{l >{\itshape}l}
$\lambda_1=1$ \qquad & controllable and observable.\\ $\lambda_2=2$&
controllable but not observable.\\ $\lambda_3=-1$& observable but not
controllable.
\end{tabular}
\end{center}
Define the control as $u=-Kx+r$.  This is \defn{linear state
feedback}, with an auxiliary input $r$.  With this controller, the
\defn{closed-loop system} becomes
\begin{equation}
\elabel{star}\begin{array}{rcl} \dot x&=&(A-BK)x+Br = A_{cl}x+Br\\
y&=&Cx
\end{array}
\end{equation}
where $A_{cl}$, the closed-loop system matrix, can be explicitly
written as
\begin{eqnarray*}
A_{cl}=\left[\begin{matrix}1 - k_1& - k_2 & - k_3\\ - k_1&2 - k_2& -
k_3\\ 0&0&-1\\\end{matrix}\right]
\end{eqnarray*}
with $K=\left[\begin{matrix}k_1& k_2& k_3\end{matrix}\right]$.  We now
study properties of the closed-loop system \eq star/, such as
controllability, observability, stability, as a function of $K$.

\head{Fact 1.} The controllability matrix
\begin{eqnarray*}
\clC =[ B \mid  A_{cl}B \mid A_{cl}^2B ]
\end{eqnarray*}
has rank $2$ for all $k_1,k_2,k_3$.  Hence,
\begin{quote}
\textit{ The rank of $\clC$ is invariant under state feedback}.
\textit{ The controllable subspace is invariant under state feedback.}
\end{quote}

\head{Fact 2.} Using the \textit{Hautus-Rosenbrock test}, one may show
that $\lambda=-1$ is the eigenvalue corresponding to the uncontrollable mode of
$A_{cl}$, independent of $K$.  The other two modes of $A_{cl}$ are controllable for
all $K$, but their values depend explicitly on $K$. Hence,
\begin{quote}
\textit{Uncontrollable modes are invariant under state feedback.}
\end{quote}

\head{Fact 3.} The observability matrix for the controlled system is
\begin{eqnarray*}
\clO=\left[\begin{matrix} C\\ CA_{cl}\\ CA_{cl}^2
\end{matrix}\right] = \left[\begin{matrix}
1&0&1\\ 1 - k_1& - k_2& - k_3-1\\ (1 - k_1)^2+k_1k_2& - k_2(3 - k_1 -
k_2)& 1+k_3(k_1+k_2) \end{matrix}\right].
\end{eqnarray*}
This matrix \textit{does not} have a fixed rank. It could
be
\begin{description}
\item \textit{one}\quad (corresponding to $k_1=1, k_2=0,k_3=-1$),
\item \textit{two}\quad (corresponding to $k_1= k_2=k_3=0$),
\textit{or}
\item \textit{three}\quad (corresponding to $k_1=k_3=0, k_2=-1$)
\end{description}

Hence,
\begin{quote}
\textit{Observability is not invariant under state feedback}.
\end{quote}

\head{Fact 4.} The eigenvalues of $A_{cl}$ can be found by solving the
characteristic equation
\begin{eqnarray*}
\det (\lambda I-A_{cl})=0 \Rightarrow (\lambda +1)(\lambda^2-(3 - k_1
- k_2)\lambda +2 - 2k_1 - k_2)=0.
\end{eqnarray*}
One eigenvalue is $\lambda=-1$, the other two depend on $k_1$ and
$k_2$.  Since the two terms $k_1+k_2$ and $2k_1+k_2$ are linearly
independent, these two eigenvalues can be chosen arbitrarily (subject
to the restriction that complex eigenvalues can appear only as
conjugate pairs) by proper choices of real scalars $k_1$ and
$k_2$. Hence,
\begin{quote}
\textit{Arbitrary pole placement is possible by linear state feedback,
in the controllable subspace.}
\end{quote}

The LTI model \eq star1/ is \defn{stabilizable} by linear state
feedback, since the uncontrollable mode is asymptotically stable.
\medskip

After gaining an understanding of the state feedback control problem
we will turn to output feedback of the form $u=-Gy + r$.  Even if an
LTI system is both controllable and observable, it is in general not
possible to assign the closed-loop poles arbitrarily by static output
feedback.  We will see shortly that all poles can be assigned
arbitrarily ({\it subject to the complex conjugate pair restriction})
if we  use \textit{dynamic} output feedback, where $G$ is a
transfer function rather than a static matricial gain.

\section{State feedback}

Returning to the regulation problem, assuming that we can measure all
of the states, the simplest controller is given by the state feedback
control law
\begin{eqnarray*}
u = - Kx,
\end{eqnarray*}
where $K$ is a $n\times m$ matrix.  This gives rise to the closed-loop
system
\begin{eqnarray*}
\dot x = (A - BK) x = A_{cl} x.
\end{eqnarray*}
To determine whether or not $x(t)\to \zero$ from any initial
condition, we must consider the eigenvalues of the closed-loop matrix
$A_{cl}$.

Consider the SISO case, where we initially assume that the pair
$(A,B)$ is in CCF (controllable canonical form):
\begin{eqnarray*}
A = \left[ \begin{array}{ccccc} 0 & 1 & 0 & \cdots & 0 \\ 0 & 0 & 1 &
\cdots & 0 \\ \vdots & \vdots &\ddots& \ddots & \vdots\\ 0 & 0 &
\cdots & 0 & 1 \\ - \alpha_n & -\alpha_{n-1}& \cdots & \cdots & -
\alpha_1 \end{array} \right], \qquad B = \left[ \begin{array}{c} 0 \\
0 \\ \vdots \\ 0 \\ 1 \end{array} \right].
\end{eqnarray*}
Recall that \balphlist
\item
The characteristic polynomial of this model is
\begin{eqnarray*}
\Delta(s) = \det (sI - A) = s^n + \alpha_1 s^{n-1} + \cdots + \alpha_n
\end{eqnarray*}
\item
$(A,B)$ is a controllable pair.
\end{list}

Writing the control as $u = - k_1 x_1 - k_2 x_2 \cdots - k_n x_n$, it
is immediate that the closed-loop system has the form
\begin{eqnarray*}
A - BK = A_{cl} = \left[ \begin{array}{ccccc} 0 & 1 & 0 & \cdots & 0
\\ 0 & 0 & 1 & \cdots & 0 \\ \vdots & \vdots &\ddots& \ddots &
\vdots\\ 0 & 0 & \cdots & 0 & 1 \\ - (\alpha_n + k_1) & -
(\alpha_{n-1} + k_2) & \cdots & \cdots & - (\alpha_1 + k_n)
\end{array} \right],
\end{eqnarray*}
so that the closed-loop characteristic polynomial is
\begin{eqnarray}
\Delta(s)= \det (sI - [A-BK]) = s^n + (\alpha_1 + k_n) s^{n-1} +
\cdots +(\alpha_n + k_1). \elabel{eq3.1.1}
\end{eqnarray}
We can assign the coefficients of this polynomial, and hence the
closed-loop poles can be assigned arbitrarily.

This analysis can be applied to any controllable SISO system.  It is
performed on a system in CCF, but we can transform to CCF for any time
invariant SISO state space model, provided the model is
controllable. The question is, what is the state transformation $P$
that is required to ensure that $\bar x = P x$ is described by a model
in CCF?  The controllability matrices $\clC$ and $\bar\clC$ are
related by the formula
\begin{eqnarray*}
\bar\clC= P \clC.
\end{eqnarray*}
Hence the matrix $P$ is given by $P = \bar{\clC} \clC^{-1}$.  The
inverse exists if the original state space model is controllable.
System eigenvalues do not change with a state transformation, so by
converting to CCF and then applying pole placement to the transformed
model, one can place the poles as desired.  This is summarized in the
following theorem.  While this result was derived for SISO systems, it
is also true in the MIMO case.

\begin{theorem}
\tlabel{pole-place} The eigenvalues of $(A-BK)$ can be placed
arbitrarily, respecting complex conjugate constraints, if and only if
$(A,B)$ is a controllable pair.
\end{theorem}

All of this is applicable if $(A,B)$ is controllable.  What do we do
if $(A,B)$ is not controllable?  Consider the Kalman controllability
canonical form introduced previously in \Figure{structure-KalmanCCF},
and given by the equations below
\begin{eqnarray*}
\dot x & = & \left[ \begin{array}{cc} A_c & A_{12} \\ 0 & A_{\barc}
\end{array} \right] x + \left[ \begin{array}{c} B_c \\ 0 \end{array}
\right] u \\ x & = & \left[ \begin{array}{c} x_c \\ x_{\barc}
\end{array} \right].
\end{eqnarray*}
We can still apply a control law of the form
\begin{eqnarray*}
u = - [K_1 \ \ K_2] \left[ \begin{array}{c} x_c \\ x_{\barc}
\end{array} \right]
\end{eqnarray*}
to obtain the closed-loop system
\begin{eqnarray}
\dot x = \left[ \begin{array}{cc} A_c - B_c K_1 & A_{12} - B_c K_2 \\
0 & A_{\barc} \end{array} \right] x \elabel{eq3.1.2}
\end{eqnarray}
From the form of the closed-loop system, we find that the eigenvalues
are dependent only upon $(A_c - B_c K_1)$ and not on $A_{\barc}$.  In
fact, the characteristic polynomial for the closed-loop system becomes
\begin{eqnarray*}
\Delta(s)= \det (sI - A_{cl}) = \underbrace{\det (sI - (A_c - B_c
K_1))}_{\parbox{.6in}{\tiny arbitrary \newline Based on $K$}}
\underbrace{\det (sI - A_{\barc})}_{\parbox{.6in}{\tiny
independent\newline of $K$}}
\end{eqnarray*}
In other words, one can place the controllable modes arbitrarily, but
the uncontrollable modes remain unchanged.

The pair $(A,B)$ is said to be \defn{stabilizable} if there exists a
$K$ such that $(A-BK)$ is a Hurwitz matrix.  From these arguments we
see that
\begin{theorem}
\tlabel{stabilizable} If the eigenvalues are distinct, then $(A,B)$ is
stabilizable if and only if the eigenvalues of the uncontrollable
modes are in the strict LHP.

For a general LTI model, $(A,B)$ is stabilizable if and only if the
eigenvalues of $ A_{\barc}$ lie in the strict LHP, where $A_{\barc}$
is the matrix defined in the KCCF.
\end{theorem}

Geometrically, \Theorem{stabilizable} implies that the model is
stabilizable if and only if the unstable subspace is contained in the
controllable subspace.

\section{Observers}\label{s:Observers}
The state feedback approach can be generalized to the situation where
only partial measurements of the state are available.  Consider again
the LTI model
\begin{eqnarray*}
\dot x & = & Ax + Bu,\\ y & = & Cx.
\end{eqnarray*}
We would like to define an observer, of the form illustrated in
\Figure{feedback-observer}, in which input-output measurements are
collected on-line to give estimates $\hatx(t)$ of the state $x(t)$.
To mimic the behavior of the system one can try
\begin{eqnarray*}
\dot\hatx(t) = A \hatx(t) + Bu(t).
\end{eqnarray*}
Defining the error as $\tilx(t) = x(t) -\hatx(t)$, this gives the
error equation $\dot \tilx(t) = A\tilx(t)$, from which we deduce that
\begin{eqnarray*}
\tilx(t) = e^{At} \tilx(0).
\end{eqnarray*}
This is a poor approach since there is no flexibility in design.  If
for example the open-loop system is unstable, then for some initial
conditions the error will not converge to zero, and may diverge to
infinity.

\begin{figure}[ht]
\ebox{.55}{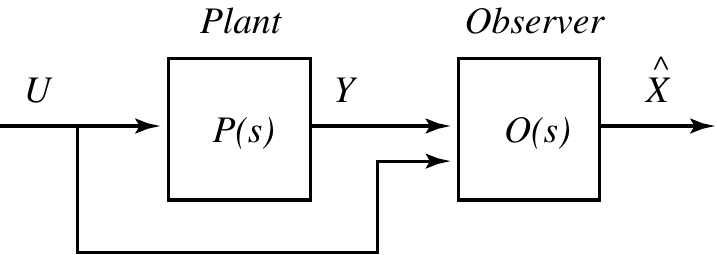}
\caption[Observer   design]{Can the observer $\clO$ be designed so that $\hat x(t)-
x(t)\to0$, as $t\to\infty$, at a predesignated rate?}
\flabel{feedback-observer}
\end{figure}

The observers we consider are of a similar form, but we adjoin an
output error term:
\begin{eqnarray*}
\dot{\hatx} = A \hatx + Bu + L (y - \haty), \qquad \hatx(0) \in\Re^n,
\end{eqnarray*}
where $\haty = C \hatx$.  For any fixed $n\times p$ matrix $L$ one
obtains
\begin{eqnarray*}
\dote= \dot x - \dot{\hatx} = Ax + Bu - A \hatx - Bu - L (Cx - C
\hatx) = (A - LC) e.
\end{eqnarray*}
To ensure that $\tilx(t)\to 0$ as $t\to\infty$, we must compute the
eigenvalues of the matrix $(A-LC)$ as a function of $L$. But note:
\begin{eqnarray*}
\eig (A-LC) = \eig (A^* - C^* L^*)
\end{eqnarray*}
We are thus exactly in the same position as when we considered pole
placement.  It follows that the eigenvalues of $(A-LC)$ can be placed
arbitrarily, provided that the matrix pair $(A^*,C^*)$ is
controllable. Based on duality, this is simply observability of the
pair $(A,C)$. Thus, we can place the observer poles arbitrarily if and
only if $(A,C)$ is observable.

This again raises the question, what do we do if this basic assumption
fails?  The pair $(A,C)$ is said to be \defn{detectable} if there
exists some $L$ such that $(A-LC)$ is a Hurwitz matrix.  This property
is the dual of stabilizability, in the sense that $(A,C)$ is
detectable if and only if $(A^*,C^*)$ is stabilizable.  We then obtain
the following dual statement to \Theorem{stabilizable}:

\begin{theorem}\tlabel{detectable}
If the eigenvalues are distinct, then $(A,C)$ is detectable if and
only if the eigenvalues corresponding to the unobservable modes lie in the strict LHP.

For a general LTI model, $(A,C)$ is detectable if and only if the
eigenvalues of $ A_{\baro}$ lie in the strict LHP, where $A_{\baro}$
is the defined in the KOCF.
\end{theorem}

A state space model is thus detectable if and only if the unobservable
subspace is contained in the stable subspace.

\section{Observer feedback}
The objective in the previous section was to construct a useful state
estimator for a state space model, but this of course was rooted in
the desired to control the system.  Suppose that we ignore that
$\hatx$ is an estimate, and we apply the observer feedback control law
\begin{eqnarray*}
u(t)=- K \hatx(t).
\end{eqnarray*}
This control law is designed based on the faith that $\tilx(t)$ will
converge to zero fast enough so that this control law is essentially
equivalent to full state feedback. To see if this faith is in vain, we
must first see if the overall system is asymptotically stable.  The
overall system is linear, so we can check this by computing the closed-loop eigenvalues.

\begin{figure}[ht]
\ebox{.75}{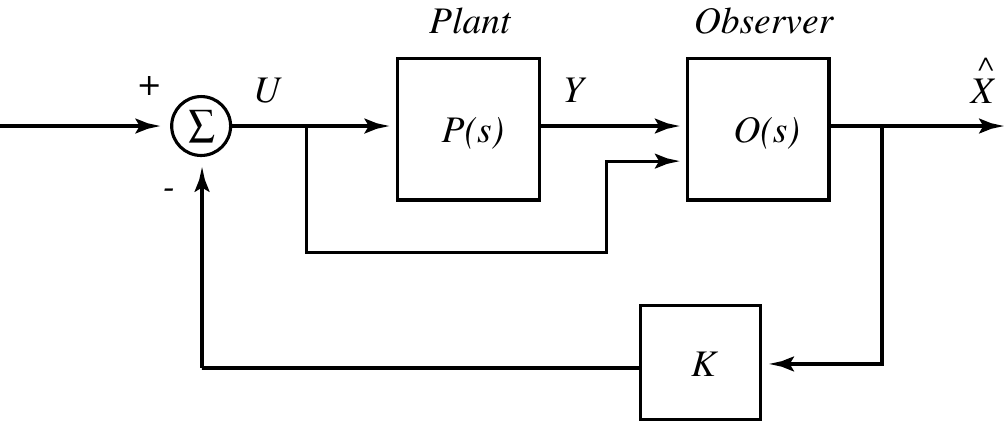}
\caption[The separation principle]{The separation principle justifies the use of state estimates
in place of the true state values in the feedback control law.
However, full state feedback is not completely equivalent to feedback
based on state estimation.}  \flabel{feedback-observerFbk}
\end{figure}

First apply the following state transformation for the overall state
\begin{eqnarray*}
\left[ \begin{array}{c} x \\ e \end{array} \right] & = &
\underbrace{\left[ \begin{array}{rr} I & 0 \\ I & -I \end{array}
\right]}_{P} \left[ \begin{array}{c} x \\ \hatx
\end{array} \right]
\qquad \qquad
\begin{array}{rcl}
e & = & x - \hatx \\ \hatx & = & x - e.
\end{array}
\end{eqnarray*}
The transformed state $ \left(\begin{matrix}x\\ e\end{matrix}\right)$
is described by the closed-loop equations
\begin{eqnarray*}
\left[ \begin{array}{c} \dot x \\ \dote \end{array} \right] & = &
\underbrace{\left[ \begin{array}{c|c}A-BK & BK \\ \hline 0 & A - LC
\end{array} \right]}_{A_{cl}}  \left[ \begin{array}{c} x \\ e
\end{array} \right]
\end{eqnarray*}
Thus, the closed-loop poles of the overall system are defined by the
characteristic polynomial
\begin{equation}
\det (sI - A_{cl}) = \det (sI - [A-BK]) \cdot \det (sI - [A-LC])
\elabel{separation}
\end{equation}
where
\balphlist
\item
$\det (sI - [A-BK]) \equiv $ state feedback eigenvalues -- arbitrary
placement if and only if $(A,B)$ is controllable.

\item
$\det (sI - [A-LC]) \equiv$ observer eigenvalues -- arbitrary
placement if and only if $(A,C)$ is observable.
\end{list}
Equation \eq separation/ is known as the \defn{separation principle}.
Its obvious consequence to control design is that the design of the
feedback gain can be conducted independently of the design of an
observer.  One must be cautious when interpreting this result however,
since we are only considering stability and the placement of closed-loop poles.  In particular, we are ignoring such issues as sensitivity
to plant uncertainty.

The separation principle allows us to separately place the state
feedback eigenvalues for good dynamic response, and the
observer-eigenvalues so they are faster than the desired closed-loop
response.  This still leaves open the question, where do we place the
observer poles?  A typical rule of thumb is that the observer poles
should be faster than the ``slowest'' state feedback pole by a factor
of $2$ to $5$.  This however is very problem specific.  In some
applications, the location of the observer poles is not very
important, so long as the poles are faster than those of the state
feedback poles.  In other examples, the observer must be designed with
care.

\section{Reduced-order (Luenberger) observers}     

Considering the state space equations \eq feedback-ss/ we see that at
least part of the state is directly observed through $y$.  Assume for
simplicity that $D=0$, and construct a state transformation $\barx =
Px$ so that $y$ is explicitly a part of the state.  This can be
accomplished if $C$ has rank $p$, since then we can define the matrix
$P$ so that
\begin{eqnarray*}
\barx = \left[ \begin{array}{c} y \\ \barx_2 \end{array} \right] =
\left[ \begin{array}{c} \frac{C\strut}{ \parbox{1.1in}{\tiny
\smallskip anything for \newline linear independence}}
\end{array}
\right] x.
\end{eqnarray*}
Then we do have a transformed model in the desired form:
\begin{eqnarray*}
\left[ \begin{array}{c} \dot{\barx}_1 \\ \dot{\barx}_2 \end{array}
\right] = \underbrace{\left[ \begin{array}{cc} A_{11} & A_{12} \\
A_{21} & A_{22} \end{array} \right]}_{PAP^{-1}} \left[
\begin{array}{c} \barx_1 \\ \barx_2 \end{array} \right] +
\underbrace{\left[ \begin{array}{c} B_1 \\ B_2 \end{array}
\right]}_{PB} u.
\end{eqnarray*}
Since by definition we have $\barx_1 = y$, the output equation is
evidently
\begin{eqnarray*}
y = [I \ \ 0] \left[ \begin{array}{c} \barx_1 \\ \barx_2 \end{array}
\right].
\end{eqnarray*}

The state estimation problem has been reduced in complexity through
this transformation: rather than construct an estimator of dimension
$n$, we only have to estimate the $(n-p)$ dimensional vector
$\barx_2$.  The next step then is to construct an observer for
$\barx_2$ alone.  To see how this can be accomplished, write the
transformed state equations as follows.  The equation for $\barx_2$
may be written in the suggestive form,
\begin{equation}
\dot{\barx}_2 = A_{22} \barx_2 + \underbrace{A_{21} y + B_2
u}_{\parbox{.8in}{\tiny known from \newline measurements}} .
\elabel{eq3.3.1}
\end{equation}
The equation for $\barx_1$ can be written
\begin{equation}
\underbrace{\doty - A_{11} y - B_1 u}_{\parbox{.8in}{ \tiny known from
\newline measurements}} =A_{12} \barx_2 \elabel{eq3.3.2}
\end{equation}
Defining $\baru = A_{21} y + B_2 u$ and $\bary=\doty - A_{11} y - B_1
u$, we obtain the state space model
\begin{eqnarray*}
\dot{\barx}_2 &=& A_{22} \barx_2 + \baru \\ \bary &=& A_{12} \barx_2
\end{eqnarray*}
We can now write down the corresponding observer to estimate
$\barx_2$. The only technicality is that $\bary$ contains the
derivative $\doty$, which strictly speaking is not known.  We will
assume for now that we know $\doty$, but through some manipulations of
the observer equations we will relax this assumption below. An
observer for $\barx_2$ takes the form
\begin{equation}
\dot{\hatx}_2 = A_{22} \hatx_2 + \baru + L (\bary - \hat{\bary})
\elabel{lue-obs}
\end{equation}
where $\hat{\bary}=A_{12} \hatx_2$.  Since $\bary = A_{12} \barx_2$,
the error equation for $e \eqdef \barx_2 - \hatx_2 $ takes the form
\begin{eqnarray*}
\dote = \dot{\barx}_2 - \dot{\hatx}_2 = (A_{22} - L A_{12}) e
\end{eqnarray*}
So we can place the observer poles arbitrarily if and only if the pair
$(A_{22}, A_{12})$ is observable.

We now show that it is unnecessary to differentiate the output
measurements.  First write the observer equations in the expanded form
\begin{eqnarray*}
\dot{\hatx}_2 = A_{22} \hatx_2 + A_{21} y + B_2 u + L (\doty - A_{11}
y - B_1 u - A_{12} \hatx_2) .
\end{eqnarray*}
These equations can be represented through the block diagram

\ebox{.8}{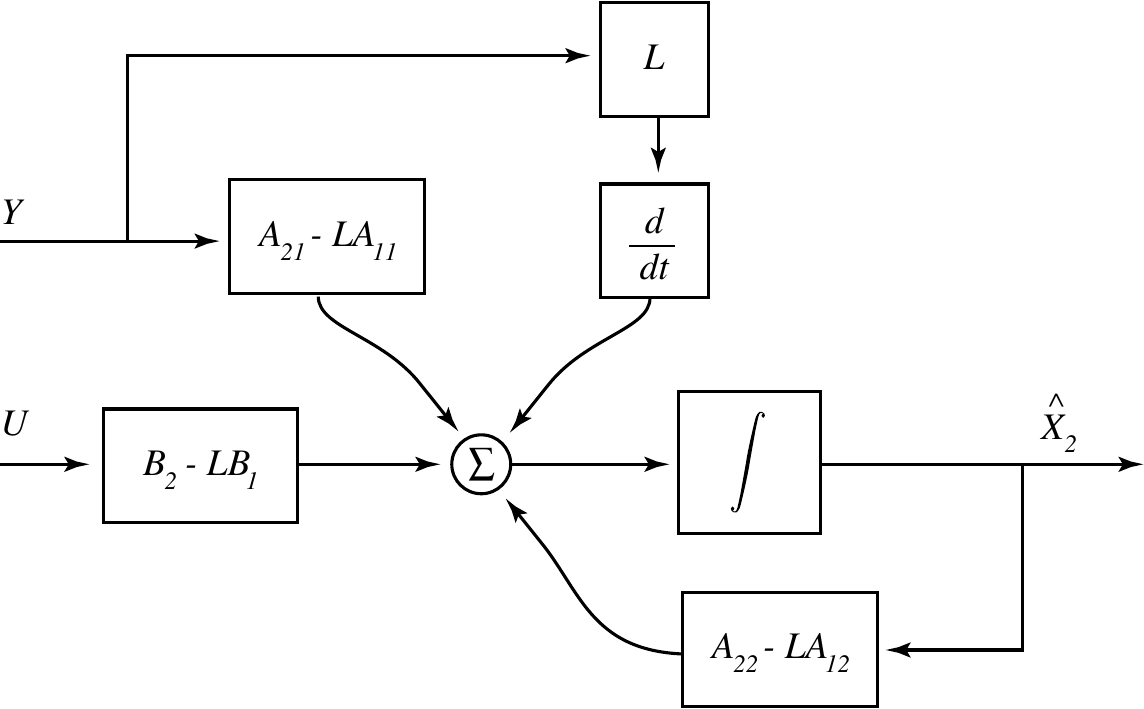}

It is evident in this diagram that the term $L\dot y$ is immediately
integrated.  To eliminate the derivative, we simply cancel the
derivative and integral to obtain the equivalent block diagram

\ebox{.8}{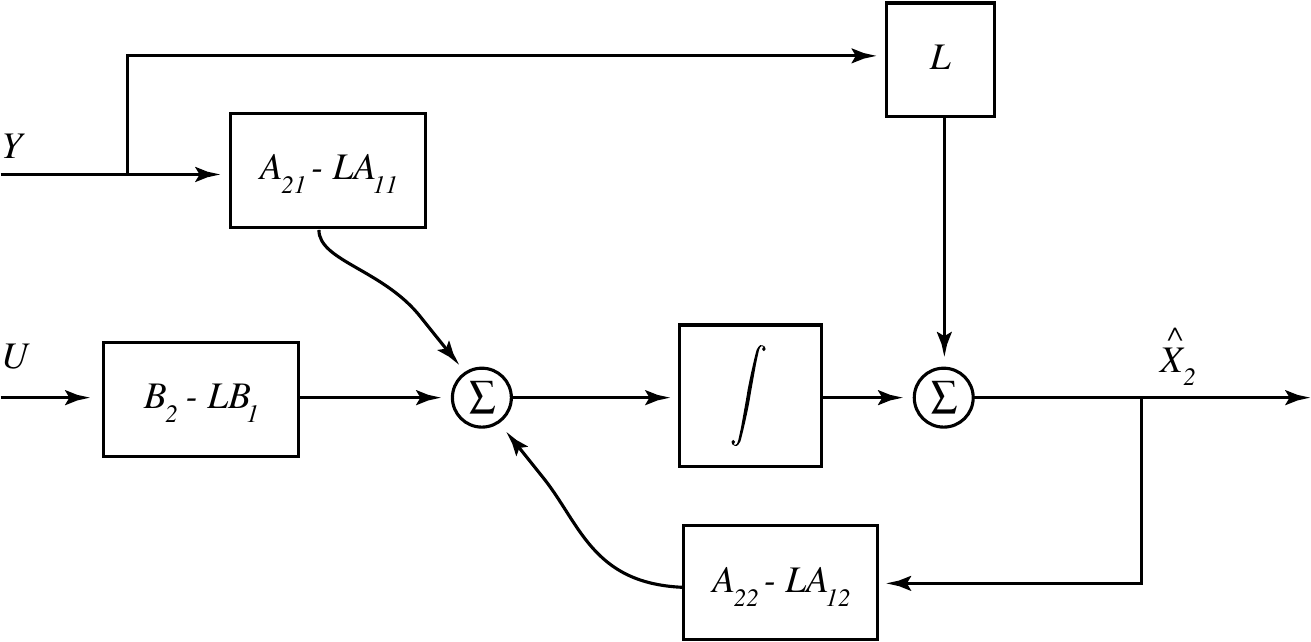}

We then obtain an estimate $\hatx$ of the original state through

\ebox{.35}{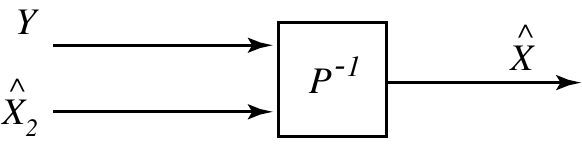}

In conclusion, provided $(A_{22}, A_{12})$ is an observable pair, the
state can be estimated based upon an $(n-p)$-dimensional observer
using the measurements $u,y$, and it is not necessary to differentiate
these measurements.  This is important, since if the measurement $y$
is noisy, its derivative may be much worse.

If we define
\begin{eqnarray*}
\barA & = & \left[ \begin{array}{rr} A_{11} & A_{12} \\ A_{21} &
A_{22} \end{array} \right] \\ \barC & = & [I\ \ 0]
\end{eqnarray*}
then by the Hautus-Rosenbrock test we know that $(\barA, \barC)$ is
observable if and only if
\begin{eqnarray*}
\rank \left[ \begin{array}{c} sI - \barA \\ \barC \end{array} \right]
= n \qquad \mbox{for any } s\in\Co .
\end{eqnarray*}
This matrix may be written as
\begin{eqnarray*}
\left[ \begin{array}{rr} sI - A_{11} & -A_{12} \\ -A_{21} & sI -
A_{22} \\ I & 0 \end{array} \right].
\end{eqnarray*}
The first $p$ columns are automatically independent because of the
identity matrix $I$.  Because of the $0$, the matrix is full rank if
the last $n-p$ columns are independent of each other.  That is, we
must have
\begin{eqnarray*}
\rank \left[ \begin{array}{c} -A_{12} \\ sI - A_{22} \end{array}
\right] = (n-p).
\end{eqnarray*}
This condition holds if $(A_{22}, A_{12})$ is observable, since this
is the Hautus test on the reduced complexity observer problem.  We
conclude that the condition for observability of the reduced order
observer is equivalent to observability of $(\barA, \barC)$, which is
equivalent to observability of the original pair $(A,C)$.

\begin{ex}
Consider the magnetically suspended ball, whose linearized and
normalized state space model is given by
\begin{eqnarray*}
\dotx = \begin{bmatrix} 0 & 1 \\ 1 & 0\end{bmatrix} x +
\begin{bmatrix} 0 \\ 1 \end{bmatrix} u,\qquad y = \begin{bmatrix}1 & 0
\end{bmatrix} x.
\end{eqnarray*}
Full state feedback is simply proportional-derivative (PD) feedback:
\begin{eqnarray*}
u=-Kx = -K_1 y -K_2 \doty.
\end{eqnarray*}
If the derivative cannot be measured directly, then there are several
possibilities: \balphlist
\item
The second state $x_2=\dot y$ could be obtained by differentiating
$y$, but this will greatly amplify any high-frequency noise.

\item
An approximate derivative could be obtained through a filter:
\begin{eqnarray*}
\dot{\hat{x}}_2 + \ell \hat{x}_2 = \ell \dot y,
\end{eqnarray*}
or in the frequency domain
\begin{eqnarray*}
\dot{\hat{X}}_2(s) = \frac{\ell s}{s + \ell} Y(s),
\end{eqnarray*}
where $\ell$ is a large constant.  The trouble with this approach is
that it does not make use of the system model, and for this reason one
cannot expect that $\hat{x}_2(t) - \dot y(t)\to 0$ as $t\to 0$.

\item
A reduced order observer can be constructed to estimate $\dot y$.  If
there is noise, the observer can be designed to obtain a good trade
off between tracking of $\dot y$, and rejection of this noise.  In the
noise free case we have that $\hat{x}_2(t) - \dot y(t)\to 0$, as long
as the observer gain is chosen to give a stable observer.
\end{list}

The reduced order observer is defined by the equations
\begin{eqnarray*}
\dot{\hat{x}}_2 = x_1 + u + \ell(\dot x_1 - {\hat{x}}_2).
\end{eqnarray*}
Taking transforms gives a formula which is similar to the estimator
(b) given above:
\begin{eqnarray*}
\dot{\hat{X}}_2(s) = \frac{\ell s +1}{\ell+ s } X_1(s) +
\frac{1}{\ell+ s } U(s).
\end{eqnarray*}
However, this observer does make use of system information, and hence
does exhibit desirable convergence properties.

The controller can be written in the frequency domain as
\begin{eqnarray*}
U(s) = -K_1 X_1(s) -K_2 \frac{\ell s +1}{\ell+ s } X_1(s) -K_2
\frac{1}{\ell+ s } U(s).
\end{eqnarray*}
Solving for $U$, it is readily seen that for large $\ell$ this is a
lead compensator.  That is, the above expression may be written as
\begin{eqnarray*}
U(s) = - g \frac{s-z}{s-p} Y(s)
\end{eqnarray*}
where $p < z < 0$. In conclusion, we see that in simple models,
observer based feedback generalizes classical control design.
\end{ex}

\begin{matlab}
\item
In Matlab there are two methods for assigning eigenvalues.  These can
be used to obtain the state feedback gain $K$, or in observer design
by placing the eigenvalues of $(A^* - C^* L^*)$ to compute the
observer gain $L$.

\item[ACKER] This command is named after Ackerman, who derived a
formula for explicit solution of the feedback gain $K$ using the
controller canonical form.  This is only suitable for SISO plants.

\item[PLACE] A more numerically robust algorithm which can be used for
MIMO systems.  This command cannot be used with repeated roots.
\end{matlab}

\begin{exercises}
\item
For the LTI model
\begin{eqnarray*}
\dot x = \left[ \begin{matrix} -1 & 2\\ 0 & 3\end{matrix} \right]x +
\left[ \begin{matrix} 1 \\ 2 \end{matrix} \right] u
\end{eqnarray*}
\balphlist
\item
Find the controllable subspace.
\item
Show that the system is unstable, and compute the unstable modes and
corresponding eigenvector.
\item
Is the system stabilizable?  Explain carefully using (a) and (b).  If
the system is stabilizable, find a stabilizing state feedback control
law.
\end{list}

\item
Consider the SISO, LTI system
\begin{eqnarray*}
\dot x = \left[\begin{matrix}1 \ 1 \\ 0\
3\end{matrix}\right] x + \left[\begin{matrix}1 \\
1\end{matrix}\right] u;\qquad y= [ 7\ \ 8] x
\end{eqnarray*}
\balphlist
\item
There is a matrix $P$ such that $(\bar A,\bar B)=(PAP^{-1}, PB)$ is in
controllable canonical form. Compute $P$ and $P^{-1}$ using the
formula $P = \bar \clC \clC^{-1}$.

\item
Compute by hand a feedback control law $\bar K=[\bar k_1,\bar k_2]$
which places the eigenvalues of $\bar A - \bar B \bar K$ at $-2$
(twice).

\item
Let $K= \bar K P$. Verify that the eigenvalues of $A-BK$ are at $-2$.
\end{list}

\item
Given the system
\begin{eqnarray*}
\dot x = \left[\begin{matrix} -1 & 4 \\ 1&
-1\end{matrix}\right] x + \left[\begin{matrix}2\\ 1
\end{matrix}\right]u,
\end{eqnarray*}
\balphlist
\item
Show that when $u\equiv \zero $ the origin is unstable.

\item
Determine the stable subspace $\Sigma_s$, the span of the eigenvectors
corresponding to unstable modes, $\Sigma_u$, and the controllable
subspace $\Sigma_c$.  Plot all on the same graph.

\item
Now choose the feedback control
\begin{eqnarray*}
u & = & (1\; \; \beta) x + r
\end{eqnarray*}
where $\beta$ is a scalar, and $r$ is a reference input.  Show that
the closed-loop system will be asymptotically stable for some value(s)
of $\beta$.

\item
What are the closed-loop poles?  Do different choices of $\beta$
affect the stable mode?
\end{list}

\item
Consider the single input/single output second-order linear
time-invariant system
\begin{eqnarray*}
\dot x &=& A x - \left[\begin{smallmatrix} 0\\
a\end{smallmatrix}\right] u \\ y &=& (1,b) x
\end{eqnarray*}
where $a$ and $b$ are scalar real parameters.

It is known that for some constant $c>0$,
\begin{eqnarray*}
A^T \left[\begin{smallmatrix} 1\\ c\end{smallmatrix}\right] =
\left[\begin{smallmatrix} -2\\ -2c\end{smallmatrix}\right], \qquad
\hbox{ and } \quad A \left[\begin{smallmatrix} 1\\
1\end{smallmatrix}\right] = \left[\begin{smallmatrix} -2\\
-2\end{smallmatrix}\right],
\end{eqnarray*}
where $A^T$ denotes the transpose of $A$.  Also, it is known that
$\lambda$, one of the eigenvalues of $A$, is strictly positive.

\balphlist
\item
Compute the left eigenvector {\it and} the right eigenvector of $A$
corresponding to the positive eigenvalue $\lambda$. {\it Hint: left
eigenvectors are reciprocal basis vectors for the basis consisting of
right eigenvectors}.

\item
For what values of $(a,b)$ is the system observable?

\item
For what values of $(a,b)$, if any, is this system BIBO (bounded-input
bounded-output) stable?  \end{list}

Explain your answers using appropriate system structural concepts.  It
may help to draw pictures of left and right eigenvectors: Think
geometrically!

\item
This problem is a follow-up to \Exercise{CtbleMode} of
\Chapter{Controllability}.  Consider again the LTI model
\begin{eqnarray*}
\dot x &=& \left[\begin{matrix}-3&1&0\\ 0&-3&0\\
0&0&4\\\end{matrix}\right] x + \left[\begin{matrix}0\\ 1 \\ 0
\\\end{matrix}\right] u \\ y &=& \left[1 \ 0 \ 1\right] x.
\end{eqnarray*}
Suppose now that the initial condition of the model is not known, and
that $x(t)$ is not observed directly.  Is it possible to choose a
control so that $y(t) = t e^{-3t}$ for $t>1$?

\item
The system $\dot x = A x + B u$; $y= C x$ is controlled using static
output feedback $u=-H y + v$.  Show that the resulting closed-loop
system
\begin{eqnarray*}
\dot x = (A-BHC)x + B v;\qquad y= C x
\end{eqnarray*}
has the same controllability/observability properties as the original
system.  \textit{Hint: use the Hautus-Rosenbrock test}.

\item
In this problem you will see that the feedback gain $K$ may not be
unique for multivariate models.  Consider the state space model
defined by the pair of matrices
\begin{eqnarray*}
A= \left[ \begin{array}{ccc} 0& 1& 0\\ 0& 1& 1\\ 0& 0& 0 \end{array}
\right] \qquad B=\left[\begin{array}{cc} 0& 0\\ 1& 0\\ 0& 1
\end{array}  \right]\qquad C=\left[\begin{array}{ccc} 1 & 0 & 0
\end{array}  \right]
\end{eqnarray*}
\balphlist
\item
Find vectors $w^a,w^b\in\Re^2$ such that each pair $(A,B_a)$,
$(A,B_b)$ is controllable, where
\begin{eqnarray*}
B_a = B w^a \qquad B_b= B w^b.
\end{eqnarray*}
Try to pick the vectors so that the two $3\times 1$ matrices $B_a,B_b$
are significantly different.

\item
Now, find two controllers of the form
\begin{eqnarray*}
u_a = - w^a K_a x \qquad u_b = - w^b K_b x
\end{eqnarray*}
so that in each case, the closed-loop poles are placed at $-1$ and
$-1\pm j$.  These designs should be performed by placing the poles of
$A-B_a K_a$, $A-B_bK_b$, respectively.

\item
Provide simultations of the step response $u(t)=-Kx(t) + (1, 1)^T$
with zero initial conditions for the two designs. Include in your
plots both $ x_1(t) $ and the two dimensional input $u(t)$, for $0\le
t\le 1$.
\end{list}

\item
A schematic of an active suspension system is illustrated below.
\ebox{.8}{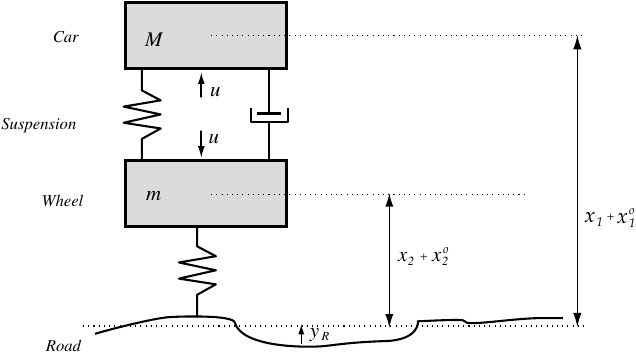} A state space model in scientific
units is given below, where $v_i= x_i'$.  The states $x_i$ have been
normalized so that their nominal values are $0$.
\begin{eqnarray*}
\frac{d}{dt} \left[ \begin{matrix} x_1\\ x_2\\ v_1\\ v_2\end{matrix}
\right] = \left[ \begin{matrix}0&0&1&0\\ 0&0&0&1 \\-10 &10 &-2 &2 \\
60 &-660 &12 &-12\end{matrix} \right] \left[ \begin{matrix} x_1\\
x_2\\ v_1\\ v_2\end{matrix} \right] +\left[ \begin{matrix} 0 \\ 0 \\
3.34 \\ -20\end{matrix} \right] u + \left[ \begin{matrix} 0 \\ 0 \\ 0
\\ 600\end{matrix} \right] y_R
\end{eqnarray*}

\balphlist
\item
Obtain an open-loop simulation of $x_1$ and $x_2$ for $0\le t \le 5 $
s., with $x_1(0)=0.5\, m$; $x_2(0)=0\, m$; $v_1(0)=v_1(0)=0$;
$u(t)=y_R(t)\equiv 0$.

\item Repeat(a), but let $y_R$ be a square wave of amplitude 0.2
m, and fundamental frequency of $0.2$ Hz.  You may take zero
initial conditions.

\item
Obtain a state feedback control law for a comfortable ride, but make
sure that the car does not bottom-out - Given the normalized state
variables, this means that $x_1(t)-x_2(t) > - 0.5m$.  You should
obtain a critically damped response for $X_1(s)/Y_R(s)$, and you
should also speed up the settling time to approximately half the open
loop response.

Show simulations under the conditions of (a) and (b) above.

\item
Full state feedback is not feasible, since it is not practical to
measure the distance from the car to the road.  Letting $y=x_1-x_2$
denote the actual measurement available to the controller, is the
resulting system observable?

\item
Repeat (c) with a full order observer, where the measurement is taken
to be $y$.
\end{list}

\item
The linearized (and normalized) magnetically suspended ball is
described by
\begin{eqnarray*}
\dot x= \left(\begin{matrix}0&1\\ 1&0\\\end{matrix}\right)x+ \left(
\begin{matrix}0\\ 1\end{matrix}\right)u.
\end{eqnarray*}

\balphlist
\item
Show that the system is unstable with $u=0$.
\item
Explain why it is possible to place the poles of the system at
arbitrary locations (with the restriction of conjugate pairs) by
linear state feedback.
\item
Find a state feedback which would place the poles of the closed-loop
system at $-1\pm j1$.
\item
Simulate your controller for the nonlinear model $\ddot y =
1-y^2/u^2$, with the nominal values $u_0=y_0=1$, and experiment with
initial conditions to find a region of asymptotic stability for the
controlled system.
\end{list}

\item\hwlabel{control-ball} In the previous problem, the ball position
$x_1$ can be measured using a photocell, but the velocity $x_2$ is
more difficult to obtain. Suppose, therefore, that the output is
$y=x_1$.  \balphlist
\item
Design a full-order observer having eigenvalues at $-5, -6$, and use
the observer feedback to produce closed-loop eigenvalues at $-1\pm j1,
-5, -6$.

\item
Simulate your controller for the nonlinear model $\ddot y =
1-y^2/u^2$, as in the previous problem, and again experiment with
initial conditions to find a region of asymptotic stability.
\end{list}

\item You will now construct a \textit{reduced-order} observer for
the previous model. \balphlist \item Repeat Exercise~10 using a
 reduced-order observer to yield closed-loop eigenvalues
at $-1 \pm j1$ and $-5$. \item Letting $z=\dot y$ in (a), compute
the transform $\hat Z(s)$, and show that it is approximately equal
to $sY(s)$ when the observer gain is large.
\end{list}

\item Below is a drawing of a cart of mass $M$ with a uniform
stick of mass $m$ pivoted on top: \ebox{.25}{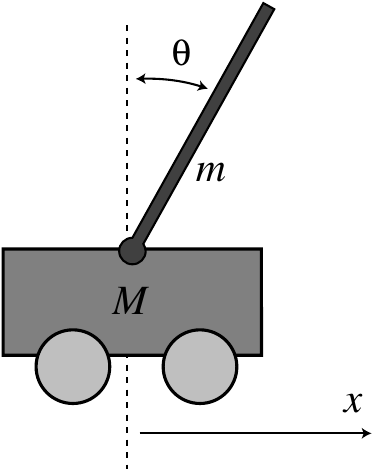} In
appropriate dimensionless units, the equations of motion may be
written as
\begin{eqnarray*}
\ddot\theta=\theta + u\;,\;\quad \ddot x= -\beta\theta - u,
\end{eqnarray*}
where $\beta \eqdef\frac{3}{4}[m/(M+m)]$ is a parameter of the
system, and $u$ is the torque applied to the wheels of the cart by
an electric motor.  We wish to find a linear feedback control that
will balance the stick (i.e., keep $\theta\approx 0$) and keep the
cart near $x=0$. To do this, find the gains $k_1, k_2, k_3$ and
$k_4$ in the state-variable feedback
\begin{eqnarray*}
u=k_1\theta +k_2\dot\theta +k_3 x+k_4\dot x
\end{eqnarray*}
such that the closed-loop system has a double pole at $s=-1$ and a
pair of complex poles at $s=-1\pm j1$

\item
On the line connecting the center of the earth to the center of the
moon, there is a so-called \textit{libration point} where the pull of
the earth on a satellite (in an orbit about the earth with the same
period as the moon's orbit) exactly equals the pull of the moon plus
the centrifugal force.  The dynamic equations for small deviations in
position away from the libration point can be shown to be:
\begin{eqnarray*}
\ddot x-2\omega\dot y-9\omega^2x&=&0 \\ \ddot y+2\omega\dot
x+4\omega^2y&=&u
\end{eqnarray*}
where $x\eqdef$ radial position perturbation, $y\eqdef$ azimuthal
position perturbation, $u=F/m\omega^2$ (control exerted by a small
reaction engine), $F\eqdef$ engine thrust in the $y$ direction,
$m\eqdef$ satellite mass, and $\omega \eqdef2\pi /29$ rad/day.
\balphlist
\item
With $u=0$, show that the equilibrium point $x\equiv y \equiv 0$ is
unstable.

\item
To stabilize the position, one can use state feedback
\begin{eqnarray*}
u=k_1x+k_2\dot x+k_3y+k_4\dot y .
\end{eqnarray*}
Show that it is possible to stabilize the system with a control law of
this form.  Determine the constants $k_1,...,k_4$ such that the
closed-loop system has poles at $s=-3\omega, s=-4\omega$, and
$s=(-3\pm j3)\omega$.
\end{list}

\item
Given the LTI system: $\dot x=Ax+Bu,\;\;y=Cx$, where
\begin{eqnarray*}
A= \left(\begin{matrix}0&0&-2\\ 1&0&1\\ 0&1&2\\\end{matrix}\right),
\quad B= \left(\begin{matrix}0&\beta\\1&0\\ 0&1\\\end{matrix}\right),
\quad C= \left(\begin{matrix}0&0&1\\\end{matrix}\right).
\end{eqnarray*}
\balphlist
\item
Determine the range of values of the scalar real parameter $\beta$ for
which we can place $2$ poles arbitrarily, using static output feedback
(as usual, we assume that if the two poles are complex then they
should be in conjugate pairs).

\item
for $\beta=1$, determine static output feedback so that two of the
eigenvalues of the closed-loop system will be $\lambda_{1,2}=-1\pm
j1$.
\end{list}

\item
Consider the linear system with transfer function
\begin{eqnarray*}
P(s) = \frac{s+28}{(s+27)(s+29)}
\end{eqnarray*}
We would like to find a control law which places the poles at $(-29,
-28)$, so that the zero is canceled.

\balphlist
\item
Place the system in modal canonical form, and apply the feedback
control law $u=-K_m x +r$. Compute $K_m$ for the desired closed-loop
pole locations.

\item
Place the system in observable canonical form, and apply the feedback
control law $u=-K_o x +r$.  Compute $K_o$ for the desired closed-loop
pole locations.

\item
Why do you think the gains are larger in the second case.
\end{list}

\item
Design a first-order controller of the form
\begin{eqnarray*}
G(s)=K\frac{(s+a)}{(s+b)}
\end{eqnarray*}
($K, a, b$ free parameters) for the plant with transfer function
\begin{eqnarray*}
P(s)=\frac{{s+28}}{{(s+27)(s+29)}}
\end{eqnarray*}
to place poles at $(-30.5, -31, -31.5)$. Repeat this for poles at
$(-31, -32, -33)$. Note that the two controllers are quite different,
even though the closed-loop poles differ only slightly.

\item
Design a first-order controller of the form
\begin{eqnarray*}
G(s)=K\frac{{(s+a)}}{{(s+b)}}
\end{eqnarray*}
($K, a, b$ free parameters) for the plant with transfer function
\begin{eqnarray*}
P(s)=\frac{{s-\alpha}}{{(s+4)}}
\end{eqnarray*}
to place poles at $-10\pm 3j$, where $\alpha>0$ is a scalar parameter.
Plot the magnitude of the sensitivity function $S = \frac{1}{P G}$ for
several values of $\alpha$ ranging from $1$ to $100$.  For what values
of $\alpha$ is the controlled system most sensitive to plant
uncertainty?

\item
Below is a diagram of a flexible structure found in the undergraduate
controls lab at the University of Illinois
\ebox{.3}{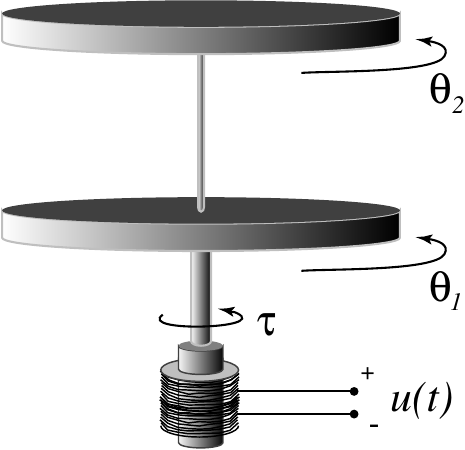} An accurate state space model can be
found, which is defined by the following matrices:
\begin{equation}
\hspace{-.4cm} A =\left[ \begin{matrix} 0 &1 &0 &0\\ -500.3513 &
-2.626 & 500.3513 &0.3972\\ 0 & 0 &0 & 1\\ 558.964 &0 & -558.964
&-2.1267
\end{matrix}\right];
\qquad B = \left[ \begin{matrix}0\\ 299.535\\ 0\\ 0\end{matrix}\right]
\elabel{torsion-ss}
\end{equation}
The input is the voltage to a DC motor, and the states of this model
are
\begin{eqnarray*}
x_1 = \theta_1, x_2 = \dot\theta_1, x_3 = \theta_2, x_4 =
\dot\theta_2.
\end{eqnarray*}
\balphlist
\item
Make a Simulink model for this system and plot the open-loop response
$\theta_2(t)$ to a non-zero initial condition with $\theta_1(0)\neq
\theta_2 (0)$.

\item
Design a full state feedback controller $u=-Kx + g r$ to place the
poles in the region $-10 \le \sigma \le -5$, and choose the gain $g$
so that the DC gain of the closed-loop system from $r$ to $\theta_2$
is one.  Plot step responses for your controlled system, with $r =
\pi/2$ radians.  In one step response, use zero initial conditions, in
the other, choose the same initial conditions as in (a).  Design $K$
so that the overshoot is less than 20\%.

\item
Make a Simulink model for an observer, assuming only $\theta_1$ is
available through direct measurements.  Combine the observer, your
full state feedback controller, and plant as on page 129 of the
lecture notes.  After designing appropriate observer gains, obtain two
step responses of the controlled system, as in (b).

\item
Repeat (c) with a reduced order observer, assuming only $\theta_1$ is
directly measured.
\end{list}
\end{exercises}

\chapter{Tracking and Disturbance Rejection}
\section{Internal model principle}\label{s:IMP}

Tracking and disturbance rejection are two of the basic goals in
control design.  This chapter addresses each of these issues by
applying the state space theory developed so far.

Assume that one wishes to make the output $y$ track a constant
reference input $r$.  One approach is to modify the state feedback
control law
\begin{eqnarray*}
u = - Kx + Nr
\end{eqnarray*}
where $N$ is chosen based upon the DC gain of the system.  Using the
final value theorem
\begin{eqnarray*}
y (\infty) = \lim_{s \to 0} s Y (s)
\end{eqnarray*}
a formula for $N$ is easily obtained.  This is a reasonable approach
if the DC gain of the plant is known.  However, the scaling matrix $N$
does not account for disturbances and parameter variations that will
hamper tracking.

To obtain a solution which is robust to DC disturbances we explicitly
model the unknown constant disturbance as follows:
\begin{eqnarray*}
\dot x & = & Ax + Bu + Ew \\ y & = & Cx + Fw,
\end{eqnarray*}
where $ w\in\Re^n $ is a constant. Our goal then is to make $e(t) = y
(t) - r \to 0$, as $t\to\infty$, regardless of the value of the
disturbance $w$, or the reference input $r$.  The approach taken
allows us to achieve these three objectives simultaneously: \balphlist
\item
Asymptotic tracking: $y(t) \to r$ as $t \to \infty$.
\item
Complete insensitivity to $w$.
\item
Tuned transient response through pole placement.
\end{list}
Conceptually, the approach taken is to note that the exogenous signal
$ \left(\begin{matrix}r\\ w\end{matrix}\right)$ is generated by a
state space model
\begin{eqnarray*}
\dot z &=& A_m z\qquad z(0) \in \Re^q; \\ \left(\begin{matrix}r \\
w\end{matrix}\right) &=& C_m z
\end{eqnarray*}
where in this special case of constant disturbances and a constant
reference input, the matrix $A_m$ is equal to zero.  This
reference/disturbance model is simply an integrator, and based upon
this we incorporate an integrator in the control law.  Thus, the
control law we adopt is of the general form
\begin{eqnarray*}
u = -K_1 x -K_2 \eta,\qquad \dot \eta = e = y-r.
\end{eqnarray*}
This procedure is known as the \defn{internal model principle}.
\notes{need a chapter! reference wonham `88 and others}

The internal model principle is most easily understood when viewed in
the frequency domain.  The controller transfer function (from $x$ to
$u$) possesses a pole at the origin, which in the frequency domain is
equivalent to demanding infinite gain at DC.  A look at the
sensitivity function then shows that sensitivity with respect to plant
uncertainty or disturbances at DC will be zero, provided the plant
itself does not possess a zero at DC.  However in this course we
remain in the time domain, and so our approach is to show that a
desired equilibrium is asymptotically stable.

First note that since the external inputs $r$ and $w$ are constant,
assuming stability one can expect that all signals will converge to
\textit{some} constant values, which will form an equilibrium for the
controlled system.  If the integrated error $\eta$ converges in this
sense, then the error itself $e$ will necessarily converge to $\zero$.
To make this precise, consider the closed-loop system equations, given
by
\begin{eqnarray*}
\left[ \begin{array}{c} \dot x \\ \doteta \end{array} \right] =
\left[\begin{array}{cc} A - BK_1 & -BK_2 \\ C & 0 \end{array} \right]
\left[\begin{array}{c} x \\ \eta \end{array} \right] +
\left[\begin{array}{c} 0 \\ -I \end{array} \right] r +
\left[\begin{array}{c} E \\ F \end{array} \right] w.
\end{eqnarray*}
At an equilibrium we have $ \left(\begin{matrix}\dotx\\
\doteta\end{matrix}\right) = \zero$, from which we conclude that
\begin{eqnarray*}
\zero = \doteta = Cx - r + Fw = y - r .
\end{eqnarray*}
Thus, we will have asymptotic disturbance rejection and asymptotic
tracking if the controlled system is stable.  To determine if this is
the case, we must consider the eigenvalues of the matrix
\begin{eqnarray*}
A_{cl} = \left[ \begin{array}{cc} A - BK_1 & -BK_2 \\ C & 0
\end{array}\right].
\end{eqnarray*}
Writing the closed-loop matrix as
\begin{eqnarray*}
\left[ \begin{array}{cc} A-BK_1 & -BK_2 \\ C & 0 \end{array} \right] =
\underbrace{\left[ \begin{array}{cc} A & 0 \\ C & 0 \end{array}
\right]}_{\tilA} - \underbrace{\left[ \begin{array}{c} B \\ 0
\end{array} \right]}_{\tilB} \underbrace{[K_1 \ K_2]}_{\tilK}
\end{eqnarray*}
we conclude that we can achieve arbitrary pole placement if and only
if the pair $(\tilA, \tilB)$ is controllable.  If this is the case
then the two gains $(K_1,K_2)= \tilK$ can be designed simultaneously
using a single \textit{Matlab} command.

We use the Hautus-Rosenbrock test to better understand the
controllablility of $(\tilA, \tilB)$.  We have for any complex $s$,
\begin{eqnarray*}
\rank (sI - \tilA \mid \tilB) = \rank \left[ \begin{array}{cc|c} sI -
A & 0 & B \\ -C & sI & 0 \end{array} \right].
\end{eqnarray*}
For controllability, we must have for all $s$,
\begin{equation}
\rank \left[ \begin{array}{cc|c} sI - A & 0 & B \\ -C & sI & 0
\end{array} \right] = n+p.
\elabel{integrate-Hautus}
\end{equation}
That is, all $(n+p)$ rows must be linearly independent.

First suppose that $s \neq 0$. Then the bottom $p$ rows (given by
$[-C\ sI \mid 0]$) are linearly independent because of the presence of
the rank $p$ matrix $sI$. The remaining upper $n$ rows will be
linearly independent if $(A, B)$ is controllable (the $0$ term does
not affect the rank of these rows). Thus if $(A,B)$ is controllable,
the rank condition is met for $s\neq 0$.  For $s = 0$, the rank
condition \eq integrate-Hautus/ becomes
\begin{eqnarray*}
\rank \left[ \begin{array}{rrr} -A & 0 & B \\ -C & 0 & 0 \end{array}
\right] = n + p.
\end{eqnarray*}
To summarize, the following two conditions are equivalent to
controllability of the pair $(\tilA,\tilB)$: \balphlist
\item
$(A,B)$ is controllable
\item
$\displaystyle \rank \left[ \begin{array}{rr} -A & B \\ -C & 0
\end{array} \right] = n + p.  $
\end{list}
For $m\ge p$, condition (b) will be true if and only if no
transmission zero exists at $s = s_0 = 0$ (see \Page{zero}).

Note that if the rank condition (b) fails then there exists a non-zero
pair $(x_0,u_0)\in \Re^{n+p}$ such that
\begin{eqnarray*}
\left[ \begin{array}{rr} A & B \\ C & 0 \end{array} \right] \left(
\begin{array}{r} x_0 \\ u_0 \end{array} \right) = \left(
\begin{array}{r} \zero \\ \zero \end{array} \right).
\end{eqnarray*}
If the constant control $u\equiv u_0$ is applied when the initial
condition is $x_0$ then we have
\begin{eqnarray*}
\dot x = A x_0 + B u_0 =\zero.
\end{eqnarray*}
That is, the state $x_0$ is an equilibrium when this constant input is
applied.  The output equation is $y=C x_0 = \zero$, so that the
initial state/input signal is invisible to the output.  It is this
complete lack of response to a DC input that must be ruled out if
integral control, or any other approach to DC tracking is to be
successful.

If the states are not available, we can still use an observer.  If
there are disturbances then the estimates may not converge, but one
can still show that the unique equilibrium will satisfy $e=0$.  In
\Exercise{observer-disturbance} you will see how the estimator behaves
in the presence of a disturbance.

\begin{ex}
Consider the first order model
\begin{eqnarray*}
\dot x & = & -2x + u + w \\ y & = & x.
\end{eqnarray*}
To regulate this system to track a constant reference $r$, the control
law will be of the form illustrated below.
\begin{center}
\ebox{1}{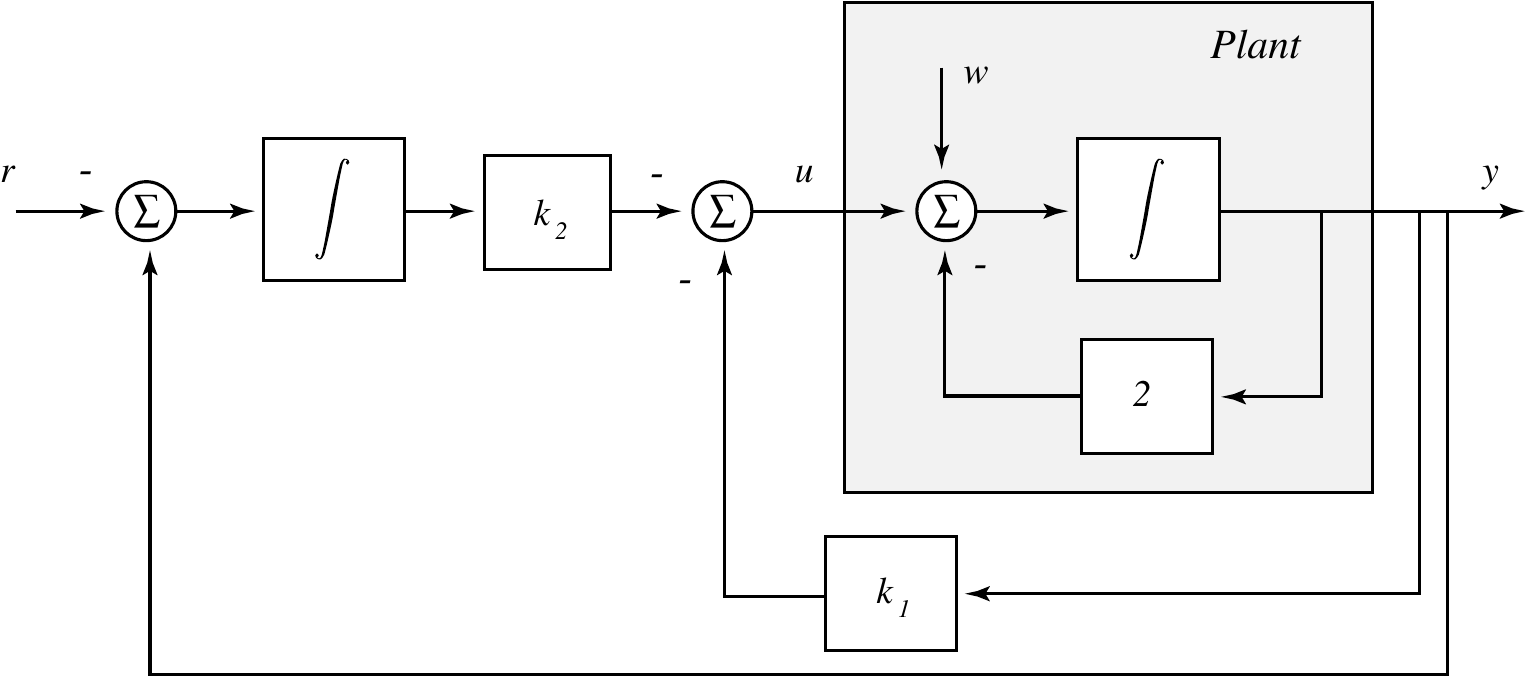}
\end{center}
To place the closed-loop poles we consider the closed-loop system
matrix
\begin{eqnarray*}
A_{cl}=\tilA - \tilB \tilK = \left[ \begin{array}{cc} -2 -k_1 & -k_2
\\ 1 & 0 \end{array} \right].
\end{eqnarray*}
For example, to place the closed-loop poles at $-2 \pm j2$ we set
\begin{eqnarray*}
\det [sI - (\tilA - \tilB \tilK)] = s^2 + (2 + k_1) s + k_2 = s^2 + 4s
+ 8.
\end{eqnarray*}
This gives $k_1 = 2, k_2 = 8$, so that we obtain the ``PI'' controller
\begin{eqnarray*}
u = -2y - 8 \int_0^t e \, d t.
\end{eqnarray*}
\end{ex}

The internal model principle can be extended to the more general
situation where the disturbance and reference signal are the sum of
periodic signals of known frequency.  If this is the case, then we may
again find a matrix $A_m$ such that
\begin{eqnarray*}
\dot z &=& A_m z\qquad z(0) \in \Re^q; \\ \left(\begin{matrix}r\\
w\end{matrix}\right) &=& C_m z,
\end{eqnarray*}
where the eigenvalues of the matrix $A_m$ lie on the $j\omega$-axis.
To design a controller one replicates this model to define the signal
$\eta$:
\begin{eqnarray*}
\dot w &=& A_m w + B_1 e \qquad w(0) \in \Re^q; \\ \eta &=& C_1 w.
\end{eqnarray*}
As before, the control can then be defined as $u= -K_1 x  -K_2\eta$.

\section{Transfer function approach}
When we use an observer together with state feedback to control a
plant, we are actually designing a form of dynamic feedback
compensation, similar to what is obtained through loop shaping in a
frequency domain design.  We show here that equivalent dynamic
compensation may be conducted entirely in the Laplace transform domain
by manipulating transfer functions.  We concern ourselves only with
the SISO case, and for simplicity we treat only the regulator problem.
Integral control or reduced order observer based feedback can be
analyzed in the transform domain in the same manner.

\begin{figure}[ht]
\ebox{.75}{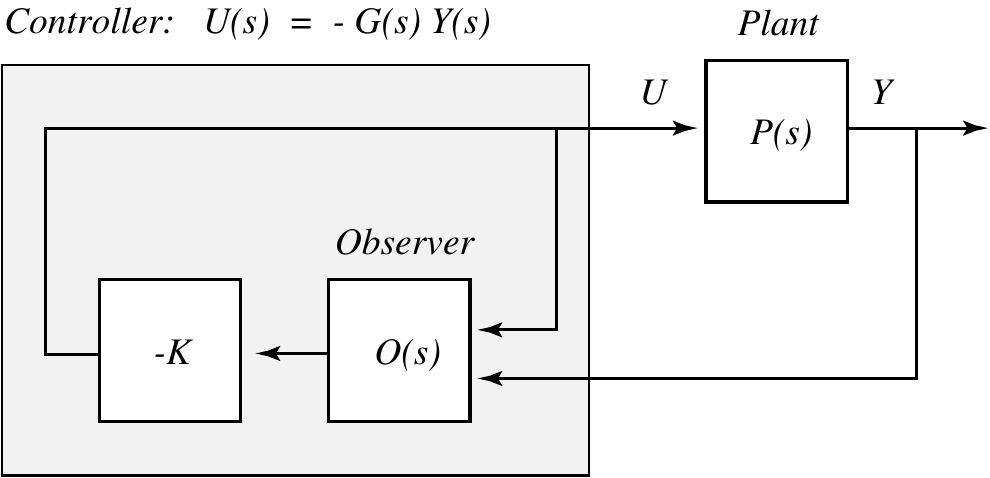}
\caption[Observer/state feedback]{The observer/state feedback combination can be regarded as
dynamic compensation} \flabel{feedback-observerTF}
\end{figure}

The regulator problem using a full order observer is illustrated in
\Figure{feedback-observerTF}.  To compute the controller transfer
function $G(s)$, we first construct the transfer function from $Y$ to
$\hat X$.  The state estimates are described by the equation
\begin{eqnarray*}
\dot{\hatx}= A \hatx + Bu + L(y - C \hatx) = (A - BK - LC) \hatx + Ly
\end{eqnarray*}
where the second equation follows from the control law
\begin{eqnarray*}
u = - K \hatx.
\end{eqnarray*}
Taking transforms, we see that
\begin{eqnarray*}
U(s) = - K \hat X(s) = - K (sI - [A - BK - LC])^{-1} L Y(s).
\end{eqnarray*}
This results in the formula
\begin{eqnarray}
G (s) = K (sI - [A - BK - LC])^{-1} L. \elabel{eq3.5.1}
\end{eqnarray}
Since there is no direct feedthrough term (no ``$D$'' matrix) in the
formula \eq eq3.5.1/, the transfer function $G$ is strictly proper.

This transfer function can be directly designed, without consideration
of a state space model.  To see this, first write the plant and
controller transfer functions as a ratio of polynomials:
\begin{eqnarray*}
\hbox{\it Plant transfer function} = P(s) =\frac{b (s)}{a (s)} \\
\hbox{\it Controller transfer function} = G(s) = \frac{n (s)}{d (s)}.
\end{eqnarray*}
If the control $U(s) = - G(s) Y(s) + R(s)$ is applied, then the closed
loop transfer function becomes
\begin{eqnarray*}
P_{cl}(s) =\frac{Y(s)}{R(s)} = \frac{\frac{b (s)}{a(s)}}{1+\frac{n
(s)}{d(s)} \frac{b(s)}{a(s)}} = \frac{b(s) d(s)}{a(s)d(s)+b(s)n(s)}.
\end{eqnarray*}
The poles of the closed-loop system are determined by the
corresponding characteristic equation
\begin{eqnarray*}
\Delta_{cl}(s) = a (s) d (s) + b (s) n (s) = 0.
\end{eqnarray*}
Thus, if we let $\alpha_c$, $\alpha_o$ denote the polynomials which
define the state feedback poles, and observer poles, respectively:
\begin{eqnarray*}
\left. \begin{array}{rcl} \alpha_c (s) & = & \det (sI - (A-BK)) \\
\alpha_o (s) & = & \det (sI - (A-LC))
\end{array} \right\} = \hbox{known,}
\end{eqnarray*}
then the polynomials $n(s),d(s)$ must satisfy the following
\defn{Diophantine equation}:
\begin{equation}
a (s) d (s) + b (s) n (s) = \alpha_c(s) \alpha_o (s). \elabel{eq3.5.2}
\end{equation}

We assume that both $a (s)$ and $d (s)$ are monic, so that each of
these polynomials has degree $n$.  We have already noted that $G$ is
strictly proper, which means that the degree of $n(s)$ is strictly
less than $n$.  The product $\alpha_c (s) \alpha_o (s)$ should also be
monic.  With this information, the unknown polynomials $d (s)$ and $n
(s)$ can be computed.  However, some conditions must be imposed.  When
we treated this problem using state space arguments we found that the
system must be controllable and observable to compute the gains $K$
and $L$.  This is equivalent to minimality, which means that the
transfer function $P(s) = C(Is - A)^{-1} B$ cannot have any common
pole-zero pairs, or that the polynomials $a$ and $b$ do not have
common roots.

\notes{should add a theorem}
\begin{ex}
Consider the transfer function description of the magnetically
suspended ball
\begin{eqnarray*}
P(s) = \frac{1}{s^2 -1}.
\end{eqnarray*}
There can be no pole zero cancellations since $b(s) = 1$, and hence
any two dimensional state space model is minimal.  We conclude that
the degree of $n(s)$ is less than or equal to $1$, and the degree of
$d(s)$ is equal to $2$.  So, the controller can be written
\begin{eqnarray*}
U(s) = \frac{n_1 s + n_0}{s^2 + d_1 s + d_0} Y(s)
\end{eqnarray*}
From equation \eq eq3.5.2/ we must solve the following equation to
place the two state feedback poles, and the two observer poles:
\begin{eqnarray*}
(s^2 -1)(s^2 + d_1 s + d_0) + (1) (n_1 s + n_0) = (s - p_1)(s - p_2)(s
- \ell_1)(s - \ell_2),
\end{eqnarray*}
which gives four equations and four unknowns through equating of the
coefficients.  These equations are called the \defn{Sylvester
equations}.
\end{ex}

Based on our criteria for solving the Diophantine equation, it follows
that for the Sylvester equations to be independent and solvable for
any set of observer poles and state feedback poles, the plant transfer
function cannot have common pole-zero pairs.  For example, consider
the plant transfer function
\begin{eqnarray*}
P(s) = \frac{s+1}{(s+1)(s+2)}.
\end{eqnarray*}
The corresponding Diophantine equation is
\begin{eqnarray*}
(s+1)(s+2) d (s) +(s+1) n (s) =\alpha_c(s) \alpha_o (s)
\end{eqnarray*}
Since the left hand side of this equation vanishes when $s=-1$,
obviously this equation cannot be solved if neither $\alpha_c$ nor
$\alpha_o$ have a root at $-1$. Pole placement is not possible in this
example because the transfer function $P$ does contain a common
pole-zero pair.

\notes{given title of section, should add more discussion on integral
control and sensitivity}

\begin{exercises}
\item
Consider the servo motor with transfer function
\begin{eqnarray*}
Y(s) = \frac{1}{s(1+s)} U(s),
\end{eqnarray*}
where $u$ is the voltage input to the motor, and $y$ is its position.
In this exercise you will compare closed-loop sensitivity with and
without integral control.  \balphlist
\item
Design a first order dynamic compensator
\begin{eqnarray*}
U(s)= k \frac{s-z}{s-p} (g R(s)-Y(s))
\end{eqnarray*}
so that the closed-loop poles are at $-3$, $-2\pm j$.  Choose $g$ so
that the closed-loop transfer function $Y/R$ has unity gain at DC.

\item
Plot the step response, and a Bode plot of the sensitivity transfer
function $S_1$.
\end{list}

\item
Consider the LTI system
\begin{equation} \elabel{DistEx}
\begin{array}{rcl} \dot x_1 &=&x_2 \\ \dot x_2&=& u + w \end{array}
\qquad y = x_1,
\end{equation}
where $w$ is a constant but unknown disturbance.  Rather than
introduce integral control, if we could obtain estimates of $w$, then
we could design a control law to cancel the effect of this
disturbance.

One way to attempt to do just this is to consider $x_3 = w$ as an
additional state variable:
\begin{eqnarray*}
\begin{array}{rcl} \dot x_1 &= &x_2\\  \dot x_2&=& u + x_3 \\  \dot
x_3&=& 0\end{array} \qquad y = x_1 .
\end{eqnarray*}
\balphlist
\item
Carefully show that this third-order ``augmented'' system is
observable.

\item
Design a full-order observer for this system.  Draw the corresponding
all-integrator diagram of the observer in full detail, including all
integrators.

\item
Based upon this observer, design a controller for this system to place
the closed-loop poles at $-1,-2,-3, -4, -5$.  Explain carefully why
the poles will be placed, and the state estimates will converge.

\item
Is it possible to design the observer-based feedback control law in
(c) so that the closed-loop system matrix is Hurwitz?
\end{list}

\item \hwlabel{observer-disturbance} Consider again the LTI system
given in \eq DistEx/.  \balphlist
\item
Design a feedback compensator of the form $u = -K_1 \hatx - K_2 \eta$
to regulate the output to zero.  Take $\dot\eta = y$, and $\hatx_2$ as
an estimate of $x_2$ based on a reduced order observer,
\textit{assuming $w\equiv 0$}.

\item
With this controller, compute $\hatx_2^\infty =
\lim_{t\to\infty}\hatx_2(t)$, as a function of $w$.  \textit{Note}:
you do not need to compute the entire transfer function.

\item
Provide simulations of your controller using \textit{Simulink} for
various initial conditions, and values of $w$.
\end{list}

\item
An electromagnetic suspension system for a magnetic ball is modeled by
\begin{eqnarray*}
m\ddot y=mg-ci^2/y^2
\end{eqnarray*}
where $m$ is the mass of the ball ($=10^{-3}$ kg.), $g=9.8 {\rm m.}/
{\rm sec}^2$, $c=9.8\times 10^{-9}$ newton-meter$^2$/amp$^2$, $y$ is
the distance between the magnet and the ball, and $i$ is the current
through the electromagnet. Suppose that transducers are available for
measuring $x_1=y, x_2=\dot y$. It is desired that in the steady state,
$x_1\rightarrow v$, where $v$ is a constant reference input.

To accomplish this, it is proposed to feed back the integral of the
error $x_1-v$ in addition to the states $x_1$ and $x_2$. That is, the
control $i$ is to be implemented as
\begin{eqnarray*}
i=k_0x_0+k_1x_1+k_2x_2
\end{eqnarray*}
where $x_0$ is the integral of $x_1-v$, and $k_0, k_1, k_2$ are
constant gains, chosen by the designer as $k_0=100$, $k_1=1200$,
$k_2=1$. Find the equilibria of the closed-loop system when $v=0.025$
m, and determine whether or not they are asymptotically stable.
\end{exercises}

\chapter{Control Design Goals}\clabel{goals}
\section{Performance}
Control can be broadly defined as the use of feedback to obtain
satisfactory performance in the presence of uncertainty.
\defn{Feedback} is simply the act of incorporating measurements as
they become available to define the inputs to the plant.  The
definition of \defn{performance} depends of course on the nature of
the particular control application, but there are some universal
performance issues.  Among these are:
\begin{description}
\item[(i)] \textit{Stability} of the controlled system is of course
the most basic goal. Anyone who has held a microphone too close to its
speaker system knows that feedback can result in instability, but when
properly designed, feedback can prevent unstable behavior.

\item[(ii)] Frequently one desires reasonably accurate \defn{tracking}
of some desired trajectory.  If the trajectory is a constant reference
value, this is known as the \defn{regulation} problem. In the
production of chemicals, particularly in the pharmacuetical industry,
yields must be controlled within stringent specifications.  In the
automatic control of an aircraft in flight, again one objective is to
keep the flight at a level altitude and a constant speed.

\item[(iii)] Another time domain specification concerns \textit{speed
of response}.  For example, if the pilot of an aircraft desires a
change in altitude, this should occur quickly, but not at the expense
of the comfort of the passengers.

\item[(iv)] Rarely do we have available a perfect model of the plant
to be controlled.  Even if an accurate model is available, it may be
convenient to linearize the nonlinear model, which can introduce
significant inaccuracies.  Such inaccuracies combined with unforeseen
influences on the plant are collectively known as disturbances, and we
must always consider \textit{robustness} to such disturbances when
designing a control system.

\end{description}
All of these issues are typically addressed in the design of a control
system in the face of physical and economic constraints, and in most
current industries, time constraints as well!

One cannot expect good performance if the physical system is poorly
designed.  In particular, it is important from the start that the
measurement devices be reasonably precise, and that actuation - the
final synthesis of the control action - be performed accurately.
However, in this course we assume that these issues were considered in
the construction of the plant, so that what is left for the control
engineer can be largely completed with a pen, paper, and a computer.

\section{Measurements}
This course concerns state space models, which in many of our examples
are naturally constructed based upon signals which have physical
significance, and which can be measured if sensors are adjoined to the
plant.  For instance, in the pendubot illustrated in
\Figure{mod-pendubot}, the states are the link angles, and their
velocities.  The link angles can be measured with high precision using
an optical encoder.  The velocities are more difficult to measure, but
they can be estimated using finite differences.  If, as in this
example, an accurate state space model is available together with
accurate state measurements, then control based upon pole placement or
optimal control is frequently simple and very effective.

However, there are times when full state information is not feasible.
For instance, in many applications the state space model is
constructed by first finding a frequency domain model (i.e.\ a Bode
plot) based upon input-output measurements.  A state space model can
then be constructed so that its frequency response accurately matches
this empirical Bode plot.  In this case, the state has no physical
significance, and the only information available for controlling the
system may be the input and output.  In \Section{Observers} we
circumvent this difficulty by constructing \textit{observers} to
estimate the state as input-output measurements are collected.  In
principle at least, one can then treat the estimates as if they are
perfect measurements of the state, and then apply one of the state
feedback control techniques.  This gives a mathematically elegant
approach to control design which is also easily analyzed.  In
particular, one can easily verify that the goals (i)--(iii) outlined
above will be met if the observer and state feedback controller are
appropriately designed.

The fourth goal, robustness, can be easily forgotten in a state space
design, and this can lead to absurd control designs!  To see this,
consider once again the pendubot shown in \Figure{mod-pendubot}.  With
the input $u$ taken as the torque applied to the lower link, and with
the output $y$ taken to be the lower link angle, the linearized plant
in the vertical position is both controllable and observable.
Theoretically then, it is easy to design a controller which will keep
the pendubot balanced in the vertical position based upon measurements
of the lower angle alone.  Physically, this seems very unlikely, and
in fact such a controller has never been successfully implemented on
this plant.  The reason for this divergence between theory and
practice is that when evaluating a controller one must also consider
robustness to plant uncertainty.  Using a frequency domain analysis,
one can more easily see that any controller for this plant based upon
the measurements $u,y$ alone will be highly sensitive to model
inaccuracy.  An insensitive design is possible if both link angles are
used in the feedback compensator.

\section{Robustness and sensitivity}
The Bode sensitivity function is a measure of the sensitivity of the
closed-loop transfer function to a deviation in the plant model.  For
the control system illustrated in \Figure{feedback-feedback}, the
closed-loop transfer function may be expressed
\begin{eqnarray*}
T = \frac{Y}{R} = \frac{GP}{1+GP}
\end{eqnarray*}
If the transfer function $P$ changes to form a new transfer function
$P'$, the controller will not function precisely as it was intended,
since the transfer function $T$ will undergo a corresponding change.
Rather than consider absolute changes, it is more appropriate to
analyze the affect of the relative change in the plant
\begin{eqnarray*}
\frac{\Delta P}{P},
\end{eqnarray*}
where $\Delta P = P'-P$ denotes the difference between the true
transfer function and the nominal one.  For small changes, this can be
interpreted as a differential.  The \defn{sensitivity} is then the
ratio of this differential, and the corresponding differential in $T$:
\begin{equation} \elabel{sensitivity-def} S \eqdef \frac{d T/T}{dP/P}
= \frac{P}{T} \frac{dT}{dP}.
\end{equation}
For the error feedback configuration given in
\Figure{feedback-feedback}, the derivative can be computed to give
\begin{eqnarray*}
S = \frac{1}{1+PG} = \frac{1}{1+L}
\end{eqnarray*}
where $L (s)$ is the \defn{loop gain} $ P (s) G (s)$.  The quantity
$1+L(s)$ is known as the \defn{return difference}.

In the MIMO case we can again write for the error feedback
configuration
\begin{eqnarray*}
Y = \left[ (I + PG)^{-1} PG\right] R,
\end{eqnarray*}
and the sensitivity function is then defined to be the matrix-valued
function of $\omega$,
\begin{eqnarray*}
S(j\omega) = (I + P(j\omega)G(j\omega))^{-1}.
\end{eqnarray*}

\begin{figure}[ht]
\ebox{.65}{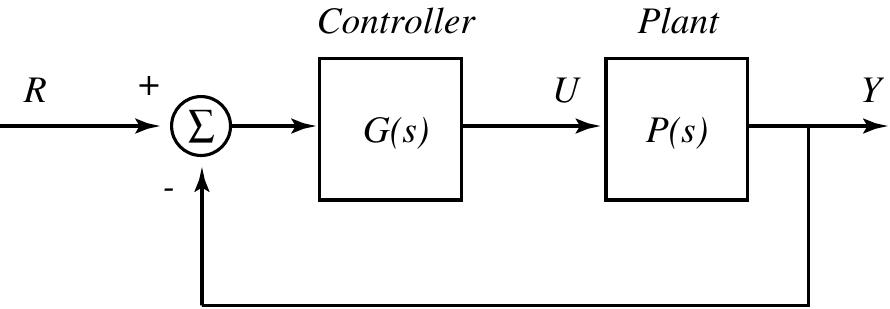}
\caption{An error feedback configuration} \flabel{feedback-feedback}
\end{figure}

We can conclude from this formula that if we want a system which is
insensitive to changes in the plant, then we should make the loop gain
$L (s)$ as large as possible. For SISO systems the size of the
sensitivity function is measured by means of the magnitude of its
frequency response.  We might require for example that
\begin{eqnarray*}
| S (j \omega) | < 1, \qquad \omega \ge 0,
\end{eqnarray*}
which does correspond to a ``large'' loop gain.  In making the
sensitivity small however, we must keep all of the previous design
goals in mind, and there are also other constraints which are
absolute.  For instance, for the error feedback configuration we will
always have
\begin{eqnarray*}
S + T = 1.
\end{eqnarray*}
Moreover, if the transfer function $L$ is relative degree at least
two, and all poles lie in the left half plane in $\Co$, we have the
constraint
\begin{eqnarray*}
\int_0^\infty \log | S(j\omega) | d \omega = 0 .
\end{eqnarray*}
In this case it is impossible to make $| S (j \omega) | < 1$ for all
$\omega$, since then $\log | S(j\omega) |$ is always negative.

\section{Zeros and sensitivity}
It is clear that the structure of the open-loop poles will play an
important role in the design of an effective controller.  What is less
obvious at first glance is the importance of the open-loop zeros.

To illustrate the importance of zeros, we again take the pendubot
introduced in \Chapter{chap1.1}.  The pendubot can be linearized in
the upright position illustrated in \Figure{mod-pendubot}. With the
input $u$ equal to the applied torque, and the output $y$ equal to the
lower link angle, the resulting state space model is defined by the
following set of matrices:
\begin{eqnarray}
A &=& \left(\begin{matrix} 0 & 1.0000 & 0 & 0\\ 51.9243 & 0 & -13.9700
& 0\\ 0 & 0 & 0 & 1.0000\\ -52.8376 & 0 68.4187 & 0 & 0
\end{matrix} \right) \; B = \left(\begin{matrix}
0\\ 15.9549\\ 0\\ -29.3596\end{matrix}\right) \nonumber \\ C &=&
\left(\begin{matrix} 1 & 0 & 0 & 0\end{matrix} \right) \qquad D = 0.
\elabel{CDpend}
\end{eqnarray}
Using \textit{Matlab}, the corresponding transfer function model is
found to be
\begin{eqnarray*}
P(s) = \frac{Y(s)}{U(s)} = 15.9549 \frac{(s-6.5354)(s+6.5354)}{
(s-9.4109)(s+9.4109)(s-5.6372)(s+5.6372)}
\end{eqnarray*}
This is of the general form
\begin{eqnarray*}
P(s) = k \frac{(s-\gamma)(s+\gamma)}{
(s-\alpha)(s+\alpha)(s-\beta)(s+\beta)},
\end{eqnarray*}
where $0<\alpha <\gamma<\beta$.

Because $P$ contains no common pole-zero pairs, the state space model
must be both controllable and observable, i.e., minimal.  In this
chapter we will show that it is therefore possible to design a fourth
order dynamic compensator $U(s) = - G(s) Y(s)$ which stabilizes this
model at $x=\zero$.  We will find in \Chapter{optimal} that a set of
closed-loop poles which will perform well using full state feedback
are
\begin{equation}
p_1 = -36 \quad p_2 = -2.9 \quad p_3 = -5.3 + 0.84j \quad p_4 = -5.3 -
0.84j. \elabel{openPendPoles}
\end{equation}

The compensation will introduce four poles, which we will place at
\begin{eqnarray*}
p_5 = -3.9932 \quad p_6 = -7.4167 \quad p_7 = -8.4862 \quad p_8 =
-11.3012.
\end{eqnarray*}
All of these poles are based upon an ``optimal'' design.  By examining
a root locus plot it may be shown that any compensator which
stabilizes this system will have poles in the right half plane in
$\Co$.  In this example, we find that the compensator which places
these poles is given by
\begin{eqnarray*}
G(s) = 24,424 \frac{(s+ 9.4)(s+5.64)(s- 7.41)}{(s+60.24 +63.5j) (s +
60.24 -63.5j) (s- 46.37) (s+ 6.54)}.
\end{eqnarray*}
As expected, this compensator has an unstable pole at $s=46.37$.

\begin{figure}[th]
\ebox{.95}{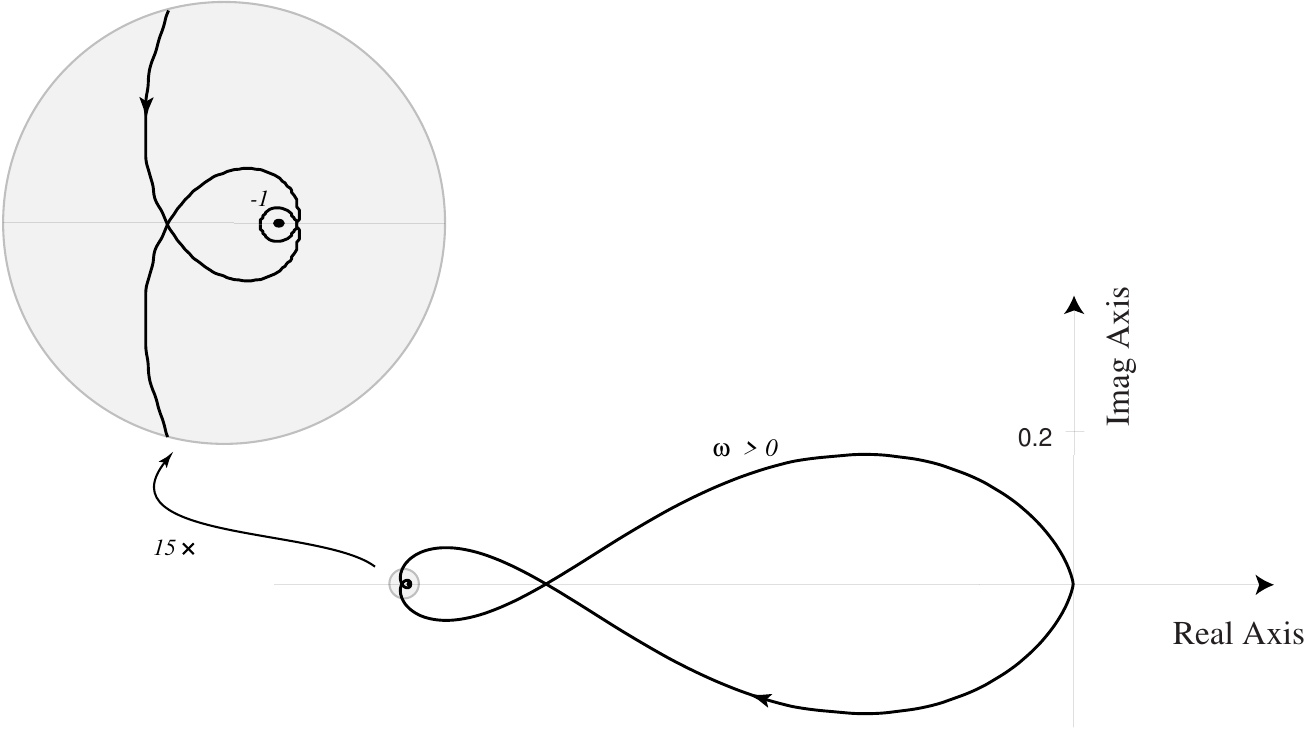}
\caption[A Nyquist plot for the Pendubot]{The Nyquist plot of $G(j\omega)P(j\omega)$.  The compensated
open-loop system with transfer function $GP$ possesses two poles in
the right half plane.  The three encirclements of $-1$ ensure that the
controlled system is asymptotically stable, from the Nyquist Stability
Criterion, but this stability is very fragile.}
\flabel{feedback-nyquist}
\end{figure}

\begin{figure}[ht]
\ebox{.7}{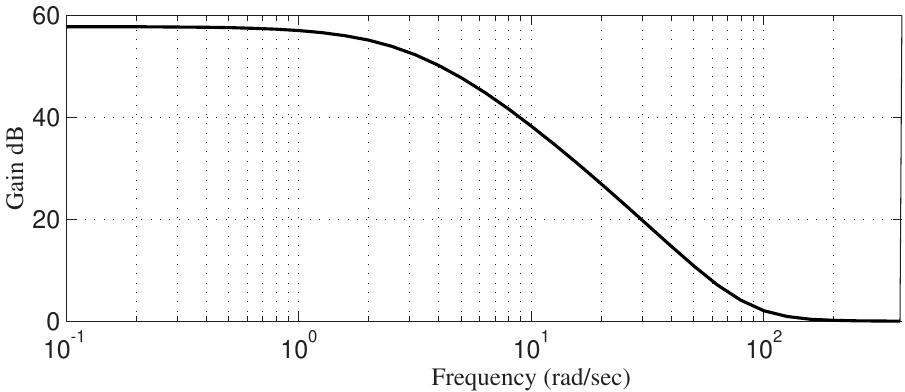}
\caption[Sensitivity function for the Pendubot]{The magnitude plot of the frequency response for the
sensitivity function $S=1/(1+GP)$.  The controlled system is highly
sensitive to plant uncertainty at low frequencies. }
\flabel{feedback-bode}
\end{figure}

A Nyquist plot of the loop transfer function $L(s)=G(s) P(s)$ is given
in \Figure{feedback-nyquist}. If the plant is changed slightly, then
the number of encirclements of $-1$ will change, resulting in an
unstable closed-loop system. A Bode plot of $S=1/(1+GP)$ is shown in
\Figure{feedback-bode}.  Here it is seen that the controlled system is
highly sensitive to plant uncertainty in the low frequency range.  By
applying higher feedback gain, it is possible to reduce sensitivity in
this range, but then one will find high sensitivity at higher
frequencies.

One way to explain these difficulties is the presence of zeros in the
right half plane, which impose strong limits on achievable
sensitivity. For the zero free model with transfer function
\begin{eqnarray*} P(s) = \frac{1}{
(s-\alpha)(s+\alpha)(s-\beta)(s+\beta)},
\end{eqnarray*}
it is is easy to obtain a low sensitivity design.  Another explanation
is again common sense - it is difficult to imagine that this plant can
be controlled using only information at the lower link.  By adding
another sensor, an insensitive design may be obtained using state
space methods.

For a SISO transfer function $P(s) = b(s)/a(s)$, the poles are the
roots of $a$, and the zeros are the roots of $b$, so that we have
$P(s_0)=0$ for a zero $s_0$.  In the MIMO case, we define the
characteristic polynomial as $\Delta(s)=\det(Is-A)$, and the roots of
$\Delta$ are then the poles of $P$.  We define zeros similarly, but it
would be far too strong to suppose that $P(s_0)=\zero$ at a zero
$s_0$.  For the state space model
\begin{eqnarray*}
\dot x = Ax + Bu, \quad y = Cx, \quad x_0 = \zero
\end{eqnarray*}
we can take Laplace transforms to obtain
\begin{eqnarray*}
(sI - A) X (s) = BU (s), \qquad Y(s) = Cx (s).
\end{eqnarray*}
Combining these equations in one matrix expression gives
\begin{eqnarray*}
\underbrace{\left[ \begin{array}{cc} sI-A & -B \\ C & \zero
\end{array} \right]}_{(n+p) \times (n+m)} \left[ \begin{array}{c} X
(s) \\ U (s) \end{array} \right] = \left[ \begin{array}{c} \zero \\ Y
(s) \end{array} \right].
\end{eqnarray*}
Assume that there are as many plant inputs as outputs $(p=m)$, so that
the plant is \defn{square}.  A natural definition of a zero can be
formulated as follows: At $s_0$, can I find a fixed $U(s_0)$ and
$X(s_0)$ so that the resulting $Y(s_0)$ is zero?  That is,
\begin{eqnarray*}
\left[ \begin{array}{cc} s_0I - A & -B \\ C & 0 \end{array} \right]
\left[ \begin{array}{c} X (s_0) \\U (s_0) \end{array} \right] = \left[
\begin{array}{c} \zero \\ \zero \end{array} \right].
\end{eqnarray*}
This is possible if and only if the matrix above is rank deficient:
\begin{eqnarray*}
\rho\left(\left[ \begin{array}{cc} s_0I - A & -B \\ C & 0 \end{array}
\right] \right) < n+m.
\end{eqnarray*}
Provided that $s_0$ is not also a pole of $P$, this is equivalent to
rank deficiency of the transfer function $P$ at $s_0$:
\begin{eqnarray*}
\rho\left(P(s_0) \right)=\rho\left(C(Is_0-A)^{-1} B \right) < m,
\end{eqnarray*}
which does generalize the definition given in the SISO case. The
complex number $s_0$ is called a \label{zero} \defn{transmission
zero}.  If all of the transmission zeros lie within the open left
half plane in $\Co$ then the model is called \defn{minimum phase}.
In the multivariate case, the location of the transmission zeros
can play an important role in the achievable sensitivity of the
controlled plant, just as in the SISO case.

\begin{exercises}
\item
For the Pendubot described by the state space model \eq CDpend/, you
will design two controllers.  The design can be done using the {\tt
PLACE} commands within \textit{Matlab}, and the formula
\eq eq3.5.1/.
\balphlist
\item
Design a fifth order dynamic compensator $U=-GY$ of the form
\begin{eqnarray*}
G(s) = \frac{n(s)}{sd(s)},
\end{eqnarray*}
where $d$ is a fourth order, monic polynomial.  The dominant closed
loop poles should lie to the left of the line $\hbox{Re}(s) = -10$ in
the complex plane.  Obtain a Bode plot of the sensitivity function for
your design.

\item
You will now design a compensator based on the two angle measurements.
We then have
\begin{eqnarray*}
U(s) =\frac{n_1(s)}{sd_1(s)} Y_1(s) +\frac{n_2(s)}{sd_2(s)} Y_2(s)
\end{eqnarray*}
where both $d_1$ and $d_2$ are fourth order polynomials.  See if you
can find an observer based dynamic feedback controller of this form
which has significantly improved sensitivity.  Obtain a Bode plot of
the sensitivity function for your design.
\end{list}
\end{exercises}


\part{Optimal Control}
\chapter{Dynamic Programming and the HJB Equation} \clabel{optimal}
\section{Problem formulation}
This chapter concerns optimal control of dynamical systems.  Most of
this development concerns linear models with a particularly simple
notion of optimality.  However, to understand the most basic concepts
in optimal control, and not become lost in complex notation, it is
most convenient to consider first the general model given in nonlinear
state space form
\begin{equation}
\dot x = f (x, u, t), \qquad x (t_0) = x_0. \elabel{optFormulation}
\end{equation}
The cost $V$ of a particular control input $u$ is defined by the
following integral
\begin{eqnarray}
V(u) = \int_{t_0}^{t_1} \ell (x, u, \tau) \, d\tau + m (x({t_1}))
\elabel{cost}
\end{eqnarray}
where \balphlist
\item
${t_1}$ is the \textit{final time} of the control problem.

\item $\ell$ is a scalar-valued function of $x$, $u$, and $t$

\item $m$ is a function of $x$. It is called the \textit{terminal
penalty function}.
\end{list}
We assume that $x_0$, $t_0$, and ${t_1}$ are known, fixed values, and
$x ({t_1})$ is free.  Our goal is to choose the control $u_{[t_0,
{t_1}]}$ to minimize $V$.

\notes{how do you determine matrices? what is the benefit of quadratic
cost?}

A case which is typical in applications is where $m$ and $\ell$ are
quadratic functions of their arguments.  For an LTV model, the system
description and cost are then given by
\begin{equation}
\begin{aligned}
\dot x &= A (t) x + B (t) u, \qquad x (t_0) = x_0 \\ V (u) &=
\int_{t_0}^{t_1} (x^T Q (t) x + u^T R (t) u)\, dt + x^T (t_1) M x
(t_1)
\end{aligned}
\elabel{LQR}
\end{equation}
where $M$, $Q$ and $R$ are positive semidefinite matrix-valued
functions of time. These matrices can be chosen by the designer to
obtain desirable closed-loop response.  The minimization of the
quadratic cost $V$ for a linear system is known as the \defn{linear
quadratic regulator} (LQR) problem.

We will study the LQR problem in detail, but first we develop some
general results for the nonlinear state space model \eq
optFormulation/.

\section{Hamilton-Jacobi-Bellman equations}
The \defn{value function} $V^\optchar=V^\optchar (x_0, t_0)$ is
defined to be the minimum value of $V$ over all controls.  This is a
function of the two variables $x$ and $t$ which can be written
explicitly as
\begin{eqnarray*}
V^\optchar (x,t) = \min_{u_{[t,{t_1}]}} \Bigl[ \int_t^{t_1} \ell
(x(\tau), u(\tau), \tau) \, d\tau + m(x({t_1})) \Bigr].
\end{eqnarray*}
Under very general conditions, the value function satisfies a partial
differential equation (PDE) known as the \defn{Hamilton-Jacobi-Bellman (HJB)
equation}.  To derive this result, let $x$ and $t$ be an arbitrary
initial time and initial state, and let $t_m$ be an intermediate time,
$t< t_m < t_1$.  Assuming that $x(\tau)$, $t\le \tau\le t_1$, is a
solution to the state equations with $x(t)=x$, we must have
\begin{eqnarray*}
V^\optchar (x,t) & = & \min_{u_{[t,{t_1}]}} \Bigl[ \int_t^{t_m}
\ell (x(\tau), u(\tau), \tau) \, d\tau + \int_{t_m}^{t_1} \ell
(x(\tau), u(\tau), \tau) \, d\tau + m (x({t_1})) \Bigr] \\ & = &
\min_{u_{[t, t_m]}} \Bigl[ \int_t^{t_m} \ell (x(\tau), u(\tau),
\tau) \, d\tau + \underbrace{\min_{u_{[t,{t_1}]}}
\Bigl(\int_{t_m}^{t_1} \ell (x(\tau), u(\tau), \tau) \, d\tau + m
(x({t_1}))}_{V^\optchar (x(t_m),t_m)} \Bigr) \Bigr].
\end{eqnarray*}
This gives the functional equation
\begin{eqnarray}
V^\optchar (x,t) = \min_{u_{[t,{t_m}]}} \Bigl[ \int_{t}^{t_m} \ell
(x(\tau), u(\tau), \tau) \, d\tau + V^\optchar (x(t_m), t_m) \Bigr].
\elabel{OP}
\end{eqnarray}
As a consequence, the optimal control over the whole interval has the
property illustrated in \Figure{optimal-traj}: If the optimal
trajectory passes through the state $x_m$ at time $x(t_m)$ using the
control $u^\optchar$, then the control $u^\optchar_{[t_m,{t_1}]}$ must
be optimal for the system starting at $x_m$ at time $t_m$.  If a
better $u^{*}$ existed on $[t_m,{t_1}]$, we would have chosen it.
This concept is called the \defn{principle of optimality}.

\begin{figure}[ht]
\ebox{.75}{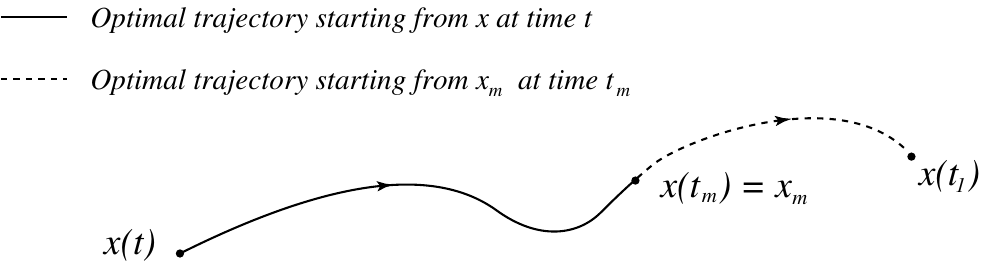}
\caption[If a better control existed...]{If a better control existed on $[t_m,{t_1}]$, we would have
chosen it. } \flabel{optimal-traj}
\end{figure}

By letting $t_m$ approach $t$, we can derive a partial differential
equation for the value function $V^\optchar$.  Let $\Delta t$ denote a
small positive number, and define
\begin{eqnarray*}
t_m & = & t + \Delta t \\ x_m & = & x (t_m) = x (t + \Delta t)= x (t)
+ \Delta x.
\end{eqnarray*}
Assuming that the value function is sufficiently smooth, we may
perform a Taylor series expansion on $V^\optchar$ using the optimality
equation \eq OP/ to obtain
\begin{eqnarray*}
V^\optchar (x,t) = \min_{u_{[t,{t_m}]}} \left\{\ell (x(t), u (t), t)
\Delta t + V^\optchar (x,t) + \frac{\partial V^\optchar }{\partial x}
(x(t),t)\Delta x + \frac{\partial V^\optchar }{\partial t}(x(t),t)
\Delta t \right\} \\
+ p (\Delta t).
\end{eqnarray*}
Dividing through by $\Delta t$ and recalling that $x(t) = x$ then
gives
\begin{eqnarray*}
0 = \min_{u_{[t,{t_m}]}} \left\{\ell (x, u (t), t) \frac{\Delta t
}{\Delta t} + \frac{\partial V^\optchar }{\partial x}
(x,t)\frac{\Delta x }{\Delta t} + \frac{\partial V^\optchar }{
\partial t} (x,t)\frac{\Delta t }{\Delta t} \right\}.
\end{eqnarray*}
Letting $\Delta t \to 0$, the ratio $\Delta x /\Delta t$ can be
replaced by a derivative to give
\begin{eqnarray*}
0 = \min_{u} \left[ \ell (x,u,t) + \frac{\partial V^\optchar}{\partial
x} (x,t)\dot x (t) \right] + \frac{\partial V^\optchar }{\partial t}
\end{eqnarray*}
where
\begin{eqnarray*}
\frac{\partial V^\optchar }{\partial x} = \left[ \frac{\partial
V^\optchar }{\partial x_1} \dots \frac{\partial V^\optchar }{\partial
x_n} \right] = (\nabla_x V^\optchar)^T.
\end{eqnarray*}
Thus, we have obtained the following partial differential equation
which the value function must satisfy if it is smooth.  The resulting
equation \eq eq4.2.1/ is the \defn{Hamilton-Jacobi-Bellman (HJB)
equation},
\begin{equation}
- \frac{\partial V^\optchar }{\partial t} \, (x,t) = \min_u \left[
\ell(x,u,t) + \frac{\partial V^\optchar}{\partial x} (x,t) f(x,u,t)
\right]. \elabel{eq4.2.1}
\end{equation}
The terminal penalty term gives a boundary condition for the HJB
equation
\begin{eqnarray}
V^\optchar (x(t_1),{t_1}) = m (x(t_1)) .
\elabel{terminal}
\end{eqnarray}

The term in brackets in \eq eq4.2.1/ is called the \defn{Hamiltonian},
and is denoted $H$:
\begin{eqnarray}
H (x, p, u, t) \eqdef \ell(x,u,t) + p^T f(x,u,t), \elabel{eq4.2.2}
\end{eqnarray}
where $p=\nabla_x V^\optchar$.  We thus arrive at the following
\begin{theorem}
\tlabel{HJB} If the value function $V^\optchar$ has continuous partial
derivatives, then it satisfies the following partial differential
equation
\begin{eqnarray*}
- \frac{\partial V^\optchar }{\partial t}\, (x,t) = \min_u H (x,
\nabla_x V^\optchar\, (x,t),u,t),
\end{eqnarray*}
and the optimal control $u^\optchar(t)$ and corresponding state
trajectory $x^\optchar(t)$ must satisfy \begin{equation} \min_u H
(x^\optchar(t), \nabla_x V^\optchar\, (x^\optchar(t),t),u,t) = H
(x^\optchar(t), \nabla_x V^\optchar\,
(x^\optchar(t),t),u^\optchar(t),t).  \elabel{HJBcontrolDef}
\end{equation}
Conversely, if there exists a value function $V^\optchar (x,t)$
solving \eq eq4.2.1/ subject to \eq terminal/, and a control
$u^\optchar$ obtained from \eq HJBcontrolDef/, then $u^\optchar$
is the optimal controller minimizing \eq cost/ for the system \eq
optFormulation/ and $V^\optchar (x_0, t_0)$ is the minimum value
of \eq cost/.  \qed
\end{theorem}

An important consequence of \Theorem{HJB} is that the optimal control
can be written in state feedback form $u^\optchar(t)
=\baru(x^\optchar(t),t)$, where the function $\baru$ is defined
through the minimization in \eq HJBcontrolDef/.

\proof {\bf Proof of Theorem 10.2.1} The first part of the statement
of the theorem has already been proved through the derivation that led
to \eq eq4.2.1/. For the second part of the statement (the
sufficiently part), let us first note that for an arbitrary control $u$, from \eq eq4.2.1/:
\begin{eqnarray*}
- \frac{\partial V^\optchar }{\partial t} \, (x,t) & = &
 \ell(x,u^\optchar,t) + \frac{\partial V^\optchar}{\partial x} (x,t) f(x,u^\optchar,t)\\
& \leq & \ell(x,u,t) + \frac{\partial V^\optchar}{\partial x} (x,t) f(x,u,t)
\end{eqnarray*}
which is equivalent to
\begin{equation}
0 = \ell (x, u^\optchar, t) + \ddt  V^\optchar (x,y) \leq
\ell (x,u,t) +\ddt  V^\optchar (x,t)
\elabel{diamond}
\end{equation}
where the total derivative (with respect to $t$) is evaluated on
the trajectory generated by \eq optFormulation/, with $u =
u^\optchar$ on the left hand-side and arbitrary $u$ on the right
hand-side.  Now integrate both side of the inequality in \eq
diamond/ from $t_0$ to $t_1$, the left hand-side over the
trajectory of \eq optFormulation/ with $u = u^\optchar$, and the
right hand-side also over the trajectory of \eq optFormulation/
but for an arbitrary $u$.  Further use \eq terminal/ to arrive at
\begin{eqnarray*}
V^\optchar (x_0, t_0) & = & m (x^\optchar (t_1)) +
\int_{t_0}^{t_1} \ell (x^\optchar, u^\optchar, t) \, dt  \\
& \leq & m (x (t_1)) + \int_{t_0}^{t_1} \ell (x, u`, t) \, dt
\end{eqnarray*}
where $x^\optchar$ is the trajectory from \eq optFormulation/
corresponding to $u^\optchar$.  The inequality above shows that
$u^\optchar$ is indeed optimal, and $V^\optchar (x_0, t_0)$ is the
optimal value of the cost.
\qed

\begin{ex}
Consider the simple integrator model, with the polynomial cost
criterion
\begin{eqnarray*}
\dot x = u \qquad V(u) = \int_0^{t_1} [ u^2 + x^4 ] \, dt
\end{eqnarray*}
Here we have $f(x,u,t)=u$, $\ell(x,u,t) = u^2 + x^4$, and $m(x,{t_1})
\equiv 0$.  The Hamiltonian is thus
\begin{eqnarray*}
H(x,p,u,t) = p u +u^2 + x^4,
\end{eqnarray*}
and the HJB equation becomes
\begin{eqnarray*}
- \frac{\partial V^\optchar }{\partial t} = \min_u\big\{\frac{\partial
V^\optchar }{\partial x} u +u^2 + x^4 \big\}.
\end{eqnarray*}

Minimizing with respect to $u$ gives
\begin{eqnarray*}
u^\optchar = -\half \frac{\partial V^\optchar }{\partial x}\, (x,t),
\end{eqnarray*}
which is a form of state feedback.  The closed-loop system has the
appealing form
\begin{eqnarray*}
\dot x^\optchar (t) = -\half \frac{\partial V^\optchar }{\partial x}
(x^\optchar(t),t)
\end{eqnarray*}
This equation shows that the control forces the state to move in the
direction in which the ``cost to go'' $V^\optchar$ decreases.

Substituting the formula for $u^\optchar$ back into the HJB equation
gives the PDE
\begin{eqnarray*}
- \frac{\partial V^\optchar }{\partial t} \, (x,t) =
-\frac{1}{4}\left(\frac{\partial V^\optchar }{\partial x} \, (x,t)
\right)^2 + x^4,
\end{eqnarray*}
with the boundary condition $V^\optchar(x,{t_1})=0$.  This is as far
as we can go, since we do not have available methods to solve a PDE of
this form.  If a solution is required, it may be found numerically.
However, a simpler set of equations is obtained in the limit as
${t_1}\to \infty$.  This simpler problem is treated in
\Exercise{InfHorHJBexample} below.
\end{ex}

We now leave the general nonlinear model and concentrate on linear
systems with quadratic cost.  We will return to the more general
problem in \Chapter{minimum}.

\section{A solution to the LQR problem}
For the remainder of this chapter we consider the LQR problem whose
system description and cost are given in \eq LQR/.  For this control
problem, $x_0$, $t_0$ and ${t_1}$ are given, and $x ({t_1})$ is
free. To ensure that this problem has a solution we assume that $R$ is
strictly positive definite $(R>0)$.

To begin, we now compute the optimal control $u^\optchar$ by solving
the HJB partial differential equation.  The Hamiltonian for this
problem is given by
\begin{eqnarray*}
H(x,p,u,t) = \ell + p^T f = x^T Qx + u^T Ru + p^T (Ax + Bu).
\end{eqnarray*}
To minimize $H$ with respect to $u$ we compute the derivative
\begin{eqnarray*}
\nabla_u H = 0 + 2 Ru + B^T p.
\end{eqnarray*}
Setting this equal to zero, we see that the optimal control is given
by
\begin{equation}
u^\optchar = - \frac{1}{2} R^{-1} B^T p \elabel{eq4.2.4}
\end{equation}
where we recall that $p =\nabla_x V^\optchar $.  Hence the controlled
state evolves according to the equation
\begin{eqnarray*}
\dot x\, (t) = A(t) x(t) - \frac{1}{2} B R^{-1} B^T \nabla_x
V^\optchar\, (x,t).
\end{eqnarray*}
Since the matrix $B R^{-1} B^T$ is positive definite, we see once
again that the control tends to force $x$ in a direction in which
$V^\optchar$ is decreasing.

Equation \eq eq4.2.4/ shows that the optimal control may be written in
state feedback form.  However, we cannot compute the optimal control
law unless we can compute the value function.  Substituting the
formula for the optimal control into the HJB equation gives
\begin{eqnarray*} - \frac{\partial V^\optchar }{\partial t} &=& H(x,
\nabla_x V^\optchar,u^\optchar, t) \\ &=& {x }^T Qx + \frac{1}{4}
(\nabla_x V^\optchar)^T BR^{-1} RR^{-1} B^T \nabla_x V^\optchar \\
&&\qquad + (\nabla_x V^\optchar)^TAx - \frac{1}{2} (\nabla_x
V^\optchar)^T BR^{-1} B^T (\nabla_x V^\optchar).
\end{eqnarray*}
This yields the partial differential equation
\begin{equation}
- \frac{\partial V^\optchar }{\partial t} = {x }^T Qx + (\nabla_x
V^\optchar)^T Ax^\optchar - \frac{1}{4} (\nabla_x V^\optchar)^T
BR^{-1} B^T (\nabla_x V^\optchar). \elabel{lqrHJB}
\end{equation}
To solve this PDE, start at time ${t_1}$ when the value function is
known:
\begin{eqnarray*}
V^\optchar (x,{t_1}) = x^T M x.
\end{eqnarray*}
The terminal $V^\optchar$ is quadratic, which suggests that we try a
quadratic solution for $V^\optchar$ over the entire interval
\begin{eqnarray*}
V^\optchar (x,t) = x^T P (t) x
\end{eqnarray*}
where $P$ is a positive semidefinite, $n\times n$ matrix.  If
$V^\optchar$ is of this form then we have
\begin{eqnarray*}
\frac{\partial V^\optchar}{\partial t} = x^T \dotP x \qquad \nabla_x
V^\optchar = 2 P (t) x.
\end{eqnarray*}
Thus, the PDE \eq lqrHJB/ becomes
\begin{eqnarray*}
-x^T \dotP x = x^T Qx + 2x^T PAx - x^T PBR^{-1} B^T Px.
\end{eqnarray*}
This can be equivalently expressed
\begin{eqnarray*}
-x^T \dotP x = x^T Qx + x^T (PA + A^T P)x - x^T PBR^{-1} B^T Px.
\end{eqnarray*}
Since this equation must hold for \textit{any} $x$, the matrix $P$
must satisfy the following matrix differential equation:
\begin{equation}
- \dotP = Q + PA + A^T P - PBR^{-1} B^T P. \elabel{eq4.2.5}
\end{equation}
This is an \textit{ordinary} differential equation for $P$, but it is
time-varying, and nonlinear due to the term which is quadratic in $P$.
Adding the boundary condition
\begin{eqnarray*}
V (x, t_1) = x^T M x = x^T (t_1) Px (t_1)
\end{eqnarray*}
gives
\begin{equation}
P ({t_1}) = M. \elabel{eq4.2.6}
\end{equation}
Equation \eq eq4.2.5/ is called the \defn{Riccati Differential
Equation} (RDE), and its boundary condition is given by \eq
eq4.2.6/. Since the boundary condition is given at the final time, the
solution to this differential equation may be viewed as
\textit{starting} at time $t_1$, and then traveling backwards in time
until time $t_0$.  The Riccati equation possesses a unique solution,
provided that the matrices $A$, $B$, $Q$ and $R$ are piecewise
continuous in $t$.  Further note that this unique solution is
symmetric, because $P^T$ satisfies the same differential equation (to
see this simply take transpose of both side of \eq eq4.2.5/).

The solution to the LQR problem is now summarized in the following
\begin{theorem}\tlabel{LQRsolution}
For the LQR optimal control problem \eq LQR/, assuming that $A$, $B$,
$Q$ and $R$ are piecewise continuous in $t$, and that $R(t)>0$ for all
$t$, \balphlist
\item
The optimal control $u^\optchar$ is given in linear state feedback
form $u^\optchar(t) = -K (t) x^\optchar(t)$, where $K (t)$ is a matrix
of time-varying gains:
\begin{eqnarray*}
u^\optchar = - \half R^{-1} B^T \nabla_x V^\optchar = - R^{-1}(t)
B^T(t) P(t)x^\optchar (t) = -K (t) x^\optchar (t).
\end{eqnarray*}

\item
The matrix-valued function $P (t)$ is positive semidefinite for each
$t$, and is defined by the RDE \eq eq4.2.5/ with boundary condition
\eq eq4.2.6/.

\item
The value function $V^\optchar$ is quadratic:
\begin{eqnarray*}
V^\optchar (x_0, t_0) = x_0^T P (t_0) x_0.
\end{eqnarray*}
\end{list}
\qed
\end{theorem}

\section{The Hamiltonian matrix}\label{s:HamMatrix}

We now show how the Riccati equation can be solved by computing the
solution to a \textit{linear} ODE.  In addition to providing a
complete solution to the LQR problem, this will provide a solution to
the infinite time horizon optimal control problem where
${t_1}=\infty$.  Consider the $2n$ dimensional LTV model, where the
state $(X(t),Y(t))'$ is a $2n\times n$ matrix-valued function of $t$.
\begin{eqnarray}
\begin{bmatrix} \dotX (t) \\ \dotY (t) \end{bmatrix} & = &
\underbrace{\left[ \begin{array}{cc} A & -BR^{-1} B^T \\ -Q & -A^T
\end{array}
\right]}_{2n \times 2n} \underbrace{
\begin{bmatrix} X (t) \\Y (t) \end{bmatrix}}_{2n \times n},
\qquad \begin{bmatrix} X ({t_1}) \\ Y ({t_1}) \end{bmatrix} =
\begin{bmatrix} I \\ M \end{bmatrix}
\elabel{eq4.2.8}
\end{eqnarray}
The $2n \times 2n$ system-matrix of \eq eq4.2.8/ is called the
\defn{Hamiltonian matrix}, denoted by $\underH$.

\begin{theorem}
\tlabel{HamiltonianLQR} The solution to the RDE \eq eq4.2.5/ with
boundary condition \eq eq4.2.6/ is given by the formula
\begin{eqnarray*}
P (t) = Y (t) X^{-1} (t), && t_0\le t\le t_1.
\end{eqnarray*}
\end{theorem}

\proof To prove the theorem, observe that the desired initial
conditions are satisfied:
\begin{eqnarray*}
\underbrace{Y ({t_1})}_{M} \underbrace{X^{-1} ({t_1})}_{I} = M = P
({t_1}).
\end{eqnarray*}
We next verify that the Riccati differential equation is also
satisfied.  Using the product rule and the matricial quotient rule
$\frac{d}{dt} (X^{-1}) = - X^{-1} \dotX X^{-1}$, we have
\begin{eqnarray*}
\frac{d}{dt} (YX^{-1}) = Y \frac{d}{dt} (X^{-1}) + \dotY X^{-1} = -
YX^{-1} \dotX X^{-1} + \dotY X^{-1}.
\end{eqnarray*}
Substituting the definitions of $\dot X$ and $\dot Y$ then gives
\begin{eqnarray*}
\frac{d}{dt} (YX^{-1}) &=& - YX^{-1} (AX - BR^{-1} B^T Y) X^{-1} + (-
QX - A^T Y) X^{-1} \\ &=& - (YX^{-1}) A + (YX^{-1}) BR^{-1} B^T
(YX^{-1}) - Q - A^T (YX^{-1})
\end{eqnarray*}
which is precisely Riccati equation for $P = YX^{-1}$.  Note that this
immediately tells us that $YX^{-1}$ is symmetric and postive
semidefinite.  \qed

In the LTI case we may extend this analysis further to explicitly
solve the RDE.  We first require the following result, which shows
that the eigenvalues of the Hamiltonian matrix possess symmetry about
the complex axis in $\Co$.  \begin{lemma} \tlabel{H-evalues} For an
LTI model, if $\lambda$ is an eigenvalue of $\underH$, so is $-
\lambda$.
\end{lemma}

\proof Define the $2n\times 2n$ matrix $J$ as
\begin{eqnarray*}
J \eqdef \left[ \begin{array}{rr} 0 & I \\ -I & 0\end{array} \right],
\end{eqnarray*}
where each element is $n \times n$.  The inverse of $J$ is equal to
\begin{eqnarray*}
J^{-1} = \left[ \begin{array}{rr} 0 & -I \\ I & 0 \end{array}\right] =
J^T.
\end{eqnarray*}
Using the definition of $\underH$ gives
\begin{eqnarray*}
J \underH J^T = - \underH^T.
\end{eqnarray*}
By the definition of similarity transformations, we have for any
eigenvalue $\lambda$ of $\underH$,
\begin{eqnarray*}
0 = \det (\lambda I - \underH) = \det (\lambda I - J \underH J^{-1}).
\end{eqnarray*}
Moreover, for this special case
\begin{eqnarray*}
\det (\lambda I - J \underH J^{-1}) = \det (\lambda I + \underH^T)
=(-1)^n \det ((-\lambda) I - J \underH J^{-1}).
\end{eqnarray*}
These two equations combined show that when $\lambda$ is equal to an
eigenvalue of $\underH$ then so is $-\lambda$, as claimed.  \qed

We will see below that the matrix $\underH$ cannot have eigenvalues on
the $j\omega$ axis, provided the model is stabilizable and detectable.
Consider the special case where the eigenvalues of $\underH$ are
distinct.  In this case we can solve the RDE by solving a linear
differential equation.  We first diagonalize the Hamiltonian matrix:
\begin{eqnarray*}
U^{-1} \underH U = \left[ \begin{array}{rr} \Lambda_s & 0 \\ 0 &
-\Lambda_s \end{array} \right]
\end{eqnarray*}
where the diagonal elements of $\Lambda_s$ lie in the strict left hand
plane in $\Co$.  Note that if the eigenvalues of $\underH$ are
distinct then, by \Lemma{H-evalues}, none can lie on the $j\omega$
axis.  The matrix of eigenvectors $U$ can be written
\begin{eqnarray*}
U = \left[ \begin{array}{c} U_{11} \\ U_{21} \end{array} \left|
\begin{array}{c} U_{12} \\ U_{22} \end{array} \right. \right]
\end{eqnarray*}
where $U_{11}, U_{21}$ correspond to the stable eigenvalues. The
solution to the differential equation \eq eq4.2.8/ can be computed in
this case:
\begin{eqnarray}
\begin{bmatrix} X (t) \\ Y (t) \end{bmatrix} & = & U
\left[ \begin{array}{cc} e^{\Lambda_s (t-{t_1})} & 0 \\ 0 &
e^{-\Lambda_s (t-{t_1})} \end{array} \right] U^{-1} \begin{bmatrix} I
\\ M
\end{bmatrix}
\elabel{eq4.2.10}
\end{eqnarray}
This is a special case of the solution to an LTI model
\begin{eqnarray*}
x (t) = \phi (t - {t_1}) x ({t_1})
\end{eqnarray*}
using the state transition matrix $\phi$.

With considerable algebra, equation \eq eq4.2.10/ can be solved to
give closed-form expressions for $X (t)$ and $Y (t)$, which can then
be combined to form $P (t, {t_1})$:
\begin{equation}
P (t,{t_1}) = [ U_{21} + U_{22} e^{-\Lambda_s (t-{t_1})} G e^{-
\Lambda_s (t-{t_1})} ] [U_{11} + U_{12} e^{-\Lambda_s (t-{t_1})} G
e^{- \Lambda_s (t-{t_1})} ]^{-1} \elabel{eq4.2.11}
\end{equation}
where
\begin{eqnarray*}
G = - [U_{22} - MU_{12}]^{-1} [U_{21} - MU_{11}].
\end{eqnarray*}
Hence the RDE can be solved explicitly for an LTI model.

\section{Infinite horizon  regulator}
\label{s:infHor}

We now consider the infinite-horizon or steady-state control problem.
We assume that the model is LTI, and since the time horizon is
infinite, the cost only depends upon the initial condition of the
model:
\begin{equation}
\begin{aligned}
\dot x & = Ax + Bu, \quad x (0) = x_0 \\ V(x) & = \int_0^\infty
(x^T(t) Qx(t) + u^T(t) Ru(t))\, dt
\end{aligned}
\elabel{LQRss}
\end{equation}
There is no need to include a terminal penalty function because the
assumptions we impose will imply that $x(t)\to \zero$ as $t\to\infty$.

Our approach is to consider an optimal control problem on the finite
interval $[0,t_1]$, and then let $t_1\to\infty$.  Let $P(t,t_1)$
denote the solution to the RDE with time horizon $t_1$, so that the
optimal control may be written $u(t)=-K(t,t_1)x(t)=-R^{-1}B^T P(t,t_1)
x(t)$.  Suppose that $P(t,t_1)\to \barP$ as $t_1\to\infty$, where
$\barP$ is an $n\times n$ matrix independent of time.  Then the
limiting control law $u(t)=-K x(t)=-R^{-1}B^T \barP x(t)$ is optimal
for the infinite horizon control problem.  To see this, observe that
for any control $u$ on $[0,\infty)$, since $x^T P(0,t_1) x$ is equal
to the optimal cost for the finite horizon problem,
\begin{eqnarray*}
x^T P(0,t_1) x & \le & \int_0^{t_1} (x^T(t) Qx(t) + u^T(t) Ru(t))\, dt \\
& \le &
\int_0^\infty (x^T(t) Qx(t) + u^T(t) Ru(t))\, dt.
\end{eqnarray*}
Letting $V^\optchar$ denote the optimal cost for the infinite horizon
problem, it follows that $x^T P(0,t_1) x \le V^\optchar(x)$ for all
$x$ and $t_1$, and hence
\begin{eqnarray*}
\lim_{t_1\to\infty} x^T P(0,t_1) x \le V^\optchar(x).
\end{eqnarray*}
This gives a lower bound on the cost $V^\optchar(x)$, and
 this lower
bound is attained with the limiting controller.

Under certain conditions we can show that the limit does exist, and
thereby solve the infinite horizon control problem.

\begin{ex}
Consider the simple model
\begin{eqnarray*}
\begin{aligned}
\dot x & = x + u \\ V(u) & = \int_0^{t_1} (x^2 + u^2) \, dt + 5
x(t_1)^2.
\end{aligned}
\end{eqnarray*}
The RDE becomes
\begin{eqnarray*}
- \frac{d}{dt} P(t,t_1) = 2P(t,t_1) + 1 - P(t,t_1)^2,\qquad P(t_1,t_1)
= 5.
\end{eqnarray*}
A solution to this differential equation is shown in
\Figure{simpleARE}.  These plots suggest that for any fixed $t$,
$P(t,t_1)$ converges to a constant $\bar P$ as $t_1\to\infty$. If this
is the case, then the limit $\bar P$ must be an equilibrium for the
ARE.  That is,
\begin{eqnarray*}
0 = 2\bar P + 1 - {\bar P}^2.
\end{eqnarray*}
Solving this equation gives $\barP = 1+\sqrt{2}$, which is consistent
with the limiting value shown in the plots.  We will see that the
limiting control law $u=-(1+\sqrt{2}) x$ is optimal for the infinite
horizon cost with $t_1=\infty$.
\end{ex}

\begin{figure}[ht]
\ebox{.6}{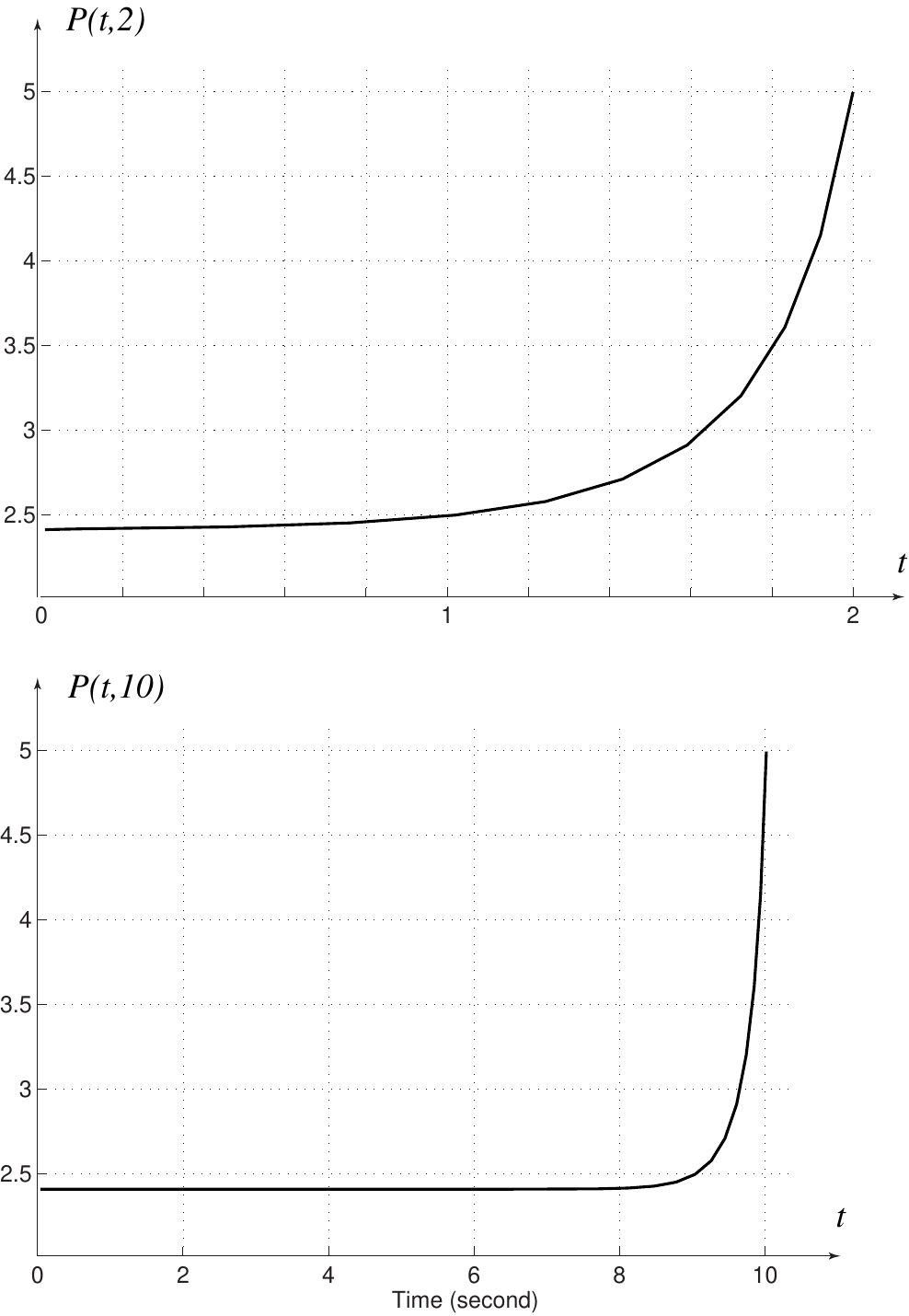}
\caption[Solutions to the ARE for a scalar model]{Solutions to the ARE for a scalar model with $t_1=2$ and
$t_1=10$.  It appears that for any fixed $t$, $P(t,t_1)$ converges to
a constant $\bar P$ as $t_1\to\infty$.}  \flabel{simpleARE}
\end{figure}
\clearpage

Using the formula \eq eq4.2.11/ it is easy to show that the limiting
matrix $\barP$ does exist in the special case where the eigenvalues of
the Hamiltonian are distinct.  Because the matrix $\Lambda_s$ is
Hurwitz, in this case it follows from \eq eq4.2.11/ that
\begin{eqnarray}
\lim_{{t_1} \to \infty} P (t,{t_1}) = U_{21} U_{11}^{-1} \eqdef \barP
\elabel{eq4.2.12}
\end{eqnarray}
where $\barP$ is an $n\times n$ symmetric matrix which does not depend on
$t$. The optimal control for \eq LQRss/ is given by the limiting
formula
\begin{eqnarray*}
u^\optchar(t) = -R^{-1} B^T \barP x^\optchar (t).
\end{eqnarray*}
Note that $u^\optchar$ is given by a state feedback control law of the
form $u^\optchar=-Kx$, which is a form of pole placement.  We will see
in \Section{RDE} that the location of the closed-loop poles possess
desirable properties.

\notes{why is U11 invertible?}

We now use the Hamiltonian matrix to characterize the closed-loop
poles under optimal control. From the eigenvalue equation $\underH U =
U\Lambda$ we have
\begin{equation}
\underH \begin{bmatrix} U_{11} \\ U_{21} \end{bmatrix} = \left[
\begin{array}{cc} A & -BR^{-1} B^T \\ -Q & -A^T \end{array} \right]
\begin{bmatrix} U_{11} \\ U_{21} \end{bmatrix} = \begin{bmatrix}
U_{11} \\ U_{21} \end{bmatrix} \Lambda_s. \elabel{ssH}
\end{equation}
The upper set of equations defined by \eq ssH/ may be written $
U_{11}\Lambda_s = A U_{11} - BR^{-1} B^T U_{21} $, and this can be
transformed to obtain
\begin{eqnarray*}
U_{11}\Lambda_s & = & A U _{11} - BR^{-1} B^T U_{21} (U_{11}^{-1}
 U_{11}) \\ &=&(A - BR^{-1} B^T U_{21} U_{11}^{-1}) U_{11}.
\end{eqnarray*}
Since $K=R^{-1} B^T \bar P$ and $\bar P=U_{21} U_{11}^{-1}$, this
gives the formula
\begin{eqnarray*}
(A - BK) U_{11} = U_{11} \Lambda_s.
\end{eqnarray*}
The matrix $U_{11}$ is thus a modal matrix for the optimal closed-loop
system matrix $A_{cl}=(A-BK)$.  That is,
\begin{theorem}
\tlabel{InfiniteHamiltonian} Assume that the controlled system is
asymptotically stable, so that the closed-loop poles determined by
\begin{eqnarray*}
\Delta(s) =\det[Is - (A-B R^{-1} B^T \bar P)]
\end{eqnarray*}
lie in the strict left half plane.  Then, \balphlist
\item
The stable eigenvalues $\Lambda_s$ for $\underH$ are also the optimal
closed-loop poles.

\item
The columns of $U_{11}$ are the eigenvectors for $A_{cl}=A-B R^{-1}
B^T \bar P$.
\end{list}
\qed
\end{theorem}

Using the Hamiltonian we now give a quadratic algebraic equation which
the matrix $\barP$ must satisfy.  Considering both the first and
second rows of \eq ssH/ gives
\begin{eqnarray}
A U_{11} - BR^{-1} B^T U_{21} & = & U_{11} \Lambda_s
\elabel{eq4.2.14a} \\ -Q U_{11} - A^T U_{21} & = & U_{21} \Lambda_s
\elabel{eq4.2.14b}
\end{eqnarray}
Multiplying the first equation by $(U_{21} U_{11}^{-1})$ on the left
hand side and by $U_{11}^{-1}$ on the right gives
\begin{eqnarray*}
(U_{21} U_{11}^{-1}) A U_{11} U_{11}^{-1} - (U_{21} U_{11}^{-1})
 BR^{-1} B^T U_{21} U_{11}^{-1} = (U_{21} U_{11}^{-1}) U_{11}
 \Lambda_s U_{11}^{-1}
\end{eqnarray*}
or
\begin{eqnarray}
\barP A - \barP BR^{-1} B^T \barP & = & U_{21} \Lambda_s U_{11}^{-1}
\elabel{eq4.2.15}
\end{eqnarray}
Multiplying \eq eq4.2.14b/ by $U_{11}^{-1}$ on the right we obtain
\begin{eqnarray*}
-Q U_{11} U_{11}^{-1} - A^T U_{21} U_{11}^{-1}= U_{21} \Lambda_s
U_{11}^{-1}
\end{eqnarray*}
or \begin{equation} -Q -A^T \barP = U_{21} \Lambda_s U_{11}^{-1}.
\elabel{eq4.2.16}
\end{equation}
Combining \eq eq4.2.15/ and \eq eq4.2.16/ gives $ \barP A - \barP
BR^{-1} B^T \barP = - Q - A^T \barP$, or
\begin{eqnarray}
A^T \barP + \barP A + Q - \barP BR^{-1} B^T \barP = 0 \elabel{ARE}
\end{eqnarray}
This final formula is known as the \defn{algebraic Riccati equation
(ARE)}.  It contains $  n (n+1)/2 $ equations to solve for the components
of $\barP$, since $\bar{P}$ is symmetric.

\begin{ex}
Consider the system/cost combination
\begin{eqnarray*}
\dot x & = & x + u \\ V & = & \int_0^\infty u^2\, dt.
\end{eqnarray*}
We have $R=1$, which is strictly positive, and $Q=0\ge 0$, so our
conditions on the cost function are satisfied.  However, it is obvious
that the optimal control is to set $u$ equal to zero so that $V$ is
also zero. The resulting closed-loop system is unstable, so this is
obviously a poor design!  Clearly, some additional conditions on the
cost function must be imposed.
\end{ex}

\begin{ex}
The difficulty in the previous example is that the cost is not a good
measure of the internal stability of the model.  This is reminiscent
of detectability, a connection that we further explore in two more
examples.

Consider first the LTI model with quadratic cost
\begin{equation}
\begin{aligned}
\dot x & = \begin{bmatrix}-3 & -2 \\ 1 &0 \end{bmatrix}x +
\begin{bmatrix} 0 \\ 1 \end{bmatrix} u; \\
V(u) & = \int_0^\infty
\left[ (x_1+x_2)^2 + u^2 \right] \, dt
\end{aligned}
\elabel{model-detect}
\end{equation}
In this example
\begin{eqnarray*}
A = \begin{bmatrix}-3 & -2 \\ 1 &0 \end{bmatrix}, \qquad B =
\begin{bmatrix} 0 \\ 1 \end{bmatrix}, \qquad Q=\begin{bmatrix}1&1 \\ 1
& 1 \end{bmatrix}\ge 0, \qquad R=1.
\end{eqnarray*}
The matrix $A$ is Hurwitz, with eigenvalues at $-2,-1$.

Letting $y=Cx = x_1 +x_2$, the cost becomes $V(u)=\int (y^2 + u^2) \,
dt$.  Consider a very special case where the initial condition $x(0)$
belongs to the unobservable subspace $\Sigma_{\bar o}$.  What then is
the optimal control?  Recall that if $u\equiv 0$ then the output
becomes
\begin{eqnarray*}
y(t) = C e^{At} x(0), \qquad t\ge 0.
\end{eqnarray*}
Since $x(0)\in \Sigma_{\bar o}$ we then have $y\equiv 0$. Obviously
then, this is the optimal control since it yields $V(u)=0$!  It
follows from the definition of the optimal cost that
\begin{eqnarray*}
V^\circ(x) = 0 \hbox{ for } x\in \Sigma_{\bar o}.
\end{eqnarray*}

This can be seen explicitly for this example by solving the ARE, and
computing $\Sigma_{\bar o}$.  The observability matrix for this
example is
\begin{eqnarray*}
\clO = \begin{bmatrix} C \\ CA
\end{bmatrix} = \begin{bmatrix} 1&1 \\ -2 & -2  \end{bmatrix}.
\end{eqnarray*}
Hence the unobservable subspace becomes
\begin{eqnarray*}
\Sigma_{\bar o} = \clN(\clO) = \Span \left\{\begin{pmatrix} 1 \\ -1
\end{pmatrix}\right\}.
\end{eqnarray*}
Using the {\tt lqr} command in Matlab, the solution to the ARE is
approximately
\begin{eqnarray*}
P = 0.24 \begin{bmatrix}1&1 \\ 1 & 1 \end{bmatrix}
\end{eqnarray*}
Hence we do indeed find that $V^\optchar(x) = x^T P x = 0$ for
$x\in\Sigma_{\bar o}$.

Setting $u=0$ may be reasonable for a stable system such as \eq
model-detect/.  Consider however the reversed situation where the
matrix $A$ is replaced by $-A$:
\begin{eqnarray*}
\dot x = \begin{bmatrix}3 & 2 \\ -1 &0 \end{bmatrix}x +
\begin{bmatrix} 0 \\ 1 \end{bmatrix} u.
\end{eqnarray*}
The new state matrix has eigenvalues at $+2, +1$.  Hence for any
non-zero initial condition the state explodes as $t\to\infty$ when the
input is set to zero.  For this system can we again conclude that a
zero input is optimal?  For this cost criterion, the answer is
\textit{yes}, for precisely the same reason as before.  The
observability matrix for this model is
\begin{eqnarray*}
\clO = \begin{bmatrix} C \\ CA \end{bmatrix} = \begin{bmatrix} 1&1 \\
2 & 2 \end{bmatrix},
\end{eqnarray*}
so that the unobservable subspace is again equal to $\Sigma_{\bar o} $
is again equal to the span of the single vector $(1,-1)^T$.  If
$x(0)\in \Sigma_{\bar o}$, and $u\equiv 0$, then $y=Cx=x_1+x_2$ is
identically zero.

We conclude that $V(u) = 0$, even though the state $x(t)\to\infty$ if
$x(0)\neq 0$ lies in the unobservable subspace.  This unstable
behavior of the optimized system is a consequence of a lack of
\textit{detectability}.
\end{ex}

For a general LQR problem, suppose that the matrix $Q$ is factored as
$Q= C^T C$ for a $p\times n$ matrix $C$. As in the previous example,
we then \textit{define} $y=C x$ so that the cost function becomes
\begin{eqnarray*}
V = \int_0^\infty (x^T Qx + u^T Ru)\, dt= \int_0^\infty (y^T y + u^T
Ru)\, dt.
\end{eqnarray*}
If the pair $(A,C)$ is \textit{observable}, then for any control $u$,
and any initial state $x(0)\neq \zero$, the cost $V$ is strictly
positive.  To see this, note that if $u$ is identically zero, then by
observability we must have $\int x^T Qx\, dt = \int |y|^2 \, dt >0$.
If $u$ is not identically zero, then $V>0$ since $R>0$. Since for an
initial condition $x_0$ the value function may be expressed
$V^\optchar(x_0) = x_0^T\barP x_0$, it follows that the matrix $\barP$
is positive definite ($\barP>0$). From the previous example it is
clear that $P$ will be singular if the system is not observable. If
the system is not detectable, then the last example shows that we
might expect disaster.  These results can be refined to give the
following

\begin{theorem}
\tlabel{LQR} For the LTI system \eq LQRss/ with $R > 0$ and $Q=C^T
C\ge 0$,

\begin{description}
\item[(i)] If $(A,B)$ is stabilizable, and $(A,C)$ is detectable
then \balphlist \item There exists a positive semi-definite
solution $\barP$ to the ARE \eq ARE/, which is unique in the class
of positive semi-definite matrices.

\item The closed-loop system matrix $A_{cl} = A - BR^{-1} B^T \barP$
is Hurwitz.

\item The infinite horizon optimal control is $u^\optchar = - R^{-1}
B^T \barP x$, and the infinite horizon optimal cost is $V^\optchar
(x_0) = x_0^T \barP x_0$.

\item If $P (t, t_1)$ is the solution to the RDE \eq eq4.2.5/ subject
to $P (t_1, t_1) = 0$, then $\lim_{t_1 \to \infty} P (t, t_1) =
\barP$.
\end{list}

\item[(ii)] If $(A,B)$ is stabilizable, and $(A,C)$ is observable, then
(a)--(d) hold, and in addition the matrix $\barP$ is positive
definite.
\end{description}
\qed
\end{theorem}

\begin{ex}
Consider the multivariate model
\begin{eqnarray*}
\dot x & = & \left[ \begin{array}{rr} 0 & -1 \\ 0 & 0 \end{array}
\right] x + \left[ \begin{array}{rr} 1 & 0 \\ 0 & 1 \end{array}
\right] u \\ V & = & \int_0^\infty (4x_1^2 + 4x_1 x_2 + x_2^2 + u_1^2
+ u_2^2)\, dt.
\end{eqnarray*}
From these equations we may deduce that the weighting matrices are
\begin{eqnarray*}
Q & = & \left[ \begin{array}{rr} 4 & 2 \\ 2 & 1 \end{array} \right],
\qquad R = \left[ \begin{array}{rr} 1 & 0 \\ 0 & 1 \end{array}
\right].
\end{eqnarray*}
The matrix $Q$ has rank one, and $Q \geq 0$. It follows from the
spectral representation of a positive definite matrix that there is a
$1\times 2$ matrix $C$ such that $Q= C^T C$.  In fact $C= (2,1)$, so
we define $y= Cx= 2x_1 + x_2$.  Thus, the controllability and
observability matrices become
\begin{eqnarray*}
\clC = \left[ \begin{array}{rrrr} 1 & 0 & 0 &-1 \\ 0 & 1 & 0 &
0\end{array} \right] \qquad \clO = \left[ \begin{array}{rr} 2 & 1 \\ 0
& -2 \end{array} \right].
\end{eqnarray*}
Since both of these matrices have rank 2, the conditions of
\Theorem{LQR}~(ii) are satisfied.

We now check directly to see if the controlled system is stable by
computing the closed-loop system poles.  The Hamiltonian matrix for
this model is
\begin{eqnarray*}
\underH & = & \left[ \begin{array}{rr|rr} 0 & -1 & -1 & 0 \\ 0 & 0 & 0
& -1 \\ \hline -4 & -2 & 0 & 0 \\ -2 & -1 & 1 & 0 \end{array} \right].
\end{eqnarray*}
The eigenvalues of $\underH$ are found by solving the characteristic
equation
\begin{eqnarray*}
\det [sI - \underH] & = & s^4 - 5s^2 + 4 \\ & = & (s-2)(s+2)(s-1)(s+1).
\end{eqnarray*}
Note that the roots of this equation possess the symmetry required by
\Lemma{H-evalues}.  The eigenvectors $\{w^i\}$ of the stable modes
$\{\lambda^i\}$ are
\begin{eqnarray*}
\lambda_1 = -1, \quad w^1=\begin{bmatrix} 1 \\ 0 \\ 2 \\ 0
\end{bmatrix}, \qquad \lambda_2 = - 2, \quad w^2=\left[
\begin{array}{r} 1 \\ -1 \\ 2 \\ -1
\end{array} \right].
\end{eqnarray*}
Thus, the matrices $\{U_{ij} \}$ may be computed:
\begin{eqnarray*}
\left[ \begin{array}{c} U_{11} \\ \hline U_{21} \end{array} \right] =
\left[ \begin{array}{c} \begin{array}{rr} 1 & 1 \\ 0 & -1 \\
\end{array} \\ \hline \begin{array}{rr} 2 & 2 \\ 0 & -1 \end{array}
\end{array} \right].
\end{eqnarray*}
Applying \eq eq4.2.12/ we can compute the solution to the ARE:
\begin{eqnarray*}
\barP = \left[ \begin{array}{rr} 2 & 2 \\ 0 & -1 \end{array} \right]
\left[ \begin{array}{rr} 1 & 1 \\ 0 & -1 \end{array} \right]^{-1} =
\left[ \begin{array}{rr} 2 & 0 \\ 0 & 1 \end{array} \right],
\end{eqnarray*}
and the optimal feedback gain matrix is
\begin{eqnarray*}
K = R^{-1} B^T \barP = \left[ \begin{array}{rr} 2 & 0 \\ 0 & 1
\end{array} \right].
\end{eqnarray*}

The resulting optimal control
\begin{eqnarray*}
u^\optchar = \begin{pmatrix} -2x_1 \\ -x_2 \end{pmatrix}.
\end{eqnarray*}
places the closed-loop eigenvalues at $(-2,-1)$, which are the stable
eigenvalues of $ \underH$.  The closed-loop eigenvectors are contained
in the columns of $U_{11}$.
\end{ex}

\begin{ex}
To show how the ARE can be solved directly without constructing the
Hamiltonian, consider the model
\begin{equation}
\begin{aligned}
\dot x & = \left[ \begin{array}{rr} 0 & 1 \\ 0 & -1 \end{array}
\right] x + \begin{bmatrix} 0 \\ 1 \end{bmatrix} u; \qquad y = x_1 \\
V & = \int_0^\infty (x_1^2 + u^2)\, dt
\end{aligned}
\elabel{AREex}
\end{equation}

In this example we have $R=1$, and $Q = C^T C$, where $C=(1\ 0)$.  The
controllability and observability matrices become
\begin{eqnarray*}
\clC = \left[ \begin{array}{rr} 0 & 1 \\ 1 & -1 \end{array} \right]
\qquad \clO = \left[ \begin{array}{rr} 1 & 0 \\ 0 & 1 \end{array}
\right].
\end{eqnarray*}
Again we see that the conditions of \Theorem{LQR}~(ii) are satisfied,
and it follows that the solution $\barP$ to the ARE will be positive
definite. We now verify this directly.

The ARE for this model is
\begin{eqnarray*}
A^T \barP + \barP A - \barP BR^{-1} B^T \barP + Q = 0.
\end{eqnarray*}
Writing $\barP$ as $ \barP = \left[ \begin{array}{rr} \ell & m \\ m &
n \end{array} \right]$, this gives the equation
\begin{eqnarray*}
\left[ \begin{array}{rr} 0 & 0 \\ \ell - m & m - n \end{array} \right]
+ \left[ \begin{array}{rr} 0 & \ell - m \\ 0 & m-n \end{array} \right]
- \left[ \begin{array}{rr} m^2 & mn \\ mn & n^2 \end{array} \right] +
\left[ \begin{array}{rr} 1 & 0 \\ 0 & 0 \end{array} \right] = 0.
\end{eqnarray*}
Examining the $(1,1)$ entry of each matrix gives
\begin{eqnarray*}
-m^2 = -1 \Rightarrow m = +1 \ \hbox{(since $\barP\ge 0$).}
\end{eqnarray*}
Examining the $(1,2)$ entry of the ARE we see that $\ell - 1- n =0$,
and from the $(2,2)$ entry $n^2 +2n -2=0$.  These equations result in
the unique positive definite solution
\begin{eqnarray*}
\barP = \left[ \begin{array}{rr} \sqrt{3} & 1 \\ 1 & - 1 +\sqrt{3}
\end{array} \right]
\end{eqnarray*}
The optimal feedback control law is $ u^\optchar = -R^{-1}B^T P x = -
[1 \ \sqrt{3}-1] x$.

To find the closed-loop poles, we must compute the eigenvalues of
\begin{eqnarray*}
A_{cl} = A-BK = \left[ \begin{array}{rr} 0 & 1 \\ 0 & -1 \end{array}
\right] - \left[ \begin{array}{rr} 0 & 0 \\ 1 & \sqrt{3}-1 \end{array}
\right] = \left[ \begin{array}{rr} 0 & 1 \\ -1 & -\sqrt{3} \end{array}
\right].
\end{eqnarray*}
The characteristic equation is $\Delta(s) = s^2 +2\frac{\sqrt{3}}{2} s
+ 1 = 0$, which gives
\begin{eqnarray*}
p_i = -\half(\sqrt{3}\pm j).
\end{eqnarray*}
The closed-loop response will be somewhat under-damped, but of course
it is stable.
\end{ex}

\section{Return difference equation}
\label{s:RDE} The LQR solution has some striking and highly desirable
properties when viewed in the frequency domain. These results have
many applications to design and analysis of linear control systems. In
particular, \balphlist
\item
For a SISO system and certain MIMO systems, a generalization of the
root locus method can be used to find the optimal state feedback
controller.

\item
For an LQR design based upon full state feedback, the closed-loop
system is insensitive to plant uncertainty, regardless of the choice
of weighting matrices $Q$ and $R$.  In particular, the gain margin is
always infinite, and the phase margin at least 60\%.
\end{list}
The conclusion (b) is surprising, given the purely time domain
approach that we have taken. In the end this serves as perhaps the
greatest motivation for using LQR methods in control design.

For the LTI model with infinite horizon cost criterion expressed in
Equation \eq LQRss/, we have the following identity for any $s\in\Co$,
known as the \defn{return difference equation}:
\begin{eqnarray}
&& R + B^T (-sI -A)^{-T} Q (sI - A)^{-1} B \nonumber \\
&& \qquad = [I + K (-sI - A)^{-1} B]^T R
[I + K (sI - A)^{-1} B]. \elabel{ReturnDiff}
\end{eqnarray}
The vector $K$ is the optimal feedback gain
\begin{eqnarray*}
K = R^{-1} B^T \barP.
\end{eqnarray*}
The left hand side of \eq ReturnDiff/ involves only the open-loop
system, while the right hand side depends on the feedback gain $K$.
The return difference equation is derived by simply expanding both
sides of the equation, and substituting the definition of $K$ and the
algebraic Riccati equation.

The return difference equation appears complex, but after the matrices
are given proper interpretation it has a relatively simple form.
First, factor the state weighting matrix as $Q = C^T C$.  By defining
$y=Cx$ we obtain the controlled model illustrated in
\Figure{optimal-blockRD}.  With $P(s) = C(Is-A)^{-1} B$, the left hand
side of \eq ReturnDiff/ becomes
\begin{eqnarray*}
R + P^T (-s) P (s).
\end{eqnarray*}
From the feedback from $u = -Kx$ we define the \defn{loop transfer
function} described in
\Chapter{goals} as $L(s) = K(Is-A)^{-1}B$, so that the  right hand
side of \eq ReturnDiff/ becomes
\begin{eqnarray*}
[I + L (-s)]^T R [I + L (s)].
\end{eqnarray*}
Hence the return difference equation can be expressed succinctly as
\begin{equation}
R + P^T (-s) P (s) = [I + L (-s)]^T R [I + L (s)].
\elabel{ReturnDiff2}
\end{equation}
On the complex axis we have $H^T(-j\omega) = H^*(j\omega)$ for any
transfer function $H$, so that \eq ReturnDiff2/ may be written
\begin{equation}
R + P^* (j\omega) P (j\omega) = [I + L (j\omega)]^* R [I + L
(j\omega)].  \elabel{ReturnDiff3}
\end{equation}

\begin{figure}[ht]
\ebox{.75}{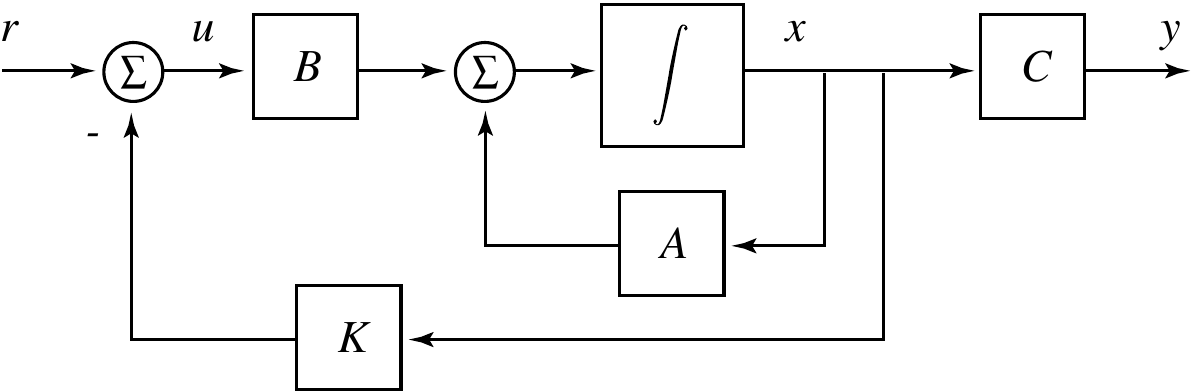}
\caption[The loop transfer function]{The loop transfer function of the LTI system under full state
feedback $u=-Kx + r$ is $L(s) = K(Is-A)^{-1} B$.}
\flabel{optimal-blockRD}
\end{figure}

In the SISO case with $R=r>0$ a scalar, the return difference equation
becomes
\begin{equation}
r + P(-s) P (s) = r [1+L(-s)][1+L (s)].  \elabel{eq5.1.7}
\end{equation}
The poles of the closed-loop system are found by solving the equation
\begin{eqnarray*}
1 + L (s) = 0.
\end{eqnarray*}
In view of \eq eq5.1.7/, the return difference equation gives the
closed-loop poles as a function of the weighting parameter $r$ and the
open-loop transfer function $P(s)$ through the formula
\begin{equation}
1 + \frac{1}{r} P (-s) P (s) = 0. \elabel{eq5.1.8}
\end{equation}
The set of all solutions to \eq eq5.1.8/, with $r$ ranging from $0$ to
$\infty$, is called the \defn{symmetric root locus}.  For fixed $r$,
the roots of this equation are the eigenvalues of the Hamiltonian
matrix $\underH$ with $R=r$.  We have seen in
\Section{infHor} that the stable (left hand plane) eigenvalues of the
Hamiltonian matrix are the poles of the closed-loop system.

\notes{this should be a formal example}

As an illustration consider the general second order system with
transfer function $P(s) = (s-z)/(s^2 + 2\zeta s +1)$, where $z<0$, and
$\zeta>0$. The symmetric root locus is illustrated in
\Figure{optimal-SRL2}.  Note the symmetry of the poles and zeros with
respect to the $j\omega$ axis.  It is seen here that when $r\approx
\infty $, the closed-loop poles are approximately equal to the open
loop poles of the system.  For small $r$, one closed-loop pole
converges to $\infty$ along the negative real axis, and the other
converges to the zero $z$.

\begin{figure}[ht]
\ebox{.85}{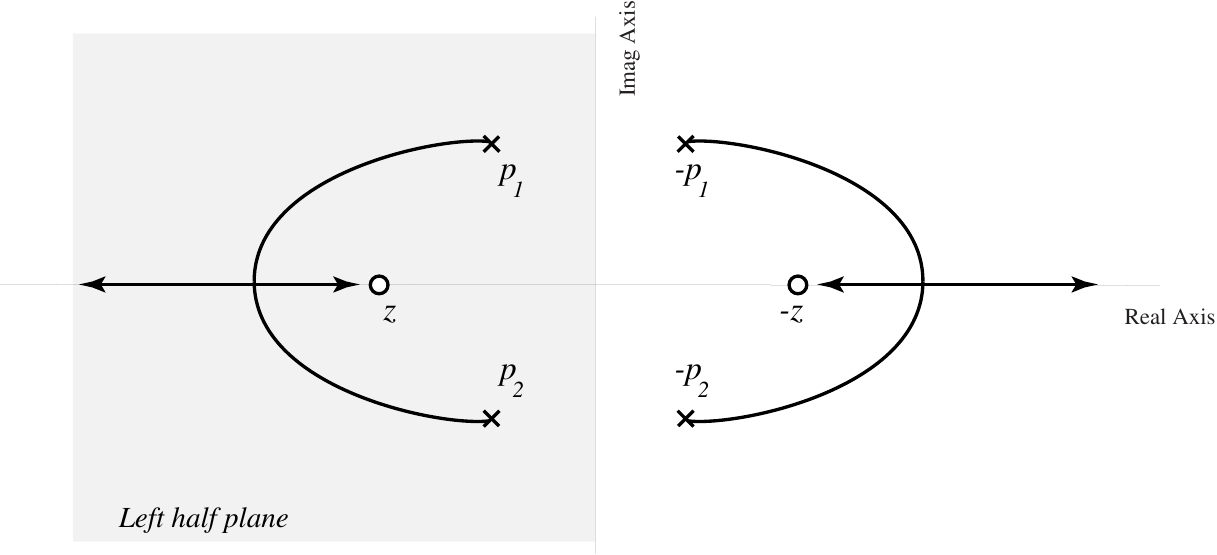}
\caption[The symmetric root locus]{The symmetric root locus for the system with transfer
function $P(s) = (s-z)/(s^2 + 2\zeta s +1)$. The roots of $1 +
\frac{1}{r} P (-s) P (s) = 0$ give the optimal closed-loop poles,
where $r$ is the control weighting in the LQR cost. As $r$ varies from
$0$ to $\infty$, this gives the symmetric root locus.}
\flabel{optimal-SRL1}
\end{figure}

\begin{ex}
We now revisit \eq AREex/ to show how a control design may be
accomplished using the symmetric root locus. The point is that
everything on the left hand side of the return difference equation is
known.  The desired closed-loop feedback poles are unknown, but are
defined through the right hand side of this equation.  In the SISO
case, this can be exploited to solve for the feedback gain $K$
directly, without solving a Riccati equation.

We must first compute the transfer function $P(s)$.  From the state
equations we have
\begin{eqnarray*}
s X_1(s) &=& X_2(s) \\ s X_2(s) &=& - X_2(s) + U(s) \\ Y(s) &=& X_1(s).
\end{eqnarray*}
Combining these equations gives $s^2 Y(s) = -s Y(s) + U(s)$, or
\begin{eqnarray*}
P(s) = \frac{Y(s)}{U(s)} = \frac{1}{s(s+1)}.
\end{eqnarray*}
To compute the symmetric root locus we solve for the roots of the
equation
\begin{eqnarray*}
0=1+r^{-1}\Bigl[ \frac{1}{s(s+1)} \Bigr]
\Bigl[\frac{1}{-s(-s+1)}\Bigr] =1+r^{-1}\Bigl[ \frac{1}{s^4 - s^2}
\Bigr].
\end{eqnarray*}

\begin{figure}[ht]
\ebox{.75}{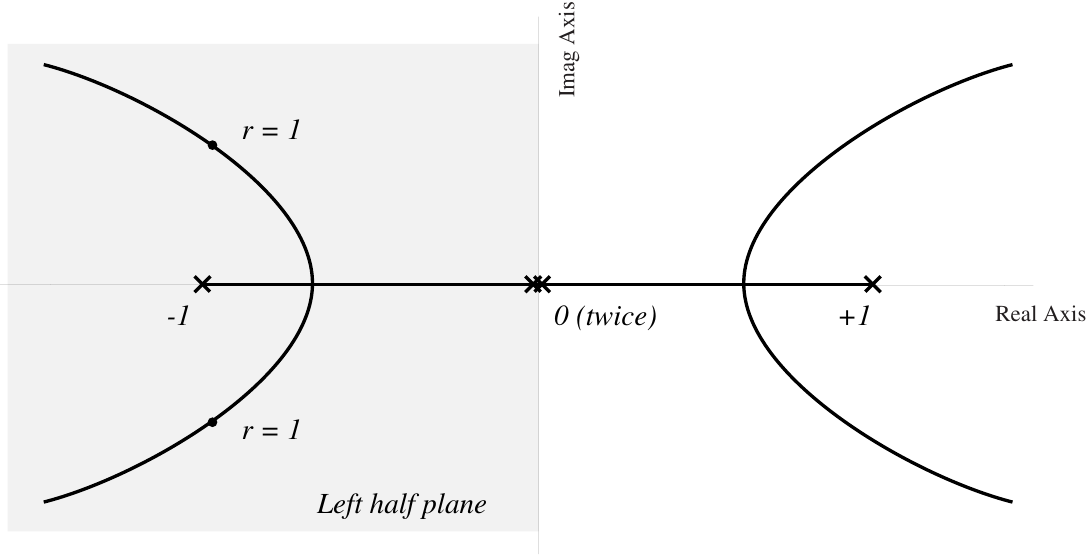}
\caption[The symmetric root locus]{The symmetric root locus for the model $\ddot y =-\doty + u$,
with cost criterion $V=\int y^2 + ru^2\, dt$. A more highly damped
closed-loop response can be obtained by increasing $r$ below the
nominal value of $1$, while decreasing $r$ will result in greater
closed-loop bandwidth. } \flabel{optimal-SRL2}
\end{figure}

A sketch of the symmetric root locus is provided in
\Figure{optimal-SRL2}.  When $r=1$ we have already seen that the
optimal control is $ u^\optchar = - [1,\ \sqrt{3}-1] x$, and that the
closed-loop poles are $ -\half(\sqrt{3}\pm j)$. From
\Figure{optimal-SRL2} it is evident that a more highly damped closed
loop response can be obtained by increasing $r$, while reducing $r$
will result in greater closed-loop bandwidth.
\end{ex}

The return difference equation leads to an elegant approach to pole
placement.  The resulting symmetric root locus gives the designer a
`good' set of poles to choose from, which can then be narrowed down
based upon other constraints, such as closed-loop bandwidth or time
domain specifications.

The most striking implication of the return difference equation is
perhaps its application to sensitivity analysis. In the SISO case, we
have seen that by letting $s = j \omega$ we obtain
\begin{eqnarray*}
| 1 + L (j \omega) |^2 = 1 + \frac{1}{r} | P (j \omega) |^2.
\end{eqnarray*}
This implies \defn{Kalman's inequality},
\begin{eqnarray*}
 | 1 + L (j \omega) |^2 \geq 1, \qquad \omega \in\Re.
\end{eqnarray*}
If we define the sensitivity function as $S= 1/(1+L)$, then we also
have
\begin{eqnarray}
|S(j \omega) | \leq 1, \qquad \omega \in\Re.  \elabel{eq5.1.9}
\end{eqnarray}
It follows from \eq sensitivity-def/ that the controlled system is
insensitive to plant uncertainty at all frequencies.  A typical plot
of the loop transfer function $L$ is shown in
\Figure{optimal-nyquist1}.  We see from this plot that the controlled
system has a gain margin of infinity and a phase margin of $\pm 60$
degrees.  These results hold for \textit{any} full state LQR design.

\begin{figure}[ht]
\ebox{.75}{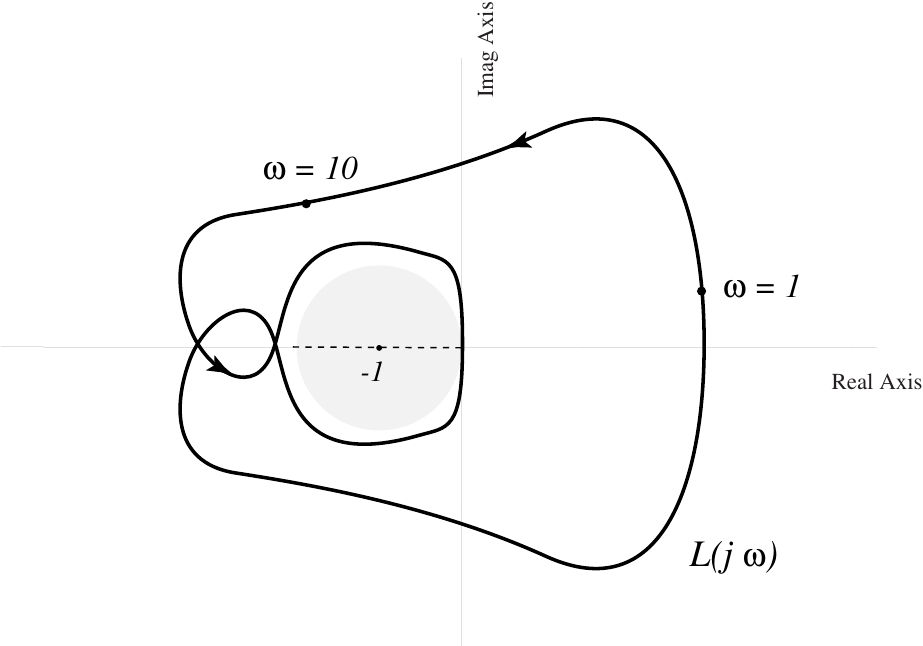}
\caption[The loop transfer function]{The loop transfer function $L(j\omega)=K(Ij\omega - A)^{-1}
B$ for the LTI system under optimal control cannot enter the shaded
disk centered at $-1$.}  \flabel{optimal-nyquist1}
\end{figure}

\begin{ex}
As an illustration of the frequency domain properties of a state space
model under optimal control, we now return to the Pendubot model
linearized in the upright position.  To control the system we define
the cost as
\begin{eqnarray*}
V=\int_0^\infty \left[ |x(t)|^2 + ru(t)^2 \right]\, dt.
\end{eqnarray*}
where $x=(q_1,\dot q_1,q_2,\dot q_2)^T$.  We therefore have $Q = C^T
C$, where $Q=C=I$.  Since the matrix $Q$ has rank $n$, so that $C$ is
square, the transfer function defined as $P(s) = C(Is-A)^{-1} B$ is an
$n\times 1$ vector-valued function of $s$.  The return difference
equation for the Pendubot with this cost criterion becomes
\begin{eqnarray*}
r+\sum_{i=1}^4 P_i(-s)P_i(s) = r[1 + L(-s)][1 + L(s)],
\end{eqnarray*}
so that a generalization of the symmetric root locus is still
possible, since the closed-loop poles are still equal to the roots of
$1+L(s)$.  For any $i$, the poles of $P_i$ are equal to the
eigenvalues of $A$.  These are simply the open-loop poles, which are
roughly located at $\{-9.5, -5.5, 5.5, 9.5\}$. The numerator of $P_i$
depends upon $i$, but the sum on the left hand side is readily
computable using \textit{Matlab}.  The resulting symmetric root locus
is shown in \Figure{optimal-pendubotSRL}.
\begin{figure}[th]
\ebox{.85}{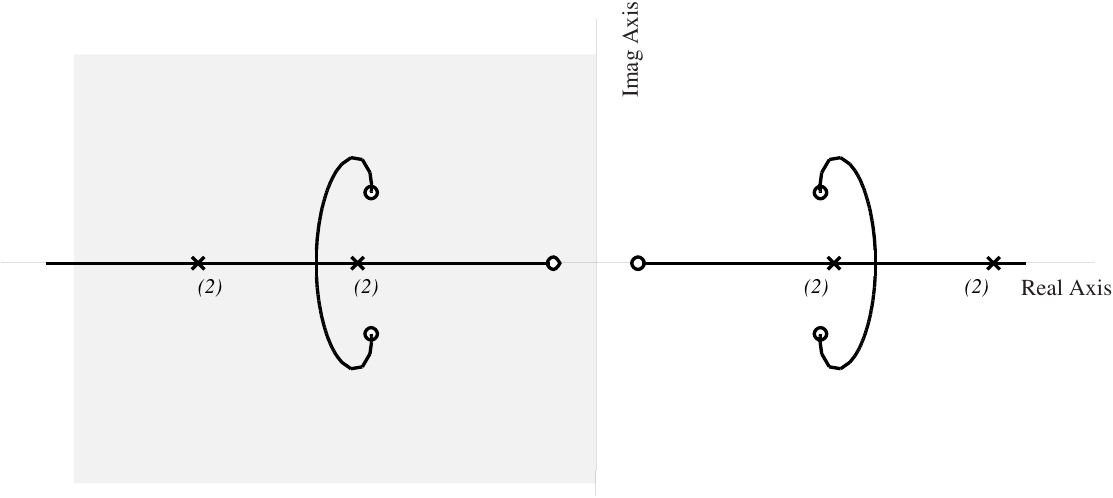} 
\caption[The symmetric root
locus for the Pendubot]{The symmetric root
locus for the Pendubot with $Q=I$.  Each of the poles and zeros in
this diagram are indicated with crosses and zeros, respectively.
Each pole appears twice because of the symmetry of the open-loop
transfer function.  } \flabel{optimal-pendubotSRL}

\ebox{.75}{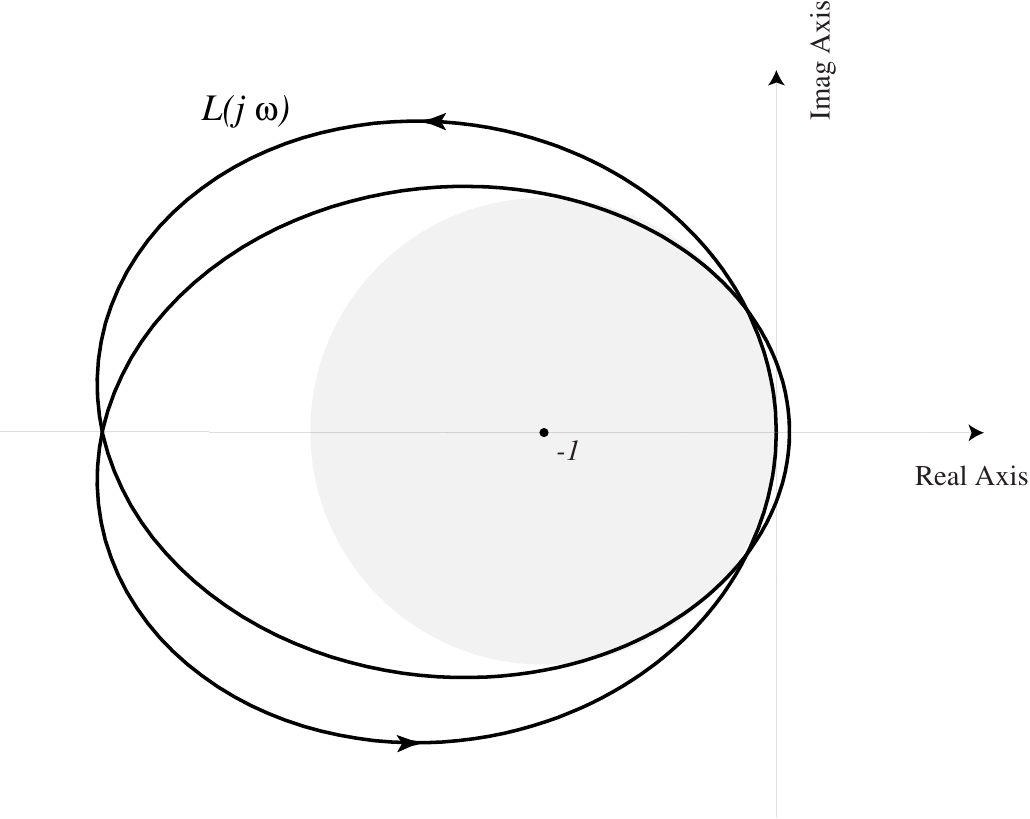}
\caption[Nyquist plot for the Pendubot]{The Nyquist plot for the Pendubot controlled with full state
feedback, based on an LQR design.  If all of the states can be
measured, a far more robust design is possible compared with what is
found using a single output measurement.}
\flabel{optimal-PendFullStateL}
\end{figure}
\end{ex}

For $r=1$ the closed-loop poles are placed at \eq openPendPoles/, as
described in the beginning of \Chapter{goals}.  We saw then that when
this controller was combined with an observer based on measurements at
the lower link alone, the controlled system is highly sensitive to
modeling errors.  For instance, the gain margin is found to be less
than $10^{-2}$. With full state feedback we know that the controlled
system is far more robust.  A Nyquist plot of the loop transfer
function $L(j\omega)$ is shown in \Figure{optimal-PendFullStateL}. We
see in this figure that, as expected, the gain margin with this design
is infinite, and that the gain can also be reduced by one half, and
the system will remain stable.

\begin{summary}
This chapter has provided only a brief survey on the LQR problem.  For
a more detailed discussion, see the book \cite{andmoo90}.  Some
historic papers on the topic are:
\begin{quote}
R. E. Kalman, ``Contributions to the Theory of Optimal Control,'' {\it
Bol. Soc. Matem. Mex.}, (1960), pp. 102--119.

R. E. Kalman, ``When is a Linear Control System Optimal?'', {\it
Trans. ASME Ser. D: J. Basic Eng.}, Vol.~86 (March 1964), pp.~1--10.
\end{quote}
Another good source for survey/tutorial papers on the linear quadratic
theory is the {\it December~1971} special issue of the {\it
IEEE~Transactions~on~Automatic~Control}.
\end{summary}

\begin{matlab}
\item[LQG] Solves the LQR problem using the data $A$, $B$, $Q$, and
$R$.

\item[LQG2] A different numerical algorithm for solving the LQR
problem.

\item[ARE] Gives solution to the algebraic Riccati equation.

\item[RLOCUS] Used together with the convolution command
\textbf{CONV}, this can be used to graph the symmetric root locus.
\end{matlab}

\begin{exercises}

\item
\hwlabel{InfHorHJBexample} The purpose of this exercise is to review
the derivation of the HJB equation, and to provide an introduction to
a general theory of infinite horizon optimal control.  Consider the
time invariant nonlinear system
\begin{eqnarray*}
\dot x = f(x,u), \qquad x(0) = x_0\in\Re^n.
\end{eqnarray*}
We assume that there is an equilibrium at $\zero\in\Re^n$, in the
sense that $f(\zero,\zero)=\zero$.  Assuming that the goal is to force
the state to zero, consider the cost criterion
\begin{eqnarray*}
J(u) = \int_0^\infty \ell(x,u)\, dt,
\end{eqnarray*}
where $\ell\colon\Re^n\times\Re^m\to\posRe$ is smooth, with
$\ell(\zero,\zero)=0$.  We typically also assume that $\ell$ is
positive definite, as defined in \Section{LyapunovDirect}.  We let
$J^\optchar(x)$ denote the minimum of $J$ over all controls, when the
initial condition of the model is $x$.  \balphlist
\item
Derive a ``principle of optimality'' for this problem by considering
an optimal trajectory $(x^\optchar,u^\optchar)$, as done in the
derivation of the HJB equations.

\item Using (a), derive a differential equation that $J^\optchar$
must satisfy.

\item
Apply your solution in (b) to solve the optimal control problem
defined by the simple integrator model with polynomial cost defined
below:
\begin{eqnarray*}
\dot x = u \qquad J(u) = \int_0^\infty [ u^2 + x^4 ] \, dt
\end{eqnarray*}
\end{list}

\item
Solve \textit{by hand} the ARE with
 $$A= \left[ \begin{matrix}0 &
1\\ 0 & 0\end{matrix}\right],\; B= \left[ \begin{matrix}0\\
1\end{matrix}\right], \; Q= \left[ \begin{matrix}1 & 0\\ 0 &
0\end{matrix}\right], \; R=1.$$

\item
For the double integrator $\ddot y = u$, what is the minimal value of
the cost
\begin{eqnarray*}
V(u)=\int_0^\infty y(t)^2 + u(t)^2\, dt,
\end{eqnarray*}
over all controls $u(\cdot)$, when $y(0)= 1$, $\dot y(0) = 0$? Explain
why the minimum exists.

\item
For the scalar time-invariant plant $\dot x=ax(t) +u(t)$, and under
the quadratic performance index
\begin{eqnarray*}
V=\int_0^{\infty}[x^2(t)+ru^2(t)]dt\;,
\end{eqnarray*}
\balphlist
\item Obtain the optimal feedback control as a function of the scalar
parameters $a$ and $r$;

\item Study the eigenvalues of the optimal closed-loop system as
\begin{description}
\item[(i)] $r\to \infty$ (expensive control)
\item[(ii)] $r\to 0$ (cheap control)
\end{description}
\end{list}

\item
The linearized (and normalized) magnetically suspended ball is
described by $\dot x_1=x_2,\quad \dot x_2=x_1+u$, where the initial
state $x(0)$ is specified.  Answer the following three questions by
constructing the Hamiltonian matrix. You may check your work using
Matlab.  \balphlist
\item
Obtain the feedback control $u(x)$ that minimizes the infinite-horizon
quadratic performance index
\begin{eqnarray*}
V=\int_0^{\infty}(x_1^2+2x_1x_2+x_2^2+4u^2)\,dt.
\end{eqnarray*}

\item Determine the eigenvalues of the closed-loop system matrix
$A_{cl}$.

\item Find the minimum value of $V$ as a function of $x(0)$.  Explain
why in this example we can have $V(x_0)=0$ even though $x_0\neq\zero$.
\end{list}

\item
The following is the so-called \textit{generalized regulator problem}:
The system equation is
\begin{eqnarray*}
\dot x(t)=A(t)x(t)+b(t)u(t)+c(t),
\end{eqnarray*}
where $c(\cdot)$ is an $n$-dimensional vector function defined on the
finite interval $[t_0, t_1]$, with continuous entries; and the
performance index is
\begin{eqnarray*}
V=x^T(t_1)Q_fx(t_1)+\int_{t_0}^{t_1}[x^T(t)R_1(t)x(t)+2x^T(t)R_{12}(t)u(t)
+u^T(t)R_2(t)u(t)]\, dt
\end{eqnarray*}
where all matrices have continuous entries, and
\begin{eqnarray*}
R_2(\cdot)>0,\qquad Q_f\geq 0,\qquad \left[
\begin{matrix}R_1(\cdot)&R_{12}(\cdot)\\
R_{12}^T(\cdot)&R_2(\cdot)\end{matrix}\right]\geq 0\;.
\end{eqnarray*}
\balphlist
\item
Using the \textit{Hamilton-Jacobi-Bellman equation}, show that the
optimal control for this problem can be written as
\begin{eqnarray*}
u^\optchar(t)=\gamma^\optchar(x,t)=-R_2^{-1}(t)\bigl[
[B^T(t)P(t)+R_{12}^T(t)]x+B^T(t) k(t)\bigr],
\end{eqnarray*}
where $P(\cdot)$ is a matrix function ($\geq 0$) and $k(\cdot)$ is a
vector function. Obtain the expressions (differential equations)
satisfied by $P$ and $k$.
\item
Obtain an expression for the minimum value of $V$.
\end{list}

\item In the previous problem, assume that all matrices and the vector
$c(\cdot)$ are time-invariant, and $t_1\to\infty$. Obtain an
expression for the steady-state optimal feedback control under
appropriate controllability and observability assumptions. Is this a
stabilizing control?

\item
The Pendubot arranged in the vertical position may be approximately
described by the linear system $\dot x = A x + B u$; $y = x_1$ where
\begin{eqnarray*}
A = \left[ \begin{matrix} 0 & 1.0000 & 0 & 0\\ 42.0933 & 0 & -4.3893 &
0\\ 0 & 0 & 0 & 1.0000\\ -29.4557 & 0 & 43.3215 & 0
\end{matrix}\right] \qquad B = \left[\begin{matrix} 0\\ 6.2427\\ 0\\
-9.7741\end{matrix}\right].
\end{eqnarray*}
There exists a feedback control law $u=-K(r) x$ which minimizes the
quadratic performance index
\begin{eqnarray*}
V(u) = \int_0^\infty \left( y^2 + \alpha {\dot y}^2 + r
u^2\,\right) d\tau
\end{eqnarray*}
Plot the resulting closed-loop poles as $r$ is varied from $0$ to
$\infty$, and discuss your results. How does the introduction of $\dot
y$ in the cost functional affect the root locus?  Note that this will
require a \textit{generalization} of the symmetric root locus derived
in these notes. Include the derivation of this result with your plots.

\item
LQR theory can also be extended to problems where also the derivative
of the control variable appears in the performance index
(\textit{i.e.} we have a soft constraint not only on the control
variable but also on its derivative). Specifically, the LQR problem
with cost
\begin{eqnarray*}
V(u)=\int_0^{\infty}[x^T(t)Qx(t)+u^T(t)R_1u(t)+ \dot u^T(t)R_2\dot
u(t)]\, dt
\end{eqnarray*}
subject to $\dot x=Ax+Bu$ can be solved by introducing the additional
state variable, $z=u$, and treating $v=\dot u$ as the control
variable.  \balphlist
\item
By following this procedure, assuming that $Q\geq 0$, $R_1\geq 0$,
$R_2>0$, and $x(0)=x_0$, $u(0)=u_0$ are specified, obtain the general
solution of the above by simply applying the LQR theory developed in
these notes.

\item
Can you formulate a generalization of the symmetric root locus for
this more general problem?

\item
After obtaining the general solution, apply it to the specific
problem:
\begin{eqnarray*}
\dot x=\left[\begin{matrix}0&1\\ 0&0\end{matrix}\right]x
+\left[\begin{matrix}0\\ 1\\\end{matrix}\right]u ;\quad
x_1(0)=x_2(0)=1, u(0)=1
\end{eqnarray*}
\begin{eqnarray*}
V(u)=\int_0^{\infty}[x_1^2(t)+ \rho \dot u^2(t)]\,dt.
\end{eqnarray*}
Provide a root locus plot for $0<\rho<\infty$, and compute the optimal
control for $\rho=2$.

\textit{Note:} For this numerical example, your solution will be in
the form \quad $v\equiv\dot u=k_1x_1+k_2x_2+k_0u$, where $k_1, k_2,
k_0$ are scalar quantities, whose values will have to be determined as
part of the solution process.
\end{list}

\item This problem deals with optimal control problems in which the
closed-loop system is required to have a prescribed degree of
stability.  Optimal control is an expression of the tradeoff between
the desire to place closed-loop poles deep into the left-half-plane
and the desire not to use excessive control (because of possible
actuator saturation) or have large excursions of the state
trajectories (because of modeling uncertainty).  Sometimes we want to
guarantee that the closed-loop poles are stable by a ``margin'' of at
least $\alpha$.  That is, the poles lie in the region $\Omega_\alpha =
\{s : {\rm Re}(s)\le -\alpha\}$.

To this end, consider the plant \quad $\dot{x} = Ax + Bu$, and the
cost functional
\begin{eqnarray*}
V(u) = \int_0^\infty e^{2\alpha t} \left(x^T Q x + u^T R u\right)
\,dt, \qquad \alpha > 0
\end{eqnarray*}
where $Q \geq 0$ and $R >0$. Observe that if $V(u^\optchar) < \infty$
then $\sqrt{Q} x^\optchar(t)$ must decay faster than $e^{-\alpha t}$
which, under appropriate conditions, means the closed-loop poles lie
in the left-half plane as desired.  In this problem you will
rigorously establish this assertion.

\balphlist
\item
Set $z(t) = e^{\alpha t} x(t)$ and $v(t) = e^{\alpha t} u(t)$, and denote the resulting as a function $v$ by $\tilde{V} (v)$.  Compute $\tilde{V} (v)$ and the corresponding state equations.

\item
Carefully write down a complete solution to this optimal control
problem including necessary and sufficient conditions under which
solutions exist and conditions under which these solutions are
stabilizing.  In the case where the eigenvalues are distinct,
interpret your conditions using the Hautus-Rosenbrock criteria.

\item
Show that the closed-loop system has all its poles inside the left
hand plane as prescribed under appropriate conditions on the open-loop
system.

\item
With $A=B=Q=R=x_0=1$ find the optimal control law and the optimal
cost, and observe that the poles are indeed placed within
$\Omega_\alpha$.
\end{list}

\item
Consider again the flexible structure defined in \eq torsion-ss/.  The
input is the voltage to a DC motor, and the states of this model are
\begin{eqnarray*}
x_1 = \theta_1, x_2 = \dot\theta_1, x_3 = \theta_2, x_4 =
\dot\theta_2.
\end{eqnarray*}
\balphlist
\item
Plot the loop transfer function with your previous design.  You should
use the Matlab command {\tt NYQUIST}, as follows
\begin{eqnarray*}
\hbox{\tt nyquist(A-.00001*eye(4),B,K,0) }
\end{eqnarray*}
where $K$ is the feedback gain.  This will give a plot of $L(j\omega)
= K(Ij\omega - A +\epsilon I)^{-1} B$.  The term $\epsilon $ is needed
since the open-loop system has a pole at the origin.  Does your plot
enter the region dist$(s,-1) < 1$?

\item
Plot the optimal closed-loop poles for the cost function
\begin{eqnarray*}
V(u) =\int_0^\infty \theta_2^2 + \beta (\dot\theta_2)^2 + \rho u^2 \,
d\tau
\end{eqnarray*}
where $\rho$ varies from $0$ to $\infty$.  Try $\beta = 0$, and then
some $\beta >1$.

\item
Obtain an optimal design for which the dominant poles satisfy $-10 \le
\sigma \le -5$, and plot the loop transfer function with this design.
\end{list}

\item
Consider the input-output system described in the frequency domain by
\begin{eqnarray*}
Y(s) = \frac{s-1}{s^2 + s +1} U(s) + W(s)
\end{eqnarray*}
where $y$ is the measurement, $u$ is the input, and $w$ is a
disturbance.

Two control designs will be obtained below by setting $w\equiv 0$.  In
practice however, $w\not= 0$, and hence its affect on the output must
be understood.  \balphlist
\item
Obtain a state space model in observable canonical form: Let
$x_1=y-w$, and $\dot x_2 = - x_1 -u$, so that $y=x_1 + w$.

\item
Derive a state feedback control law $u=-Kx$ for an approximately
critically damped response which minimizes $V(u)=\int y^2 + ru^2\,
d\tau$ for some $r$. Your closed-loop poles should be constrained in
magnitude between $1$ and $5$.

Sketch the symmetric root locus by hand - use Matlab to verify your
plot, and to determine an appropriate feedback gain.

What is the frequency response $\frac{Y}{W}(j\omega)$ for the closed-loop system when the state $x$ is perfectly observed?

\item
Design a reduced order observer, and combine this with the control law
obtained above so that now $u= - K_1 y -K_2 \hat x_2$. Note that the
observer will have to use the approximation $y\approx x_1$ since
neither $x_1$ nor $x_2$ are directly measurable.

Your observer pole should be at least three times as large in
magnitude than the magnitude of the largest closed-loop pole obtained
in (b).

Plot the frequency response $\frac{Y}{W}(j\omega)$ using the {\tt
BODE} command, and discuss your findings.
\end{list}

\item
In this problem you will design a controller for a model of a Boeing
747 aircraft.

A linearized model for the longitudinal motion of this plant when
operating shortly after take-off, at sea level, with a speed of 190
miles per hour is given by
\begin{eqnarray*}
\frac{d}{dt} \left[ \begin{matrix}\delta V \\ \delta \alpha \\ q \\
\delta \theta \\\end{matrix}\right] = \left[ \begin{matrix} -0.0188 &
11.5959 & 0 & -32.2 \\ -0.0007 & -0.5357 & 1 & 0\\ 0.000048 & -0.4944
& -0.4935 & 0\\ 0&0&1&0\\\end{matrix}\right] \left[
\begin{matrix}\delta V \\ \delta \alpha \\ q \\ \delta \theta
\\\end{matrix}\right] + \left[ \begin{matrix}0\\ 0\\ -0.5632\\ 0
\\\end{matrix}\right] \delta e
\end{eqnarray*}
where $\delta V$ is the deviation from the nominal speed of 190 mph,
$\delta \alpha $ is the \textit{angle of attack}, which is the
longitudinal angle of the plane with respect to the longitudinal
velocity vector, $\delta \theta$ is the \textit{pitch angle}, which is
the longitudinal angle of the plane with respect to the ground, and
$q=\dot\delta\theta$.  The input $\delta e$ is the deviation of the
elevators with respect to some nominal value.  Of course, the throttle
position is also a significant control, but we will hold the throttle
value constant in this design.

The altitude of the plane $h$ is the variable that we would like to
control.  Linearizing, we find that a reasonable approximation to the
rate of change of $h$ is given by
\begin{eqnarray*}
\dot\delta h = 290(\delta \theta - \delta\alpha)
\end{eqnarray*}
\balphlist
\item
Adjoin the variable $h$ to the state space model above to obtain a
fifth order system with state $x$.  Obtain an LQR design to obtain a
feedback control $u=-Kx + r$ for this system using symmetric root
locus.  Plot the step response of both the altitude $h$ and the pitch
angle $\theta$ of the controlled system.  Try to obtain a rise time of
250 seconds with reasonable undershoot.

Note that the response is always slow due to one closed-loop pole close to the origin.

\item
Can you obtain a better controller with a more general choice of $Q$
and $R$?  For your best design, where are the resulting closed-loop
poles?

\item
Using your best choice of $K$ in (b), improve the response by
modifying the controller as follows: Pass the reference signal $r$
through a first order filter $G(s)$. Hence, in the frequency domain,
your controller is of the form
\begin{eqnarray*}
U(s) = -K X(s) + k\frac{s-z}{s-p} R(s).
\end{eqnarray*}
You should place the zero $z$ at the location of the slow pole
(approximately), and place the pole $p$ further from the origin.  This
should cancel the effect of the slow closed-loop pole in the closed-loop response. The gain $k$ should be chosen to achieve good steady
state tracking for d.c.\ reference inputs.

Again plot the step response of both the altitude $h$ and the pitch
angle $\theta$ of the controlled system.  Can you choose $k$, $z$ and
$p$ to improve your response over that obtained in (b)?
\end{list}
\end{exercises}

\chapter{An Introduction to the Minimum Principle}
\clabel{minimum} We now return to the general nonlinear state space
model, with general cost criterion of the form
\begin{equation}
\begin{aligned}
\dot x(t) & = f (x(t), u(t), t), \quad x (t_0) = x_0\in\Re^n \\ V(u) &
= \int_{t_0}^{t_1} \ell (x(t), u(t), t) \,dt + m (x ({t_1}))
\end{aligned}
\elabel{min-opt}
\end{equation}
We have already shown under certain conditions that if the input $u$
minimizes $V(u)$, then a partial differential equation known as the
HJB equation must be satisfied.  One consequence of this result is
that the optimal control can be written in state feedback form through
the derivative of the value function $V^\optchar$.  The biggest
drawback to this result is that one must solve a PDE which in general
can be very complex.  Even for the relatively simple LQR problem, the
solution of the HJB equations required significant ingenuity.

The minimum principle is again a method for obtaining necessary
conditions under which a control $u$ is optimal.  This result is based
upon the solution to an ordinary differential equation in $2n$
dimensions.  Because it is an ordinary differential equation, in many
instances it may be solved even though the HJB equation is
intractable.  Unfortunately, this simplicity comes with a price. The
solution $u^\optchar$ to this ordinary differential equation is in
open-loop form, rather than state space form.  The two view points,
 the Minimum Principle and the HJB equations, each have value in
nonlinear optimization, and neither approach can tell the whole story.

\section{Minimum Principle and the HJB equations}
In a course at this level, it is not possible to
give a complete proof of the Minimum Principle.  We can however give
some heuristic arguments to make the result seem plausible, and to
gain some insight.  We initially consider the optimization problem \eq
min-opt/ where $x_0$, $t_0$, ${t_1}$ are fixed, and $x(t_1)$ is
free. Our first approach is through the HJB equation
\begin{eqnarray*}
- \frac{\partial }{\partial t} V^\optchar (x,t) & = & \min_{u} (\ell
(x,u,t) + \frac{\partial V^\optchar}{\partial x} f (x,u,t)) \\
u^\optchar (t) & = & \mbox{arg}\min_{u} (\ell (x,u,t) + \frac{\partial
V^\optchar}{\partial x} f (x,u,t)).
\end{eqnarray*}

For $x,t$ fixed, let $ \bar{u} (x,t)$ denote the value of $u$ which
attains the minimum above, so that
$u^\optchar=\baru(x^\optchar(t),t)$.  In the derivation below we
assume that $\baru$ is a smooth function of $x$. This assumption is
false in many models, which indicates that another derivation of the
Minimum Principle which does not rely on the HJB equation is required
to create a general theory.  With $\baru$ so defined, the HJB equation
becomes
\begin{eqnarray*}
- \frac{\partial }{\partial t} V^\optchar (x,t)= \ell (x, \bar{u}
(x,t),t) + \frac{\partial V^\optchar}{\partial x} (x,t) f (x, \bar{u}
(x,t), t).
\end{eqnarray*}
Taking partial derivatives of both sides with respect to $x$ gives the
term $\frac{\partial V^\optchar}{\partial x}$ on both sides of the
resulting equation:
\begin{eqnarray*}
\begin{aligned}
- \frac{\partial^2 }{\partial x \partial t} V^\optchar (x,t) & =
\frac{\partial\ell }{\partial x} (x, \bar{u} (x,t), t) \\ & +
\frac{\partial^2 V^\optchar }{\partial^2 x} (x,t)f (x, \bar{u} (x,t),
t) \\ & + \frac{\partial V^\optchar}{\partial x}(x,t) \frac{\partial
f}{\partial x} (x, \bar{u} (x,t), t) \\ & +
\underbrace{\frac{\partial}{\partial u} \left(\ell + \frac{\partial
V^\optchar}{\partial x}f \right) (x, \bar{u} (x,t), t)}_{\mbox{
\textit{derivative vanishes at} } \bar{u}} \frac{\partial \bar{u}
(x,t)}{\partial x}.
\end{aligned}
\end{eqnarray*}
As indicated above, since $\bar{u}$ is the unconstrained minimum of
$H=\ell + \frac{\partial V^\optchar}{\partial x} f$, the partial
derivative with respect to $u$ must vanish.

This PDE holds for any state $x$ and time $t$ - they are treated as
independent variables.  Consider now the optimal trajectory
$x^\optchar(t)$ with optimal input
$u^\optchar(t)=\baru(x^\optchar(t),t)$.  By a simple substitution we
obtain
\begin{eqnarray}
0 = \frac{\partial^2 }{\partial x \partial t} V^\optchar (x^\optchar
(t), t) & + & \frac{\partial \ell }{\partial x} (x^\optchar (t),
u^\optchar (t), t) \nonumber \\ & + & \frac{\partial^2 V^\optchar }{
\partial x^2}(x^\optchar (t), t) f (x^\optchar (t), u^\optchar (t),t)
\nonumber \\ & + & \frac{\partial V^\optchar}{\partial x}(x^\optchar
(t), t) \frac{\partial f }{\partial x} (x^\optchar (t), u^\optchar
(t), t).  \elabel{HJBode}
\end{eqnarray}
We can now convert this PDE into an ODE.  The main trick is to define
\begin{eqnarray*}
p (t) \eqdef \frac{\partial V^\optchar}{\partial x} (x^\optchar (t),
t).
\end{eqnarray*}
This function is interpreted as the \textit{sensitivity of the cost
with respect to current state}, and takes values in $\Re^n$.  The
derivative of $p$ with respect to $t$ can be computed as follows:
\begin{eqnarray*}
\frac{d}{dt} p (t) & = & \frac{\partial^2 V^\optchar}{\partial x^2}
(x^\optchar (t), t) \dot{x}^\optchar (t) + \frac{\partial^2
V^\optchar}{\partial x \partial t} (x^\optchar (t), t) \\ & = &
\frac{\partial^2 V^\optchar}{\partial x^2} (x^\optchar (t), t) f
(x^\optchar (t), u^\optchar (t), t) + \frac{\partial^2
V^\optchar}{\partial x \partial t} (x^\optchar (t), t).
\end{eqnarray*}
The two mixed partial terms in this equation are also included in \eq
HJBode/.  Combining these two equations, we thereby eliminate these
terms to obtain
\begin{equation}
\begin{aligned}
0 = \frac{d}{dt} p (t) & + \frac{\partial \ell }{\partial x}
(x^\optchar (t), u^\optchar (t), t) \\ & + \underbrace{\frac{\partial
V^\optchar}{\partial x} (x^\optchar (t), t) }_{p^T (t)} \frac{\partial
f}{\partial x} (x^\optchar (t), u^\optchar (t), t).
\end{aligned}
\elabel{Min0}
\end{equation}
From the form of the Hamiltonian
\begin{eqnarray*}
H(x,p,u,t) = p^T f (x, u, t) + \ell (x, u, t),
\end{eqnarray*}
the differential equation \eq Min0/ may be written
\begin{eqnarray*}
\dot{p} (t) = - \nabla_x H (x^\optchar (t),p(t),u^\optchar (t), t)
\end{eqnarray*}
with the boundary condition
\begin{eqnarray*}
p (t_1) = \nabla_x V^\optchar (x^\optchar (t_1), t_1) = \nabla_x m
 (x^\optchar (t_1), t_1).
\end{eqnarray*}

This is not a proof, since for example we do not know if the value
function $V^\optchar$ will be smooth. However, it does make the
following result seem plausible.  For a proof see \cite{lue69}.
\begin{theorem}[Minimum Principle] \tlabel{min} Suppose that $t_1$ is
fixed, $x (t_1)$ is free, and suppose that $u^\optchar$ is a solution
to the optimal control problem \eq min-opt/.  Then

\balphlist
\item
There exists a \textit{costate vector} $p (t)$ such that for $t_0\le
t\le t_1$,
\begin{eqnarray*}
u^\optchar(t) = \argmin_{u} H (x^\optchar (t), p (t),u, t).
\end{eqnarray*}

\item
The pair $(p,x^\optchar)$ satisfy the 2-point boundary value problem:
\begin{equation}
\begin{aligned}
\dot{x}^\optchar (t) & = \phantom{-} \nabla_p H (x^\optchar
(t),p(t),u^\optchar (t), t) \qquad \Bigl(= f (x^\optchar (t),
u^\optchar (t), t) \Bigr) \\ \dotp (t) & = - \nabla_x H (x^\optchar
(t),p(t),u^\optchar (t), t)
\end{aligned}
\elabel{2ptBVP}
\end{equation}
with the two boundary conditions
\begin{eqnarray*}
x (t_0) = x_0; \qquad p (t_1) = \nabla_x m (x (t_1), t_1).
\end{eqnarray*}
\end{list}
\qed
\end{theorem}

\section{Minimum Principle and Lagrange multipliers${}^\star$}
Optimal control involves a functional minimization which is similar in
form to ordinary optimization in $\Re^m$ as described in a second year
calculus course.  For \textit{unconstrained} optimization problems in
$\Re^m$, the main idea is to look at the derivative of the function
$V$ to be minimized, and find points in $\Re^m$ at which the
derivative is zero, so that $\nabla V(x) = \zero$.  Such an $x$ is
called a \textit{stationary point} of the optimization problem.  By
examining all stationary points of the function to be minimized one
can frequently find among these the optimal solution. The
\textit{calculus of variations} is the infinite dimensional
generalization of ordinary optimization in $\Re^m$.  Conceptually, it
is no more complex than its finite dimensional counterpart.

In more advanced calculus courses, solutions to \textit{constrained}
optimization problems in $\Re^m$ are addressed using \defn{Lagrange
multipliers}.  Suppose for example that one desires to minimize the
function $V(x)$ subject to the constraint $g(x) = \zero$, $x\in\Re^m$,
where $g\colon\Re^m\to\Re^d$.  Consider the new cost function
$\hatV(x,p) = V(x) + p^T g(x)$, where $p\in\Re^d$.  The vector $p$ is
known as a Lagrange multiplier.  The point of extending the state
space in this way is that we can now solve an unconstrained
optimization problem, and very often this gives a solution to the
constrained problem.  Indeed, if $x^0,p^0$ is a stationary point for
$\hatV$ then we must have
\begin{eqnarray}
\zero &=& \nabla_x \hatV \,(x^0,p^0) = \nabla V(x^0) + \nabla
g(x^0){p^0} \elabel{LM1} \\ \zero &=& \nabla_p \hatV \,(x^0,p^0) =
g(x^0). \elabel{LM2}
\end{eqnarray}
Equation \eq LM1/ implies that the gradient of $V$ is a linear
combination of the gradients of the components of $g$ at $x_0$, which
is a necessary condition for optimality under very general
conditions. \Figure{optimal-multipliers} illustrates this with $m=2$
and $d=1$.  Equation \eq LM2/ is simply a restatement of the
constraint $g =\zero$.  This Lagrange multiplier approach can also be
generalized to infinite dimensional problems of the form \eq min-opt/,
and this then gives a direct derivation of the Minimum Principle which
does not rely on the HJB equations.

\begin{figure}[ht]
\ebox{.85}{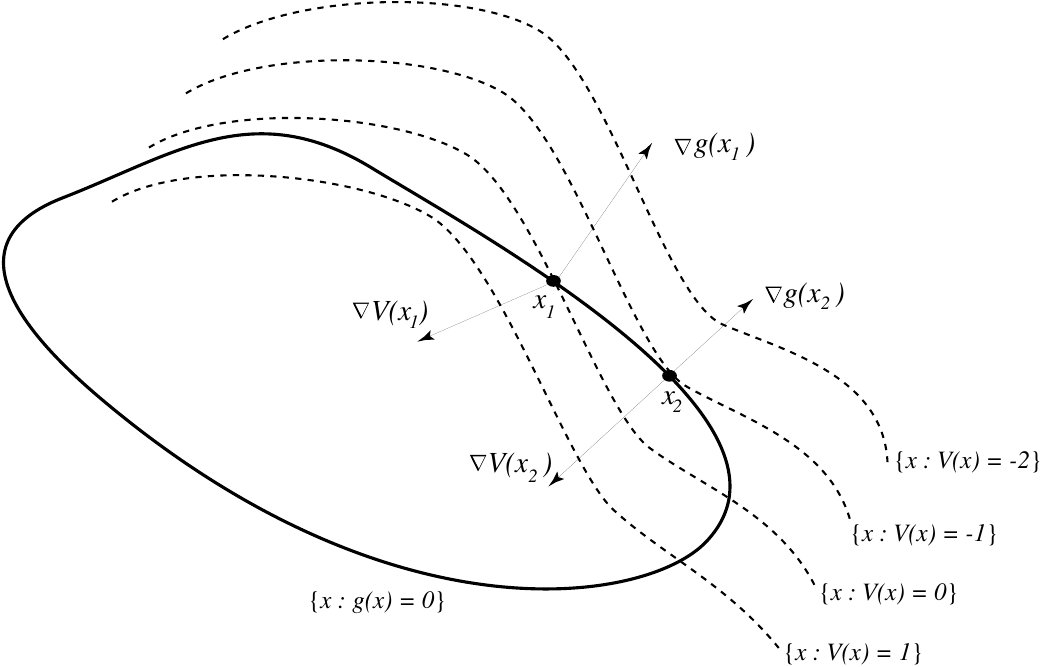}
\caption[An optimization problem on $\Re^2$ with a single constraint]{An optimization problem on $\Re^2$ with a single constraint
$g=0$. The point $x_1$ is not optimal: By moving to the right along
the constraint set, the function $V(x)$ will decrease.  The point
$x_2$ is a minimum of $V$ subject to the constraint $g(x) = 0$ since
at this point we have $V(x_2) = -1$.  $V$ can get no lower - for
example, $V(x)$ is never equal to $-2$ since the level set $\{x : V(x)
= -2\}$ does not intersect the constraint set $\{x : g(x) = 0\}$.  At
$x_2$, the gradient of $V$ and the gradient of $g$ are parallel.}
\flabel{optimal-multipliers}
\end{figure}

To generalize the Lagrange multiplier approach we must first
generalize the concept of a stationary point.  Suppose that $F$ is a
functional on $D^r[t_0,t_1]$.  That is, for any function $z\in
D^r[t_0,t_1]$, $F(z)$ is a real number.  For any $\eta \in
D^r[t_0,t_1]$ we can define a directional derivative as follows:
\begin{eqnarray*}
D_\eta F\, (z) = \lim_{\epsilon \to 0} \frac{F(z+\epsilon
\eta)-F(z)}{\epsilon },
\end{eqnarray*}
whenever the limit exists.  The function $z(t,\epsilon) =
z(t)+\epsilon \eta(t)$ may be viewed as a perturbation of $z$, as
shown in \Figure{optimal-variation}

\begin{figure}[ht]
\ebox{.55}{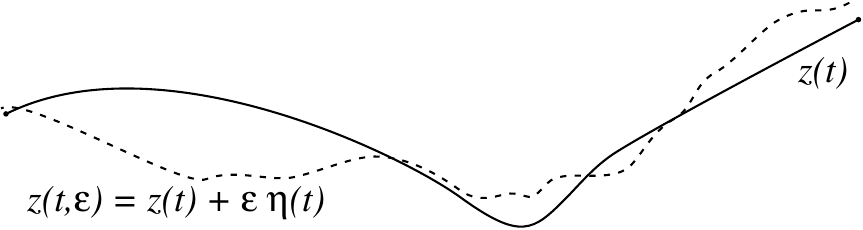}
\caption{A perturbation of the function $z\in D[t_0,t_1]$.}
\flabel{optimal-variation}
\end{figure}

We call $z_0$ a \defn{stationary point} of $F$ if $D_\eta F\, (z_0) =
0$ for \textit{any} $\eta\in D^r[t_0,t_1]$.  If $z_0$ is a minimum of
$F$, then a perturbation cannot decrease the cost, and hence we must
have for any $\eta$,
\begin{eqnarray*}
F(z_0+\epsilon \eta)\ge F(z_0),\qquad \epsilon \in\Re.
\end{eqnarray*}
From the definition of the derivative it then follows that an optimal
solution must be a stationary point, just as in the finite dimensional
case!

The problem at hand can be cast as a constrained functional
minimization problem,
\begin{eqnarray*}
\hbox{Minimize} && V(x,u) = \int \ell\, d\tau + m \\ \hbox{Subject to}
 && \dotx -f = \zero, \qquad x \in D^n[t_0,t_1],\ u \in D^m[t_0,t_1].
\end{eqnarray*}
To obtain necessary conditions for optimality for constrained problems
of this form we can extend the Lagrange multiplier method.  We outline
this approach in five steps:

\head{Step 1.}  Append the state equations to obtain the new cost
functional
\begin{eqnarray*}
\hatV(x,u) = \int_{t_0}^{t_1} \ell \,dt + m (x({t_1})) +
\int_{t_0}^{t_1} p^T (f - \dot x) \,dt.
\end{eqnarray*}
The \textit{Lagrange multiplier vector} $p$ lies in
$D^n[t_0,t_1]$. The purpose of this is to gain an expression for the
cost in which $x$ and $u$ can be varied independently.

\head{Step 2.}  Use integration by parts to eliminate the derivative
of $x$ in the definition of $\hatV$:
\begin{eqnarray*}
\int_{t_0}^{t_1} p^T \dot x \,dt & = & \int_{t_0}^{t_1} p^T dx =
\left.\left(p^T x \right)\right|_{t_0}^{t_1} - \int_{t_0}^{t_1}
\dotp^T x \, dt.
\end{eqnarray*}

\head{Step 3.}  Recall the definition of the Hamiltonian,
\begin{eqnarray*}
H(x,p,u,t) \eqdef \ell(x, u, t) + p^T f(x, u, t).
\end{eqnarray*}
Combining this with the formulas given in the previous steps gives
\begin{eqnarray}
\hatV(x,u) &=& \int_{t_0}^{t_1} H(x,p,u,t) \,dt + \int_{t_0}^{t_1}
\dotp^T x \,dt \nonumber \\ && + p^T (t_0) x (t_0) - p^T ({t_1}) x
({t_1}) + m (x({t_1})). \elabel{eq5.1.1}
\end{eqnarray}

\head{Step 4.}  Suppose that $u^\optchar$ is an optimal control, and
that $x^\optchar$ is the corresponding optimal state trajectory.  The
Lagrange multiplier theorem asserts that the pair
$(x^\optchar,u^\optchar)$ is a stationary point of $\hatV$ for some
$p^\optchar$.  Perform perturbations of the optimal control and state
trajectories to form
\begin{equation}
u (t, \delta) = u^\optchar (t) + \delta \psi(t) \qquad x (t, \epsilon)
= x^\optchar (t) + \epsilon \eta (t). \elabel{eq5.1.2}
\end{equation}
Since we are insisting that $x(t_0)=x_0$, we may assume that
$\eta(t_0)=\zero$.  Consider first variations in $\epsilon$, with
$\delta$ set equal to zero. Letting $\hatV (\epsilon) =
\hatV(x^\optchar + \epsilon\eta,u^\optchar) $, we must have
\begin{eqnarray*}
\frac{d}{d\epsilon}\hatV(\epsilon) = 0.
\end{eqnarray*}
Using \eq eq5.1.1/ to compute the derivative gives
\begin{eqnarray}
0 = \int_{t_0}^{t_1} \frac{\partial}{\partial x} H (x^\optchar,p,
u^\optchar, t) \eta(t) \,dt + \int_{t_0}^{t_1} \dotp^T \eta \,dt
\nonumber \\ - p^T ({t_1}) \eta ({t_1}) + p^T(t_0) \eta(t_0) +
\frac{\partial}{\partial x} m(x^\optchar({t_1}))\eta({t_1}).
\elabel{eq5.1.3}
\end{eqnarray}

Similarly, by considering perturbations in $u^\optchar$ we obtain for
any $\psi\in D^m[t_0,t_1]$,
\begin{eqnarray}
\zero = \int_{t_0}^{t_1} \frac{\partial}{\partial u} H (x^\optchar,p,
u^\optchar, t) \psi(t) \,dt. \elabel{eq5.1.3b}
\end{eqnarray}
This simpler expression is obtained because only the first term in \eq
eq5.1.1/ depends upon $u$.

\head{Step 5.}  We can choose $\eta(t)$ \textit{freely} in \eq
eq5.1.3/. From this it follows that
\begin{eqnarray*}
\frac{\partial H }{\partial x} + \dotp^T = \zero^T \Rightarrow \dotp =
- \nabla_x H
\end{eqnarray*}
and since $\eta(t_0)=\zero$, and $\eta(t_1)$ is free,
\begin{eqnarray*}
-p^T ({t_1}) + \frac{\partial m }{\partial x } ({t_1}) = \zero^T
\Rightarrow p ({t_1}) = \nabla_x m (x({t_1})).
\end{eqnarray*}

Similarly, by \eq eq5.1.3b/ we have
\begin{eqnarray*}
\frac{\partial H }{\partial u} = \zero^T \Rightarrow \nabla_u H =
\zero.
\end{eqnarray*}
In fact, if $u$ is to be a minimum of $\hat V$, then in fact it must
minimize $H$ pointwise.  These final equations then give the Minimum
Principle \Theorem{min}.

From this proof it is clear that many generalizations of the Minimum
Principle are possible.  Suppose for instance that the final state
$x(t_1)=x_1$ is specified.  Then the perturbation $\eta$ will satisfy
$\eta(t_1) = \zero$, and hence using \eq eq5.1.3/, it is impossible to
find a boundary condition for $p$.  None is needed in this case, since
to solve the $2n$-dimensional coupled state and costate equations, it
is enough to know the initial and final conditions of $x^\optchar$.

\section{The penalty approach${}^\star$}
A third heuristic approach to the Minimum Principle involves relaxing
the hard constraint $\dotx - f = \zero$, but instead impose a large,
yet ``soft'' constraint by defining the cost
\begin{eqnarray*}
\hatV (x,u)= \int_{t_0}^{t_1} \ell(x(t),u(t),t) \,dt
+\frac{k}{2}\int_{t_0}^{t_1} |\dot x(t) - f(x(t),u(t),t)|^2 \,dt + m
(x({t_1})).
\end{eqnarray*}
The constant $k$ in this equation is assumed large, so that $\dot x(t)
- f(x(t),u(t),t)\approx \zero$.

We assume that a pair $(x_k,u_k)$ exists which minimizes $\hatV_k$.
Letting $(x^\optchar,u^\optchar)$ denote a solution to the original
optimization problem, we have by optimality,
\begin{eqnarray*}
\hatV_k (x_k,u_k) \le \hatV_k (x^\optchar,u^\optchar) = V^\optchar.
\end{eqnarray*}
Assuming $\ell$ and $m$ are positive, this gives the following uniform
bound
\begin{eqnarray*}
\int_{t_0}^{t_1} |\dot x(t) - f(x(t),u(t),t)|^2 \,dt \le \frac{2}{k}
V^\optchar.
\end{eqnarray*}
Hence, for large $k$, the pair $(x_k,u_k)$ will indeed approximately
satisfy the differential equation $\dot x = f$.

If we perturb $x_k$ to form $x_k + \epsilon \eta$ and define
$\hatV(\epsilon) = \hatV(x_k + \epsilon \eta, u_k)$ then we must have
$d/d\epsilon\, \hatV(\epsilon) =0$ when $\epsilon =0$.  Using the
definition of $\hat V$ gives
\begin{eqnarray*}
\hatV(\epsilon) &=& \int_{t_0}^{t_1} \ell(x_k(t) + \epsilon
\eta(t),u_k(t),t) \,dt  + m (x_k({t_1} + \epsilon \eta(t_1)) \\
&& \qquad + \frac{k}{2}\int_{t_0}^{t_1} |\dot x_k(t) +
\epsilon \dot\eta(t) - f(x_k(t) + \epsilon \eta(t),u_k(t),t)|^2 \,dt.
\end{eqnarray*}
The derivative of this expression with respect to $\epsilon$ can be
computed as follows:
\begin{eqnarray*}
\begin{aligned}
0=\frac{d}{d\epsilon} \hatV\, (0) &= \int_{t_0}^{t_1} \frac{d}{dx}
\ell(x(t),u(t),t) \eta(t) \,dt \\ \quad & + k \int_{t_0}^{t_1} (\dot
x(t) - f(x(t), u(t),t))^T [ \dot\eta -\frac{\partial} {\partial x}f(
x(t), u(t),t) \eta(t) ] \,dt \\ & + \frac{\partial}{\partial x} m (
x({t_1})) \eta(t_1).
\end{aligned}
\end{eqnarray*}
To eliminate the derivative term $\dot\eta$ we integrate by parts, and
this gives the expression
\begin{equation}
\begin{aligned}
0 &= \int_{t_0}^{t_1} \Bigl\{\frac{d}{dx} \ell(x_k(t),u_k(t),t) \,dt +
p_k(t)^T \frac{\partial}{\partial x}f(x_k(t), u_k(t),t) + \frac{d}{dt}
\Bigl(p_k(t)^T\Bigr) \Bigr\} \eta(t) \, dt \\ & -p_k(t_1) +
\frac{\partial}{\partial x} m (x_k({t_1}) \eta(t_1),
\end{aligned}
\elabel{penaltyHam}
\end{equation}
where we have set $p_k(t) = - k (\dot x_k(t) - f(x_k(t), u_k(t),t))$.
Since $\eta$ is \textit{arbitrary}, we see that
\begin{eqnarray*}
\zero = \frac{d}{dt} p_k (t)^T + \frac{\partial \ell }{\partial x}
(x_k(t), u_k (t), t) + p_k(t)^T \frac{\partial }{\partial x} f(x_k
(t), u_k (t), t),
\end{eqnarray*}
and we again obtain the boundary condition
\begin{eqnarray*}
p_k(t_1) = \frac{\partial}{\partial x} m (x({t_1}).
\end{eqnarray*}
Considering perturbations in $u$ gives the equation $\nabla_u H\, (x_k
(t), p_k(t),u_k (t), t) = 0$, which is a weak form of the Minimum
Principle for the perturbed problem.

By letting $k\to \infty$, this penalty approach can be used to prove
both the complete Minimum Principle and the Lagrange Multiplier
Theorem.

\notes{leaving the half out of the cost really causes problems here!}

\section{Application to LQR}
The LQR problem is a good test case to see if the Minimum Principle is
a useful tool for the construction of optimal policies.  Consider
again the general LTI model with quadratic cost:
\begin{eqnarray*}
\dot x & = & Ax + Bu \\ V & = & \int_{t_0}^{t_1} (x^T Qx + u^T Ru)
\,dt + x^T (t_1) M x (t_1).
\end{eqnarray*}
The Hamiltonian is written
\begin{eqnarray*}
H = x^T Q x + u^T Ru + p^T (Ax + Bu)
\end{eqnarray*}
and hence the control can be computed through
\begin{eqnarray*}
\nabla_u H = 0 = 2 Ru + B^T p \Longrightarrow u = - \half R^{-1} B^T
p.
\end{eqnarray*}
This then gives the first set of differential equations:
\begin{eqnarray}
\dot x = Ax + Bu = Ax - \half BR^{-1} B^T p. \elabel{eq5.1.4}
\end{eqnarray}
Through the expression $\nabla_x H = 2Qx + A^T p$, we find that the
derivative of $p$ is
\begin{eqnarray}
\dotp = - \nabla_x H = -2Qx - A^T p. \elabel{eq5.1.5}
\end{eqnarray}
Equations \eq eq5.1.4/ and \eq eq5.1.5/ form the coupled set of
differential equations
\begin{eqnarray*}
\begin{bmatrix} \dot x \\ \dotp \end{bmatrix}
 & = & \left[ \begin{array}{cc} A & - \half BR^{-1} B^T \\ -2Q & -A^T
\end{array} \right] \begin{bmatrix} x \\ p \end{bmatrix}
\end{eqnarray*}
with boundary conditions
\begin{eqnarray*}
x (t_0) & = & x_0 \\ p ({t_1}) & = & 2M x ({t_1}).
\end{eqnarray*}
If we scale $p$ we again arrive at the Hamiltonian defined in
\Section{HamMatrix}: With $\lambda = \half p$,
\begin{equation}
\begin{bmatrix} \dot x \\ \dotlambda \end{bmatrix}
= \left[ \begin{array}{cc} A & -BR^{-1} B^T \\ -Q & -A^T
\end{array} \right]   \begin{bmatrix} x \\ \lambda
\end{bmatrix} \underH \left[\begin{matrix} \dot x \\ \dot \lambda
\end{matrix}\right], \elabel{Ham}
\end{equation}
with $ u^\optchar(t) = - R^{-1} B^T \lambda (t) $.

This ODE can be solved using the ``sweep method''.  Suppose that
$\lambda (t) = P (t) x(t)$.  Then
\begin{eqnarray*}
\dotlambda = \dotP x + P\dot x = \dotP x + P (Ax + Bu),
\end{eqnarray*}
and substituting $ u^\optchar = -R^{-1} B^T \lambda = - R^{-1} B^T Px$
gives
\begin{eqnarray*}
\dotlambda = \dotP x + P (Ax - BR^{-1} B^T Px).
\end{eqnarray*}
From \eq Ham/ we also have
\begin{eqnarray*}
\dot \lambda = -Qx -A^T \lambda = - Qx -A^T Px.
\end{eqnarray*}
Equating the two expressions for $\dot \lambda$ gives
\begin{eqnarray*}
- Qx - A^T Px = \dot P x + PAx - PBR^{-1} B^T x
\end{eqnarray*}
which yields the RDE given in \eq eq4.2.5/.  The boundary condition
for $\lambda$ is $\lambda =\half p = \half m'\, (x_1)$, which from the
definitions of $p$ and $m$ implies that
\begin{eqnarray*}
P (t_1) = M.
\end{eqnarray*}
Solving for $P(t)$ gives $\lambda (t)$, which in turn gives $p (t)$.

\section{Nonlinear examples}
We now solve some nonlinear problems where it is possible to obtain an
explicit solution to the coupled state-costate equations given in the
Minimum Principle.  At the same time, we also give several extensions
of this result.

Our first extension involves constraints on the input. Suppose that
$\clU$ is a subset of $\Re^m$, and that we require that $u(t)\in\clU$
for all $t$.  Since this is a hard constraint, \Theorem{min} is not
directly applicable, but we have the following simple extension:

\begin{theorem}[Minimum Principle with constraints]
\tlabel{minConstraints} Suppose that $x (t_1)$ is free, $t_1$ is
fixed, and suppose that $u^\optchar$ is a solution to the optimal
control problem \eq min-opt/ under the constraint $u\in\clU$.  That
is,
\begin{eqnarray*}
V^\optchar = \min \{V(u) : u(t)\in\clU \ \hbox{for all $t$} \}.
\end{eqnarray*}
We then have \balphlist
\item
There exists a costate vector $p (t)$ such that
\begin{eqnarray*}
u^\optchar(t)= \argmin_{u\in\clU} H (x^\optchar (t), p (t), u,t).
\end{eqnarray*}

\item
The pair $(p,x^\optchar)$ satisfy the 2-point boundary value problem
\eq 2ptBVP/, with the two boundary conditions
\begin{eqnarray*}
x (t_0) = x_0; \qquad p (t_1) = \frac{\partial}{\partial x} m (x
(t_1), t).
\end{eqnarray*}
\end{list}
\qed
\end{theorem}

\begin{ex}
To illustrate \Theorem{minConstraints} consider the control of a
bilinear system, defined by the differential equation
\begin{eqnarray*}
\dot{x} = ux, \qquad x (0) = x_0 > 0, \ 0 \leq u(t) \leq 1.
\end{eqnarray*}
Suppose that the goal is to make $x$ large, while keeping the
derivative of $x$ small on average. Then a reasonable cost criterion
is
\begin{eqnarray*}
V (u)= \int_0^{t_1} \dotx (\tau) - x (\tau) \, d\tau = \int_0^{t_1} [
u (\tau) - 1] x (\tau) \, d\tau.
\end{eqnarray*}
We assume that $t_1 \geq 1$ is fixed, and that $x (t_1)$ is free.

The Hamiltonian becomes
\begin{eqnarray*}
H(x,p,u,t) & = & p f + \ell = p (ux) + [u-1] x \\ & = & x \{u (p+1) -
1\}
\end{eqnarray*}
and by the Minimum Principle, the optimal control takes the form
\begin{eqnarray*}
u^\optchar (t)=\argmin_{u\in\clU} H (x^\optchar (t), p (t), u,t).
\end{eqnarray*}
Since $x^\optchar (t) > 0$ for all $t$, the minimization leads to
\begin{eqnarray*}
u^\optchar (t) & = &
\begin{cases}
1 & \hbox{if } p (t) + 1 < 0 \\ 0 & \hbox{if } p (t) + 1 > 0 \\
\hbox{unknown} & \hbox{if } p (t) = - 1
\end{cases}
\end{eqnarray*}

The costate variable and state variable satisfy the differential
equations
\begin{eqnarray*}
\dot{p} & = & - \frac{\partial H }{\partial x} = - (p+1) u^\optchar +
1 \\ \dot{x^\optchar} & = & \frac{\partial H}{\partial p} = u^\optchar
x^\optchar
\end{eqnarray*}
with the boundary conditions
\begin{eqnarray*}
x^\optchar (0) = x_0,\qquad p (t_1) = \frac{\partial m }{\partial x} =
0.
\end{eqnarray*}

How do we solve these coupled differential equations to compute
$u^\optchar$?  First note that if $t = t_1$ then $p (t) + 1 = 1 > 0$
so $u^\optchar (t) = 0$.  By continuity, if $t \approx t_1$, then $p
(t) > 0$, so we still must have $u^\optchar (t)= 0$.  Thus $\dot{p}
(t) = + 1$ for $t \approx t_1$, which implies
\begin{eqnarray*}
p (t) = t - t_1, \qquad t \approx t_1.
\end{eqnarray*}
As we move backwards in time from time $t_1$, the costate variable
will continue to have this form up until the first time that $p (t) +
1$ changes sign.  This will occur when $0 = t - t_1 + 1$, or $ t = t_1
- 1$.

For $t < t_1-1$ we then have $p (t) + 1 < 0$ so that $u^\optchar (t) =
1$.  For these time values, the costate vector is described by the
equation
\begin{eqnarray*}
\dot{p} = - (p+1) + 1 = - p.
\end{eqnarray*}
Since $p (t_1 -1) = -1$, we obtain
\begin{eqnarray*}
p (t) = -e^{-t + t_1 -1}, t < t_1 -1.
\end{eqnarray*}
The costate trajectory is sketched in \Figure{optimal-bilinear}.

\begin{figure}[ht]
\ebox{1}{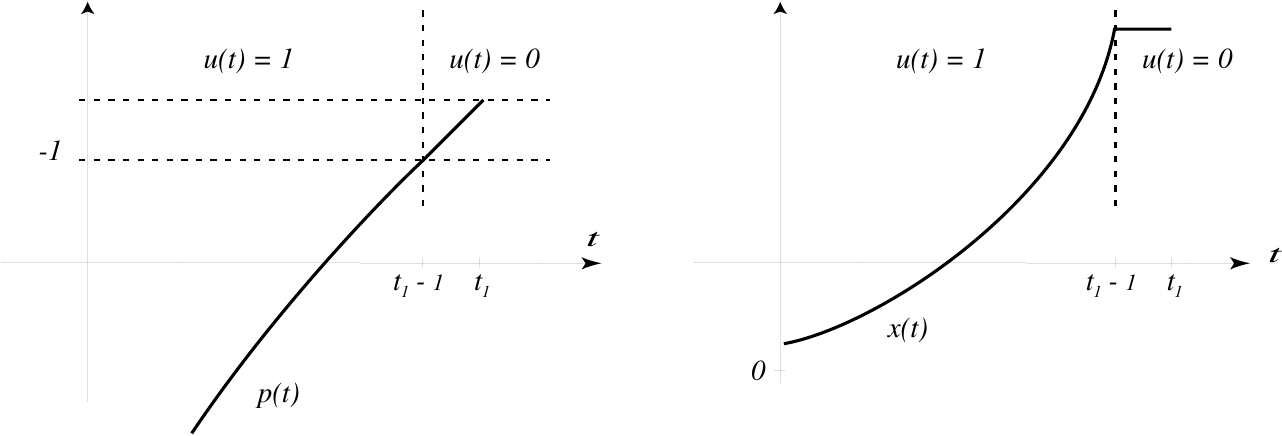}
\caption[Optimal state
trajectory for the bilinear model]{The costate trajectory, optimal control, and optimal state
trajectory for the bilinear model. } \flabel{optimal-bilinear}
\end{figure}

We conclude that the optimal control is \defn{bang-bang}.  That is, it
takes on only two values, which are the upper and lower limits:
\begin{eqnarray*}
u (t) =\left\{\begin{array}{ll} 1 & t < t_1 -1 \\ 0 & t > t_1
-1. \end{array} \right.
\end{eqnarray*}
The time $t_1-1$ is called the \textit{switching point}.  The optimal
controlled state trajectory is defined by the differential equation
\begin{eqnarray*}
\dot{x}^\optchar(t) = \begin{cases} x^\optchar(t) & t< t_1 - 1 \\ 0 &
t > t_1 - 1 \end{cases}
\end{eqnarray*}
so the resulting state trajectory is expressed
\begin{eqnarray*}
x^\optchar (t) = \left\{
\begin{array}{ll}
x_0 e^t & t < t_1 -1 \\ x_0 e^{t_1 -1} & t > t_1 -1.
\end{array} \right.
\end{eqnarray*}
This is also illustrated in \Figure{optimal-bilinear}.
\end{ex}

The next example we treat has a hard constraint on the \textit{final
state value}, so that $t_0$, $t_1$, $x (t_0) = x_0$, and $x (t_1) =
x_1$ are all fixed and prespecified.  In this situation the Minimum
Principle is the same, but the boundary conditions for differential
equations must be modified.

\begin{theorem}[Minimum Principle with final value constraints]
\tlabel{minEnd} Suppose that $t_0$, $t_1$, $x (t_0) = x_0$, and $x
(t_1) = x_1$ are prespecified, and suppose that $u^\optchar$ is a
solution to the optimal control problem \eq min-opt/, subject to these
constraints.  Then \balphlist
\item
There exists a costate vector $p (t)$ such that
\begin{eqnarray*}
u^\optchar (t) = \argmin_{u} H (x^\optchar (t), p (t), u,t).
\end{eqnarray*}

\item
The pair $(p,x^\optchar)$ satisfy the 2-point boundary value problem
\eq 2ptBVP/, with the two boundary conditions
\begin{eqnarray*}
x (t_0) = x_0; \qquad x(t_1)=x_1.
\end{eqnarray*}
\end{list}
\qed
\end{theorem}

\begin{ex}
To illustrate the application of \Theorem{minEnd} we consider the LQR
problem where the terminal state has been specified. Note that in
general we cannot make the state {\it stay} at the terminal state
$x_1$.  The cost criterion is
\begin{eqnarray*}
V & = & \half \int_0^{t_1} (x^T Q x + u^T R u) dt \qquad R > 0, Q \geq
0 \\ H & = & p (Ax + Bu) + \half (x^T Q x + u^T Ru).
\end{eqnarray*}
The Minimum Principle again implies that the optimal control has the
form
\begin{eqnarray*}
u^\optchar (t)= -R^{-1} B^T p (t),
\end{eqnarray*}
and that the costate vector satisfies the differential equation
\begin{eqnarray*}
\dot{p} (t) = - \frac{\partial H }{\partial x} = - p^T A - x^T Q
\end{eqnarray*}

Thus we again arrive at the linear differential equation based on the
Hamiltonian:
\begin{eqnarray*}
\begin{pmatrix} \dot{x} \\ \dot{p}  \end{pmatrix}
= \left(\begin{array}{ll} A & - BR^{-1} B^T \\ -Q & -A^T \end{array}
\right) \begin{pmatrix} x \\ p \end{pmatrix}; \qquad
\begin{array}{rcl} x (t_0) & = & x_0 \\ x (t_1) & = & x_1
\end{array}.
\end{eqnarray*}
Let $\psi (t, \tau)$ denote the state transition matrix, so that
\begin{eqnarray*}
\frac{d}{dt} \psi (t, \tau) = \underH(t) \psi (t, \tau); \qquad\psi
(t,t) = I.
\end{eqnarray*}
This is a $2n\times 2n$ matrix which we decompose as
\begin{eqnarray*}
\psi (t, \tau) = \left[ \begin{array}{ll} \psi_{11} (t, \tau) &
\psi_{12} (t, \tau) \\ \psi_{21} (t, \tau) & \psi_{22} (t, \tau)
\end{array} \right].
\end{eqnarray*}

To solve for the optimal control, we will compute the initial
condition $p (t_0) = p_0$.  This will allow us to compute $p(t)$ for
all $t$, and hence also the optimal control expressed in terms of $p$
above as $u^\optchar= -R^{-1} B^T p $.  The unknown term $p_0$ can be
expressed in terms of the other known quantities as follows:
\begin{eqnarray*}
x_1 = x (t_1) = \psi_{11} (t_1, t_0) x_0 + \psi_{12} (t_1, t_0) p_0
\end{eqnarray*}
\textit{Assuming} that $\psi_{12} (t_1, t_0)$ is invertible, this
gives the formula
\begin{eqnarray*}
p_0 = \psi_{12} (t_1, t_0)^{-1} \{x_1 - \psi_{11} (t_1, t_0) x_0 \}.
\end{eqnarray*}
For all $t$ we then have
\begin{eqnarray*}
p (t) = \psi_{21} (t, t_0) x_0 + \psi_{22} (t, t_0) p_0.
\end{eqnarray*}

The optimal control is
\begin{eqnarray*}
u^\optchar (t) = - R^{-1} B^T p (t).
\end{eqnarray*}
and the optimal state trajectory becomes
\begin{eqnarray*}
x^\optchar (t) = \psi_{11} (t, t_0) x_0 + \psi_{12} (t, t_0) p_0,
\qquad t \geq t_0.
\end{eqnarray*}
\end{ex}

Our last extension concerns the case where the terminal time is not
fixed.  At the same time, suppose that \textit{some} components of $x
(t_1)$ are specified.  We assume that $t_0$ and $x (t_0) = x_0$ are
fixed.

\begin{theorem}[Minimum Principle with free terminal time]
\tlabel{minEndTime} Suppose that $t_0$ and $x (t_0) = x_0$ are fixed,
and that for some index set $I\subset \{1,\dots,n\}$,
\begin{eqnarray*}
x_i(t_1) = x_{1i},\quad \mbox{ if } i\in I.
\end{eqnarray*}
Suppose that $u^\optchar$ is a solution to the optimal control problem
\eq min-opt/, subject to these constraints.  Then

\balphlist
\item
There exists a costate vector $p (t)$ such that
\begin{eqnarray*}
u^\optchar (t) = \argmin_{u} H (x^\optchar (t), p (t), u,t).
\end{eqnarray*}

\item
The pair $(p,x^\optchar)$ satisfy the 2-point boundary value problem
\eq 2ptBVP/, with the boundary conditions
\begin{eqnarray*}
x_i (t_i) & = & x_{1i}, \quad i\in I \\ p_j (t_1) & = & \frac{\partial
m}{\partial x_j} (x^\optchar (t_1), t_1), \quad i\in I^c.
\end{eqnarray*}
The unspecified terminal time $t_1$ satisfies
\begin{equation}
\frac{\partial m }{\partial t} (x^\optchar (t_1), t) + H (x^\optchar
(t_1),p (t_1),u^\optchar (t_1), t_1) = 0.  \elabel{termTime}
\end{equation}
\end{list}
\qed
\end{theorem}

The most common example in which the terminal time is free is a
minimum time problem.  As an example consider a single input linear
state space model
\begin{eqnarray*}
\dot{x} = Ax + bu.
\end{eqnarray*}
We wish to find the input $u^\optchar$ which drives $x$ from $x (0) =
x_0$ to $x (t_1) = x_1$ in \textit{minimum time}, under the constraint
$|u (t)| \leq 1$ for all $t$.  Hence the cost criterion we consider is
\begin{eqnarray*}
V (u) = t_1 = \int_0^{t_1} 1\, dt,
\end{eqnarray*}
so that $\ell\equiv 1$, and $m\equiv 0$.

The Hamiltonian for this problem is
\begin{eqnarray*}
H = \ell + p^Tf = 1 + p^T A x + p^Tbu.
\end{eqnarray*}
Minimizing over all $H$ again gives a bang-bang control law
\begin{eqnarray*}
u^\optchar (t) = \left\{\begin{array}{ll} 1 & p(t)^T b < 0 \\ -1 & p
(t)^T b > 0 \end{array} \right.
\end{eqnarray*}
If $A$, $b$ are time invariant and the eigenvalues $\{\lambda_i\}$ of
$A$ are real, distinct, and negative, then since the modes are all
decreasing, the quantity $p (t)^T b$ changes sign at most $n-1$
times. Hence, the number of switching times for the optimal control is
also bounded above by $n-1$.  \notes{give other boundary conditions?}

The costate equation is
\begin{eqnarray*}
\dot{p} = - \nabla_x H = -A^T p.
\end{eqnarray*}
Using the final time boundary condition \eq termTime/, or $\left. H
\right|_{t=t_1} = 0$, we obtain
\begin{eqnarray*}
1 + p^T (t_1) A x^\optchar (t_1) + p^T (t_1) bu^\optchar (t_1) = 1 +
p^T (t_1) A x^\optchar (t_1) - |p^T (t_1) b| = 0.
\end{eqnarray*}
We now show how these equations may be solved using a second order
example.

\begin{ex}
Consider the special case where the model is the double integrator
$\ddot y = u$, which is described by the state space model with
\begin{eqnarray*}
A = \left(\begin{array}{cc} 0 & 1 \\ 0 & 0 \end{array} \right) \qquad
b = \begin{pmatrix} 0 \\ 1 \end{pmatrix}
\end{eqnarray*}
so that $ \dot{x}_1 = x_2$, $ \dot{x}_2 = u$.  The optimal control has
the form $u^\optchar = -\sgn(p^T b) = -\sgn (p_2 (t))$. From the
costate equations
\begin{eqnarray*}
\begin{pmatrix}\dot p_1 \\ \dot p_2 \end{pmatrix}  = -A^T
\begin{pmatrix}\dot p_1 \\ \dot p_2 \end{pmatrix}
\end{eqnarray*}
we find that for constants $c_1$ and $c_2$,
\begin{eqnarray*}
p_1 (t) = c_1 \qquad p_2 (t) = -c_1 t + c_2.
\end{eqnarray*}
It follows that $p_2 (t)$ changes sign on $(0, t_1)$ \textit{at most
once}.

\begin{figure}[ht]
\ebox{.5}{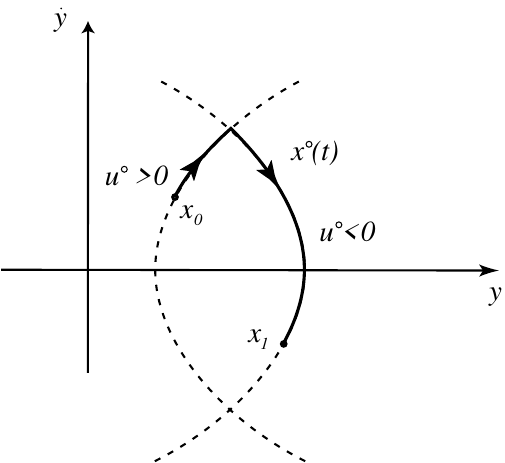}
\caption[Optimal trajectories for the minimum time problem]{The optimal trajectories for the minimum time problem in this
two dimensional example follow quadratic paths, and may switch paths
at most once.  A typical optimal path is illustrated as a union of two
dark line segments in this drawing.}  \flabel{optimal-minTime}
\end{figure}

If $p_2(t)<0$ then $ u^\optchar (t) = + 1$, so that by solving the
differential equations
\begin{eqnarray*}
\dot{x}_1^\optchar = x_2^\optchar, \qquad \dotx_2^\optchar = 1,
\end{eqnarray*}
we find that $x_1^\optchar (t) =\half (x_2^\optchar (t))^2 + K_1$ for
some constant $K_1$.  If $ p_2 (t) > 0$ so that $u^\optchar (t) = -
1$, then for another constant, $x_1^\optchar (t) = - \half
(x_2^\optchar (t))^2 + K_2$.  Hence the optimal state trajectories
follow the quadratic paths illustrated in \Figure{optimal-minTime}.
\end{ex}

\begin{exercises}
\item
Consider the system/performance index
\begin{eqnarray*}
\dot x = u \qquad V (u) = \frac{q}{2} x^2(T) + \frac{1}{2} \int_0^T
u^2\, dt,
\end{eqnarray*}
with $x(0) = x_0$, $T$ given, and $q\ge 0$.  Find $u^\optchar$ using
the Minimum Principle. First consider the case where $\infty>q>0$.
Then, obtain the control for ``$q=\infty$''.  That is, impose the
constraint that $x(T)=0$.  Do the controllers converge as
$q\to\infty$?

\item
A simplified model for a reactor is given by the bilinear system $\dot
x=xu$, where $x$ is the concentration, and $u$ is the rate, taken as
the control variable.  Suppose that $x(0)=0.5$, and the terminal time
is $t_1=1$. Using the Minimum Principle, find the control function
that minimizes the performance index
\begin{eqnarray*}
V =[x(1)-1]^2+\int_0^1u^2(t)\, dt.
\end{eqnarray*}

\item
The temperature of a room is described by \quad $\dot x =ax+u\;,$\quad
where $x$ is the room temperature, $a$ is a negative constant, and $u$
is the heat input rate. Suppose that $x(0)=x_0$ is given and it is
desired to find a control function $u(\cdot)$ so that $x(t_1)=q$,
where $q>x_0$, and the performance index
\begin{eqnarray*}
V=\int_0^{t_1}u^2(t)\,dt
\end{eqnarray*}
is minimized. Note that, in this problem, the state trajectory
$x(\cdot)$ is specified at two end points, while the terminal time
$t_1$ is free.  Using the Minimum Principle, determine the optimal
open-loop control function and the corresponding value of $t_1$, in
terms of $a$, $x_0$ and $q$.
\end{exercises}


\phantomsection
\addcontentsline{toc}{chapter}{Index} 

\begin{theindex}

  \item adjoint, \hyperpage{38}
  \item adjoint equation, \hyperpage{62}
  \item algebraic Riccati equation (ARE), \hyperpage{181}
  \item asymptotically stable, \hyperpage{68}

  \indexspace

  \item bang-bang, \hyperpage{215}
  \item basis, \hyperpage{26}
  \item bounded, \hyperpage{82}
  \item bounded input/bounded output stable, \hyperpage{82}

  \indexspace

  \item characteristic equation, \hyperpage{33}
  \item characteristic polynomial, \hyperpage{33}
  \item closed-loop system, \hyperpage{121}
  \item controllability grammian, \hyperpage{93}
  \item controllability matrix, \hyperpage{92}
  \item controllable, \hyperpage{90}
  \item controllable canonical form, \hyperpage{11}
  \item controllable subspace, \hyperpage{97}

  \indexspace

  \item detectable, \hyperpage{127}
  \item dimension, \hyperpage{26}
  \item Diophantine equation, \hyperpage{151}
  \item dual, \hyperpage{38}, \hyperpage{110}

  \indexspace

  \item eigenvalue, \hyperpage{32}
  \item eigenvector, \hyperpage{32}
  \item equilibrium, \hyperpage{67}
  \item equilibrium state, \hyperpage{4}

  \indexspace

  \item Feedback, \hyperpage{155}
  \item field, \hyperpage{23}
  \item fundamental matrix, \hyperpage{56}

  \indexspace

  \item generalized eigenvectors, \hyperpage{34}
  \item globally asymptotically stable, \hyperpage{68}
  \item Grammian, \hyperpage{37}

  \indexspace

  \item Hamilton-Jacobi-Bellman (HJB) equation, \hyperpage{168}, 
		\hyperpage{170}
  \item Hamiltonian, \hyperpage{170}
  \item Hamiltonian matrix, \hyperpage{175}
  \item Hermitian, \hyperpage{76}
  \item Hurwitz, \hyperpage{69}

  \indexspace

  \item image, \hyperpage{29}
  \item internal model principle, \hyperpage{146}
  \item invariant, \hyperpage{75}

  \indexspace

  \item Jordan form, \hyperpage{35}

  \indexspace

  \item Kalman Controllability Canonical Form, \hyperpage{111}
  \item Kalman Observability Canonical Form, \hyperpage{112}
  \item Kalman's inequality, \hyperpage{191}

  \indexspace

  \item Lagrange multipliers, \hyperpage{206}
  \item linear, \hyperpage{29}
  \item linear operator, \hyperpage{29}
  \item linear quadratic regulator, \hyperpage{168}
  \item linear state feedback, \hyperpage{121}
  \item linearity, \hyperpage{1}
  \item linearly independent, \hyperpage{25}
  \item loop gain, \hyperpage{158}
  \item loop transfer function, \hyperpage{187}
  \item Lyapunov equation, \hyperpage{61}, \hyperpage{77}, 
		\hyperpage{84}
  \item Lyapunov function, \hyperpage{84}

  \indexspace

  \item minimal, \hyperpage{113}
  \item minimum phase, \hyperpage{162}
  \item modal form, \hyperpage{14}
  \item modes, \hyperpage{14}

  \indexspace

  \item natural basis, \hyperpage{27}
  \item norm, \hyperpage{36}
  \item nullspace, \hyperpage{30}

  \indexspace

  \item observability grammian, \hyperpage{108}
  \item observability matrix, \hyperpage{108}
  \item observable, \hyperpage{107}
  \item observable canonical form, \hyperpage{13}
  \item orthogonal, \hyperpage{37}
  \item orthonormal, \hyperpage{37}

  \indexspace

  \item Partial fraction expansion, \hyperpage{13}
  \item Peano-Baker series, \hyperpage{58}
  \item performance, \hyperpage{155}
  \item positive definite, \hyperpage{45}, \hyperpage{72}, 
		\hyperpage{76}
  \item positive semi-definite, \hyperpage{45}, \hyperpage{77}
  \item principle minors, \hyperpage{77}
  \item principle of optimality, \hyperpage{169}

  \indexspace

  \item range, \hyperpage{29}
  \item rank, \hyperpage{29}
  \item realization, \hyperpage{10}
  \item reciprocal, \hyperpage{38}
  \item regulation, \hyperpage{155}
  \item return difference, \hyperpage{158}
  \item return difference equation, \hyperpage{187}
  \item Riccati Differential Equation, \hyperpage{174}

  \indexspace

  \item semigroup property, \hyperpage{49}, \hyperpage{57}
  \item sensitivity, \hyperpage{158}
  \item separation principle, \hyperpage{129}
  \item similar, \hyperpage{32}
  \item square, \hyperpage{161}
  \item stabilizable, \hyperpage{123}, \hyperpage{125}
  \item stable in the sense of Lyapunov, \hyperpage{67}
  \item stable subspace, \hyperpage{79}
  \item state transition matrix, \hyperpage{48}, \hyperpage{57}
  \item stationary point, \hyperpage{85}, \hyperpage{208}
  \item strictly proper, \hyperpage{13}
  \item subspace, \hyperpage{25}
  \item Sylvester equations, \hyperpage{152}
  \item symmetric root locus, \hyperpage{188}

  \indexspace

  \item tracking, \hyperpage{155}
  \item transmission zero, \hyperpage{162}

  \indexspace

  \item value function, \hyperpage{168}
  \item vector space, \hyperpage{24}

\end{theindex}

\end{document}